\documentclass[12pt]{article}
\usepackage{graphics,graphicx}
\usepackage{amssymb,amsmath,amsopn,amsfonts}
\usepackage[dvips]{color}
\usepackage{colordvi,multicol}
\usepackage{epsf}
\newfam\msbfam
\font\tenmsb=msbm10 \textfont\msbfam=\tenmsb \font\sevenmsb=msbm7
\scriptfont\msbfam=\sevenmsb \font\fivemsb=msbm5
\scriptscriptfont\msbfam=\fivemsb

\def\th#1{\vspace{1mm}\noindent{\bf #1}\quad}

\def\proof{\vspace{1mm}\noindent{\it Proof}\quad}

\numberwithin{equation}{section}

\def\bc{\begin{center}}
\def\ec{\end{center}}
\def\no{\noindent}
\def\hang{\hangindent\parindent}
\def\textindent#1{\indent\llap{\qquad #1\ \ \enspace}\ignorespaces}
\def\ref{\par\hang\textindent}

\begin{document}

\title{ {\bf Approximating three-dimensional Navier-Stokes equations driven by space-time white noise
\thanks{Research supported in part  by NSFC (No.11301026)}\\} }
\author{  {\bf Rongchan Zhu}$^{\mbox{a}}$, {\bf Xiangchan Zhu}$^{\mbox{b},}$\thanks{Corresponding author}
\date{}
\thanks{E-mail address:
zhurongchan@126.com(R. C. Zhu), zhuxiangchan@126.com(X. C. Zhu)}\\\\
$^{\mbox{a}}$Department of Mathematics, Beijing Institute of Technology, Beijing 100081,  China\\
$^{\mbox{b}}$School of Science, Beijing Jiaotong University, Beijing 100044, China  }

\maketitle

\noindent {\bf Abstract}

In this paper we study approximations to 3D Navier-Stokes (NS) equation driven by space-time white noise by  paracontrolled distribution proposed in [GIP13]. A solution theory for this equation has been developed recently in [ZZ14] based on regularity structure theory and paracontrolled distribution. In order to make the approximating equation converge to 3D NS equation driven by space-time white noise, we should subtract some drift terms in approximating equations. These drift terms, which come from renormalizations in the solution theory,  converge to the solution multiplied by some constant depending on approximations.

\vspace{1mm}
\no{\footnotesize{\bf 2000 Mathematics Subject Classification AMS}:\hspace{2mm} 60H15, 82C28}
 \vspace{2mm}

\no{\footnotesize{\bf Keywords}:\hspace{2mm}  stochastic Navier-Stokes equation, regularity structure, paracontrolled distribution, space-time white noise, renormalisation}

\section{Introduction}
In this paper, we study approximations to 3D Navier-Stokes equation driven by space-time white noise:
 Recall that the Navier-Stokes equations describe the time
evolution of an incompressible fluid and are given by
\begin{equation}\aligned d u= &(\nu \Delta u-P(u\cdot \nabla u))dt +PdW \\u(0)=&u_0, \quad div u = 0\endaligned\end{equation}
where $u(x, t)\in \mathbb{R}^3$ denotes the value of the velocity field at time $t$ and position
$x$, $P$ denotes the Leray projection, and $W$ is an external force field acting on
the fluid. We will consider the case when $x\in \mathbb{T}^3$, the three-dimensional torus.
Our mathematical model for the driving force $W$ is a cylindrical Wiener process on $L^2(\mathbb{T}^3;\mathbb{R}^3)$.

In one dimensional case, approximations to general stochastic partial differential equations driven by space-time white noise have been very well studied (see
[Gy98, Gy99, DG01, HMW14] and the reference therein). Now in this paper we consider the approximations to 3D NS equation driven by space-time white noise.

Random Navier-Stokes equations have been studied in many articles (see e.g. [FG95], [HM06], [De13], [DD02], [RZZ14] and the reference therein). For two dimensional case: existence and uniqueness of the  solutions have been obtained if the noisy forcing term is white in time and
colored (or white) in space. For three dimensional case,  existence of Markov solutions  for stochastic 3D Navier-Stokes equations driven by trace-class noise has been obtained in [FR08], [DD03].

The difficulty in treating (1.1) lies in the lack of spatial regularity of its solution and the nonlinear term cannot be well defined in the classical sense if driven by space-time white noise. Local existence and uniqueness of the solutions to (1.1) has been established recently in [ZZ14] based on the regularity structure theory introduced by Martin Hairer in [Hai14] and the paracontrolled distribution proposed by Gubinelli, Imkeller and  Perkowski in [GIP13].

In the theory of regularity structures,  the right objects, e.g. regularity structure that could
possibly take the place of Taylor polynomials can be constructed. The regularity can also be endowed with a model, which is a concrete way of associating every distribution to the abstract regularity structure. Multiplication, differentiation, the living space of the solutions,  and the convolution with singular kernel can be defined on this regularity structure and then the equation has been lifted on the regularity structure. On this regularity structure, the fixed point argument can be applied to obtain  local existence and uniqueness of the solutions. Furthermore, we can go back to the real world with the help of another central tool of the theory the reconstruction operator $\mathcal{R}$. If $W$ is a smooth  process, $\mathcal{R}u$ coincides with the classic solution of the equation.

The theory of paracontrolled distribution combines the idea of Gubinelli's controlled rough path [Gub04] and Bony's paraproduct [Bon84], which is defined by the following: Let $\Delta_jf$ be the jth Littlewood-Paley block of a distribution $f$, define
$$\pi_<(f,g)=\pi_>(g,f)=\sum_{j\geq-1}\sum_{i<j-1}\Delta_if\Delta_jg, \quad\pi_0(f,g)=\sum_{|i-j|\leq1}\Delta_if\Delta_jg.$$ Formally $fg=\pi_<(f,g)+\pi_0(f,g)+\pi_>(f,g)$. Observing that if $f$ is regular $\pi_<(f,g)$ behaves like $g$ and is the only term in the Bony's paraproduct not raising the regularities, the authors in [GIP13]  consider paracontrolled ansatz of the type
$$u=\pi_<(u',g)+u^\sharp,$$ where $\pi_<(u',g)$ represents the "bad-term" in the solution, $g$ is some distribution we can handle and $u^\sharp$ is regular enough to define the multiplication required. Then to make sense of the product of $uf$ we only need to define $gf$.

It is natural to try to use these two methods to prove the approximation result. A typical example of approximations to (1.1) is that if $u_\varepsilon$ solves
$$\aligned d u_\varepsilon= &(\nu \Delta_\varepsilon u_\varepsilon-P\sum_{j=1}^3D^\varepsilon_j(u_\varepsilon u_\varepsilon^j))dt +PdW, \endaligned\eqno(1.2)$$
where $$\Delta_\varepsilon F(x)=\sum_{l=1}^3\frac{F(x_1,...,x_l+\varepsilon,...,x_3)-2F(x)+F(x_1,...,x_l-\varepsilon,...,x_3)}{\varepsilon^2},$$
and $$D^\varepsilon_j F(x)=\frac{F(x_1,...,x_j+\varepsilon,...,x_3)-F(x_1,...,x_j,...,x_3)}{\varepsilon},$$
then we want to prove that $u_\varepsilon$ converges to $u$ as $\varepsilon\downarrow0$.
In this paper we use the paracontrolled distribution method to prove the result.

To make sense of (1.1) we need to do renormalizations for some ill-defined term and by some cancelations due to symmetry these renormalized terms disappear in (1.1). However, we have to modify the equation (1.2) by subtracting some drift terms from renormalizations (see (1.3) below), since there is no cancelations for the approximation equations. If not, the limiting process $\bar{u}=\lim_{\varepsilon\rightarrow0}u_\varepsilon$ turns out not to be a solution of (1.1). Instead it solves a similar equation with an additional drift term. This extra term depends on the specific choice of approximations and it can be calculated explicitly. This situation is similar to the case in [HM12,  HMW14], where they consider the one dimensional Burgers-like equations by using rough path theory (see also [HW13]).

One motivation for this paper is to illustrate how  to obtain concrete approximation results by using the paracontrolled distribution method. This is
particularly interesting for the KPZ equation [KPZ86]
$$\partial_th=\partial^2_x h+\lambda (\partial_x h)^2-\infty+\xi,$$
where $\xi$ denotes space-time white noise and “$\infty$” denotes an “infinite constant” that
needs to be subtracted in order to make sense of the diverging term $(\partial_x h)^2$ (see [Hai13]). This
equation is an important model for surface growth and it is related to the KPZ universality class (see e.g. [Cor11] and the references
therein).
The present article illustrate how one can obtain approximation
results for 3D Navier-Stokes equation driven by space-time white noise, which has similar singularity as the KPZ equation. The approximation result for KPZ equation will be studied in our future work.

The first difficulty we encounter is that Schauder estimate for $\Delta_\varepsilon$ defined above, which is essential for the proof,  does not hold any more. Here we have to modify the operator such that the approximated heat semigroup $e^{t\Delta_\varepsilon}$ only involve the projection onto a finite-dimensional subspace (see Example 1.1 (i)), which is natural
in the context of numerical approximations. However this will cause another problem. In general we need some differentiability for the corresponding Fourier multiplier (see Mihlin multiplier theorem) to prove this result. In one dimensional case we can  overcome this difficulty  by applying Marcinkiewicz multiplier theorem to control the bounded variation norm of multiplier (see [HMW14]). But in three dimensional case, the Marcinkievicz multiplier theorem also needs some differentiability for the corresponding Fourier multiplier and cannot be used here.  However, under our assumptions we can view the Fourier multiplier as a smooth function multiply some $L^p$ Fourier multiplier, $p>1$,  and prove the Schauder estimate.

The second difficulty lies in how to prove the commutator estimate for the approximated heat semigroup $e^{t\Delta_\varepsilon}$ since the original proof based on Taylor expansion, which requires some differentiability for the corresponding Fourier multiplier. Here we prove the commutator estimate  for a special case (see Lemmas 3.5 and 3.6). This can be used in the case that we approximate  $W$ by only
keeping finite dimensional Fourier modes which is also natural in the context of numerical approximations.

\textbf{Framework and main result}

For $\varepsilon>0$  we consider the following approximating stochastic PDEs given by
$$\aligned d u^{\varepsilon,i_0}=&(\Delta_\varepsilon u^{\varepsilon,i_0}-\frac{1}{2}\sum_{i,j=1}^3P^{i_0 i}D_j^\varepsilon(u^{\varepsilon,i} u^{\varepsilon,j}+\sum_{i_1=1}^3(C^{\varepsilon,i,i_1,j}+\tilde{C}^{\varepsilon,i,i_1,j}+C^{\varepsilon,j,i_1,i}+\tilde{C}^{\varepsilon,j,i_1,i}) u^{\varepsilon,i_1}))dt\\&+\sum_{i=1}^3P^{i_0 i}H_\varepsilon dW^i
\\u^\varepsilon(0)=&Pu^\varepsilon_0,\endaligned\eqno(1.3)$$
for $i_0=1,2,3$.
Here as (1.1) $W$ is a cylindrical Wiener process on $L^2(\mathbb{T}^3;\mathbb{R}^3)$, $P$ is the Leray projection.
The operators $\Delta_\varepsilon, D^\varepsilon_j,H_\varepsilon$ are given by their action in Fourier space
$$\widehat{\Delta_\varepsilon v}(k)=-|k|^2f(\varepsilon k)\hat{v}(k),$$
$$\widehat{D_j^\varepsilon v}(k)=k^jg(\varepsilon k^j)\hat{v}(k),$$
$$\widehat{H_\varepsilon W}(k)= h(\varepsilon k)\hat{W}(k).$$
Throughout the paper we assume that there exists some $L_0>0$ such that
$$f(x)=\left\{\begin{array}{ll}\tilde{f}(x),&\ \ \ \ \textrm{ if } |x^1|\leq L_0,|x^2|\leq L_0,|x^3|\leq L_0\\\infty&\ \ \ \ \textrm{ otherwise,}\end{array}\right.$$
where $\tilde{f}:\mathbb{R}^3\rightarrow[0,\infty)$ is even satisfying $\tilde{f}(0)=1$ and for $|x^1|\leq 3L_0,|x^2|\leq 3L_0,|x^3|\leq3L_0$, $|D_k\tilde{f}(x)|\lesssim \frac{1}{|x|^{|k|-1}}+C$ for $|k|\leq5$, $\tilde{f}(x)\geq c_f>0$.  Moreover we assume that $$g(x)=\frac{e^{ax\imath}-e^{-bx\imath}}{(a+b)x}$$ for some $a,b\geq0, a+b>0$
and
$h(x)$ is a bounded even function and is continuously differentiable on $\{|\cdot|\leq \bar{L}_0\}$ for some $\bar{L}_0\leq L_0/2$ satisfying ${h}(0)=1$, and $\textrm{supp} h\subset \{x:|x|\leq L_0/2\}$.

As mentioned before here we should subtract some drift term in (1.2) and for $i,i_1,j=1,2,3$,  $$\aligned C^{\varepsilon,i,i_1,j}(t)=&\frac{1}{\varepsilon}(2\pi)^{-3}\sum_{i_2,i_3=1}^3\sum_{k_2\in\mathbb{Z}^3\backslash\{0\}}(1-e^{-2|k_2|^2tf(\varepsilon k_2)}) (\cos(a\varepsilon k_2^{i_2})-\cos(b\varepsilon k_2^{i_2}))\frac{h(\varepsilon k_2)^2}{8(a+b)|k_2|^4f(\varepsilon k_2)^2}\\&\hat{P}^{ii_1}(k_2)\hat{P}^{i_2i_3}(k_{2})\hat{P}^{ji_3}(k_{2})\\\rightarrow &\Lambda:=\sum_{i_2,i_3=1}^3(2\pi)^{-3}\int_{\mathbb{R}^3} (\cos ax^{i_2}-\cos bx^{i_2})\frac{h(x)^2}{8(a+b)|x|^4f(x)^2}\hat{P}^{ii_1}(x)\hat{P}^{i_2i_3}(x)\hat{P}^{ji_3}(x)dx,\endaligned$$
as $\varepsilon\rightarrow0$ and
$$\aligned \tilde{C}^{\varepsilon,i,i_1,j}(t)=&\frac{1}{\varepsilon}(2\pi)^{-3}\sum_{i_2,i_3=1}^3\sum_{k_2\in\mathbb{Z}^3\backslash\{0\}}(1-e^{-2|k_2|^2tf(\varepsilon k_2)}) (\cos(a\varepsilon k_2^{i_1})-\cos(b\varepsilon k_2^{i_1}))\frac{h(\varepsilon k_2)^2}{8(a+b)|k_2|^4f(\varepsilon k_2)^2}\\&\hat{P}^{ii_2}(k_2)\hat{P}^{i_2i_3}(k_{2})\hat{P}^{ji_3}(k_{2})\\\rightarrow &\Lambda_1:=\sum_{i_2,i_3=1}^3(2\pi)^{-3}\int_{\mathbb{R}^3} (\cos ax^{i_1}-\cos bx^{i_1})\frac{h(x)^2}{8(a+b)|x|^4f(x)^2}\hat{P}^{ii_2}(x)\hat{P}^{i_2i_3}(x)\hat{P}^{ji_3}(x)dx,\endaligned$$
as $\varepsilon\rightarrow0$. Here $\hat{P}^{ii_1}(k)=\delta_{ii_1}-\frac{k^ik^{i_1}}{|k|^2}$ for $k\in\mathbb{R}^3\backslash\{0\}$. Note that $\Lambda, \Lambda_1$ are indeed well-defined by assumptions and depend on the specific choice of approximation. Before we proceed, we list some of the most common examples of discretizations that do fit our framework which is similar to examples in [HM12].

\vskip.10in
\th{Example 1.1} i)\emph{Finite difference discretization.} In this case we take $$\Delta_\varepsilon F(x)=\sum_{l=1}^3\frac{F(x_1,...,x_l+\varepsilon,...,x_3)-2F(x)+F(x_1,...,x_l-\varepsilon,...,x_3)}{\varepsilon^2}.$$
We also approximate $\Delta_\varepsilon$ and $W$ by only
keeping finite Fourier modes.
This corresponds to the choice $$\tilde{f}(x)=\frac{4}{|x|^2}(\sin^2\frac{x_1}{2}+\sin^2\frac{x_2}{2}+\sin^2\frac{x_3}{2}), \quad h(x)=1_{\{|x|\leq L_0/2\}}.$$

ii)\emph{Galerkin discretization}. In this case, we approximate $\Delta$ and $W$ by only
keeping finite dimensional Fourier modes. This corresponds to the choice
$$\tilde{f}(x)=1,\quad h(x)=1_{\{|x|\leq L_0/2\}}.$$
\vskip.10in
\th{Remark 1.2}  Note that $$D^\varepsilon_j F(x)=\frac{F(x_1,...,x_j+a\varepsilon,...,x_3)-F(x_1,...,x_j-b\varepsilon,...,x_3)}{(a+b)\varepsilon}$$ for some $a,b\geq0, a+b>0$.

\vskip.10in
The main result of this paper is the following theorem:
\vskip.10in
\th{Theorem 1.3} Let $z\in (1/2,1/2+\delta_0)$ with $0<\delta_0<1/2$ and $u_0\in \mathcal{C}^{-z}$. Let $(u,\tau)$ be the unique maximal solution of the following equation
$$du^i=\Delta u^i+\sum_{i_1=1}^3P^{ii_1}dW^{i_1}-\frac{1}{2}\sum_{i_1=1}^3P^{ii_1}(\sum_{j=1}^3D_j(u^{i_1}u^{j}))dt\quad u(0)=u_0,\eqno(1.4)$$
and let for $\varepsilon\in(0,1)$ the function $u^\varepsilon$ be the unique maximal solution to (1.3). If the initial data satisfies $u_0^\varepsilon-u_0\rightarrow0 \textrm{ in }\mathcal{C}^{-z}$ then there exists a sequence of random time  $\tau_L$ such that  $\lim_{L\rightarrow\infty}\tau_L=\tau$ and $$\sup_{t\in[0,\tau_L]}\|u^\varepsilon-u\|_{-z}\rightarrow0\quad \textrm{ in probability}, \quad \textrm{ as } \varepsilon\rightarrow0 .$$
\vskip.10in
\th{Remark 1.4} i) For the definition of $\mathcal{C}^{-z}$ and norm $\|\cdot\|_{-z}$ see Section 2.

ii) Indeed by a modification we can prove that there exists $r>0$ and  a sequence of random time  $\tau_\varepsilon$ satisfying  $\lim_{\varepsilon\rightarrow0}\tau_\varepsilon=\tau$  such that
$$P(\sup_{t\in[0,\tau_\varepsilon]}\|u^\varepsilon-u\|_{-z}>\varepsilon^r)\rightarrow0, \textrm{ as } \varepsilon\rightarrow0.$$
This gives the convergence rate. We suspect this to be the true rate of convergence and we omit it here for simpicity.
\vskip.10in
\textbf{Structure of the paper} This paper is organized as follows. In Section 2, we recall some basic notions and results in paracontrolled distribution method. In Section 3 we first prove some estimates for the approximated operators. Then we construct  solutions to (1.3) and (1.4) and prove uniform bounds by a similar argument as [ZZ14]. Finally at the end of this section we give the proof of our main result. In Section 4  convergence of  terms involving the renormalized terms is proved.
\section{Besov spaces and paraproduct}

In the following we recall the definitions and some properties of Besov spaces and paraproducts. For a general introduction to these theories we refer to [BCD11, GIP13].
First we introduce the following notations. The space of real valued infinitely differentiable functions of compact support is denoted by $\mathcal{D}(\mathbb{R}^d)$ or $\mathcal{D}$. The space of Schwartz functions is denoted by $\mathcal{S}(\mathbb{R}^d)$. Its dual, the space of tempered distributions is denoted by $\mathcal{S}'(\mathbb{R}^d)$. If $u$ is a vector of $n$ tempered distributions on $\mathbb{R}^d$, then we write $u\in \mathcal{S}'(\mathbb{R}^d,\mathbb{R}^n)$. The Fourier transform and the inverse Fourier transform are denoted by $\mathcal{F}u$ and $\mathcal{F}^{-1}u$.

 Let $\chi,\theta\in \mathcal{D}$ be nonnegative radial functions on $\mathbb{R}^d$, such that

i. the support of $\chi$ is contained in a ball and the support of $\theta$ is contained in an annulus;

ii. $\chi(z)+\sum_{j\geq0}\theta(2^{-j}z)=1$ for all $z\in \mathbb{R}^d$.

iii. $\textrm{supp}(\chi)\cap \textrm{supp}(\theta(2^{-j}\cdot))=\emptyset$ for $j\geq1$ and $\textrm{supp}(\theta(2^{-i}\cdot))\cap \textrm{supp}(\theta(2^{-j}\cdot))=\emptyset$ for $|i-j|>1$.

We call such $(\chi,\theta)$ dyadic partition of unity, and for the existence of dyadic partitions of unity see [BCD11, Proposition 2.10]. The Littlewood-Paley blocks are now defined as
$$\Delta_{-1}u=\mathcal{F}^{-1}(\chi\mathcal{F}u)\quad \Delta_{j}u=\mathcal{F}^{-1}(\theta(2^{-j}\cdot)\mathcal{F}u).$$

For $\alpha\in \mathbb{R}$, the H\"{o}lder-Besov space $\mathcal{C}^\alpha$ is given by $\mathcal{C}^\alpha=B^\alpha_{\infty,\infty}(\mathbb{R}^d,\mathbb{R}^n)$, where for $p,q\in [1,\infty]$ we define
$$B^\alpha_{p,q}(\mathbb{R}^d,\mathbb{R}^n)=\{u=(u^1,...,u^n)\in\mathcal{S}'(\mathbb{R}^d,\mathbb{R}^n):\|u\|_{B^\alpha_{p,q}}=\sum_{i=1}^n(\sum_{j\geq-1}(2^{j\alpha}\|\Delta_ju^i\|_{L^p})^q)^{1/q}<\infty\},$$
with the usual interpretation as $l^\infty$ norm in case $q=\infty$. We write $\|\cdot\|_{\alpha}$ instead of $\|\cdot\|_{B^\alpha_{\infty,\infty}}$ in the following for simplicity.

We point out that everything above and everything that follows can be applied to distributions on the torus (see [S85, SW71]). More precisely, let $\mathcal{S}'(\mathbb{T}^d)$ be the space of distributions on $\mathbb{T}^d$. Therefore, Besov spaces on the torus with general indices $p,q\in[1,\infty]$ are defined as
$$B^\alpha_{p,q}(\mathbb{T}^d,\mathbb{R}^n)=\{u\in\mathcal{S}'(\mathbb{T}^d,\mathbb{R}^n):\|u\|_{B^\alpha_{p,q}}=\sum_{i=1}^n(\sum_{j\geq-1}(2^{j\alpha}\|\Delta_ju^i\|_{L^p(\mathbb{T}^d)})^q)^{1/q}<\infty\}.$$
 We  will need the following Besov embedding theorem on the torus (c.f. [GIP13, Lemma 41]):
\vskip.10in
 \th{Lemma 2.1} Let $1\leq p_1\leq p_2\leq\infty$ and $1\leq q_1\leq q_2\leq\infty$, and let $\alpha\in\mathbb{R}$. Then $B^\alpha_{p_1,q_1}(\mathbb{T}^d)$ is continuously embedded in $B^{\alpha-d(1/p_1-1/p_2)}_{p_2,q_2}(\mathbb{T}^d)$.
\vskip.10in

 Now we recall the following paraproduct introduced by Bony (see [Bon81]). In general, the product $fg$ of two distributions $f\in \mathcal{C}^\alpha, g\in \mathcal{C}^\beta$ is well defined if and only if $\alpha+\beta>0$. In terms of Littlewood-Paley blocks, the product $fg$ can be formally decomposed as
 $$fg=\sum_{j\geq-1}\sum_{i\geq-1}\Delta_if\Delta_jg=\pi_<(f,g)+\pi_0(f,g)+\pi_>(f,g),$$
 with $$\pi_<(f,g)=\pi_>(g,f)=\sum_{j\geq-1}\sum_{i<j-1}\Delta_if\Delta_jg, \quad\pi_0(f,g)=\sum_{|i-j|\leq1}\Delta_if\Delta_jg.$$
We also use the notation for $j\geq0$
$$S_jf=\sum_{i\leq j-1}\Delta_if.$$
We will use without comment that $\|\cdot\|_\alpha\leq\|\cdot\|_\beta$ for $\alpha\leq\beta$, that $\|\cdot\|_{L^\infty}\lesssim \|\cdot\|_\alpha$ for $\alpha>0$, and that $\|\cdot\|_\alpha\lesssim\|\cdot\|_{L^\infty}$ for $\alpha\leq0$. We will also use that $\|S_ju\|_{L^\infty}\lesssim 2^{-j\alpha}\|u\|_\alpha$ for $\alpha<0, j\geq0$ and $u\in \mathcal{C}^\alpha$, where $\|\cdot\|_\alpha$ denotes the norm in $\mathcal{C}^\alpha, \alpha\in\mathbb{R}$.

\vskip.10in
 The basic result about these bilinear operations is given by the following estimates:
\vskip.10in
 \th{Lemma 2.2}(Paraproduct estimates, [Bon 81, GIP13, Lemma 2]) For any $\beta\in \mathbb{R}$ we have
 $$\|\pi_<(f,g)\|_\beta\lesssim \|f\|_{L^\infty}\|g\|_\beta\quad f\in L^\infty, g\in \mathcal{C}^\beta,$$
 and for $\alpha<0$ furthermore
 $$\|\pi_<(f,g)\|_{\alpha+\beta}\lesssim \|f\|_{\alpha}\|g\|_\beta\quad f\in \mathcal{C}^\alpha, g\in \mathcal{C}^\beta.$$
 For $\alpha+\beta>0$ we have
 $$\|\pi_0(f,g)\|_{\alpha+\beta}\lesssim \|f\|_{\alpha}\|g\|_\beta\quad f\in \mathcal{C}^\alpha, g\in \mathcal{C}^\beta.$$

\vskip.10in
 The following basic commutator lemma is important for our use:
\vskip.10in
 \th{Lemma 2.3}([GIP13, Lemma 5]) Assume that $\alpha\in (0,1)$ and $\beta,\gamma\in \mathbb{R}$ are such that $\alpha+\beta+\gamma>0$ and $\beta+\gamma<0$. Then for smooth $f,g,h,$ the trilinear operator
 $$C(f,g,h)=\pi_0(\pi_<(f,g),h)-f\pi_0(g,h)$$ allows for the bound
 $$\|C(f,g,h)\|_{\alpha+\beta+\gamma}\lesssim\|f\|_\alpha\|g\|_\beta\|h\|_\gamma.$$
 Thus, $C$ can be uniquely extended to a bounded trilinear operator from $\mathcal{C}^\alpha\times \mathcal{C}^\beta \times \mathcal{C}^\gamma$ to $ \mathcal{C}^{\alpha+\beta+\gamma}$.

\vskip.10in

Now we recall the following commutator estimate from [ZZ14, Lemma 3.4, Lemma 3.6].
\vskip.10in
\th{Lemma 2.4} Let $u\in \mathcal{C}^\alpha$ for some $\alpha<1$ and $v\in \mathcal{C}^\beta$ for some $\beta\in \mathbb{R}$. Then for  every $k,l=1,2,3$ we have
$$\|P^{kl}\pi_{<}(u,v)-\pi_<(u,P^{kl}v)\|_{\alpha+\beta}\lesssim \|u\|_\alpha\|v\|_\beta,$$
where  $P$ is the Leray projection.

\vskip.10in

\th{Lemma 2.5} Let $u\in \mathcal{C}^\alpha$ for some $\alpha\in \mathbb{R}$. Then we have for every $k,l=1,2,3$
$$\|P^{kl}u\|_{\alpha}\lesssim \|u\|_\alpha,$$
where  $P$ is the Leray projection.
\vskip.10in
Now we recall the following estimate for heat semigroup $P_t:=e^{t\Delta}$.
\vskip.10in
\th{Lemma 2.6}([GIP13, Lemma 47]) Let $u\in \mathcal{C}^\alpha$ for some $\alpha\in \mathbb{R}$. Then for every $\delta\geq0$
$$\|P_tu\|_{\alpha+\delta}\lesssim t^{-\delta/2}\|u\|_\alpha.$$

\vskip.10in

\th{Lemma 2.7} ([CC13, Lemma A.1]) Let $u\in \mathcal{C}^\alpha$ for some $\alpha<1$ and $v\in \mathcal{C}^\beta$ for some $\beta\in \mathbb{R}$. Then for $\delta\geq \alpha+\beta$
$$\|P_t\pi_{<}(u,v)-\pi_<(u,P_t v)\|_{\delta}\lesssim t^{\frac{\alpha+\beta-\delta}{2}}\|u\|_\alpha\|v\|_\beta.$$
\vskip.10in

\th{Lemma 2.8} ([CC13, Lemma 2.5])  Let $u\in \mathcal{C}^{\alpha+\delta}$ for some $\alpha\in \mathbb{R},\delta>0$. Then for every $ t\geq0$
$$\|(P_t-I)u\|_{\alpha}\lesssim  t^{\delta/2}\|u\|_{\alpha+\delta}.$$

\section{Proof of the main result}

\subsection{Estimates for the approximated operators}

In the following without loss of generality we assume $\textrm{supp} \theta(2^{-j}\cdot)\subset \{\xi:2^{j-1}\leq|\xi|\leq2^{j+1}\}$. 
As we mentioned in introduction, in general we need some differentiability for the corresponding Fourier multiplier (see Mihlin multiplier theorem) to prove this result. In one dimensional case we can  overcome this difficulty  by applying Marcinkiewicz multiplier theorem to control the bounded variation norm of multiplier (see [HMW14]). But in three dimensional case, the Marcinkievicz multiplier theorem also needs some differentiability for the corresponding Fourier multiplier and cannot be used here.  However, under our assumptions we can view the Fourier multiplier as a smooth function multiply some $L^p$ Fourier multiplier, $p>1$,  and prove the Schauder estimate.
Now we introduce the following notations and operators from [Tri78, Section 2.11.2]. Let
$$q_k=\{x||x_j|\leq 2^k,\textrm{ for }j=1,2,3\}, \qquad k=0,1,2,...,$$
and $$Q_0=q_0,\quad Q_k=q_k-q_{k-1}, \qquad k=0,1,2,....$$
The characteristic function $\chi_k$ of $Q_k$ is a multiplier on $L^p,1<p<\infty$.
\vskip.10in
First  we have the following equivalent norm, which  makes the estimates of the approximated operators much easier:
\vskip.10in
\th{Lemma 3.1} Let $1<p<\infty,$ $s\in\mathbb{R}$. $\|\{\mathcal{F}^{-1}\chi_\cdot\mathcal{F}u\}\|_{l_\infty^s(L^p)}:=\sup_j2^{js}\|\mathcal{F}^{-1}\chi_j\mathcal{F}u\|_{L^p}$ is an equivalent norm of $B^s_{p,\infty}$.

\proof A homogeneity argument shows that there exists a constant $C$ independent of $j$ such that for $j\geq0$
$$\|\mathcal{F}^{-1}\chi_j\mathcal{F}u\|_{L^p}\leq C\|u\|_{L^p},$$
and $$\|\mathcal{F}^{-1}\theta(2^{-j}\cdot)\mathcal{F}u\|_{L^p}\leq C\|u\|_{L^p}.$$
 Then we have for $j\geq0$
 $$2^{js}\|\mathcal{F}^{-1}\theta(2^{-j}\cdot)\mathcal{F}u\|_{L^p}\leq 2^{js}\sum_{k=j-1}^{j+2}\|\mathcal{F}^{-1}\theta(2^{-j}\cdot)\chi_k\mathcal{F}u\|_{L^p}\lesssim \|\{\mathcal{F}^{-1}\chi_\cdot\mathcal{F}u\}\|_{l_\infty^s(L^p)},$$
and for $j=-1$,
$$\|\mathcal{F}^{-1}\chi\mathcal{F}u\|_{L^p}\leq \|\mathcal{F}^{-1}\chi\chi_0\mathcal{F}u\|_{L^p}\lesssim\|\{\mathcal{F}^{-1}\chi_\cdot\mathcal{F}u\}\|_{l_\infty^s(L^p)};$$
on the other hand, for $j\geq1$
$$2^{js}\|\mathcal{F}^{-1}\chi_j\mathcal{F}u\|_{L^p}\leq 2^{js}\sum_{k=j-1}^{j+1}\|\mathcal{F}^{-1}\chi_j\theta_k\mathcal{F}u\|_{L^p}\lesssim \|u\|_{B^s_{p,\infty}},$$
and for $j=0$,
$$\|\mathcal{F}^{-1}\chi_0\mathcal{F}u\|_{L^p}\leq\|\mathcal{F}^{-1}\chi_0(\chi+\theta)\mathcal{F}u\|_{L^p}\lesssim \|u\|_{B^s_{p,\infty}}.$$
Thus the result follows.
$\hfill\Box$
\vskip.10in

\th{Lemma 3.2} Let $u\in \mathcal{C}^\alpha$ for some $\alpha\in \mathbb{R}$. Then for every $\delta\geq0,\kappa>0, t>0$
$$\sup_{\varepsilon\in(0,1)}\|P^\varepsilon_tu\|_{\alpha+\delta-\kappa}\lesssim t^{-\delta/2}\|u\|_\alpha.$$

\proof Without loss of generality we suppose that $L_0>1$. First we consider the operator on $\mathbb{R}^3$. For  $N\in\mathbb{N}$ and $2^{-N}L_0\leq\varepsilon<2^{-N+1}L_0,0<j\leq N, p>1,$ we obtain
$$\aligned &\|\mathcal{F}^{-1}\chi_j\mathcal{F}P^\varepsilon_tu\|_{L^p}=\|\mathcal{F}^{-1}\chi_j\varphi^\varepsilon\mathcal{F}u\|_{L^p}\\\leq& \|\mathcal{F}^{-1}\chi_j\varphi^\varepsilon\phi_1(2^{-j}\cdot)\mathcal{F}u\|_{L^p}+\|\mathcal{F}^{-1}\chi_j\varphi^\varepsilon\phi_2(2^{-j}\cdot)\mathcal{F}u\|_{L^p}+\|\mathcal{F}^{-1}\chi_j\varphi^\varepsilon\phi_3(2^{-j}\cdot)\mathcal{F}u\|_{L^p}.\endaligned$$
Here $$\varphi^\varepsilon(\xi)=e^{-t|\xi|^2f(\varepsilon\xi)},$$ and $$(\phi_1(2^{-j}\cdot)+\phi_2(2^{-j}\cdot)+\phi_3(2^{-j}\cdot))\varphi^\varepsilon\chi_j=\varphi^\varepsilon\chi_j$$ with $\phi_i$ smooth and $\textrm{supp}\phi_i\subset \{\frac{1}{4}\leq |\xi^i|\leq \frac{3}{2},|\xi^l|\leq \frac{3}{2},l\neq i\}$. Moreover, we have for $0<j\leq N$ $$\chi_j(\xi)\phi_i(2^{-j}\xi)\varphi^\varepsilon(\xi)=\chi_j(\xi)\phi_i(2^{-j}\xi)\tilde{\varphi}^\varepsilon(\xi)1_{\{|\xi^l|\leq \varepsilon^{-1}L_0,l=1,2,3\}},$$
with $\tilde{\varphi}^\varepsilon(\xi)= e^{-t|\xi|^2\tilde{f}(\varepsilon\xi)}$,  which implies that
$$\|\mathcal{F}^{-1}\chi_j\varphi^\varepsilon\phi_i(2^{-j}\cdot)\mathcal{F}u\|_{L^p}\lesssim \|\mathcal{F}^{-1}\tilde{\varphi}^\varepsilon\phi_i(2^{-j}\cdot)\mathcal{F}\mathcal{F}^{-1}\chi_j\mathcal{F}u\|_{L^p}\leq \|\mathcal{F}^{-1}\tilde{\varphi}^\varepsilon\phi_i(2^{-j}\cdot)\|_{L^1}\|\mathcal{F}^{-1}\chi_j\mathcal{F}u\|_{L^p},$$
where we used that $1_{\{|\xi^l|\leq \varepsilon^{-1}L_0,l=1,2,3\}}$ is an $L^p$-multiplier.
By calculation on $\mathbb{R}^3$ we get that for $\delta\geq0$ $$\aligned&\|\mathcal{F}^{-1}(\tilde{\varphi}^\varepsilon\phi_i(2^{-j}\cdot))\|_{L^1}=\|\mathcal{F}^{-1}(\tilde{\varphi}^\varepsilon(2^{j}\cdot)\phi_i)\|_{L^1}
\lesssim \|(1-\Delta)^2(\tilde{\varphi}^\varepsilon(2^j\cdot)\phi_i(\cdot))\|_{L^1}\\\lesssim& \sum_{0\leq|k|\leq 4}2^{j|k|}\|(D_k\tilde{\varphi}^\varepsilon)(2^j\cdot)|_{\cdot\in\textrm{supp} \phi_i}\|_{L^\infty}\lesssim \sum_{0\leq|k|\leq 4}2^{j|k|}\frac{1}{2^{j|k|}(2^j\sqrt{t})^{\delta}}\lesssim (2^j\sqrt{t})^{-\delta},\endaligned$$
which yields that for $0<j\leq N$
$$\|\mathcal{F}^{-1}\chi_j\mathcal{F}P^\varepsilon_tu\|_{L^p}\lesssim (2^j\sqrt{t})^{-\delta}\|\mathcal{F}^{-1}\chi_j\mathcal{F}u\|_{L^p}. $$
Here in the second inequality we used $\phi_i$ has compact support and $\inf_{\xi\in\textrm{supp}\phi_i}\tilde{f}(\varepsilon2^j\xi)\geq c>0$ and $|D_k\tilde{f}(\xi)|\lesssim \frac{1}{|\xi|^{|k|-1}}+C$ for any multiindices $k$ satisfying $|k|\leq4$ to deduce that for every $\delta\geq0$ $|D_k\tilde{\varphi}^\varepsilon(\xi)|\lesssim \frac{1}{|\xi|^{|k|+\delta}t^{\delta/2}}$.

For $j>N$ we have that
$$\|\mathcal{F}^{-1}\chi_j\mathcal{F}P^\varepsilon_tu\|_{L^p}=0.$$

For $j=0$ by the estimate above and the Mihlin multiplier theorem we obtain $$\|\mathcal{F}^{-1}\chi_0\varphi^\varepsilon\mathcal{F}u\|_{L^p}= \|\mathcal{F}^{-1}\tilde{\varphi}^\varepsilon\mathcal{F}\mathcal{F}^{-1}\chi_0\mathcal{F}u\|_{L^p}\lesssim \|\mathcal{F}^{-1}\chi_0\mathcal{F}u\|_{L^p}.$$
Since all the constants above are independent of $\varepsilon$,  we have on $\mathbb{R}^3$
$$\sup_{\varepsilon\in(0,1)}\|P^\varepsilon_tu\|_{B^{\alpha+\delta}_{p,\infty}} \lesssim t^{-\delta/2}\|u\|_{B^{\alpha}_{p,\infty}}.$$
By the theory in [SW71] we know that the above calculations also hold on $\mathbb{T}^3$.
Thus on $\mathbb{T}^3$ we have for $p$ large enough by Lemma 2.1 $$\sup_{\varepsilon\in (0,1)}\|P^\varepsilon_tu\|_{\alpha+\delta-\kappa}\lesssim\sup_{\varepsilon\in(0,1)}\|P^\varepsilon_tu\|_{B^{\alpha+\delta}_{p,\infty}} \lesssim t^{-\delta/2}\|u\|_{B^{\alpha}_{p,\infty}} \lesssim t^{-\delta/2}\|u\|_\alpha.$$
$\hfill\Box$
\vskip.10in
\th{Lemma 3.3} Let $u\in \mathcal{C}^{\alpha+\eta}$ for some $\alpha\in \mathbb{R},0<\eta<1$. Then  for every $\delta\geq0,\kappa>0,\varepsilon\in(0,1)$
$$\|(P^\varepsilon_t-P_t)u\|_{\alpha+\delta-\kappa}\lesssim \varepsilon^\eta t^{-\delta/2}\|u\|_{\alpha+\eta}.$$

\proof  Without loss of generality we assume that $L_0>1$. First we consider the operator $\bar{P}_t^\varepsilon-P_t=\mathcal{F}^{-1}(\bar{\varphi}^\varepsilon-\varphi)\mathcal{F}$ on $\mathbb{R}^3$ with
$$\bar{\varphi}^\varepsilon(\xi)=\left\{\begin{array}{ll}\varphi(\xi)&\ \ \ \ \textrm{ for } |\xi_i|\leq 1/2,\\\varphi^\varepsilon(\xi)&\ \ \ \ \textrm{ otherwise}.\end{array}\right.$$Here $$\varphi^\varepsilon(\xi)=e^{-t|\xi|^2f(\varepsilon\xi)},\qquad\varphi(\xi)=e^{-t|\xi|^2}.$$
 For $N\in\mathbb{N}$ and $L_02^{-N}\leq\varepsilon<2^{-N+1}L_0,0<j\leq N, p>1$ we obtain that
$$\aligned &\|\mathcal{F}^{-1}\chi_j\mathcal{F}(\bar{P}^\varepsilon_t-P_t)u\|_{L^p}=\|\mathcal{F}^{-1}\chi_j(\bar{\varphi}^\varepsilon-\varphi)\mathcal{F}u\|_{L^p}\\\leq& \|\mathcal{F}^{-1}\chi_j(\bar{\varphi}^\varepsilon-\varphi)\phi_1(2^{-j}\cdot)\mathcal{F}u\|_{L^p}+\|\mathcal{F}^{-1}\chi_j(\bar{\varphi}^\varepsilon-\varphi)\phi_2(2^{-j}\cdot)\mathcal{F}u\|_{L^p}+\|\mathcal{F}^{-1}\chi_j(\bar{\varphi}^\varepsilon-\varphi)\phi_3(2^{-j}\cdot)\mathcal{F}u\|_{L^p},\endaligned$$
where $$(\phi_1(2^{-j}\cdot)+\phi_2(2^{-j}\cdot)+\phi_3(2^{-j}\cdot))\chi_j=\chi_j,\quad j>0$$ with $\phi_i$ smooth and $\textrm{supp}\phi_i\subset \{\frac{1}{4}\leq |\xi^i|\leq \frac{3}{2},|\xi^l|\leq \frac{3}{2},l\neq i\}$. Moreover, we have for $0<j\leq N, i=1,2,3,$ that $$\chi_j(\xi)\phi_i(2^{-j}\xi)\bar{\varphi}^\varepsilon(\xi)=\chi_j(\xi)\phi_i(2^{-j}\xi)\tilde{\varphi}^\varepsilon(\xi) 1_{\{|\xi^l|\leq \varepsilon^{-1}L_0,l=1,2,3\}},$$
with $\tilde{\varphi}^\varepsilon(\xi)=e^{-t|\xi|^2\tilde{f}(\varepsilon\xi)}$,  which implies that
$$\aligned &\|\mathcal{F}^{-1}\chi_j(\bar{\varphi}^\varepsilon-\varphi)\phi_i(2^{-j}\cdot)\mathcal{F}u\|_{L^p}\\\lesssim& \|\mathcal{F}^{-1}(\tilde{\varphi}^\varepsilon-\varphi)\phi_i(2^{-j}\cdot)\mathcal{F}\mathcal{F}^{-1}\chi_j\mathcal{F}u\|_{L^p}
+\|\mathcal{F}^{-1}\varphi1_{\{|\cdot^l|> \varepsilon^{-1}L_0,\textrm{ for some } l\}}\phi_i(2^{-j}\cdot)\mathcal{F}\mathcal{F}^{-1}\chi_j\mathcal{F}u\|_{L^p}\\\lesssim & \|\mathcal{F}^{-1}(\tilde{\varphi}^\varepsilon-\varphi)\phi_i(2^{-j}\cdot)\|_{L^1}\|\mathcal{F}^{-1}\chi_j\mathcal{F}u\|_{L^p}+1_{j=N}\|\mathcal{F}^{-1}\varphi\phi_i(2^{-N}\cdot)\mathcal{F}\mathcal{F}^{-1}\chi_N\mathcal{F}u\|_{L^p}\\\lesssim & \|\mathcal{F}^{-1}(\tilde{\varphi}^\varepsilon-\varphi)\phi_i(2^{-j}\cdot)\|_{L^1}\|\mathcal{F}^{-1}\chi_j\mathcal{F}u\|_{L^p}+1_{j=N}
\varepsilon^\eta2^{N\eta}(2^N\sqrt{t})^{-\delta}\|\mathcal{F}^{-1}\chi_N\mathcal{F}u\|_{L^p}.\endaligned$$
Here in the last inequality we used $|\varepsilon 2^N|\gtrsim 1$.
By calculation we get that $$\aligned&\|\mathcal{F}^{-1}(\tilde{\varphi}^\varepsilon-\varphi)\phi_i(2^{-j}\cdot)\|_{L^1}=\|\mathcal{F}^{-1}(\tilde{\varphi}^\varepsilon-\varphi)(2^{j}\cdot)\phi_i\|_{L^1}
\lesssim \|(1-\Delta)^2((\tilde{\varphi}^\varepsilon-\varphi)(2^j\cdot)\phi_i(\cdot))\|_{L^1}\\\lesssim& \sum_{0\leq|k|\leq 4}\|2^{j|k|}(D_k\tilde{\varphi}^\varepsilon-D_k\varphi)(2^j\cdot)|_{\xi\in\textrm{supp} \phi_i}\|_{L^\infty}\lesssim \sum_{0\leq|k|\leq 4}2^{j|k|}\frac{\varepsilon^\eta2^{j\eta}}{2^{j|k|}(2^j\sqrt{t})^{\delta}}\lesssim \varepsilon^\eta2^{j\eta}(2^j\sqrt{t})^{-\delta}.\endaligned$$
Here in the third inequality we used $\inf_{\xi\in\textrm{supp}\phi_i}\tilde{f}(\varepsilon2^j\xi)\geq c>0$  and $
|D_k\tilde{f}(\xi)|\lesssim \frac{1}{|\xi|^{|k|-1}}+C$ for any multiindices $k$ satisfying $|k|\leq4$ to deduce that for every $\delta\geq0,0<\eta<1$ $|D_k(\tilde{\varphi}^\varepsilon-\varphi)(\xi)|\leq \frac{(\varepsilon|\xi|)^\eta}{|\xi|^{|k|+\delta}t^{\delta/2}}$. Now we deduce that for $0<j\leq N$
$$\|\mathcal{F}^{-1}\chi_j\mathcal{F}(\bar{P}^\varepsilon_t-P_t)u\|_{L^p}\lesssim\varepsilon^\eta2^{j\eta}(2^j\sqrt{t})^{-\delta}
\|\mathcal{F}^{-1}\chi_j\mathcal{F}u\|_{L^p}.$$

For $j>N$ we have that
$$\|\mathcal{F}^{-1}\chi_j\mathcal{F}(\bar{P}^\varepsilon_t-P_t)u\|_{L^p}=\|\mathcal{F}^{-1}\chi_j\varphi\mathcal{F}u\|_{L^p}\lesssim \varepsilon^\eta2^{j\eta}(2^j\sqrt{t})^{-\delta}\|\mathcal{F}^{-1}\chi_j\mathcal{F}u\|_{L^p},$$
where we used $|\varepsilon 2^j|\gtrsim1$.

For $j=0$ we obtain that $$\aligned &\|\mathcal{F}^{-1}\chi_0(\bar{\varphi}^\varepsilon-\varphi)\mathcal{F}u\|_{L^p}\\=& \|\mathcal{F}^{-1}(\tilde{\varphi}^\varepsilon-\varphi)(\phi_1+\phi_2+\phi_3)\tilde{\chi}_0\mathcal{F}\mathcal{F}^{-1}\chi_0\mathcal{F}u\|_{L^p}\\\lesssim &
\|\mathcal{F}^{-1}(\tilde{\varphi}^\varepsilon-\varphi)(\phi_1+\phi_2+\phi_3)\mathcal{F}\mathcal{F}^{-1}\chi_0\mathcal{F}u\|_{L^p}\\\lesssim & \sum_{i=1}^3\|\mathcal{F}^{-1}(\tilde{\varphi}^\varepsilon-\varphi)\phi_i\|_{L^1}\|\mathcal{F}^{-1}\chi_0\mathcal{F}u\|_{L^p}\\\lesssim & \varepsilon^\eta(\sqrt{t})^{-\delta}\|\mathcal{F}^{-1}\chi_0\mathcal{F}u\|_{L^p},\endaligned$$
where $\tilde{\chi}_0=\chi_0-1_{\{|\xi^i|\leq1/2, i=1,2,3\}}$ and in the last inequality we used similar argument as the case $0<j\leq N$.

By the theory in [SW71] we know that the above calculations also hold on $\mathbb{T}^3$.
Since on $\mathbb{T}^3$ operator $P^\varepsilon_t-P_t$ only depends on the value of $\varphi^\varepsilon-\varphi$ on $\mathbb{Z}^3$ which equals to $\bar{\varphi}^\varepsilon-\varphi$ on $\mathbb{Z}^3$ we obtain that on $\mathbb{T}^3$ for $p>1$ large enough $$\|(P^\varepsilon_t-P_t)u\|_{\alpha+\delta-\kappa}\lesssim\|(P^\varepsilon_t-P_t)u\|_{B^{\alpha+\delta}_{p,\infty}} \lesssim \varepsilon^\eta t^{-\delta/2}\|u\|_{B^{\alpha+\eta}_{p,\infty}} \lesssim \varepsilon^\eta t^{-\delta/2}\|u\|_{\alpha+\eta}.$$
$\hfill\Box$
\vskip.10in
The next result concerns time regularity of the approximated heat semigroup. We want to emphasize that $P_t^\varepsilon$ is not strongly continuous at $0$. However, under our assumptions  $P_t^\varepsilon$ has nice time continuity properties for $t>0$.
\vskip.10in

\th{Lemma 3.4} Let $u\in \mathcal{C}^{\alpha+\delta}$ for some $\alpha\in \mathbb{R},\delta>0$. Then for every $\kappa>0,\varepsilon\in (0,1), t>s>0$
$$\|(P^\varepsilon_t-P^\varepsilon_s)u\|_{\alpha-\kappa}\lesssim  (t-s)^{\delta/2}\|u\|_{\alpha+\delta}.$$

\proof Without loss of generality we assume that $L_0>1$. First we consider the operator $\bar{P}_t^\varepsilon-\bar{P}_s^\varepsilon=\mathcal{F}^{-1}(\bar{\varphi}_t^\varepsilon-\bar{\varphi}_s^\varepsilon)\mathcal{F}$ on $\mathbb{R}^3$ with
$$\bar{\varphi}^\varepsilon_t(\xi)=\left\{\begin{array}{ll}\varphi_t(\xi)&\ \ \ \ \textrm{ for } |\xi_i|\leq 1/2,\\\varphi^\varepsilon_t(\xi)&\ \ \ \ \textrm{ otherwise}.\end{array}\right.$$Here $$\varphi_t^\varepsilon(\xi)=e^{-t|\xi|^2f(\varepsilon\xi)},\qquad\varphi_t(\xi)=e^{-t|\xi|^2}.$$
For $N\in\mathbb{N}$ and  $2^{-N}L_0\leq\varepsilon<2^{-N+1}L_0,0<j\leq N$ we obtain that
$$\aligned &\|\mathcal{F}^{-1}\chi_j\mathcal{F}(\bar{P}^\varepsilon_t-\bar{P}_s^\varepsilon)u\|_{L^p}=\|\mathcal{F}^{-1}\chi_j(\bar{\varphi}^\varepsilon_t-\bar{\varphi}^\varepsilon_s)\mathcal{F}u\|_{L^p}\\\leq& \|\mathcal{F}^{-1}\chi_j(\bar{\varphi}^\varepsilon_t-\bar{\varphi}^\varepsilon_s)\phi_1(2^{-j}\cdot)\mathcal{F}u\|_{L^p}+\|\mathcal{F}^{-1}\chi_j(\bar{\varphi}^\varepsilon_t-\bar{\varphi}^\varepsilon_s)\phi_2(2^{-j}\cdot)\mathcal{F}u\|_{L^p}+\|\mathcal{F}^{-1}\chi_j(\bar{\varphi}^\varepsilon_t-\bar{\varphi}^\varepsilon_s)\phi_3(2^{-j}\cdot)\mathcal{F}u\|_{L^p}.\endaligned$$
where
$$(\phi_1(2^{-j}\cdot)+\phi_2(2^{-j}\cdot)+\phi_3(2^{-j}\cdot))\chi_j=\chi_j$$ with $\phi_i$ smooth and $\textrm{supp}\phi_i\subset \{\frac{1}{4}\leq |\xi^i|\leq \frac{3}{2},|\xi^l|\leq \frac{3}{2},l\neq i\}$. Moreover, we have that for $0<j\leq N$ $$\chi_j(\xi)\phi_i(2^{-j}\xi)\bar{\varphi}^\varepsilon_t(\xi)=\chi_j(\xi)\phi_i(2^{-j}\xi)\tilde{\varphi}^\varepsilon_t(\xi)1_{\{|\xi^l|\leq \varepsilon^{-1}L_0,l=1,2,3\}},$$
with $\tilde{\varphi}^\varepsilon_t(\xi)=e^{-t|\xi|^2\tilde{f}(\varepsilon\xi)}$, which implies that
$$\aligned&\|\mathcal{F}^{-1}\chi_j(\bar{\varphi}^\varepsilon_t-\bar{\varphi}^\varepsilon_s)\phi_i(2^{-j}\cdot)\mathcal{F}u\|_{L^p}\\\lesssim& \|\mathcal{F}^{-1}(\tilde{\varphi}^\varepsilon_t-\tilde{\varphi}^\varepsilon_s)\phi_i(2^{-j}\cdot)\mathcal{F}\mathcal{F}^{-1}\chi_j\mathcal{F}u\|_{L^p}\\\leq& \|\mathcal{F}^{-1}(\tilde{\varphi}^\varepsilon_t-\tilde{\varphi}^\varepsilon_s)\phi_i(2^{-j}\cdot)\|_{L^1}\|\mathcal{F}^{-1}\chi_j\mathcal{F}u\|_{L^p}.\endaligned$$
By calculation we get that $$\aligned&\|\mathcal{F}^{-1}(\tilde{\varphi}^\varepsilon_t-\tilde{\varphi}^\varepsilon_s)\phi_i(2^{-j}\cdot)\|_{L^1}=\|\mathcal{F}^{-1}(\tilde{\varphi}^\varepsilon_t-\tilde{\varphi}^\varepsilon_s)(2^{j}\cdot)\phi_i\|_{L^1}
\lesssim \|(1-\Delta)^2((\tilde{\varphi}^\varepsilon_t-\tilde{\varphi}^\varepsilon_s)(2^j\cdot)\phi_i(\cdot))\|_{L^1}\\\lesssim& \sum_{0\leq|k|\leq 4}2^{j|k|}\|(D_k\tilde{\varphi}^\varepsilon_t-D_k\tilde{\varphi}^\varepsilon_s)(2^j\cdot)|_{\cdot\in\textrm{supp} \phi_i}\|_{L^\infty}\lesssim \sum_{0\leq|k|\leq 4}2^{j|k|}\frac{(t-s)^{\delta/2}2^{j\delta}}{2^{j|k|}}\lesssim (t-s)^{\delta/2}2^{j\delta}.\endaligned$$
Here in the third inequality we used $\inf_{\xi\in\textrm{supp}\phi_i}\tilde{f}(\varepsilon2^j\xi)\geq c>0$, $|1-e^{-(t-s)\tilde{f}(\varepsilon \xi)|\xi|^2}|\leq (t-s)^{\delta/2}|\xi|^\delta$  and $|D_k\tilde{f}(\xi)|\lesssim \frac{1}{|\xi|^{|k|-1}}+C$ for
for any multiindices $k$ satisfying $|k|\leq4$ to deduce that for any $\delta\geq0$ $|D_k(\tilde{\varphi}^\varepsilon_t-\tilde{\varphi}_s^\varepsilon)(\xi)|\lesssim \frac{(t-s)^{\delta/2}|\xi|^\delta}{|\xi|^{|k|}}$.
For $j>N$ we have that
$$\|\mathcal{F}^{-1}\chi_j\mathcal{F}(P^\varepsilon_t-P_s^\varepsilon)u\|_{L^p}=0.$$
For $j=0$ we obtain that $$\aligned &\|\mathcal{F}^{-1}\chi_0(\bar{\varphi}_t^\varepsilon-\bar{\varphi}_s^\varepsilon)\mathcal{F}u\|_{L^p}\\\leq& \|\mathcal{F}^{-1}(\varphi_t-\varphi_s)1_{\{|\xi_i|\leq1/2, i=1,2,3\}}\bar{\chi}_0\mathcal{F}\mathcal{F}^{-1}\chi_0\mathcal{F}u\|_{L^p}+\|\mathcal{F}^{-1}(\tilde{\varphi}_t^\varepsilon-\tilde{\varphi}_s^\varepsilon)(\phi_1+\phi_2+\phi_3)\tilde{\chi}_0\mathcal{F}\mathcal{F}^{-1}\chi_0\mathcal{F}u\|_{L^p}\\\lesssim &
\|\mathcal{F}^{-1}(\varphi_t-\varphi_s)\bar{\chi}_0\mathcal{F}\mathcal{F}^{-1}\chi_0\mathcal{F}u\|_{L^p}+\|\mathcal{F}^{-1}(\tilde{\varphi}^\varepsilon_t-\tilde{\varphi}_s^\varepsilon)(\phi_1+\phi_2+\phi_3)\mathcal{F}\mathcal{F}^{-1}\chi_0\mathcal{F}u\|_{L^p}\\\lesssim & (t-s)^{\delta/2}\|\mathcal{F}^{-1}\chi_0\mathcal{F}u\|_{L^p}+\sum_{i=1}^3\|\mathcal{F}^{-1}(\tilde{\varphi}_t^\varepsilon-\tilde{\varphi}_s^\varepsilon)\phi_i\|_{L^1}\|\mathcal{F}^{-1}\chi_0\mathcal{F}u\|_{L^p}\\\lesssim &  (t-s)^{\delta/2}\|\mathcal{F}^{-1}\chi_0\mathcal{F}u\|_{L^p},\endaligned$$
where $\tilde{\chi}_0=\chi_0-1_{\{|\xi^i|\leq1/2, i=1,2,3\}}$ and $\bar{\chi}_0$ is smooth with supp$\bar{\chi}_0\subset \{|\xi|\leq1\}$ and $1_{\{|\xi^i|\leq1/2, i=1,2,3\}}\bar{\chi}_0=1_{\{|\xi^i|\leq1/2, i=1,2,3\}}$. Here in the second inequality we used $1_{\{|\xi^i|\leq1/2, i=1,2,3\}}$ is an $L^p$ multiplier and in the last inequality we used similar argument as the case $0<j\leq N$.

By the theory in [SW71] we know that the above calculations also hold on $\mathbb{T}^3$. Since on $\mathbb{T}^3$ the operator $P^\varepsilon_t-P_s^\varepsilon$ only depends on the value of $\varphi^\varepsilon_t-\varphi_s^\varepsilon$ on $\mathbb{Z}^3$,  which equals to $\bar{\varphi}^\varepsilon_t-\bar{\varphi}_s^\varepsilon$ on $\mathbb{Z}^3$, we obtain that on $\mathbb{T}^3$ for $p$ large enough by Lemma 2.1
$$\|(P^\varepsilon_t-P_s^\varepsilon)u\|_{\alpha-\kappa}\lesssim\|(P^\varepsilon_t-P_s^\varepsilon)u\|_{B^{\alpha}_{p,\infty}} \lesssim  (t-s)^{\delta/2}\|u\|_{B^{\alpha+\delta}_{p,\infty}} \lesssim (t-s)^{\delta/2}\|u\|_{\alpha+\delta}.$$
$\hfill\Box$
\vskip.10in
The following two lemmas state the commutator estimate for the approximated semigroup in a particular case, which can be used in the case that we approximate  $W$ by only
keeping finite dimensional Fourier modes.
\vskip.10in
\th{Lemma 3.5} Let $u\in \mathcal{C}^\alpha$ for some $\alpha<1$ and $v\in \mathcal{C}^\beta$ for some $\beta\in \mathbb{R}$ satisfying $\textrm{supp} \mathcal{F}v\subset \{\xi:|\xi|\leq L_0/(2\varepsilon)\}$. Then for $\delta\geq \alpha+\beta$
$$\|P_t^\varepsilon\pi_{<}(u,v)-\pi_<(u,P_t^\varepsilon v)\|_{\delta}\lesssim t^{\frac{\alpha+\beta-\delta}{2}}\|u\|_\alpha\|v\|_\beta.$$

\proof Since $\textrm{supp} \mathcal{F}v\subset \{\xi:|\xi|\leq L_0/(2\varepsilon)\}$, by the definition of  $\Delta_j$ we have $\Delta_jv\neq0$ only for $j\leq [\log_2 (L_0/\varepsilon)]$,
which implies that
$$P_t^\varepsilon\pi_{<}(u,v)-\pi_<(u,P_t^\varepsilon v)=\sum_{j=-1}^{[\log_2 (L_0/\varepsilon)]}(P_t^\varepsilon(S_{j-1}u\Delta_jv)-S_{j-1}uP_t^\varepsilon \Delta_jv).$$
We also have that the Fourier transform of $P_t^\varepsilon(S_{j-1}u\Delta_jv)-S_{j-1}uP_t^\varepsilon \Delta_jv$ has its support in a suitable annulus $2^j\mathcal{A}$. Let $\psi\in \mathcal{D}(\mathbb{R}^3)$ with support in an annulus $\tilde{\mathcal{A}}$ be such that $\psi=1$ on $\mathcal{A}$.

Thus we obtain that for $j\leq [\log_2 (L_0/\varepsilon)]$ $\textrm{supp} (\mathcal{F}S_{j-1}u\Delta_jv) \subset \{|\xi^i|\leq L_0/\varepsilon\}$ and $\textrm{supp}(\mathcal{F}\Delta_jv) \subset \{|\xi^i|\leq L_0/\varepsilon\}$, which imply that
$$\aligned &\|[(\psi(2^{-j}\cdot)\varphi^\varepsilon)(D),S_{j-1}u]\Delta_jv\|_{L^\infty}=\|[(\psi(2^{-j}\cdot)\tilde{\varphi}^\varepsilon)(D),S_{j-1}u]\Delta_jv\|_{L^\infty}
\\\lesssim&\sum_{\eta\in \mathbb{N}^d,|\eta|=1}\|x^\eta\mathcal{F}^{-1}(\psi(2^{-j}\cdot)\tilde{\varphi}^\varepsilon)\|_{L^1}\|\partial^\eta S_{j-1}u\|_{L^\infty}\|\Delta_jv\|_{L^\infty},\endaligned$$
 where in the last inequality we used the same argument as the proof of [CC13, Lemma A.1].

 Here $$\varphi^\varepsilon(\xi)=e^{-t|\xi|^2f(\varepsilon\xi)},\quad \tilde{\varphi}^\varepsilon(\xi)=e^{-t|\xi|^2\tilde{f}(\varepsilon\xi)},$$ $(\psi(2^{-j}\cdot)\varphi^\varepsilon)(D)u=\mathcal{F}^{-1}(\psi(2^{-j}\cdot)\varphi^\varepsilon\mathcal{F}u)$, $[(\psi(2^{-j}\cdot)\varphi^\varepsilon)(D),S_{j-1}u]$ denotes the commutator.

Now we have that

$$\aligned &\|x^\eta\mathcal{F}^{-1}(\psi(2^{-j}\cdot)\tilde{\varphi^\varepsilon})\|_{L^1}\\\leq& 2^{-j}\|\mathcal{F}^{-1}(\partial^\eta\psi)(2^{-j}\cdot)\tilde{\varphi^\varepsilon})\|_{L^1}+\|\mathcal{F}^{-1}(\psi(2^{-j}\cdot)\partial^\eta\tilde{\varphi^\varepsilon})\|_{L^1}
\\=&2^{-j}\|\mathcal{F}^{-1}(\partial^\eta\psi(\cdot)\tilde{\varphi^\varepsilon}(2^j\cdot))\|_{L^1}+\|\mathcal{F}^{-1}(\psi(\cdot)\partial^\eta\tilde{\varphi^\varepsilon}(2^j\cdot))\|_{L^1}
\\\lesssim&2^{-j}\|(1+|\cdot|^2)^{2}\mathcal{F}^{-1}(\partial^\eta\psi(\cdot)\tilde{\varphi^\varepsilon}(2^j\cdot))\|_{L^\infty}+ \|(1+|\cdot|^2)^{2}\mathcal{F}^{-1}(\psi(\cdot)\partial^\eta\tilde{\varphi^\varepsilon}( 2^j\cdot))\|_{L^\infty} \\=&2^{-j}\|\mathcal{F}^{-1}((1-\Delta)^{2}(\partial^\eta\psi(\cdot)\tilde{\varphi^\varepsilon}( 2^j\cdot)))\|_{L^\infty}+ \|\mathcal{F}^{-1}((1-\Delta)^{2}(\psi(\cdot)\partial^\eta\tilde{\varphi^\varepsilon}( 2^j\cdot)))\|_{L^\infty}
\\\lesssim&2^{-j}\|(1-\Delta)^{2}(\partial^\eta\psi(\cdot)\tilde{\varphi^\varepsilon}( 2^j\cdot))\|_{L^1}+ \|(1-\Delta)^{2}(\psi(\cdot)\partial^\eta\tilde{\varphi^\varepsilon}( 2^j\cdot))\|_{L^1}
 \\\lesssim&2^{-j}\sum_{0\leq|m|\leq 4}( 2^j)^{|m|}\frac{t^{-\mu}2^{-2j\mu}}{( 2^j)^{|m|}}+ \sum_{|m|\leq 5}( 2^j)^{|m|}\frac{t^{-\mu}2^{-2j\mu}}{( 2^j)^{|m|+1}}\\\lesssim&2^{-j}t^{-\mu}2^{-2j\mu},\endaligned$$
 where  in the fourth inequality we used $|D^m\tilde{\varphi^\varepsilon}(\xi)|\lesssim |\xi|^{-|m|}t^{-\mu}|\xi|^{-2\mu}, \mu\geq0$ for any multiindices $m$ satisfying $|m|\leq5$.
Hence we get that
$$\|[\psi(2^{-j}\cdot)\varphi^\varepsilon(D),S_{j-1}u]\Delta_jv\|_{L^\infty}\lesssim t^{\frac{\alpha+\beta-\delta}{2}}2^{j(\alpha+\beta-\delta)}2^{-j(\alpha+\beta)}\|u\|_\alpha\|v\|_\beta,$$which yields the result by the same argument as in the proof of [CC13, Lemma A.1].$\hfill\Box$
\vskip.10in
\th{Lemma 3.6} Let $u\in \mathcal{C}^\alpha$ for some $\alpha<1$ and $v\in \mathcal{C}^\beta$ for some $\beta\in \mathbb{R}$ satisfying $\textrm{supp} \mathcal{F}v\subset \{|\xi|\leq L_0/(2\varepsilon)\}$. Then for $\delta\geq \alpha+\beta, \kappa>0$
$$\|(P_t^\varepsilon-P_t)\pi_{<}(u,v)-\pi_<(u,(P_t^\varepsilon-P_t) v)\|_{\delta-\kappa}\lesssim \varepsilon^\kappa t^{\frac{\alpha+\beta-\delta}{2}}\|u\|_\alpha\|v\|_\beta.$$

\proof Since $\textrm{supp} \mathcal{F}v\subset \{\xi:|\xi|\leq L_0/(2\varepsilon)\}$, by the definition of  $\Delta_j$ we have $\Delta_jv\neq0$ only for $j\leq [\log_2 (L_0/\varepsilon)]$, which implies that
$$(P_t^\varepsilon-P_t)\pi_{<}(u,v)-\pi_<(u,(P_t^\varepsilon-P_t) v)=\sum_{j=-1}^{[\log_2 (L_0/\varepsilon)]}((P_t^\varepsilon-P_t)(S_{j-1}u\Delta_jv)-S_{j-1}u(P_t^\varepsilon-P_t) \Delta_jv).$$
We also have  the Fourier transform of $(P_t^\varepsilon-P_t)(S_{j-1}u\Delta_jv)-S_{j-1}u(P_t^\varepsilon-P_t) \Delta_jv$ has its support in a suitable annulus $2^j\mathcal{A}$. Let $\psi\in \mathcal{D}(\mathbb{R}^3)$ with support in an annulus $\tilde{\mathcal{A}}$ be such that $\psi=1$ on $\mathcal{A}$.
Thus we obtain that for $j\leq [\log_2 (L_0/\varepsilon)]$ $\textrm{supp} (\mathcal{F}S_{j-1}u\Delta_jv) \subset \{|\xi^i|\leq L_0/\varepsilon\}$ and $\textrm{supp} (\mathcal{F}\Delta_jv) \subset \{|\xi^i|\leq L_0/\varepsilon\}$, which imply that
$$\aligned &\|[(\psi(2^{-j}\cdot)(\varphi^\varepsilon-\varphi))(D),S_{j-1}u]\Delta_jv\|_{L^\infty}=\|[(\psi(2^{-j}\cdot)(\tilde{\varphi}^\varepsilon-\varphi))(D),S_{j-1}u]\Delta_jv\|_{L^\infty}
\\\lesssim&\sum_{\eta\in \mathbb{N}^d,|\eta|=1}\|x^\eta\mathcal{F}^{-1}(\psi(2^{-j}\cdot)(\tilde{\varphi}^\varepsilon-\varphi))\|_{L^1}\|\partial^\eta S_{j-1}u\|_{L^\infty}\|\Delta_jv\|_{L^\infty},\endaligned$$
where in the last inequality we used the same argument as the proof of [CC13, Lemma A.1].
 Here $$\varphi^\varepsilon(\xi)=e^{-t|\xi|^2f(\varepsilon\xi)},\quad \tilde{\varphi}^\varepsilon(\xi)=e^{-t|\xi|^2\tilde{f}(\varepsilon\xi)},\quad \varphi(\xi)=e^{-t|\xi|^2},$$ $(\psi(2^{-j}\cdot)\varphi^\varepsilon)(D)u=\mathcal{F}^{-1}(\psi(2^{-j}\cdot)\varphi^\varepsilon\mathcal{F}u)$, $[(\psi(2^{-j}\cdot)\varphi^\varepsilon)(D),S_{j-1}u]$ denotes the commutator.

Now we have that
$$\aligned &\|x^\eta\mathcal{F}^{-1}(\psi(2^{-j}\cdot)(\tilde{\varphi^\varepsilon}-\varphi))\|_{L^1}\\\leq& 2^{-j}\|\mathcal{F}^{-1}(\partial^\eta\psi)(2^{-j}\cdot)(\tilde{\varphi^\varepsilon}-\varphi))\|_{L^1}+\|\mathcal{F}^{-1}(\psi(2^{-j}\cdot)\partial^\eta(\tilde{\varphi^\varepsilon}-\varphi))\|_{L^1}
\\=&2^{-j}\|\mathcal{F}^{-1}(\partial^\eta\psi(\cdot)(\tilde{\varphi^\varepsilon}-\varphi)(2^j\cdot))\|_{L^1}+\|\mathcal{F}^{-1}(\psi(\cdot)\partial^\eta(\tilde{\varphi^\varepsilon}-\varphi)(2^j\cdot))\|_{L^1}
\\\lesssim&2^{-j}\|(1+|\cdot|^2)^{2}\mathcal{F}^{-1}(\partial^\eta\psi(\cdot)(\tilde{\varphi^\varepsilon}-\varphi)(2^j\cdot))\|_{L^\infty}+ \|(1+|\cdot|^2)^{2}\mathcal{F}^{-1}(\psi(\cdot)\partial^\eta(\tilde{\varphi^\varepsilon}-\varphi)( 2^j\cdot))\|_{L^\infty} \\=&2^{-j}\|\mathcal{F}^{-1}((1-\Delta)^{2}(\partial^\eta\psi(\cdot)(\tilde{\varphi^\varepsilon}-\varphi)( 2^j\cdot)))\|_{L^\infty}+ \|\mathcal{F}^{-1}((1-\Delta)^{2}(\psi(\cdot)\partial^\eta(\tilde{\varphi^\varepsilon}-\varphi)( 2^j\cdot)))\|_{L^\infty}
\\\lesssim&2^{-j}\|(1-\Delta)^{2}(\partial^\eta\psi(\cdot)(\tilde{\varphi^\varepsilon}-\varphi)( 2^j\cdot))\|_{L^1}+ \|(1-\Delta)^{2}(\psi(\cdot)\partial^\eta(\tilde{\varphi^\varepsilon}-\varphi)( 2^j\cdot))\|_{L^1}
 \\\lesssim&2^{-j}\sum_{0\leq|m|\leq 4}( 2^j)^{|m|}\frac{\varepsilon^\kappa t^{-\mu}2^{-2j\mu+j\kappa}}{( 2^j)^{|m|}}+ \sum_{|m|\leq 5}( 2^j)^{|m|}\frac{\varepsilon^\kappa t^{-\mu}2^{-2j\mu+j\kappa}}{( 2^j)^{|m|+1}}\\\lesssim&2^{-j}\varepsilon^\kappa t^{-\mu}2^{-2j\mu+j\kappa},\endaligned$$
 where  in the fourth inequality we used $|D^m(\tilde{\varphi^\varepsilon}-\varphi)(\xi)|\lesssim \varepsilon^\kappa|\xi|^{-|m|}t^{-\mu}|\xi|^{\kappa-2\mu}, \mu,\kappa\geq0$ for any multiindices $m$ satisfying $|m|\leq5$.
Thus we get that
$$\|[\psi(2^{-j}\cdot)\varphi^\varepsilon(D),S_{j-1}u]\Delta_jv\|_{L^\infty}\lesssim \varepsilon^\kappa t^{\frac{\alpha+\beta-\delta}{2}}2^{j(\kappa+\alpha+\beta-\delta)}2^{-j(\alpha+\beta)}\|u\|_\alpha\|v\|_\beta,$$which implies the result by the same argument as in the proof of [CC13, Lemma A.1].$\hfill\Box$
\vskip.10in
The following two lemmas concern estimates for the operator $D^\varepsilon_j$.
\vskip.10in
\th{Lemma 3.7} Let $u\in \mathcal{C}^{\alpha+1}$ for some $\alpha\in \mathbb{R}$. Then  for every $\kappa>0$
$$\sup_{\varepsilon\in(0,1)}\|D^\varepsilon_ju\|_{\alpha-\kappa}\lesssim\|u\|_{\alpha+1}.$$

\proof First we consider the operator on $\mathbb{R}^3$. For $j>-1, p>1$ we have
$$\aligned &\|\mathcal{F}^{-1}\theta(2^{-j}\cdot)\mathcal{F}D^\varepsilon_iu\|_{L^p}=\|\mathcal{F}^{-1}\theta(2^{-j}\xi)
\xi^ig(\varepsilon \xi^i)\mathcal{F}u(\xi)\|_{L^p}.\endaligned$$
It is easy to get that $|g(x)|\leq C, |g'(x)|\leq C/|x|$, which by the Mihlin multiplier theorem implies that $g$ is an $L^p(\mathbb{R}^1)$ multiplier, hence an $L^p(\mathbb{R}^3)$ multiplier.
Therefore we have
 $$\aligned&\|\mathcal{F}^{-1}\theta(2^{-j}\xi)\xi^ig(\varepsilon \xi^i)\mathcal{F}u(\xi)\|_{L^p}\lesssim\|\mathcal{F}^{-1}\theta(2^{-j}\xi)\xi^i\mathcal{F}u(\xi)\|_{L^p}=\|\mathcal{F}^{-1}\psi(2^{-j}\xi)\theta(2^{-j}\xi)\xi^i\mathcal{F}u(\xi)\|_{L^p}\\\lesssim&\|\mathcal{F}^{-1}\psi(2^{-j}\xi)\xi^i\|_{L^1}\|\mathcal{F}^{-1}\theta(2^{-j}\cdot)\mathcal{F}u\|_{L^p} \lesssim2^j\|\mathcal{F}^{-1}\theta(2^{-j}\cdot)\mathcal{F}u\|_{L^p} ,\endaligned$$
 where $\psi$ is  smooth with the support of $\psi$ contained in an annulus and $\psi\theta=\theta.$ Moreover, we obtain similar estimate for $j=-1$.  By the theory in [SW71] we know that the above calculations also hold on $\mathbb{T}^3$.
Thus on $\mathbb{T}^3$ we deduce that for $p$ large enough$$\sup_{\varepsilon\in(0,1)}\|D^\varepsilon_iu\|_{\alpha-\kappa}\lesssim\sup_{\varepsilon\in(0,1)}\|D^\varepsilon_iu\|_{B^{\alpha}_{p,\infty}} \lesssim \|u\|_{B^{\alpha+1}_{p,\infty}} \lesssim\|u\|_{\alpha+1}.$$
$\hfill\Box$
\vskip.10in
\th{Lemma 3.8} Let $u\in \mathcal{C}^{\alpha+1+\eta}$ for some $\alpha\in \mathbb{R},\eta>0$. Then  for every $\kappa>0$
$$\|(D^\varepsilon_i-D_i)u\|_{\alpha-\kappa}\lesssim \varepsilon^\eta\|u\|_{\alpha+1+\eta}.$$

\proof First we consider the operator on $\mathbb{R}^3$. Suppose that $2^{-N}\leq \varepsilon <2^{-N+1}$, $N\in \mathbb{N}$ then for $0\leq j\leq N, p>1$ we have
$$\aligned &\|\mathcal{F}^{-1}\theta(2^{-j}\cdot)\mathcal{F}(D^\varepsilon_i-D_i)u\|_{L^p}=\|\mathcal{F}^{-1}\theta(2^{-j}\xi)\varepsilon (\xi^i)^2\frac{g(\varepsilon \xi^i)-\imath}{\varepsilon \xi^i}\mathcal{F}u(\xi)\|_{L^p}.\endaligned$$
Let $M(x)=\frac{g(x)-\imath }{x}$. It is easy to get that $|M(x)|\leq C, |M'(x)|\leq C/|x|$, which by Mihlin multiplier theorem implies that $M$ is an $L^p(\mathbb{R}^1)$ multiplier, hence an $L^p(\mathbb{R}^3)$ multiplier.
Therefore
 $$\aligned&\|\mathcal{F}^{-1}\theta(2^{-j}\xi)\varepsilon (\xi^i)^2M(\varepsilon \xi^i)\mathcal{F}u(\xi)\|_{L^p}\lesssim\|\mathcal{F}^{-1}\theta(2^{-j}\xi)\varepsilon (\xi^i)^2\mathcal{F}u(\xi)\|_{L^p}\\=&\|\mathcal{F}^{-1}\psi(2^{-j}\xi)\theta(2^{-j}\xi)\varepsilon (\xi^i)^2\mathcal{F}u(\xi)\|_{L^p}\lesssim\|\mathcal{F}^{-1}\psi(2^{-j}\xi)\varepsilon (\xi^i)^2\|_{L^1}\|\mathcal{F}^{-1}\theta(2^{-j}\cdot)\mathcal{F}u\|_{L^p} \\\lesssim& \varepsilon^\eta2^{j(1+\eta)}\|\mathcal{F}^{-1}\theta(2^{-j}\cdot)\mathcal{F}u\|_{L^p}, \endaligned$$
where $\psi$ is  smooth with the support of $\psi$ contained in an annulus and $\psi\theta=\theta$ and we used $\varepsilon 2^j\lesssim1$ in the last inequality.

For $j>N$ we have
$$\aligned &\|\mathcal{F}^{-1}\theta(2^{-j}\cdot)\mathcal{F}(D^\varepsilon_i-D_i)u\|_{L^p}\\\leq& \|\mathcal{F}^{-1}\theta(2^{-j}\cdot)\mathcal{F}D^\varepsilon_iu\|_{L^p}+\|\mathcal{F}^{-1}\theta(2^{-j}\cdot)\mathcal{F}D_iu\|_{L^p}\\\lesssim& 2^j\|\mathcal{F}^{-1}\theta(2^{-j}\cdot)\mathcal{F}u\|_{L^p}\\\lesssim& \varepsilon^\eta 2^{j(1+\eta)}\|\mathcal{F}^{-1}\theta(2^{-j}\cdot)\mathcal{F}u\|_{L^p},\endaligned$$
where in the second inequality we used a similar estimate as Lemma 3.7 and $\varepsilon 2^j\geq1$ in the last inequality.
Similarly, we obtain the desired estimate for $j=-1$. By the theory in [SW71] we know that the above calculations also hold on $\mathbb{T}^3$. 
Thus on $\mathbb{T}^3$ for $p$ large enough.
$$\|(D^\varepsilon_i-D_i)u\|_{\alpha-\kappa}\lesssim\|(D^\varepsilon_i-D_i)u\|_{B^{\alpha}_{p,\infty}} \lesssim \varepsilon^\eta\|u\|_{B^{\alpha+1+\eta}_{p,\infty}} \lesssim\varepsilon^\eta\|u\|_{\alpha+1+\eta}.$$
$\hfill\Box$

\vskip.10in

\subsection{ Construction of solutions to 3D NS equation driven by space-time white noise and approximating equations}

 In this subsection we give an outline of the construction of local solutions to (1.1) and (1.3). The construction given here differs slightly from the construction presented in [ZZ14] (see Remark 3.11 in the following).

\textbf{Construction of solutions to 3D NS equation driven by space-time white noise}:
Now we write (1.1) in the form of (1.4) with $u_0\in \mathcal{C}^{-z},$
where $z\in (1/2,1/2+\delta_0)$ with $0<\delta_0<1/2$.  As  in [ZZ14] we give the definition of the solution of (1.4) as limit of solutions $\bar{u}^\varepsilon$ to the following equation:
$$d \bar{u}^{\varepsilon,i}=\Delta \bar{u}^{\varepsilon,i}dt+\sum_{i_1=1}^3P^{ii_1}H_\varepsilon dW^{\varepsilon,i_1}-\frac{1}{2}\sum_{i_1=1}^3P^{ii_1}(\sum_{j=1}^3D_j(\bar{u}^{\varepsilon,i_1} \bar{u}^{\varepsilon,j}))dt,\eqno(3.1)$$
$$\bar{u}^\varepsilon(0)=Pu_0\in \mathcal{C}^{-z}.$$
Now we split (3.1) into the following four equations:
$$d \bar{u}_1^{\varepsilon,i}=\Delta \bar{u}_1^{\varepsilon,i}dt+\sum_{i_1=1}^3P^{ii_1}H_\varepsilon dW,$$

$$d \bar{u}_2^{\varepsilon,i}=\Delta \bar{u}_2^{\varepsilon,i}dt-\frac{1}{2}\sum_{i_1=1}^3P^{ii_1}(\sum_{j=1}^3D_j(\bar{u}_1^{\varepsilon,i_1}\diamond \bar{u}_1^{\varepsilon,j}))dt\quad \bar{u}^\varepsilon_2(0)=0,$$

$$d\bar{u}_3^{\varepsilon,i}=\Delta \bar{u}_3^{\varepsilon,i}-\frac{1}{2}\sum_{i_1=1}^3P^{ii_1}(\sum_{j=1}^3D_j(\bar{u}_1^{\varepsilon,i_1}\diamond \bar{u}_2^{\varepsilon,j}+\bar{u}_2^{\varepsilon,i_1}\diamond \bar{u}_1^{\varepsilon,j}))\quad \bar{u}^\varepsilon_3(0)=0,$$
and

$$\aligned \bar{u}_4^{\varepsilon,i}(t)=&P_t(\sum_{i_1=1}^3P^{ii_1}{u}_0^{i_1}-\bar{u}_1^{\varepsilon,i}(0))-\frac{1}{2}\int_0^t
P_{t-s}\bigg{[}\sum_{i_1,j=1}^3P^{ii_1}D_j(\bar{u}_1^{\varepsilon,i_1}\diamond (\bar{u}_3^{\varepsilon,j}+\bar{u}_4^{\varepsilon,j})+(\bar{u}_3^{\varepsilon,i_1} +\bar{u}_4^{\varepsilon,i_1})\diamond \bar{u}_1^{\varepsilon,j}\\&+\bar{u}_2^{\varepsilon,i_1}\diamond \bar{u}_2^{\varepsilon,j}+\bar{u}_2^{\varepsilon,i_1}(\bar{u}_3^{\varepsilon,j}+\bar{u}_4^{\varepsilon,j})+\bar{u}_2^{\varepsilon,j}
(\bar{u}_3^{\varepsilon,i_1}+\bar{u}_4^{\varepsilon,i_1})+(\bar{u}_3^{\varepsilon,i_1}+\bar{u}_4^{\varepsilon,i_1})(\bar{u}_3^{\varepsilon,j}
+\bar{u}_4^{\varepsilon,j}))\bigg{]}ds.\endaligned\eqno(3.2)$$
Here for $i,j=1,2,3$, $$\bar{u}_1^{\varepsilon,i}=\int_{-\infty}^t\sum_{i_1=1}^3P^{ii_1}P_{t-s}H_\varepsilon dW^{i_1},$$
$$\bar{u}^{\varepsilon,i}_1\diamond \bar{u}^{\varepsilon,j}_1:=\bar{u}^{\varepsilon,i}_1\bar{u}^{\varepsilon,j}_1-\bar{C}^{\varepsilon,ij}_0,$$
 $$\bar{u}_2^{\varepsilon,j}\diamond \bar{u}_1^{\varepsilon,i}=\bar{u}_1^{\varepsilon,i}\diamond \bar{u}_2^{\varepsilon,j}:=\bar{u}_1^{\varepsilon,i} \bar{u}_2^{\varepsilon,j},$$
 $$\bar{u}_2^{\varepsilon,i}\diamond \bar{u}_2^{\varepsilon,j}:=\bar{u}_2^{\varepsilon,i} \bar{u}_2^{\varepsilon,j}-\bar{\varphi}_2^{\varepsilon,ij}(t)-\bar{C}_2^{\varepsilon,ij},$$
 $$\bar{u}_1^{\varepsilon,i}\diamond \bar{u}_3^{\varepsilon,j}=\bar{u}_3^{\varepsilon,j}\diamond \bar{u}_1^{\varepsilon,i}:=\pi_<(\bar{u}_3^{\varepsilon,j},\bar{u}_1^{\varepsilon,i})+\pi_{0,\diamond}(\bar{u}_3^{\varepsilon,j},\bar{u}_1^{\varepsilon,i})
 +\pi_>(\bar{u}_3^{\varepsilon,j},\bar{u}_1^{\varepsilon,i}),$$
 $$\bar{u}_1^{\varepsilon,i}\diamond \bar{u}_4^{\varepsilon,j}=\bar{u}_4^{\varepsilon,j}\diamond \bar{u}_1^{\varepsilon,i}:=\pi_<(\bar{u}_4^{\varepsilon,j},\bar{u}_1^{\varepsilon,i})+\pi_{0,\diamond}(\bar{u}_4^{\varepsilon,j},\bar{u}_1^{\varepsilon,i})
+\pi_>(\bar{u}_4^{\varepsilon,j},\bar{u}_1^{\varepsilon,i}),$$
 with $$\pi_{0,\diamond}(\bar{u}_3^{\varepsilon,j}, \bar{u}_1^{\varepsilon,i}):=\pi_{0}(\bar{u}_3^{\varepsilon,j}, \bar{u}_1^{\varepsilon,i})-\bar{\varphi}_1^{\varepsilon,ji}(t)-\bar{C}_1^{\varepsilon,ji},$$
 and $$\pi_{0,\diamond}(\bar{u}_4^{\varepsilon,i}, \bar{u}_1^{\varepsilon,j}):=\pi_{0}(\bar{u}_4^{\varepsilon,i}, \bar{u}_1^{\varepsilon,j}),$$
with $\bar{C}^\varepsilon_0$ is defined in section 4.2,
 $\bar{C}_1^\varepsilon$ and $\bar{\varphi}_1^\varepsilon$ are defined in Section 4.4 and  $\bar{C}_2^\varepsilon$ and $\bar{\varphi}_2^\varepsilon$ are defined in Section 4.6 and $\bar{\varphi}^\varepsilon_i$ converges to some $\bar{\varphi}_i$ with respect to the norm $\|\varphi\|=\sup_{t\in[0,T]}t^\rho|\varphi(t)|$ for every $\rho>0$ and $i=1,2$. Define $$d\bar{K}^{\varepsilon,i}=(\Delta\bar{K}^{\varepsilon,i}+\bar{u}_1^{\varepsilon,i})dt\quad \bar{K}^{\varepsilon,i}(0)=0.$$
In [ZZ14] we proved that for every $\delta>0$ $$\aligned \bar{C}^\varepsilon_W(T):=&\sup_{t\in[0,T]}(\sum_{i=1}^3\|\bar{u}_1^{\varepsilon,i}\|_{-1/2-\delta/2}+\sum_{i,j=1}^3\|\bar{u}_1^{\varepsilon,i}\diamond \bar{u}_1^{\varepsilon,j}\|_{-1-\delta/2}+\sum_{i,j=1}^3\|\bar{u}_1^{\varepsilon,i}\diamond \bar{u}_2^{\varepsilon,j}\|_{-1/2-\delta/2}\\&+\sum_{i,j=1}^3\|\bar{u}_2^{\varepsilon,i}\diamond \bar{u}_2^{\varepsilon,j}\|_{-\delta}+\sum_{i,j=1}^3\|\pi_{0,\diamond}( \bar{u}_3^{\varepsilon,i},\bar{u}_1^{\varepsilon,j})\|_{-\delta}+\sum_{i,i_1,j,j_1=1}^3\|\pi_{0,\diamond}(P^{ii_1}D_{j} \bar{K}^{\varepsilon,j},\bar{u}_1^{\varepsilon,j_1})\|_{-\delta}\\&+\sum_{i,i_1,j,j_1=1}^3\|\pi_{0,\diamond}
(P^{ii_1}D_{j}\bar{K}^{\varepsilon,i_1},\bar{u}_1^{\varepsilon,j_1})\|_{-\delta})<\infty,\endaligned$$
with
 $$\pi_{0,\diamond}(P^{ii_1}D_j\bar{K}^{\varepsilon,j}, \bar{u}_1^{\varepsilon,j_2}):=\pi_{0}(P^{ii_1}D_j\bar{K}^{\varepsilon,j}, \bar{u}_1^{\varepsilon,j_2}),$$
  $$\pi_{0,\diamond}(P^{ii_1}D_j\bar{K}^{\varepsilon,i_1}, \bar{u}_1^{\varepsilon,j_2}):=\pi_{0}(P^{ii_1}D_j\bar{K}^{\varepsilon,i_1}, \bar{u}_1^{\varepsilon,j_2}).$$
  In the following we will fix $\delta>0$ small enough such that $$\delta<\delta_0\wedge \frac{1}{4}(1-z)\wedge\frac{1}{3}(1-2\delta_0)\wedge(2z-1).$$

 Moreover by Lemmas 2.5 and 2.6 we easily deduce that for $i=1,2,3$,
$$\sup_{t\in[0,T]}(\sum_{i=1}^3\|\bar{u}_2^{\varepsilon,i}\|_{-\delta}+\sum_{i=1}^3\|\bar{u}_3^{\varepsilon,i}\|_{1/2-\delta})\lesssim \bar{C}_W^\varepsilon,$$
and $$\|\bar{K}^{\varepsilon,i}(t)\|_{\frac{3}{2}-\delta}\lesssim t^{\delta/4}\sup_{s\in[0,t]}\|\bar{u}_1^{\varepsilon,i}(s)\|_{-1/2-\delta/2}\lesssim t^{\delta/4}\bar{C}_W^\varepsilon.\eqno(3.3)$$
For every $\varepsilon>0$ by a similar argument as [ZZ14] we obtain  solutions of equation (3.2): More precisely, for each $\varepsilon\in(0,1)$ there exists $\bar{u}_4^\varepsilon$ satisfying equation (3.2) respectively before $T^\varepsilon>0$ such that $\bar{u}_4^\varepsilon\in C((0,T^\varepsilon);\mathcal{C}^{1/2-\delta_0})$ with respect to the norm
$\sup_{t\in [0,T]}t^{\frac{1/2-\delta_0+z+\kappa}{2}}\|\bar{u}_4(t)\|_{1/2-\delta_0}$ for $0<3\kappa<1/2-z+\delta_0$ and satisfies
$$\sup_{t\in [0,T^\varepsilon]}t^{\frac{1/2-\delta_0+z+\kappa}{2}}\|\bar{u}_4^\varepsilon(t)\|_{1/2-\delta_0}=\infty.$$
Indeed, by (3.2) and Lemma 2.2 and Lemmas 2.5, 2.6 one has the following estimates $$\sup_{t\in [0,T]}t^{\frac{1/2-\delta_0+z+\kappa}{2}}\|\bar{u}^\varepsilon_4(t)\|_{1/2-\delta_0}\lesssim C_\varepsilon(\|u_0\|_{-z},\bar{u}^\varepsilon_1,\bar{u}^\varepsilon_2,\bar{u}^\varepsilon_3)+T^{\frac{1/2+\delta_0-z-3\kappa}{2}}(\sup_{t\in [0,T]}t^{\frac{1/2-\delta_0+z+\kappa}{2}}\|\bar{u}^\varepsilon_4(t)\|_{1/2-\delta_0})^2,$$
where $C_\varepsilon(\|u_0\|_{-z},\bar{u}^\varepsilon_1,\bar{u}^\varepsilon_2,\bar{u}^\varepsilon_3)$ are constants depending on $\varepsilon$ and we used $z<1/2+\delta_0-3\kappa$.

In the following we will fix 
$$0<\kappa<\frac{\delta}{2}\wedge(\frac{1}{6}+\frac{\delta_0}{3}-\frac{z}{3})\wedge\frac{1-2\delta_0-3\delta}{7}\wedge\frac{1-z-4\delta}{5}.$$

Now consider the paracontrolled ansatz for $i=1,2,3$,
$$\bar{u}_4^{\varepsilon,i}=-\frac{1}{2}\sum_{i_1=1}^3P^{ii_1}\bigg(\sum_{j=1}^3D_j[\pi_<(\bar{u}_3^{\varepsilon,i_1}+\bar{u}_4^{\varepsilon,i_1},
\bar{K}^{\varepsilon,j})+\pi_<(\bar{u}_3^{\varepsilon,j}+\bar{u}_4^{\varepsilon,j},\bar{K}^{\varepsilon,i_1})]\bigg)+\bar{u}^{\varepsilon,\sharp,i}$$ with $\bar{u}^{\varepsilon,\sharp,i}(t)\in \mathcal{C}^{1/2+\beta}$ for some $\delta/2<\beta<(z+2\delta-1/2)\wedge (1/2-2\delta-3\kappa)\wedge (\frac{1}{2}-\delta_0-\delta-\kappa)\wedge(2-2z-\frac{5\delta}{2}-4\kappa)$ and $t\in(0,T^\varepsilon)$ (which can be done for fix $\varepsilon>0$ since the noise term is smooth and we note that
$$t^{\frac{1/2+\beta+z+\kappa}{2}}\|\bar{u}^\varepsilon_4(t)\|_{1/2+\beta}\lesssim C_\varepsilon(\|u_0\|_{-z},\bar{u}^\varepsilon_1,\bar{u}^\varepsilon_2,\bar{u}^\varepsilon_3)+t^{\frac{1/2+\delta_0-3\kappa-z}{2}}(\sup_{s\in [0,t]}s^{\frac{1/2-\delta_0+z+\kappa}{2}}\|\bar{u}^\varepsilon_4(s)\|_{1/2-\delta_0})^2).$$
By paracontrolled ansatz and Lemmas 2.2, 2.5,  2.6  we obtain the following estimate: for $i=1,2,3$,
$$\aligned\|\bar{u}_4^{\varepsilon,i}\|_{1/2-\delta-\kappa}\lesssim& \sum_{i_1,j=1}^3\|\bar{u}_3^{\varepsilon,i_1}+\bar{u}_4^{\varepsilon,i_1}\|_{1/2-\delta_0}\sup_{t\in[0,T]}\|\bar{u}_1^{\varepsilon,j}\|_{-1/2-\delta}+\|\bar{u}^{\varepsilon,\sharp,i}\|_{1/2+\beta}.\endaligned\eqno(3.4)$$

Then $\bar{u}^\varepsilon=\bar{u}^\varepsilon_1+\bar{u}^\varepsilon_2+\bar{u}^\varepsilon_3+\bar{u}^\varepsilon_4$ solves (3.1) if and only if $\bar{u}^{\varepsilon,\sharp}$ solves the following equation:
$$\aligned \bar{u}^{\varepsilon,\sharp,i}=&P_t(\sum_{i_1=1}^3P^{ii_1}u_0^{i_1}-\bar{u}^{\varepsilon,i}_1(0)) -\frac{1}{2}\int_0^tP_{t-s}\sum_{i_1,j=1}^3P^{ii_1}D_j(\bar{u}_2^{\varepsilon,i_1}\diamond \bar{u}_2^{\varepsilon,j}+\bar{u}_2^{\varepsilon,i_1}(\bar{u}_3^{\varepsilon,j}+\bar{u}_4^{\varepsilon,j})+\bar{u}_2^{\varepsilon,j}(\bar{u}_3^{\varepsilon,i_1}+\bar{u}_4^{\varepsilon,i_1})\\&+(\bar{u}_3^{\varepsilon,i_1}+\bar{u}_4^{\varepsilon,i_1})(\bar{u}_3^{\varepsilon,j}+\bar{u}_4^{\varepsilon,j})
+\pi_>(\bar{u}_3^{\varepsilon,i_1}+\bar{u}_4^{\varepsilon,i_1},\bar{u}_1^{\varepsilon,j})+\pi_{0,\diamond}(\bar{u}_3^{\varepsilon,i_1},\bar{u}_1^{\varepsilon,j})+\pi_{0,\diamond}(\bar{u}_4^{\varepsilon,i_1},\bar{u}_1^{\varepsilon,j})\\&+\pi_>(\bar{u}_3^{\varepsilon,j}+\bar{u}_4^{\varepsilon,j},\bar{u}_1^{\varepsilon,i_1})+\pi_{0,\diamond}(\bar{u}_3^{\varepsilon,j},\bar{u}_1^{\varepsilon,i_1})
+\pi_{0,\diamond}(\bar{u}_4^{\varepsilon,j},\bar{u}_1^{\varepsilon,i_1}))ds\\&-\frac{1}{2}\int_0^tP_{t-s} \sum_{i_1,j_1=1}^3P^{ii_1}D_{j_1}(\pi_<(\bar{u}_3^{\varepsilon,i_1}+\bar{u}_4^{\varepsilon,i_1},\bar{u}_1^{\varepsilon,j_1})+\pi_<(\bar{u}_3^{\varepsilon,j_1}
+\bar{u}_4^{\varepsilon,j_1},\bar{u}_1^{\varepsilon,i_1}))ds\\&+\frac{1}{2}\sum_{i_1=1}^3P^{ii_1}\bigg(\sum_{j_1=1}^3D_{j_1}[\pi_<(\bar{u}_3^{\varepsilon,i_1}
+\bar{u}_4^{\varepsilon,i_1},\bar{K}^{\varepsilon,j_1})+\pi_<(\bar{u}_3^{\varepsilon,j_1}+\bar{u}_4^{\varepsilon,j_1},\bar{K}^{\varepsilon,i_1})]\bigg)
\\:=&P_t(\sum_{i_1=1}^3P^{ii_1}u_0^{i_1}-\bar{u}^{\varepsilon,i}_1(0))+\int_0^tP_{t-s}\bar{\phi}^{\varepsilon,\sharp,i}ds+\bar{F}^{\varepsilon,i},\endaligned\eqno(3.5)$$
where $\bar{F}^{\varepsilon,i}$ denotes the last two terms.
\vskip.10in
\no\textbf{Estimates for $\bar{\phi}^{\varepsilon,\sharp,i}$}
\vskip.10in
\th{Lemma 3.9} For $\bar{\phi}^{\varepsilon,\sharp,i}$ defined above, the following estimate holds:
$$\aligned&\|\bar{\phi}^{\varepsilon,\sharp,i}\|_{-1-2\delta-2\kappa}
\lesssim (1+(\bar{C}^\varepsilon_W)^3)(1+\|\bar{u}^{\varepsilon,\sharp}\|_{1/2+\beta}+\|\bar{u}_4^{\varepsilon}\|_{1/2-\delta_0})
+\|\bar{u}_4^\varepsilon\|_\delta^2,\endaligned\eqno(3.6)$$

\proof By a similar argument as [ZZ14] we obtain that for $\delta\leq\delta_0<1/2-3\delta/2$
$$\aligned &\|\pi_{0,\diamond}(\bar{u}_4^{\varepsilon,i},\bar{u}_1^{\varepsilon,j})\|_{-\delta}
\lesssim  \sum_{i_1,j_1=1}^3\bigg[\|\bar{u}_3^{\varepsilon,i_1}+\bar{u}_4^{\varepsilon,i_1}
\|_{1/2-\delta_0}\bar{C}_W^\varepsilon\|\bar{u}_1^{\varepsilon,j}\|_{-1/2-\delta/2}
\\&+\|\bar{u}_3^{\varepsilon,j_1}+\bar{u}_4^{\varepsilon,j_1}\|_{1/2-\delta_0}\|\pi_{0,\diamond}(P^{ii_1}D_{j_1}
\bar{K}^{\varepsilon,i_1},\bar{u}_1^{\varepsilon,j})\|_{-\delta}+\|\bar{u}^{\varepsilon,\sharp,i}\|_{1/2+\beta}\|\bar{u}_1^{\varepsilon,j}\|_{-1/2-\delta/2}
\\&+\|\bar{u}_3^{\varepsilon,i_1}+\bar{u}_4^{\varepsilon,i_1}\|_{1/2-\delta_0}\|\pi_{0,\diamond}(P^{ii_1}D_{j_1}
\bar{K}^{\varepsilon,j_1},\bar{u}_1^{\varepsilon,j})\|_{-\delta}\bigg]
,\endaligned$$
and
$$\aligned &\|\pi_>(\bar{u}_3^{\varepsilon,i_1}+\bar{u}_4^{\varepsilon,i_1},\bar{u}_1^{\varepsilon,j})
\|_{-2\delta-\kappa}\lesssim  (\|\bar{u}_3^{\varepsilon,i_1}\|_{1/2-\delta}+\|\bar{u}_4^{\varepsilon,i_1}\|_{1/2-\delta-\kappa})
\|\bar{u}_1^{\varepsilon,j}\|_{-1/2-\delta/2}\\\lesssim & [\|\bar{u}_3^{\varepsilon,i_1}\|_{1/2-\delta}+\sum_{i_2=1}^3(\|\bar{u}_3^{\varepsilon,i_2}+\bar{u}_4^{\varepsilon,i_2}\|_{1/2-\delta_0}
\bar{C}_W^\varepsilon+\|\bar{u}^{\varepsilon,\sharp,i_1}\|_{1/2+\beta})]\|\bar{u}_1^{\varepsilon,j}\|_{-1/2-\delta/2},\endaligned$$
where in the last inequality we used (3.4).

Hence by (3.5) we get that
$$\aligned&\|\bar{\phi}^{\varepsilon,\sharp,i}\|_{-1-2\delta-2\kappa}\\\lesssim& \sum_{i_1,j_1=1}^3\bigg[\|\bar{u}_2^{\varepsilon,i_1}\diamond \bar{u}_2^{\varepsilon,j_1}\|_{-\delta}+\|\pi_{0,\diamond}( \bar{u}_3^{\varepsilon,i_1},\bar{u}_1^{\varepsilon,j_1})\|_{-\delta}+\|\bar{u}_2^{\varepsilon,i_1}\|_{-\delta}\|\bar{u}_3^{\varepsilon,j_1}+\bar{u}_4^{\varepsilon,j_1}\|_{1/2-\delta_0}+\|\bar{u}_3^{\varepsilon}+\bar{u}_4^{\varepsilon}\|^2_{\delta}
\\&+\|\bar{u}_3^{\varepsilon,i_1}\|_{1/2-\delta}
\|\bar{u}_1^{\varepsilon,j_1}\|_{-1/2-\delta/2}+\|\bar{u}_3^{\varepsilon,i_1}+\bar{u}_4^{\varepsilon,i_1}\|_{1/2-\delta_0}
\|\bar{u}_1^{\varepsilon,j_1}\|_{-1/2-\delta/2}\bar{C}_W^\varepsilon\\&+\sum_{j,i_2=1}^3\|\bar{u}_3^{\varepsilon,i_1}+\bar{u}_4^{\varepsilon,i_1}\|_{1/2-\delta_0}
\|\pi_{0,\diamond}(P^{i_2i_1}D_{j_1}\bar{K}^{\varepsilon,j_1},\bar{u}_1^{\varepsilon,j})\|_{-\delta}
\\&+\sum_{j,i_2=1}^3\|\bar{u}_3^{\varepsilon,j_1}+\bar{u}_4^{\varepsilon,j_1}\|_{1/2-\delta_0}\|\pi_{0,\diamond}(P^{i_2i_1}D_{j_1}\bar{K}^{\varepsilon,i_1},
\bar{u}_1^{\varepsilon,j})\|_{-\delta}+\|\bar{u}^{\varepsilon,\sharp,i_1}\|_{1/2+\beta}\|\bar{u}_1^{\varepsilon,j_1}\|_{-1/2-\delta/2}\bigg]
\\\lesssim &(1+(\bar{C}^\varepsilon_W)^3)(1+\|\bar{u}^{\varepsilon,\sharp}\|_{1/2+\beta}+\|\bar{u}_4^{\varepsilon}\|_{1/2-\delta_0})
+\|\bar{u}_4^\varepsilon\|_\delta^2,\endaligned$$
which implies (3.6). $\hfill\Box$
\vskip.10in

\no\textbf{Estimate for $\bar{F}^\varepsilon$}
\vskip.10in

We now turn to $\bar{F}^\varepsilon$:
$$\aligned&\|\bar{F}^\varepsilon\|_{1/2+\beta}\\\lesssim& \sum_{i_1,j_1=1}^3\|\int_0^tP_{t-s}P^{ii_1}D_{j_1}\pi_<(\bar{u}_3^{\varepsilon,i_1}(s)+\bar{u}_4^{\varepsilon,i_1}(s)-(\bar{u}_3^{\varepsilon,i_1}(t)+\bar{u}_4^{\varepsilon,i_1}(t)),\bar{u}_1^{\varepsilon,j_1}(s))ds\|_{1/2+\beta}
\\&+\sum_{i_1,j_1=1}^3\|\int_0^tP_{t-s}P^{ii_1}D_{j_1}\pi_<(\bar{u}_3^{\varepsilon,i_1}(t)+\bar{u}_4^{\varepsilon,i_1}(t),\bar{u}_1^{\varepsilon,j_1}(s))ds\\&-P^{ii_1}D_{j_1}\pi_<(\bar{u}_3^{\varepsilon,i_1}(t)+\bar{u}_4^{\varepsilon,i_1}(t),\int_0^tP_{t-s}\bar{u}_1^{\varepsilon,j_1}ds)\|_{1/2+\beta}
\endaligned$$
$$\aligned&+\sum_{i_1,j_1=1}^3\|\int_0^tP_{t-s}P^{ii_1}D_{j_1}\pi_<(\bar{u}_3^{\varepsilon,j_1}(s)+\bar{u}_4^{\varepsilon,j_1}(s)-(\bar{u}_3^{\varepsilon,j_1}(t)
+\bar{u}_4^{\varepsilon,j_1}(t)),\bar{u}_1^{\varepsilon,i_1}(s))ds\|_{1/2+\beta}
\\&+\sum_{i_1,j_1=1}^3\|\int_0^tP_{t-s}P^{ii_1}D_{j_1}\pi_<(\bar{u}_3^{\varepsilon,j_1}(t)+\bar{u}_4^{\varepsilon,j_1}(t),\bar{u}_1^{\varepsilon,i_1}(s))ds
\\&-P^{ii_1}D_{j_1}\pi_<(\bar{u}_3^{\varepsilon,j_1}(t)+\bar{u}_4^{\varepsilon,j_1}(t),\int_0^tP_{t-s}\bar{u}_1^{\varepsilon,i_1}ds)\|_{1/2+\beta}
\\=&I_1+I_2+I_3+I_4.\endaligned$$
Now we only estimate $I_1$ and $I_2$ and for $I_3, I_4$ we have similar results. By Lemma 2.7 we have $$I_2\lesssim t^{\frac{1}{4}-\frac{\delta_0+\beta+\kappa+\frac{\delta}{2}}{2}}\|\bar{u}_3^\varepsilon+\bar{u}_4^\varepsilon\|_{1/2-\delta_0}\bar{C}^\varepsilon_W,\eqno(3.7)$$
where by the condition on $\beta$ we have  $\frac{1}{4}-\frac{\beta+\frac{\delta}{2}+\delta_0+\kappa}{2}>0.$
For $I_1$  Lemmas 2.2, 2.6 yield
$$\aligned I_1\lesssim&\int_0^t(t-s)^{-1-\frac{\delta/2+\beta+\kappa}{2}}\|\bar{u}_1^\varepsilon(s)\|_{-1/2-\delta/2}\|\bar{u}_3^\varepsilon(t)+\bar{u}^\varepsilon_4(t)-\bar{u}^\varepsilon_3(s)-\bar{u}^\varepsilon_4(s)\|_{\kappa/2}ds,\endaligned$$
and we note that by Lemma 2.8 that for $t>s>0$
$$\aligned &\|\bar{u}_3^\varepsilon(t)+\bar{u}^\varepsilon_4(t)-\bar{u}^\varepsilon_3(s)-\bar{u}^\varepsilon_4(s)\|_{\kappa/2}\\\lesssim
&\|(P_{t/2}-P_{s/2})(P_{t/2}+P_{s/2})(u_0-\bar{u}_1^\varepsilon(0))\|_{\kappa/2}+\|\int_0^s(P_{t-r}-P_{s-r})\bar{G}^\varepsilon(r)dr\|_{\kappa/2}+\|\int_s^tP_{t-r} \bar{G}^\varepsilon(r)dr\|_{\kappa/2}
\\\lesssim &(t-s)^{b_0}s^{-(z+2\kappa+2b_0)/2}\|u_0-\bar{u}_1^\varepsilon(0)\|_{-z}+(t-s)^{b/2}\int_0^s(s-r)^{-(3/2+\delta/2+2\kappa+b)/2}\|\bar{G}^\varepsilon(r)\|_{-3/2-\delta/2-\kappa/2}dr
\\&+(t-s)^{b_1}(\int_s^t\|P_{t-r} \bar{G}^\varepsilon(r)\|_{\kappa/2}^{\frac{1}{1-b_1}}dr)^{1-b_1}
\\\lesssim &(t-s)^{b_0}s^{-(z+2\kappa+2b_0)/2}\|u_0-\bar{u}_1^\varepsilon(0)\|_{-z}+(t-s)^{b/2}\int_0^s(s-r)^{-(3/2+\delta/2+2\kappa+b)/2}\|\bar{G}^\varepsilon(r)\|_{-3/2-\delta/2-\kappa/2}dr
\\&+(t-s)^{b_1}(\int_s^t(t-r)^{-\frac{3/2+\delta/2+2\kappa}{2(1-b_1)}}\| \bar{G}^\varepsilon(r)\|_{-3/2-\delta/2-\kappa/2}^{\frac{1}{1-b_1}}dr)^{1-b_1},\endaligned$$
where $\delta/2+\beta+2\kappa<2b_0<2-z-2\kappa$, $\delta/2+\beta+2\kappa<b<1/2-2\kappa-\delta/2$, $\frac{1}{2}(\delta/2+\beta+2\kappa)<b_1<[1-(\delta+z+\kappa)]\wedge \frac{1}{2}(1/2-\delta/2-2\kappa)$ and $$\aligned \bar{G}^\varepsilon=&\sum_{i_1,j=1}^3P^{ii_1}D_j(\bar{u}_1^{\varepsilon,i_1}\diamond \bar{u}_2^{\varepsilon,j}+\bar{u}_2^{\varepsilon,i_1}\diamond \bar{u}_1^{\varepsilon,j}+\bar{u}_1^{\varepsilon,i_1}\diamond (\bar{u}_3^{\varepsilon,j}+\bar{u}_4^{\varepsilon,j})+(\bar{u}_3^{\varepsilon,i_1} +\bar{u}_4^{\varepsilon,i_1})\diamond \bar{u}_1^{\varepsilon,j}\\&+\bar{u}_2^{\varepsilon,i_1}\diamond \bar{u}_2^{\varepsilon,j}+\bar{u}_2^{\varepsilon,i_1}(\bar{u}_3^{\varepsilon,j}+\bar{u}_4^{\varepsilon,j})+\bar{u}_2^{\varepsilon,j}(\bar{u}_3^{\varepsilon,i_1}+\bar{u}_4^{\varepsilon,i_1})+(\bar{u}_3^{\varepsilon,i_1}+\bar{u}_4^{\varepsilon,i_1})(\bar{u}_3^{\varepsilon,j}+\bar{u}_4^{\varepsilon,j})),\endaligned$$ Moreover, by a similar argument as (3.6) one has the following estimate
$$\|\bar{G}^\varepsilon\|_{-3/2-\delta/2-\kappa/2}\lesssim(1+(\bar{C}^\varepsilon_W)^3)(1+\|\bar{u}^{\varepsilon,\sharp}\|_{1/2+\beta}+\|\bar{u}_4^{\varepsilon}\|_{1/2-\delta_0})
+\|\bar{u}_4^\varepsilon\|_\delta^2.$$
Thus we obtain that
$$\aligned I_1\lesssim &\bar{C}_W^\varepsilon\bigg(\int_0^t(t-s)^{-1-\frac{\delta/2+\beta+2\kappa}{2}+b_0}s^{-(z+2\kappa+2b_0)/2}ds\|u_0-\bar{u}_1^\varepsilon(0)\|_{-z}
\\&+\int_0^t(t-s)^{-1-\frac{\delta/2+\beta+2\kappa}{2}+b/2}\int_0^s(s-r)^{-(3/2+\delta/2+2\kappa+b)/2}\|\bar{G}^\varepsilon(r)\|_{-3/2-\delta/2-\kappa/2}drds
\\&+\int_0^t(t-s)^{-1-\frac{\delta/2+\beta+2\kappa}{2}+b_1}(\int_s^t(t-r)^{-\frac{3/2+\delta/2+2\kappa}{2(1-b_1)}}\| \bar{G}^\varepsilon(r)\|_{-3/2-\delta/2-\kappa/2}^{\frac{1}{1-b_1}}dr)^{1-b_1}ds\bigg)
\\\lesssim &\bar{C}_W^\varepsilon\bigg(\int_0^t(t-s)^{-1-\frac{\delta/2+\beta+2\kappa}{2}+b_0}s^{-(z+2\kappa+2b_0)/2}ds\|u_0^\varepsilon-\bar{u}_1^\varepsilon(0)\|_{-z}
\\&+\int_0^t(t-r)^{-\frac{3}{4}-\frac{\delta+\beta+4\kappa}{2}}\int_0^1(1-l)^{-1-\frac{\delta/2+\beta+2\kappa}{2}+b/2}l^{-(3/2+\delta/2+2\kappa+b)/2}dl
\|\bar{G}^\varepsilon(r)\|_{-3/2-\delta/2-\kappa/2}dr
\\&+(\int_0^t(t-s)^{-1-\frac{\delta/2+\beta+2\kappa}{2}+b_1}ds)^{b_1}(\int_0^t\int_0^r(t-s)^{-1-\frac{\delta/2+\beta+2\kappa}{2}+b_1}(t-r)^{-\frac{3/2+\delta/2+2\kappa}{2(1-b_1)}}\\&\| \bar{G}^\varepsilon(r)\|_{-3/2-\delta/2-\kappa/2}^{\frac{1}{1-b_1}}dsdr)^{1-b_1}\bigg)
\\\lesssim &\bar{C}^\varepsilon_Wt^{-\frac{\delta/2+\beta+z}{2}-2\kappa}\|u_0-\bar{u}_1^\varepsilon(0)\|_{-z}
+C(\bar{C}^\varepsilon_W)\int_0^t(t-r)^{-\frac{3}{4}-\frac{\delta+\beta+4\kappa}{2}}(\|\bar{u}^{\varepsilon,\sharp}\|_{1/2+\beta}+\|\bar{u}_4^{\varepsilon}\|_{1/2-\delta_0}+\|\bar{u}_4^\varepsilon\|_\delta^2)dr
\\&+C(\bar{C}^\varepsilon_W)+C(t,\bar{C}^\varepsilon_W)\big{[}\int_0^t(t-r)^{-\frac{3/2+\delta/2+2\kappa}{2(1-b_1)}}(\|\bar{u}^{\varepsilon,\sharp}\|_{1/2+\beta}+\|\bar{u}_4^{\varepsilon}\|_{1/2-\delta_0}+\|\bar{u}_4^\varepsilon\|_\delta^2)^{\frac{1}{1-b_1}}dr\big{]}^{1-b_1}.
\endaligned\eqno(3.8)$$
Combining (3.7) and (3.8) we deduce that
$$\aligned&\|\bar{F}^\varepsilon\|_{1/2+\beta}\\\lesssim& t^{-\frac{\delta/2+\beta+z}{2}-2\kappa}\|u_0-\bar{u}_1^\varepsilon(0)\|_{-z}\bar{C}^\varepsilon_W
+C(\bar{C}^\varepsilon_W)\int_0^t(t-r)^{-\frac{3}{4}-\frac{\delta+\beta+4\kappa}{2}}(\|\bar{u}^{\varepsilon,\sharp}\|_{1/2+\beta}+\|\bar{u}_4^{\varepsilon}\|_{1/2-\delta_0}+\|\bar{u}_4^\varepsilon\|_\delta^2)dr
\\&+C(\bar{C}^\varepsilon_W)+C(t,\bar{C}^\varepsilon_W)\big{[}\int_0^t(t-r)^{-\frac{3/2+\delta/2+2\kappa}{2(1-b_1)}}(\|\bar{u}^{\varepsilon,\sharp}\|_{1/2+\beta}+\|\bar{u}_4^{\varepsilon}\|_{1/2-\delta_0}+\|\bar{u}_4^\varepsilon\|_\delta^2)^{\frac{1}{1-b_1}}dr\big{]}^{1-b_1}
\\&+t^{\frac{1}{4}-\frac{\delta_0+\beta+\kappa+\frac{\delta}{2}}{2}}\|\bar{u}_3^\varepsilon+\bar{u}_4^\varepsilon\|_{1/2-\delta_0}
\bar{C}^\varepsilon_W.\endaligned\eqno(3.9)$$
A similar argument also imply that
$$\aligned&\|\bar{F}^\varepsilon\|_{\delta}\\\lesssim& \bar{C}^\varepsilon_Wt^{{\frac{1}{4}-\frac{z}{2}-\frac{\kappa}{2}-\frac{3\delta}{4}}}\|u_0-\bar{u}_1^\varepsilon(0)\|_{-z}
+C(\bar{C}^\varepsilon_W)\int_0^t(t-r)^{-\frac{1+2\delta+\kappa}{2}}(\|\bar{u}^{\varepsilon,\sharp}\|_{1/2+\beta}+\|\bar{u}_4^{\varepsilon}\|_{1/2-\delta_0}+\|\bar{u}_4^\varepsilon\|_\delta^2)dr
\\&+C(\bar{C}^\varepsilon_W)+C(t,\bar{C}^\varepsilon_W)(\int_0^t(t-r)^{-\frac{1+2\delta+2\kappa}{2}}(\|\bar{u}^{\varepsilon,\sharp}\|_{1/2+\beta}+\|\bar{u}_4^{\varepsilon}\|_{1/2-\delta_0}+\|\bar{u}_4^\varepsilon\|_\delta^2)dr)
\\&+t^{\frac{1}{4}-\frac{\delta}{4}}\|\bar{u}_3^\varepsilon+\bar{u}_4^\varepsilon\|_{\delta}\bar{C}^\varepsilon_W,\endaligned\eqno(3.10)$$
where  we use $\|\bar{u}_3^\varepsilon(t)+\bar{u}^\varepsilon_4(t)-(\bar{u}_3^\varepsilon(s)+\bar{u}^\varepsilon_4(s))\|_{-1/2+3\delta/2+\kappa}$ to control the corresponding $I_1$ of $\|\bar{F}^\varepsilon\|_\delta$.
\vskip.10in
\no \textbf{Construction of the solution}
\vskip.10in
\th{Lemma 3.10}  There exists some $T_0>0$ (independent of $\varepsilon$) with $T_0\leq T^\varepsilon$ for $\varepsilon\in(0,1)$ such that
$$\aligned \sup_{t\in[0,T_0]}(t^{\delta+z+\kappa}\|\bar{u}^{ \varepsilon ,\sharp}\|_{1/2+\beta}+t^{\frac{\delta+z+\kappa}{2}}\|\bar{u}^{ \varepsilon,\sharp}(t)\|_{\delta})\lesssim C(T_0,\bar{C}_W^{\varepsilon},\|u_0\|_{-z})
.\endaligned\eqno(3.11)$$

\proof By paracontrolled ansatz Lemma 2.2 and (3.3) one then has
$$\aligned\|\bar{u}_4^{\varepsilon,i}\|_{1/2-\delta_0}\lesssim& t^{\delta/4}\sum_{i_1,j=1}^3\|\bar{u}_3^{\varepsilon,i_1}+\bar{u}_4^{\varepsilon,i_1}\|_{1/2-\delta_0}\bar{C}^\varepsilon_W+\|\bar{u}^{\varepsilon,\sharp,i}\|_{1/2-\delta_0},\endaligned$$
which shows that for $t$ small enough (only depending on $\bar{C}^\varepsilon_W$)
$$\aligned\|\bar{u}_4^{\varepsilon,i}\|_{1/2-\delta_0}\lesssim& (\bar{C}^\varepsilon_W)^2+\|\bar{u}^{\varepsilon,\sharp,i}\|_{1/2-\delta_0}.\endaligned\eqno(3.12)$$
Similarly, we have for $t$ small enough (only depending on $\bar{C}^\varepsilon_W$)
$$\aligned\|\bar{u}_4^{\varepsilon,i}\|_{\delta}\lesssim& (\bar{C}^\varepsilon_W)^2+\|\bar{u}^{\varepsilon,\sharp,i}\|_{\delta}.\endaligned\eqno(3.13)$$
Now first we assume that $\sup_\varepsilon \bar{C}^\varepsilon_W<\infty$ and in this case we can choose $t$  not depending on $\varepsilon.$
Then  (3.6) and (3.9) yield that for $\delta+z+\kappa<1$
$$\aligned &t^{\delta+z+\kappa}\|\bar{u}^{ \varepsilon,\sharp}(t)\|_{1/2+\beta}\\\lesssim & \|u_0-\bar{u}_1^{ \varepsilon}(0)\|_{-z}(1+\bar{C}_W^\varepsilon)+t^{\delta+z+\kappa}C(\bar{C}^{ \varepsilon}_W)\int_0^t(t-s)^{-3/4-\delta-\beta/2-3\kappa/2}s^{-(\delta+z+\kappa)}s^{\delta+z+\kappa}(1+\|\bar{u}^{ \varepsilon,\sharp}\|_{1/2+\beta}\\&+\|\bar{u}_4^{ \varepsilon}\|_{1/2-\delta_0}
+\|\bar{u}_4^{ \varepsilon}\|_\delta^2)ds
+t^{\delta+z+\kappa}\bigg{(}t^{-\frac{\delta/2+\beta+z}{2}-2\kappa}\|{u}_0-\bar{u}_1^\varepsilon(0)\|_{-z}\bar{C}_W^\varepsilon
\\&+C(\bar{C}^{ \varepsilon}_W)\int_0^t(t-r)^{-\frac{3}{4}-\frac{\delta+\beta+4\kappa}{2}}(\|\bar{u}^{ \varepsilon,\sharp}\|_{1/2+\beta}+\|\bar{u}_4^{ \varepsilon}\|_{1/2-\delta_0}+\|\bar{u}_4^{ \varepsilon}\|_\delta^2)dr
\\&+C(\bar{C}^{ \varepsilon}_W)+C(t,\bar{C}^{\varepsilon}_W)\big{[}\int_0^t(t-r)^{-\frac{3/2+\delta/2+2\kappa}{2(1-b_1)}}(\|\bar{u}^{ \varepsilon,\sharp}\|_{1/2+\beta}+\|\bar{u}_4^{ \varepsilon}\|_{1/2-\delta_0}+\|\bar{u}_4^{ \varepsilon}\|_\delta^2)^{\frac{1}{1-b_1}}dr\big{]}^{1-b_1}
\\&+t^{\frac{1}{4}-\frac{\delta_0+\beta+\kappa+\frac{\delta}{2}}{2}}\|\bar{u}_3^{ \varepsilon}+\bar{u}_4^{ \varepsilon}\|_{1/2-\delta_0}\bar{C}_W^{ \varepsilon}\bigg{)}
,\endaligned\eqno(3.14)$$
where we used the condition on $\beta$ to deduce $-3/4-\delta-\beta/2-3\kappa/2>-1$, $\frac{1/2+\beta+z}{2}\leq\delta+z$ and $\delta+z+\kappa-\frac{\delta/2+\beta+z}{2}-2\kappa\geq0$.
Hence by (3.12) (3.13) for $t$ small enough we have
$$\aligned &t^{\delta+z+\kappa}\|\bar{u}^{ \varepsilon,\sharp}(t)\|_{1/2+\beta}\\\lesssim & \|u_0-\bar{u}_1^{\varepsilon}(0)\|_{-z}(1+\bar{C}_W^\varepsilon)+t^{\delta+z+\kappa}C(\bar{C}_W^{\varepsilon})\int_0^t(t-s)^{-3/4-\delta-\beta/2-3\kappa/2}s^{-(\delta+z+\kappa)}s^{\delta+z+\kappa}(\|\bar{u}^{ \varepsilon,\sharp}\|_{1/2+\beta}\\&
+\|\bar{u}^{ \varepsilon,\sharp}\|_\delta^2)ds
+C(\bar{C}_W^{\varepsilon})
\\&+C(\bar{C}_W^{\varepsilon})t^{\delta+z+\kappa}\int_0^t(t-r)^{-\frac{3}{4}-\frac{\delta+\beta+4\kappa}{2}}r^{-(\delta+z+\kappa)}r^{\delta+z+\kappa}(\|\bar{u}^{ \varepsilon,\sharp}\|_{1/2+\beta}+\|\bar{u}^{ \varepsilon,\sharp}\|_\delta^2)dr
\\&+t^{\frac{\delta+z+\kappa}{1-b_1}}\int_0^t(t-r)^{-\frac{3/2+\delta/2+2\kappa}{2(1-b_1)}}r^{-\frac{(\delta+z+\kappa)}{1-b_1}}(r^{\delta+z+\kappa}(\|\bar{u}^{ \varepsilon,\sharp}\|_{1/2+\beta}+\|\bar{u}^{ \varepsilon,\sharp}\|_\delta^2))^{\frac{1}{1-b_1}}dr
.\endaligned\eqno(3.15)$$
A similar argument as (3.15) and using (3.13), (3.10) one also has that for $t$ small enough and $0<5\kappa<1-z-4\delta$
$$\aligned &t^{\frac{\delta+z+\kappa}{2}}\|\bar{u}^{ \varepsilon,\sharp}(t)\|_{\delta}\\\lesssim & \|u_0-\bar{u}_1^{\varepsilon}(0)\|_{-z}+
t^{\frac{\delta+z+\kappa}{2}}C(\bar{C}_W^{\varepsilon})\int_0^t(t-s)^{-1/2-3\delta/2-2\kappa}s^{-(\delta+z+\kappa)}s^{\delta+z+\kappa}(\|\bar{u}^{ \varepsilon,\sharp}\|_{1/2+\beta}\\&
+\|\bar{u}^{ \varepsilon,\sharp}\|_\delta^2)ds
+C(\bar{C}_W^{\varepsilon})
\\&+C(\bar{C}_W^{\varepsilon})t^{\frac{\delta+z+\kappa}{2}}\int_0^t(t-r)^{-\frac{1+2\delta+2\kappa}{2}}r^{-(\delta+z+\kappa)}r^{\delta+z+\kappa}(\|\bar{u}^{ \varepsilon,\sharp}\|_{1/2+\beta}+\|\bar{u}_4^{ \varepsilon}\|_\delta^2)dr
.\endaligned\eqno(3.16)$$
Combining  (3.15, 3.16) we get that by Bihari's inequality there exists some $T_0$ (independent of $\varepsilon$) such that
$$\aligned \sup_{t\in[0,T_0]}(t^{\delta+z+\kappa}\|\bar{u}^{ \varepsilon ,\sharp}\|_{1/2+\beta}+t^{\frac{\delta+z+\kappa}{2}}\|\bar{u}^{ \varepsilon,\sharp}(t)\|_{\delta})\lesssim C(T_0,\bar{C}_W^{\varepsilon},\|u_0\|_{-z})
,\endaligned$$
which combining with (3.6) implies that
$$\sup_{t\in[0,T_0]}t^{\delta+z+\kappa}\|\bar{\phi}^{ \varepsilon,\sharp}\|_{-1-2\delta-2\kappa}\lesssim C(T_0,\bar{C}_W^{\varepsilon},\|u_0\|_{-z}).\eqno(3.17)$$
Hence by (3.5) (3.12), (3.17) we have that
$$\sup_{t\in [0,T_0]}t^{\frac{1/2-\delta_0+z+\kappa}{2}}\|\bar{u}^{ \varepsilon}_4(t)\|_{1/2-\delta_0}\lesssim (\bar{C}_W^{\varepsilon})^2+\|u_0\|_{-z}+C(T_0,\bar{C}_W^{\varepsilon},\|u_0\|_{-z}),$$
which implies that $T^\varepsilon\geq T_0$. Here we used $z\geq1/2+\delta/2$ and
the following estimate
$$\aligned&\|\bar{F}^\varepsilon\|_{1/2-\delta_0}\\\lesssim& t^{\frac{-z+\delta_0-\delta/2-\kappa}{2}}\|u_0-\bar{u}_1^\varepsilon(0)\|_{-z}\bar{C}_W^\varepsilon
+C(\bar{C}^\varepsilon_W)\int_0^t(t-r)^{-\frac{3/2+\delta-\delta_0+2\kappa}{2}}(\|\bar{u}^{\varepsilon,\sharp}\|_{1/2+\beta}+\|\bar{u}_4^{\varepsilon}\|_{1/2-\delta_0}+\|\bar{u}_4^\varepsilon\|_\delta^2)dr
\\&+C(\bar{C}^\varepsilon_W)+C(t,\bar{C}^\varepsilon_W)\bigg(\int_0^t(t-r)^{-\frac{3/2+\delta+2\kappa-\delta_0}{2}}
(\|\bar{u}^{\varepsilon,\sharp}\|_{1/2+\beta}+\|\bar{u}_4^{\varepsilon}\|_{1/2-\delta_0}+\|\bar{u}_4^\varepsilon\|_\delta^2)dr\bigg)
\\&+t^{\frac{1}{4}-\frac{\delta}{4}}\|\bar{u}_3^\varepsilon+\bar{u}_4^\varepsilon\|_{\frac{1}{2}-\delta_0}\bar{C}^\varepsilon_W,\endaligned$$
which can be proved by a similar argument as (3.10), where one can $\|\bar{u}_3^\varepsilon(t)+\bar{u}^\varepsilon_4(t)-(\bar{u}_3^\varepsilon(s)+\bar{u}^\varepsilon_4(s))\|_{-\delta_0+\delta/2+\kappa}$ and $z+3\kappa+3\delta<1$ to control $\|\bar{F}^\varepsilon\|_{1/2-\delta_0}$.$\hfill\Box$
\vskip.10in

Moreover by paracontrolled ansatz and (3.3) we also obtain
$$\aligned\|\bar{u}_4^{ \varepsilon,i}\|_{-z}\lesssim& t^{\delta/4}\sum_{i_1=1}^3\|\bar{u}_3^{ \varepsilon,i_1}+\bar{u}_4^{ \varepsilon,i_1}\|_{-z}\bar{C}^\varepsilon_W+\|\bar{u}^{ \varepsilon,\sharp,i}\|_{-z},\endaligned$$
which combining with (3.5) implies that for $t$ small enough (only depending on $\bar{C}^\varepsilon_W$), $t\in [0,T_0]$ and $1-z-4\delta-5\kappa>0$
$$\aligned\|\bar{u}_4^{ \varepsilon,i}(t)\|_{-z}\lesssim&  C(\bar{C}^\varepsilon_W)+\|\bar{u}^{ \varepsilon,\sharp,i}(t)\|_{-z}
\\\lesssim & C(\bar{C}^\varepsilon_W)+\|Pu_0-\bar{u}_1^\varepsilon(0)\|_{-z}+\int_0^t(t-s)^{-\frac{1+2\delta+3\kappa-z}{2}}s^{-(\delta+z+\kappa)}s^{\delta+z+\kappa}\|\bar{\phi}^{ \varepsilon,\sharp}\|_{-1-2\delta-2\kappa}ds\\&+\|\bar{F}^\varepsilon(t)\|_{-z} \\\lesssim&C(\|u_0\|_{-z},T_0,\bar{C}^\varepsilon_W)+t^{\frac{1}{4}-\frac{\delta}{4}}\|\bar{u}_4^{ \varepsilon}(t)\|_{-z}\bar{C}_W^{\varepsilon}.\endaligned$$
Here in the last inequality we used  $$\aligned&\|\bar{F}^{ \varepsilon}(t)\|_{-z}\\\lesssim& C(\bar{C}^\varepsilon_W)\int_0^t(t-s)^{-\frac{3/2+\delta/2-z+\kappa}{2}}s^{-\frac{\delta+\kappa+z}{2}}ds\sup_{s\in[0,t]}s^{\frac{\delta+\kappa+z}{2}}\|\bar{u}^{ \varepsilon}_3+\bar{u}^{ \varepsilon}_4\|_{\delta}
+t^{\frac{1}{4}-\frac{\delta}{4}}\|\bar{u}_3^{ \varepsilon}(t)+\bar{u}_4^{ \varepsilon}(t)\|_{-z}\bar{C}_W^{\varepsilon}.\endaligned$$
Hence  for  $T_0$ small enough (only depending on $\bar{C}^\varepsilon_W$) one has the following bounds
$$\aligned\sup_{t\in[0,T_0]}\|\bar{u}_4^{ \varepsilon,i}(t)\|_{-z}\lesssim&  C(\|u_0\|_{-z},T_0,\bar{C}^\varepsilon_W).\endaligned\eqno(3.18)$$
Similar arguments show that for every $a>0$ there exists a sufficiently small $T_0>0$ such that the map $( u_0,\bar{u}_1^\varepsilon,\bar{u}^\varepsilon_1\diamond \bar{u}_1^\varepsilon,\bar{u}^\varepsilon_1\diamond \bar{u}^\varepsilon_2,\bar{u}^\varepsilon_2\diamond \bar{u}^\varepsilon_2,\pi_{0,\diamond}( \bar{u}^\varepsilon_3,\bar{u}^\varepsilon_1),\pi_{0,\diamond}(PD\bar{K}^\varepsilon,\bar{u}^\varepsilon_1))\mapsto \bar{u}_4^{ \varepsilon}$ is Lipschitz continuous on the set $$\max\{\|u_0\|_{-z},\bar{C}^{\varepsilon}_W\}\leq a.$$Here we consider $\bar{u}_4^{ \varepsilon}$ with respect to the norm of $$\sup_{t\in [0,T_0]}\|\bar{u}^{ \varepsilon}_4(t)\|_{-z}.$$
Thus we  obtain that  $\bar{u}_4^\varepsilon$ restricted to $[0,T_0]$ depends in a locally Lipschitz continuous way on the data $( u_0,\bar{u}_1^\varepsilon,\bar{u}^\varepsilon_1\diamond \bar{u}_1^\varepsilon,\bar{u}^\varepsilon_1\diamond \bar{u}^\varepsilon_2,\bar{u}^\varepsilon_2\diamond \bar{u}^\varepsilon_2,\pi_{0,\diamond}( \bar{u}^\varepsilon_3,\bar{u}^\varepsilon_1),\pi_{0,\diamond}(PD\bar{K}^\varepsilon,\bar{u}^\varepsilon_1))$.

In [ZZ14, Section 4] we proved that there exists some $\eta>0$ and $\bar{u}_1\in C([0,T];\mathcal{C}^{-1/2-\delta/2})$, $\bar{u}_2\in C([0,T];\mathcal{C}^{-\delta})$, $\bar{u}_3\in C([0,T];\mathcal{C}^{\frac{1}{2}-\delta})$  such that for every $p>0$ $$E\|\bar{u}_1^{\varepsilon}-\bar{u}_1\|_{C([0,T];\mathcal{C}^{-1/2-\delta/2})}^p\lesssim \varepsilon^{\eta p},$$
$$E\|\bar{u}_2^{\varepsilon}-\bar{u}_2\|_{C([0,T];\mathcal{C}^{-\delta})}^p\lesssim \varepsilon^{\eta p}.$$
$$E\|\bar{u}_3^{\varepsilon}-\bar{u}_3\|_{C([0,T];\mathcal{C}^{1/2-\delta})}^p\lesssim \varepsilon^{\eta p}.$$
Then for  $\varepsilon_k=2^{-k}\rightarrow0$ and $\epsilon>0$
$$\sum_{k=1}^\infty P(\|\bar{u}_1^{\varepsilon_k}-\bar{u}_1\|_{C([0,T];\mathcal{C}^{-1/2-\delta/2})}>\epsilon)\leq \sum_{k=1}^\infty 2^{-k\gamma}/\epsilon<\infty,\eqno(3.19)$$
which by Borel-Cantelli lemma implies that
$\bar{u}_1^{\varepsilon_k,i}-\bar{u}_1^i\rightarrow0$ in $C([0,T];\mathcal{C}^{-1/2-\delta/2})$ a.s., as $k\rightarrow\infty$.  The results for other terms are similar. Thus we obtain that $\sup_{\varepsilon_k=2^{-k},k\in\mathbb{N}}\bar{C}_W^{\varepsilon_k}<\infty $ a.s.,  $T_0$ independent of $\varepsilon$, $\bar{u}_4:=\lim_{k\rightarrow\infty}\bar{u}_4^{\varepsilon_k}$ on $[0,T_0]$, $u=\bar{u}_1+\bar{u}_2+\bar{u}_3+\bar{u}_4$ as solution of (1.4) on $[0,T_0]$ and
$$\sup_{t\in[0,T_0]}\|\bar{u}^{\varepsilon_k}-u\|_{-z}\rightarrow0\quad a.s..$$

Now we want to extend the solution to the maximal solution such that
$$\sup_{t\in[0,\tau)}\|{u}\|_{-z}=\infty.$$
A similar argument as above implies that there exists some $T_1(C(T_0))$ (for simplicity we assume $T_1\leq T_0$) such that for every $t^*\in[0,T_0]$
$$\sup_{t\in[t^*,t^*+T_1]}\big{[}(t-t^*)^{\delta+z+\kappa}\|\bar{u}^{ \varepsilon ,\sharp}\|_{1/2+\beta}+(t-t^*)^{\frac{\delta+z+\kappa}{2}}\|\bar{u}^{ \varepsilon,\sharp}(t)\|_{\delta}\big{]}\lesssim C(T_1,\bar{C}^\varepsilon_W,C(T_0),\|u(t^*)\|_{-z}).$$
Here the only difference is that $\bar{K}^{\varepsilon,i}$ satisfies the following equation
$$d\bar{K}^{\varepsilon,i}=(\Delta\bar{K}^{\varepsilon,i}+\bar{u}_1^{\varepsilon,i})dt\quad \bar{K}^{\varepsilon,i}(t^*)=0,$$
and by a similar argument as [ZZ14] we obtain that there exists some $\eta>0$ such that for every $p>1$
$$E\sup_{r\in[0,\cdot]}\|\pi_0(PD\int_r^{\cdot}P_{\cdot-s}\bar{u}^\varepsilon_1ds,\bar{u}^\varepsilon_1(\cdot))
-\pi_0(PD\int_r^{\cdot}P_{\cdot-s}\bar{u}_1ds,\bar{u}_1(\cdot))\|_{C([0,T];\mathcal{C}^{-\delta})}^p\lesssim\varepsilon^{p\eta},$$
 which implies that similar convergence also holds for  $\pi_0(PD\bar{K}^\varepsilon,\bar{u}^\varepsilon_1)$ in this case. Here we omit superscript for simplicity.

Thus for $t^*=T_0-\frac{T_1(C(T_0))}{2}$ we obtain the following estimate
$$\aligned&\sup_{t\in[T_0,T_0+\frac{T_1}{2}]}(t^{\delta+z+\kappa}\|\bar{u}^{ \varepsilon ,\sharp}\|_{1/2+\beta}+t^{\frac{\delta+z+\kappa}{2}}\|\bar{u}^{ \varepsilon,\sharp}(t)\|_{\delta})\\\lesssim& \sup_{t\in[T_0,T_0+\frac{T_1}{2}]}((t-t^*)^{\delta+z+\kappa}\|\bar{u}^{ \varepsilon ,\sharp}\|_{1/2+\beta}+(t-t^*)^{\frac{\delta+z+\kappa}{2}}\|\bar{u}^{ \varepsilon,\sharp}(t)\|_{\delta})\\\lesssim &C(T_1,\bar{C}^\varepsilon_W,C(T_0),\|u_0\|_{-z}),\endaligned$$
where we used (3.18) in the last inequality. 
Hence  by a similar argument as above we obtain the solution $u=\lim_{k\rightarrow\infty}\bar{u}^{\varepsilon_k}$ on $[0,T_0+\frac{T_1}{2}]$. Now by iterating  the above arguments  we obtain that there exist the explosion time $\tau>0$ and the maximal solution $u=\lim_{k\rightarrow\infty}\bar{u}^{\varepsilon_k}$ on $[0,\tau)$ such that
$$\sup_{t\in[0,\tau)}\|u(t)\|_{-z}=\infty.$$
Moreover for $L\geq0$ define $\tau_L:=\inf\{ t:\|u(t)\|_{-z}\geq L\}\wedge L$ and then $\tau_L$ increasing to $\tau$ and  $\bar{\rho}_L^\varepsilon:=\inf\{ t:\bar{C}^\varepsilon_W(t)\geq L\}.$
 By (3.11) and similar argument as above
we have $$\sup_{t\in[0,\tau_L\wedge \bar{\rho}_{L_1}^\varepsilon]}(t^{\delta+z+\kappa}\|\bar{u}^{ \varepsilon ,\sharp}\|_{1/2+\beta}+t^{\frac{\delta+z+\kappa}{2}}\|\bar{u}^{ \varepsilon,\sharp}(t)\|_{\delta})\lesssim C(L,L_1),\eqno(3.20)$$
which combining with (3.12), (3.13)  imply that for $L,L_1\geq0$
$$\sup_{t\in[0,\tau_L\wedge \bar{\rho}^\varepsilon_{L_1}]}(t^{\delta+z+\kappa}\| \bar{u}^{ \varepsilon}_4(t)\|_{1/2-\delta_0}+t^{\frac{\delta+z+\kappa}{2}}\| \bar{u}^{ \varepsilon}_4(t)\|_{\delta})\lesssim C(L,L_1),\eqno(3.21)$$
 $$\sup_{t\in[0,\tau_L\wedge \bar{\rho}^{\varepsilon}_{L_1}]}\|\bar{u}^{\varepsilon}-u\|_{-z}\rightarrow^P0.$$
 In particular, for a subsequence $\varepsilon_k=2^{-k}\rightarrow0,k\in\mathbb{N}$ we have
 $$\sup_{t\in[0,\tau_L]}\|\bar{u}^{\varepsilon_k}-u\|_{-z}\rightarrow0\quad  a.s..$$
\vskip.10in
\th{Remark 3.11} The construction here is slightly different from [ZZ14] as we consider the mild form of $\bar{u}_4$ here. The advantage of this approach is that it can also be used to obtain the solution of (1.3).

\vskip.10in
\textbf{Construction of solutions to approximating equation}:
Now we also split the equation (1.3) into the following four equations:
$$d u_1^{\varepsilon,i}=\Delta_\varepsilon u_1^{\varepsilon,i}dt+\sum_{i_1=1}^3P^{ii_1}H_\varepsilon dW,$$
$$d u_2^{\varepsilon,i}=\Delta_\varepsilon u_2^{\varepsilon,i}dt-\frac{1}{2}\sum_{i_1=1}^3P^{ii_1}(\sum_{j=1}^3D_j^\varepsilon(u_1^{\varepsilon,i_1}\diamond u_1^{\varepsilon,j}))dt\quad u^\varepsilon_2(0)=0,$$
$$du_3^{\varepsilon,i}=\Delta_\varepsilon u_3^{\varepsilon,i}dt-\frac{1}{2}\sum_{i_1=1}^3P^{ii_1}(\sum_{j=1}^3D^\varepsilon_j(u_1^{\varepsilon,i_1}\diamond u_2^{\varepsilon,j}+u_2^{\varepsilon,i_1}\diamond u_1^{\varepsilon,j}))dt\quad u^\varepsilon_3(0)=0,$$
and$$\aligned u_4^{\varepsilon,i}(t)=&P_t^\varepsilon(\sum_{i_1=1}^3P^{ii_1}u_0^{\varepsilon,i_1}-u_1^{\varepsilon,i_1}(0))-\frac{1}{2}\int_0^tP_{t-s}^\varepsilon
\bigg[\sum_{i_1,j=1}^3P^{ii_1}D_j^\varepsilon(u_1^{\varepsilon,i_1}\diamond (u_3^{\varepsilon,j}+u_4^{\varepsilon,j})+(u_3^{\varepsilon,i_1} +u_4^{\varepsilon,i_1})\diamond u_1^{\varepsilon,j}\\&+u_2^{\varepsilon,i_1}\diamond u_2^{\varepsilon,j}+u_2^{\varepsilon,i_1}(u_3^{\varepsilon,j}+u_4^{\varepsilon,j})+u_2^{\varepsilon,j}(u_3^{\varepsilon,i_1}
+u_4^{\varepsilon,i_1})+(u_3^{\varepsilon,i_1}+u_4^{\varepsilon,i_1})(u_3^{\varepsilon,j}+u_4^{\varepsilon,j}))\bigg]ds.\endaligned\eqno(3.22)$$
Here for $i,j=1,2,3$,
$$u_1^{\varepsilon,i}=\int_{-\infty}^t\sum_{i_1=1}^3P^{ii_1}P^\varepsilon_{t-s}H_\varepsilon dW^{i_1} ,$$
$$u^{\varepsilon,i}_1\diamond u^{\varepsilon,j}_1:=u^{\varepsilon,i}_1u^{\varepsilon,j}_1-C^{\varepsilon,ij}_0,$$
 $$u_2^{\varepsilon,j}\diamond u_1^{\varepsilon,i}=u_1^{\varepsilon,i}\diamond u_2^{\varepsilon,j}:=u_1^{\varepsilon,i} u_2^{\varepsilon,j}+\sum_{i_1=1}^3(C^{\varepsilon,i,i_1,j}(t)+\tilde{C}^{\varepsilon,i,i_1,j}(t))u_1^{\varepsilon,i_1},$$
 $$u_2^{\varepsilon,i}\diamond u_2^{\varepsilon,j}:=u_2^{\varepsilon,i} u_2^{\varepsilon,j}-\varphi_2^{\varepsilon,ij}(t)-C_2^{\varepsilon,ij},$$
 $${u}_3^{\varepsilon,i}\diamond {u}_1^{\varepsilon,j}={u}_1^{\varepsilon,j}\diamond {u}_3^{\varepsilon,i}=\pi_<({u}_3^{\varepsilon,i},{u}_1^{\varepsilon,j})+\pi_{0,\diamond}({u}_3^{\varepsilon,i},{u}_1^{\varepsilon,j})+\pi_>({u}_3^{\varepsilon,i},
 {u}_1^{\varepsilon,j}),$$
 $${u}_4^{\varepsilon,i}\diamond {u}_1^{\varepsilon,j}={u}_1^{\varepsilon,j}\diamond {u}_4^{\varepsilon,i}=\pi_<({u}_4^{\varepsilon,i},{u}_1^{\varepsilon,j})+\pi_{0,\diamond}({u}_4^{\varepsilon,i},{u}_1^{\varepsilon,j})+\pi_>({u}_4^{\varepsilon,i},
{u}_1^{\varepsilon,j}),$$
 $$\pi_{0,\diamond}(u_3^{\varepsilon,i}, u_1^{\varepsilon,j}):=\pi_{0}(u_3^{\varepsilon,i}, u_1^{\varepsilon,j})-\varphi_1^{\varepsilon,ij}(t)-C_1^{\varepsilon,ij}+\sum_{i_1=1}^3(C^{\varepsilon,i,i_1,j}(t)+\tilde{C}^{\varepsilon,i,i_1,j}(t))u_2^{\varepsilon,i_1},$$
 and$$\pi_{0,\diamond}({u}_4^{\varepsilon,i}, {u}_1^{\varepsilon,j}):=\pi_{0}({u}_4^{\varepsilon,i}, {u}_1^{\varepsilon,j})+\sum_{i_1=1}^3(C^{\varepsilon,i,i_1,j}(t)+\tilde{C}^{\varepsilon,i,i_1,j}(t))(u_3^{\varepsilon,i_1}+u_4^{\varepsilon,i_1}),$$
  with $C^\varepsilon_0$ is defined in Section 4.2, $C^\varepsilon, \tilde{C}^\varepsilon_0$ are given in introduction and Section 4.3,
 $C_1^\varepsilon$ and $\varphi_1^\varepsilon$ are defined in Section 4.4 and  $C_2^\varepsilon$ and $\varphi_2^\varepsilon$ are defined in Section 4.6 and $\varphi^\varepsilon_i$ converges to some $\varphi_i$ with respect to $\|\varphi\|=\sup_{t\in[0,T]}t^\rho|\varphi(t)|$ for every $\rho>0$ and $i=1,2$.
Define $K^\varepsilon$ be the solution to the following equation:
$$dK^{\varepsilon,i}=(\Delta_\varepsilon K^{\varepsilon,i}+u_1^{\varepsilon,i})dt\quad K^{\varepsilon,i}(0)=0.$$
For every $\varepsilon>0$ by a similar argument as the case for $\bar{u}^\varepsilon$  we obtain  solutions of equation (3.22): More precisely, for each $\varepsilon\in(0,1)$ there exists $u_4^\varepsilon$ satisfying equation (3.22) respectively before $T_0^\varepsilon>0$ such that $u_4^\varepsilon\in C((0,T_0^\varepsilon);\mathcal{C}^{1/2-\delta_0})$ with respect to the norm
$\sup_{t\in [0,T]}t^{\frac{1/2-\delta_0+z+\kappa}{2}}\|u_4(t)\|_{1/2-\delta_0}$ for $0<3\kappa<1/2-z+\delta_0$ and satisfies
$$\sup_{t\in [0,T_0^\varepsilon]}t^{\frac{1/2-\delta_0+z+\kappa}{2}}\|u_4^\varepsilon(t)\|_{1/2-\delta_0}=\infty.$$

Now by a similar argument as the case for $\bar{u}^\varepsilon$ we can also consider the paracontrolled ansatz for $i=1,2,3$,
$$u_4^{\varepsilon,i}=-\frac{1}{2}\sum_{i_1=1}^3P^{ii_1}(\sum_{j=1}^3D^\varepsilon_j[\pi_<(u_3^{\varepsilon,i_1}+u_4^{\varepsilon,i_1},K^{\varepsilon,j})
+\pi_<(u_3^{\varepsilon,j}+u_4^{\varepsilon,j},K^{\varepsilon,i_1})])+u^{\varepsilon,\sharp,i},$$
 with $\bar{u}^{\varepsilon,\sharp,i}(t)\in \mathcal{C}^{1/2+\beta}$ for $\beta$ as in Section 3.2 and $t\in(0,T_0^\varepsilon)$.

Then $u^\varepsilon=u_1^\varepsilon+u_2^\varepsilon+u_3^\varepsilon+u_4^\varepsilon$ solves (1.3) if and only if $u^{\varepsilon,\sharp}$ solves the following equation:
$$\aligned u^{\varepsilon,\sharp,i}(t)=& P_t^\varepsilon(\sum_{i_1=1}^3P^{ii_1}u_0^{\varepsilon,i_1}-u_1^{\varepsilon,i}(0))-\frac{1}{2}\int_0^tP^\varepsilon_{t-s}\sum_{i_1,j=1}^3P^{ii_1}D^\varepsilon_j(u_2^{\varepsilon,i_1}\diamond u_2^{\varepsilon,j}+u_2^{\varepsilon,i_1}(u_3^{\varepsilon,j}+u_4^{\varepsilon,j})\\&+u_2^{\varepsilon,j}(u_3^{\varepsilon,i_1}+u_4^{\varepsilon,i_1})+(u_3^{\varepsilon,i_1}+u_4^{\varepsilon,i_1})(u_3^{\varepsilon,j}+u_4^{\varepsilon,j})
+\pi_>(u_3^{\varepsilon,i_1}+u_4^{\varepsilon,i_1},u_1^{\varepsilon,j})+\pi_{0,\diamond}(u_3^{\varepsilon,i_1},u_1^{\varepsilon,j})\\&+\pi_{0,\diamond}(u_4^{\varepsilon,i_1},u_1^{\varepsilon,j})
+\pi_>(u_3^{\varepsilon,j}+u_4^{\varepsilon,j},u_1^{\varepsilon,i_1})+\pi_{0,\diamond}(u_3^{\varepsilon,j},u_1^{\varepsilon,i_1})
+\pi_{0,\diamond}(u_4^{\varepsilon,j},u_1^{\varepsilon,i_1}))ds\\&-\frac{1}{2}\int_0^tP_{t-s}^\varepsilon \sum_{i_1,j=1}^3P^{ii_1}D^\varepsilon_{j_1}(\pi_<(u_3^{\varepsilon,i_1}+u_4^{\varepsilon,i_1},u_1^{\varepsilon,j_1})
+\pi_<(u_3^{\varepsilon,j_1}+u_4^{\varepsilon,j_1},u_1^{\varepsilon,i_1}))ds\\&+\frac{1}{2}\sum_{i_1=1}^3P^{ii_1}
(\sum_{j=1}^3D^\varepsilon_j[\pi_<(u_3^{\varepsilon,i_1}+u_4^{\varepsilon,i_1},K^{\varepsilon,j})
+\pi_<(u_3^{\varepsilon,j}+u_4^{\varepsilon,j},K^{\varepsilon,i_1})])\\:=&P^\varepsilon_t(\sum_{i_1=1}^3P^{ii_1}u_0^{\varepsilon,i_1}-u_1^{\varepsilon,i}(0)
)
+\int_0^tP_{t-s}^\varepsilon\phi^{\varepsilon,\sharp,i}ds+ F^{\varepsilon,i},\endaligned\eqno(3.23)$$
where $F^{\varepsilon}$ represents the last two terms. A similar argument as the above proof and Lemmas 3.2 3.4 3.5 and 3.7 yield  similar estimates as (3.6) (3.9) (3.10) and (3.11) (3.18) for $\varepsilon\in(0,1)$.  Here the only difference lies in the estimate of $\pi_0(u_4^\varepsilon,u_1^\varepsilon)$. This is similar as the following (3.28) and we omit it here.
Moreover define $\tau^\varepsilon_L:=\inf\{ t:\|u^\varepsilon(t)\|_{-z}\geq L\}$ and $\rho_L^\varepsilon:=\inf\{ t:C^\varepsilon_W(t)\geq L\},$
where for $T>0$, $$\aligned C^\varepsilon_W(T):=&\sup_{t\in[0,T]}\bigg(\sum_{i=1}^3\|u_1^{\varepsilon,i}\|_{-1/2-\delta/2}+\sum_{i,j=1}^3\|u_1^{\varepsilon,i}\diamond u_1^{\varepsilon,j}\|_{-1-\delta/2}+\sum_{i,j=1}^3\|u_1^{\varepsilon,i}\diamond u_2^{\varepsilon,j}\|_{-1/2-\delta/2}\\&+\sum_{i,j=1}^3\|u_2^{\varepsilon,i}\diamond u_2^{\varepsilon,j}\|_{-\delta}+\sum_{i,j=1}^3\|\pi_{0,\diamond}( u_3^{\varepsilon,i},u_1^{\varepsilon,j})\|_{-\delta}\\&+\sum_{i,i_1,j,j_1=1}^3\|\pi_{0,\diamond}(P^{ii_1}D^\varepsilon_{j} K^{\varepsilon,j},u_1^{\varepsilon,j_1})\|_{-\delta}+\sum_{i,i_1,j,j_1=1}^3\|\pi_{0,\diamond}
(P^{ii_1}D_{j}^\varepsilon K^{\varepsilon,i_1},u_1^{\varepsilon,j_1})\|_{-\delta}\bigg),\endaligned$$
with
 $$\pi_{0,\diamond}(P^{ii_1}D^\varepsilon_jK^{\varepsilon,j}, u_1^{\varepsilon,j_1}):=\pi_{0}(P^{ii_1}D^\varepsilon_jK^{\varepsilon,j}, u_1^{\varepsilon,j_1})-C_3^{\varepsilon,i,i_1,j_1,j},$$
  $$\pi_{0,\diamond}(P^{ii_1}D^\varepsilon_jK^{\varepsilon,i_1}, u_1^{\varepsilon,j_1}):=\pi_{0}(P^{ii_1}D^\varepsilon_jK^{\varepsilon,i_1}, u_1^{\varepsilon,j_1})-\tilde{C}_3^{\varepsilon,i,j,j_1,i_1},$$
  where $C_3^\varepsilon, \tilde{C}_3^\varepsilon$ are given in Section 4.6.
 In Section 4 we will prove that for $i,j, i_1,j_1=1,2,3$ and every $\delta>0$ small enough, $\bar{u}_1^{\varepsilon,i}-u_1^{\varepsilon,i}\rightarrow0$ in $C([0,T];\mathcal{C}^{-\frac{1}{2}-\frac{\delta}{2}})$
and $\bar{u}_1^{\varepsilon,i}\diamond \bar{u}_1^{\varepsilon,j}-u_1^{\varepsilon,i}\diamond u_1^{\varepsilon,j}\rightarrow0$ in $C([0,T];\mathcal{C}^{-1-\delta/2})$, $\bar{u}_1^{\varepsilon,i}\diamond \bar{u}_2^{\varepsilon,j}-u_1^{\varepsilon,i}\diamond u_2^{\varepsilon,j}\rightarrow0$ in $C([0,T];\mathcal{C}^{-1/2-\delta/2})$, $\bar{u}_2^{\varepsilon,i}\diamond \bar{u}_2^{\varepsilon,j}-u_2^{\varepsilon,i}\diamond u_2^{\varepsilon,j}\rightarrow0$ in $C([0,T];\mathcal{C}^{-\delta})$, $\pi_{0,\diamond}(\bar{u}_3^{\varepsilon,i},\bar{u}_1^{\varepsilon,j})-\pi_{0,\diamond}(u_3^{\varepsilon,i},u_1^{\varepsilon,j})\rightarrow0$ in $C([0,T];\mathcal{C}^{-\delta})$ and $\pi_{0,\diamond}(P^{ii_1}D_j\bar{K}^{\varepsilon,j},\bar{u}_1^{\varepsilon,j_1})-\pi_{0,\diamond}(P^{ii_1}D^\varepsilon_jK^{\varepsilon,j},u_1^{\varepsilon,j_1})\rightarrow0,$ $\pi_{0,\diamond}(P^{ii_1}D_j\bar{K}^{\varepsilon,i_1},\bar{u}_1^{\varepsilon,j_1})-\pi_{0,\diamond}(P^{ii_1}D^\varepsilon_jK^{\varepsilon,i_1},
u_1^{\varepsilon,j_1})\rightarrow0$ in $C([0,T];\mathcal{C}^{-\delta})$, as $\varepsilon\rightarrow0$. By a similar argument as above we obtain that for $\varepsilon\in(0,1)$
$$\sup_{t\in[0,\tau_L^\varepsilon\wedge {\rho}_{L_1}^\varepsilon]}(t^{\delta+z+\kappa}\|u^{\varepsilon ,\sharp}\|_{1/2+\beta}+t^{\frac{\delta+z+\kappa}{2}}\|u^{\varepsilon,\sharp}(t)\|_{\delta})\lesssim C(L,L_1),\eqno(3.24)$$
and
$$\sup_{t\in[0,\tau_L^\varepsilon\wedge \rho^\varepsilon_{L_1}]}(t^{\delta+z+\kappa}\| u^{\varepsilon}_4(t)\|_{1/2-\delta_0}+t^{\frac{\delta+z+\kappa}{2}}\| u^{\varepsilon}_4(t)\|_{\delta})\lesssim C(L,L_1).\eqno(3.25)$$

\subsection{Estimate for $u^\varepsilon-\bar{u}^\varepsilon$ and the proof of Theorem 1.3}

First  by Lemmas 2.5, 3.3, 3.7, 3.8  for $2\kappa<\delta$ one has that
$$\sup_{t\in[0,T]}\sum_{i=1}^3\|u_2^{\varepsilon,i}-\bar{u}_2^{\varepsilon,i}\|_{-\delta}\lesssim \varepsilon^{\kappa/2}C_W^{\varepsilon}+\sup_{t\in[0,T]}\sum_{i,j=1}^3\|\bar{u}_1^{\varepsilon,i}\diamond \bar{u}_1^{\varepsilon,j}-u_1^{\varepsilon,i}\diamond u_1^{\varepsilon,j}\|_{-1-\delta/2},$$
$$\sup_{t\in[0,T]}\sum_{i=1}^3\|u_3^{\varepsilon,i}-\bar{u}_3^{\varepsilon,i}\|_{1/2-\delta}\lesssim \varepsilon^{\kappa/2}C_W^{\varepsilon}+\sup_{t\in[0,T]}\sum_{i,j=1}^3\|\bar{u}_1^{\varepsilon,i}\diamond \bar{u}_2^{\varepsilon,j}-u_1^{\varepsilon,i}\diamond u_2^{\varepsilon,j}\|_{-1/2-\delta/2}$$
and $$\|\bar{K}^{\varepsilon,i}(t)-K^{\varepsilon,i}(t)\|_{\frac{3}{2}-\delta}\lesssim \varepsilon^{\kappa/2}\sup_{s\in[0,t]}\|u_1^{\varepsilon,i}(s)\|_{-1/2-\delta/2}+t^{\delta/4}\sup_{s\in[0,t]}\|u_1^{\varepsilon,i}(s)-\bar{u}_1^{\varepsilon,i}(s)\|_{-1/2-\delta/2}.\eqno(3.26)$$

Paracontrolled ansatz and Lemmas 2.2, 2.5, 2.6, 3.2, 3.7 and 3.8  also imply the following estimate:
$$\aligned\|u_4^{\varepsilon,i}-\bar{u}_4^{\varepsilon,i}\|_{1/2-\delta-\kappa}\lesssim& \sum_{i_1,j=1}^3(\|u_3^{\varepsilon,i_1}+u_4^{\varepsilon,i_1}-\bar{u}_3^{\varepsilon,i_1}-\bar{u}_4^{\varepsilon,i_1}\|_{1/2-\delta_0}
\sup_{t\in[0,T]}\|u_1^{\varepsilon,j}\|_{-1/2-\delta/2}\\&+\varepsilon^{\kappa/2}\|\bar{u}_3^{\varepsilon,i_1}+\bar{u}_4^{\varepsilon,i_1}\|_{1/2-\delta_0}
\sup_{t\in[0,T]}\|\bar{u}_1^{\varepsilon,j}\|_{-1/2-\delta/2}\\&+\|\bar{u}_3^{\varepsilon,i_1}+\bar{u}_4^{\varepsilon,i_1}\|_{1/2-\delta_0}
\sup_{t\in[0,T]}\|u_1^{\varepsilon,j}-\bar{u}_1^{\varepsilon,j}\|_{-1/2-\delta/2})+\|u^{\varepsilon,\sharp,i}-\bar{u}^{\varepsilon,\sharp,i}\|_{1/2+\beta}.\endaligned\eqno(3.27)$$

To deal with $\|u^{\varepsilon,i}_4-\bar{u}^{\varepsilon,i}_4\|_{1/2-\delta_0}$, we first consider $\|\phi^{\varepsilon,\sharp,i}-\bar{\phi}^{\varepsilon,\sharp,i}\|_{-1-2\delta-2\kappa}$ and $\|F^\varepsilon-\bar{F}^\varepsilon\|_{1/2+\beta}$.
\vskip.10in

\no\textbf{Estimate of $\|\phi^{\varepsilon,\sharp,i}-\bar{\phi}^{\varepsilon,\sharp,i}\|_{-1-2\delta-2\kappa}$}
\vskip.10in
First we prove the following estimate.
\vskip.10in
\th{Lemma 3.12} Let $u\in \mathcal{C}^{\alpha+1}$ for $\alpha\in\mathbb{R}$. Then for every $0<\beta_0<1, \kappa>0$ and $a,b\geq0$ with $a+b>0$ we have
$$\|u(\cdot+a\varepsilon\eta_j)-u(\cdot-b\varepsilon\eta_j)\|_{\alpha+1-\beta_0-\kappa}\lesssim \varepsilon^{\beta_0}\|u\|_{\alpha+1},$$
where for $j=1,2,3,$ $\eta_j\in \mathbb{R}^3$ and the j-th component equals to $1$ and others are zero.

\proof By a similar argument as in the proof of Lemma 3.7 we have
 $$\aligned&\|\mathcal{F}^{-1}\theta(2^{-j}\cdot)\mathcal{F}(u(\cdot+a\varepsilon\eta_j)-u(\cdot-b\varepsilon\eta_j))\|_{L^p} \lesssim2^j\varepsilon\|\mathcal{F}^{-1}\theta(2^{-j}\cdot)\mathcal{F}u\|_{L^p} ,\endaligned$$
 and
 $$\aligned&\|\mathcal{F}^{-1}\theta(2^{-j}\cdot)\mathcal{F}(u(\cdot+a\varepsilon\eta_j)-u(\cdot-b\varepsilon\eta_j))\|_{L^p} \\\lesssim&\|\mathcal{F}^{-1}\theta(2^{-j}\cdot)\mathcal{F}(u(\cdot+a\varepsilon\eta_j))\|_{L^p}+\|\mathcal{F}^{-1}\theta(2^{-j}\cdot)\mathcal{F}(u(\cdot-b\varepsilon\eta_j))\|_{L^p} \\\lesssim&\|\mathcal{F}^{-1}\theta(2^{-j}\cdot)\mathcal{F}u\|_{L^p}.\endaligned$$
 Thus on $\mathbb{T}^3$  for $p$ large enough one has that $$\|u(\cdot+a\varepsilon\eta_j)-u(\cdot-b\varepsilon\eta_j)\|_{\alpha+1-\beta_0-\kappa}\lesssim\|u(\cdot+a\varepsilon\eta_j)-u(\cdot-b\varepsilon\eta_j)\|_{B^{\alpha+1-\beta_0}_{p,\infty}} \lesssim \varepsilon^{\beta_0}\|u\|_{B^{\alpha+1}_{p,\infty}} \lesssim\varepsilon^{\beta_0}\|u\|_{\alpha+1}.$$
$\hfill\Box$
\vskip.10in
Now we prove the following estimate for $\|\phi^{\varepsilon,\sharp,i}-\bar{\phi}^{\varepsilon,\sharp,i}\|_{-1-2\delta-2\kappa}$.
\vskip.10in
\th{Lemma 3.13}  For $1/2-3\delta/2-\delta_0-7\kappa/2>0$ one has the following estimate:
$$\aligned&\|\phi^{\varepsilon,\sharp,i}-\bar{\phi}^{\varepsilon,\sharp,i}\|_{-1-2\delta-2\kappa}
\\\lesssim &(\delta C_W^{\varepsilon}+\varepsilon^{\kappa/2}(C_W^{\varepsilon}+\bar{C}^\varepsilon_W+1))(1+(C_W^{\varepsilon})^2+(\bar{C}^\varepsilon_W)^2)(1+\|\bar{u}^{\varepsilon,\sharp}\|_{1/2+\beta}+\|\bar{u}_4^{\varepsilon}\|_{1/2-\delta_0}+\|u_4^{\varepsilon}\|_{1/2-\delta_0})
\\&+\|u_4^\varepsilon-\bar{u}_4^\varepsilon\|_{1/2-\delta_0}(1+\bar{C}^\varepsilon_W+(C_W^{\varepsilon})^2)+\|u_4^\varepsilon-\bar{u}_4^\varepsilon\|_{\delta}(\|u_4^{\varepsilon}\|_{\delta}+\|\bar{u}_4^{\varepsilon}\|_{\delta})+\varepsilon^{\kappa/2}\|\bar{u}_4^\varepsilon\|_\delta^2\\&+C_W^{\varepsilon}\|u^{\varepsilon,\sharp}-\bar{u}^{\varepsilon,\sharp}\|_{1/2+\beta},\endaligned\eqno(3.28)$$
where
 $$\aligned\delta C_W^{\varepsilon}:=&\sup_{t\in[0,T]}(\sum_{i=1}^3\|u_1^{\varepsilon,i}-\bar{u}_1^{\varepsilon,i}\|_{-1/2-\delta/2}+\sum_{i,j=1}^3\|u_1^{\varepsilon,i}\diamond u_1^{\varepsilon,j}-\bar{u}_1^{\varepsilon,i}\diamond \bar{u}_1^{\varepsilon,j}\|_{-1-\delta/2}\\&+\sum_{i,j=1}^3\|u_1^{\varepsilon,i}\diamond u_2^{\varepsilon,j}-\bar{u}_1^{\varepsilon,i}\diamond \bar{u}_2^{\varepsilon,j}\|_{-1/2-\delta/2}+\sum_{i,j=1}^3\|u_2^{\varepsilon,i}\diamond u_2^{\varepsilon,j}-\bar{u}_2^{\varepsilon,i}\diamond \bar{u}_2^{\varepsilon,j}\|_{-\delta}\\&+\sum_{i,j=1}^3\|\pi_{0,\diamond}( u_3^{\varepsilon,i},u_1^{\varepsilon,j})-\pi_{0,\diamond}( \bar{u}_3^{\varepsilon,i},\bar{u}_1^{\varepsilon,j})\|_{-\delta}+\sum_{i,i_1,j,j_1=1}^3\|\pi_{0,\diamond}(P^{ii_1}D_{j}^\varepsilon K^{\varepsilon,j},u_1^{\varepsilon,j_1})\\&-\pi_{0,\diamond}(P^{ii_1}D_{j} \bar{K}^{\varepsilon,j},\bar{u}_1^{\varepsilon,j_1})\|_{-\delta}+\sum_{i,i_1,j,j_1=1}^3\|\pi_{0,\diamond}
(P^{ii_1}D^\varepsilon_{j}K^{\varepsilon,i_1},u_1^{\varepsilon,j_1})-\pi_{0,\diamond}
(P^{ii_1}D_{j}\bar{K}^{\varepsilon,i_1},\bar{u}_1^{\varepsilon,j_1})\|_{-\delta}).\endaligned$$

\proof First  we consider $\pi_{0,\diamond}(u_4^{\varepsilon,i},u_1^{\varepsilon,j})-\pi_{0,\diamond}(\bar{u}_4^{\varepsilon,i},\bar{u}_1^{\varepsilon,j})$: by the paracontrolled ansatz one has for $i,j=1,2,3$,

$$\aligned&\pi_{0,\diamond}(u_4^{\varepsilon,i},u_1^{\varepsilon,j})-\pi_{0,\diamond}(\bar{u}_4^{\varepsilon,i},\bar{u}_1^{\varepsilon,j})\\=&-\frac{1}{2}
\bigg[\pi_{0,\diamond}(\sum_{i_1,j_1=1}^3P^{ii_1}\pi_<(u_3^{\varepsilon,i_1}+u_4^{\varepsilon,i_1},D^\varepsilon_{j_1}K^{\varepsilon,j_1}),u_1^{\varepsilon,j})
-\pi_{0,\diamond}(\sum_{i_1,j_1=1}^3P^{ii_1}\pi_<(\bar{u}_3^{\varepsilon,i_1}+\bar{u}_4^{\varepsilon,i_1},D_{j_1}
\bar{K}^{\varepsilon,j_1}),\bar{u}_1^{\varepsilon,j})\bigg]\\&-\frac{1}{2}\bigg[
\pi_{0,\diamond}(\sum_{i_1,j_1=1}^3P^{ii_1}\pi_<(u_3^{\varepsilon,j_1}+u_4^{\varepsilon,j_1},D^\varepsilon_{j_1}K^{\varepsilon,i_1}),u_1^{\varepsilon,j})
-\pi_{0,\diamond}(\sum_{i_1,j_1=1}^3P^{ii_1}\pi_<(\bar{u}_3^{\varepsilon,j_1}+\bar{u}_4^{\varepsilon,j_1},D_{j_1}\bar{K}^{\varepsilon,i_1}),\bar{u}_1^{\varepsilon,j})
\bigg]\\&-\frac{1}{2}\bigg[\sum_{i_1,j_1=1}^3\pi_0(P^{ii_1}\pi_<(D_{j_1}^\varepsilon(u_3^{\varepsilon,i_1}
+u_4^{\varepsilon,i_1}),K^{\varepsilon,j_1}),u_1^{\varepsilon,j})-\sum_{i_1,j_1=1}^3\pi_0(P^{ii_1}\pi_<(D_{j_1}(\bar{u}_3^{\varepsilon,i_1}
+\bar{u}_4^{\varepsilon,i_1}),\bar{K}^{\varepsilon,j_1}),\bar{u}_1^{\varepsilon,j})\bigg]\\
&-\frac{1}{2}\bigg[\sum_{i_1,j_1=1}^3\pi_0(P^{ii_1}\pi_<(D_{j_1}^\varepsilon(u_3^{\varepsilon,j_1}+u_4^{\varepsilon,j_1}),K^{\varepsilon,i_1}),u_1^{\varepsilon,j})
-\sum_{i_1,j_1=1}^3\pi_0(P^{ii_1}\pi_<(D_{j_1}(\bar{u}_3^{\varepsilon,j_1}+\bar{u}_4^{\varepsilon,j_1}),\bar{K}^{\varepsilon,i_1}),\bar{u}_1^{\varepsilon,j})\bigg]
\\&+(\pi_0(u^{\varepsilon,\sharp,i},u_1^{\varepsilon,j})-\pi_0(\bar{u}^{\varepsilon,\sharp,i},\bar{u}_1^{\varepsilon,j}))
+\psi_1^\varepsilon+\psi_2^\varepsilon.\endaligned\eqno(3.29)$$
Here $$\aligned\psi_1^\varepsilon=&-\frac{1}{2}\sum_{i_1,j_1=1}^3\pi_0(P^{ii_1}(\pi_<(D_{j_1}^\varepsilon(u_3^{\varepsilon,i_1}
+u_4^{\varepsilon,i_1}),K^{\varepsilon,j_1}(\cdot+a\varepsilon \eta_{j_1})-K^{\varepsilon,j_1}(\cdot)),u_1^{\varepsilon,j})\\&-\frac{1}{2}\sum_{i_1,j_1=1}^3\pi_0(P^{ii_1}(\pi_<((u_3^{\varepsilon,i_1}
+u_4^{\varepsilon,i_1})(\cdot-b\varepsilon\eta_{j_1})-(u_3^{\varepsilon,i_1}+u_4^{\varepsilon,i_1}),D_{j_1}^\varepsilon K^{\varepsilon,j_1}),u_1^{\varepsilon,j}))\\:=&L_1+L_2.\endaligned$$
and $$\aligned\psi_2^\varepsilon=&-\frac{1}{2}\sum_{i_1,j_1=1}^3(\pi_0(P^{ii_1}(\pi_<(D_{j_1}^\varepsilon(u_3^{\varepsilon,j_1}+u_4^{\varepsilon,j_1}),
K^{\varepsilon,i_1}(\cdot+a\varepsilon\eta_{j_1})-K^{\varepsilon,i_1}(\cdot)),u_1^{\varepsilon,j})\\&+\pi_0(P^{ii_1}(\pi_<((u_3^{\varepsilon,j_1}
+u_4^{\varepsilon,j_1})(\cdot-b\varepsilon\eta_{j_1})-(u_3^{\varepsilon,j_1}+u_4^{\varepsilon,j_1}),D_{j_1}^\varepsilon K^{\varepsilon,i_1}),u_1^{\varepsilon,j})).\endaligned$$
For $\psi_1^\varepsilon$ we have $$\aligned\|L_1\|_{\kappa}&\lesssim \sum_{i_1,j_1=1}^3\|u_1^\varepsilon\|_{-1/2-\delta/2}\|D^\varepsilon_{j_1}(u_3^{\varepsilon,i_1}+u_4^{\varepsilon,i_1})\|_{-1/2-\delta_0-\kappa}
\|K^\varepsilon(\cdot+a\varepsilon\eta_{j_1})-K^\varepsilon(\cdot)\|_{1+\delta_0+2\kappa+\delta/2}
\\&\lesssim \|u_1^\varepsilon\|_{-1/2-\delta/2}\|u_3^\varepsilon+u_4^\varepsilon\|_{1/2-\delta_0}\varepsilon^{\beta_0}\|K^\varepsilon\|_{3/2-\delta},\endaligned\eqno(3.30)$$
where $\beta_0=(1/2-3\delta/2-\delta_0-3\kappa),1/2-3\delta/2-\delta_0>3\kappa>0 $ and we used Lemma 3.12 in the last inequality.
Moreover $$\aligned&\|L_2\|_{-\delta}
\\\leq&\sum_{i_1,j_1=1}^3(\|\pi_0(P^{ii_1}\pi_<((u_3^{\varepsilon,i_1}+u_4^{\varepsilon,i_1})(\cdot-b\varepsilon\eta_{j_1})-(u_3^{\varepsilon,i_1}+u_4^{\varepsilon,i_1}),D_{j_1}^\varepsilon K^{\varepsilon,j_1}),u_1^{\varepsilon,j})
\\&-\pi_0(\pi_<((u_3^{\varepsilon,i_1}+u_4^{\varepsilon,i_1})(\cdot-b\varepsilon\eta_{j_1})-(u_3^{\varepsilon,i_1}+u_4^{\varepsilon,i_1}),P^{ii_1}
D_{j_1}^\varepsilon K^{\varepsilon,j_1}),u_1^{\varepsilon,j})\|_{\kappa}\\&+\|\pi_0(\pi_<((u_3^{\varepsilon,i_1}+u_4^{\varepsilon,i_1})(\cdot-b\varepsilon\eta_{j_1})-(u_3^{\varepsilon,i_1}+u_4^{\varepsilon,i_1}),P^{ii_1}D^\varepsilon_{j_1}K^{\varepsilon,j_1}),u_1^{\varepsilon,j})\\&-((u_3^{\varepsilon,i_1}+u_4^{\varepsilon,i_1})(\cdot-b\varepsilon\eta_{j_1})-(u_3^{\varepsilon,i_1}+u_4^{\varepsilon,i_1}))\pi_0(P^{ii_1}D^\varepsilon_{j_1}K^{\varepsilon,j_1}
,u_1^{\varepsilon,j})\|_{\kappa}\\&+\|((u_3^{\varepsilon,i_1}+u_4^{\varepsilon,i_1})(\cdot-b\varepsilon\eta_{j_1})-(u_3^{\varepsilon,i_1}+u_4^{\varepsilon,i_1}))
\pi_{0,\diamond}
(P^{ii_1}D^\varepsilon_{j_1}K^{\varepsilon,j_1},u_1^{\varepsilon,j})\|_{-\delta}\\&+\|(u_3^{\varepsilon,i_1}+u_4^{\varepsilon,i_1})
(\cdot-b\varepsilon\eta_{j_1})-(u_3^{\varepsilon,i_1}+u_4^{\varepsilon,i_1})\|_{3\delta/2+2\kappa}|C^\varepsilon(t)|)\\\lesssim&\sum_{i_1,j_1=1}^3
\|(u_3^{\varepsilon}+u_4^{\varepsilon})(\cdot-b\varepsilon\eta_{j_1})-(u_3^{\varepsilon}+u_4^{\varepsilon})\|_{3\delta/2+2\kappa}(\|D_{j_1}^\varepsilon K^{\varepsilon,j_1}\|_{1/2-\delta-\kappa}\|u_1^\varepsilon\|_{-1/2-\delta/2}\\&+\|\pi_{0,\diamond}
(P^{ii_1}D_{j_1}^\varepsilon K^{\varepsilon,j_1},u_1^{\varepsilon,j})\|_{-\delta}+C)\\\lesssim&\varepsilon^{\beta_0}\|u_3^{\varepsilon}+u_4^{\varepsilon}\|_{1/2-\delta_0}(\| K^{\varepsilon}\|_{3/2-\delta}\|u_1^\varepsilon\|_{-1/2-\delta/2}+\sum_{i_1,j_1=1}^3\|\pi_{0,\diamond}
(P^{ii_1}D_{j_1}^\varepsilon K^{\varepsilon,j_1},u_1^{\varepsilon,j})\|_{-\delta}+C),\endaligned\eqno(3.31)$$
where we used Lemmas 2.3, 2.4 in the second inequality and Lemma 3.12 in the last inequality. The estimate for $\psi_2^\varepsilon$ can be obtained similarly.

The desired estimate  for the third term, the forth term and the fifth term on the right hand side of  (3.29) are easily obtained by Lemmas 2.2 3.7 3.8 and (3.26). We only need to consider the first  term in the right hand side of (3.29) and the desired estimates for the second term follows similarly:
$$\aligned&\pi_{0,\diamond}(P^{ii_1}\pi_<(u_3^{\varepsilon,i_1}+u_4^{\varepsilon,i_1},D^\varepsilon_{j_1}K^{\varepsilon,j_1}),u_1^{\varepsilon,j})
-\pi_{0,\diamond}(P^{ii_1}\pi_<(\bar{u}_3^{\varepsilon,i_1}+\bar{u}_4^{\varepsilon,i_1},D_{j_1}\bar{K}^{\varepsilon,j_1}),\bar{u}_1^{\varepsilon,j})
\\=&\pi_0(P^{ii_1}\pi_<(u_3^{\varepsilon,i_1}+u_4^{\varepsilon,i_1},D_{j_1}^\varepsilon K^{\varepsilon,j_1}),u_1^{\varepsilon,j})
-\pi_0(\pi_<(u_3^{\varepsilon,i_1}+u_4^{\varepsilon,i_1},P^{ii_1}
D_{j_1}^\varepsilon K^{\varepsilon,j_1}),u_1^{\varepsilon,j})\\&-(\pi_0(P^{ii_1}\pi_<(\bar{u}_3^{\varepsilon,i_1}+\bar{u}_4^{\varepsilon,i_1},D_{j_1} \bar{K}^{\varepsilon,j_1}),\bar{u}_1^{\varepsilon,j})
-\pi_0(\pi_<(\bar{u}_3^{\varepsilon,i_1}+\bar{u}_4^{\varepsilon,i_1},P^{ii_1}
D_{j_1}\bar{K}^{\varepsilon,j_1}),\bar{u}_1^{\varepsilon,j}))\\&+\pi_0(\pi_<(u_3^{\varepsilon,i_1}+u_4^{\varepsilon,i_1},P^{ii_1}D^\varepsilon_{j_1}K^{\varepsilon,j_1}),u_1^{\varepsilon,j})-(u_3^{\varepsilon,i_1}+u_4^{\varepsilon,i_1})\pi_0(P^{ii_1}D^\varepsilon_{j_1}K^{\varepsilon,j_1}
,u_1^{\varepsilon,j})\\&-(\pi_0(\pi_<(\bar{u}_3^{\varepsilon,i_1}+\bar{u}_4^{\varepsilon,i_1},P^{ii_1}D_{j_1}\bar{K}^{\varepsilon,j_1}),\bar{u}_1^{\varepsilon,j})-(\bar{u}_3^{\varepsilon,i_1}+\bar{u}_4^{\varepsilon,i_1})\pi_0(P^{ii_1}D_{j_1}\bar{K}^{\varepsilon,j_1}
,u_1^{\varepsilon,j}))\\&+(u_3^{\varepsilon,i_1}+u_4^{\varepsilon,i_1})\pi_{0,\diamond}
(P^{ii_1}D^\varepsilon_{j_1}K^{\varepsilon,j_1},u_1^{\varepsilon,j})-(\bar{u}_3^{\varepsilon,i_1}+\bar{u}_4^{\varepsilon,i_1})\pi_{0,\diamond}
(P^{ii_1}D_{j_1}\bar{K}^{\varepsilon,j_1},\bar{u}_1^{\varepsilon,j}).\endaligned$$
  Lemmas 3.7, 3.8 and Lemmas 2.3, 2.4 imply for $\delta\leq\delta_0<1/2-3\delta/2-\kappa$ that
$$\aligned &\|\pi_{0,\diamond}(P^{ii_1}\pi_<(u_3^{\varepsilon,i_1}+u_4^{\varepsilon,i_1},D^\varepsilon_{j_1}K^{\varepsilon,j_1}),u_1^{\varepsilon,j})-\pi_{0,\diamond}(P^{ii_1}\pi_<(\bar{u}_3^{\varepsilon,i_1}+\bar{u}_4^{\varepsilon,i_1},D_{j_1}\bar{K}^{\varepsilon,j_1}),\bar{u}_1^{\varepsilon,j})\|_{-\delta}\\\lesssim & \|u_3^{\varepsilon,i_1}+u_4^{\varepsilon,i_1}-(\bar{u}_3^{\varepsilon,i_1}+\bar{u}_4^{\varepsilon,i_1})\|_{1/2-\delta_0}\|K^{\varepsilon,j_1}
\|_{3/2-\delta}\|u_1^{\varepsilon,j}\|_{-1/2-\delta/2}\\&+\|\bar{u}_3^{\varepsilon,i_1}+\bar{u}_4^{\varepsilon,i_1}\|_{1/2-\delta_0}
\|u_1^{\varepsilon,j}\|_{-1/2-\delta/2}(\varepsilon^{\kappa/2}\|K^{\varepsilon,j_1}\|_{3/2-\delta}+\|K^{\varepsilon,j_1}-\bar{K}^{\varepsilon,j_1}\|_{3/2-\delta})
\\&+\|\bar{u}_3^{\varepsilon,i_1}+\bar{u}_4^{\varepsilon,i_1}
\|_{1/2-\delta_0}\|\bar{K}^{\varepsilon,j_1}\|_{3/2-\delta}\|u_1^{\varepsilon,j}-\bar{u}_1^{\varepsilon,j}\|_{-1/2-\delta/2}\\&+\|u_3^{\varepsilon,i_1}+u_4^{\varepsilon,i_1}\|_{1/2-\delta_0}\|\pi_{0,\diamond}(P^{ii_1}D^\varepsilon_{j_1}K^{\varepsilon,j_1},u_1^{\varepsilon,j})-\pi_{0,\diamond}(P^{ii_1}D_{j_1}\bar{K}^{\varepsilon,j_1},\bar{u}_1^{\varepsilon,j})\|_{-\delta}
\\&+\|u_3^{\varepsilon,i_1}+u_4^{\varepsilon,i_1}-\bar{u}_3^{\varepsilon,i_1}-\bar{u}_4^{\varepsilon,i_1}\|_{1/2-\delta_0}\|\pi_{0,\diamond}(P^{ii_1}D_{j_1}\bar{K}^{\varepsilon,j_1},\bar{u}_1^{\varepsilon,j})\|_{-\delta}
.\endaligned\eqno(3.32)$$
Hence by (3.26) (3.30-3.32) we obtain for $i,j=1,2,3,$ that
$$\aligned &\|\pi_{0,\diamond}(u_4^{\varepsilon,i},u_1^{\varepsilon,j})-\pi_{0,\diamond}(\bar{u}_4^{\varepsilon,i},\bar{u}_1^{\varepsilon,j})\|_{-\delta}
\\\lesssim  &  \sum_{i_1,j_1=1}^3\big{[}\|u_3^{\varepsilon,i_1}+u_4^{\varepsilon,i_1}-(\bar{u}_3^{\varepsilon,i_1}+\bar{u}_4^{\varepsilon,i_1})\|_{1/2-\delta_0}
\sup_{t\in[0,T]}\|u_1^{\varepsilon,j_1}\|_{-1/2-\delta/2}\|u_1^{\varepsilon,j}\|_{-1/2-\delta/2}\\&+\|\bar{u}_3^{\varepsilon,i_1}
+\bar{u}_4^{\varepsilon,i_1}\|_{1/2-\delta_0}\|u_1^{\varepsilon,j}\|_{-1/2-\delta/2}(\varepsilon^{\kappa/2}\sup_{t\in[0,T]}\|u_1^{\varepsilon,j_1}\|_{-1/2-\delta/2}
+\sup_{t\in[0,T]}\|u_1^{\varepsilon,j_1}-\bar{u}_1^{\varepsilon,j_1}\|_{-1/2-\delta/2})
\\&+\|\bar{u}_3^{\varepsilon,i_1}+\bar{u}_4^{\varepsilon,i_1}
\|_{1/2-\delta_0}\sup_{t\in[0,T]}\|\bar{u}_1^{\varepsilon,j_1}\|_{-1/2-\delta/2}\|u_1^{\varepsilon,j}-\bar{u}_1^{\varepsilon,j}\|_{-1/2-\delta/2}
\\&+\|u_3^{\varepsilon,i_1}+u_4^{\varepsilon,i_1}\|_{1/2-\delta_0}\|\pi_{0,\diamond}(P^{ii_1}D^\varepsilon_{j_1}K^{\varepsilon,j_1},u_1^{\varepsilon,j})
-\pi_{0,\diamond}(P^{ii_1}D_{j_1}\bar{K}^{\varepsilon,j_1},\bar{u}_1^{\varepsilon,j})\|_{-\delta}
\\&+\|u_3^{\varepsilon,i_1}+u_4^{\varepsilon,i_1}-\bar{u}_3^{\varepsilon,i_1}-\bar{u}_4^{\varepsilon,i_1}\|_{1/2-\delta_0}\|\pi_{0,\diamond}(P^{ii_1}D_{j_1}
\bar{K}^{\varepsilon,j_1},\bar{u}_1^{\varepsilon,j})\|_{-\delta}\\&+\|u_3^{\varepsilon,j_1}+u_4^{\varepsilon,j_1}\|_{1/2-\delta_0}
\|\pi_{0,\diamond}(P^{ii_1}D^\varepsilon_{j_1}K^{\varepsilon,i_1},u_1^{\varepsilon,j})
-\pi_{0,\diamond}(P^{ii_1}D_{j_1}\bar{K}^{\varepsilon,i_1},\bar{u}_1^{\varepsilon,j})\|_{-\delta}
\\&+\|u_3^{\varepsilon,j_1}+u_4^{\varepsilon,j_1}-\bar{u}_3^{\varepsilon,j_1}-\bar{u}_4^{\varepsilon,j_1}\|_{1/2-\delta_0}\|\pi_{0,\diamond}(P^{ii_1}D_{j_1}
\bar{K}^{\varepsilon,i_1},\bar{u}_1^{\varepsilon,j})\|_{-\delta}\\&+\|u^{\varepsilon,\sharp,i}-\bar{u}^{\varepsilon,\sharp,i}\|_{1/2+\beta}\|u_1^{\varepsilon,j}\|_{-1/2-\delta/2}+\|\bar{u}^{\varepsilon,\sharp,i}\|_{1/2+\beta}\|u_1^{\varepsilon,j}-\bar{u}_1^{\varepsilon,j}\|_{-1/2-\delta/2}
\\&+\varepsilon^{\beta_0}(\|u_3^{\varepsilon,i_1}+u_4^{\varepsilon,i_1}\|_{1/2-\delta_0}\|u_1^{\varepsilon,j}\|_{-1/2-\delta/2}
\|K^{\varepsilon,j_1}\|_{3/2-\delta}+\|\pi_{0,\diamond}(P^{ii_1}D^\varepsilon_{j_1}K^{\varepsilon,i_1},u_1^{\varepsilon,j})\|_{-\delta}\|u_3^{\varepsilon,j_1}+u_4^{\varepsilon,j_1}\|_{1/2-\delta_0}\\&
+\|\pi_{0,\diamond}(P^{ii_1}D^\varepsilon_{j_1}K^{\varepsilon,j_1},u_1^{\varepsilon,j})\|_{-\delta}\|u_3^{\varepsilon,i_1}+u_4^{\varepsilon,i_1}\|_{1/2-\delta_0}+C\|u_3^{\varepsilon,i_1}+u_4^{\varepsilon,i_1}\|_{1/2-\delta_0})\big{]}\\\lesssim& [\delta C_W^{\varepsilon}+\varepsilon^{\kappa/2}(C_W^{\varepsilon}+\bar{C}^\varepsilon_W+1)](1+(C_W^{\varepsilon})^2+(\bar{C}^\varepsilon_W)^2)
(1+\|\bar{u}_4^{\varepsilon}\|_{\frac{1}{2}-\delta_0}+\|\bar{u}^{\varepsilon,\sharp}\|_{\frac{1}{2}+\beta}
+\|u_4^{\varepsilon}\|_{\frac{1}{2}-\delta_0})\\&+\|u_4^{\varepsilon}-\bar{u}_4^{\varepsilon}\|_{\frac{1}{2}-\delta_0}[(C_W^{\varepsilon})^2+\bar{C}^\varepsilon_W]
+\|u^{\varepsilon,\sharp}-\bar{u}^{\varepsilon,\sharp}\|_{\frac{1}{2}+\beta}C_W^{\varepsilon}.\endaligned\eqno(3.33)$$
Here we also used Lemma 2.2 to obtain the estimates for  the third and the fourth terms on the right hand side of (3.29).
Now we turn to the term $\pi_>(u_3^{\varepsilon,i_1}+u_4^{\varepsilon,i_1},u_1^{\varepsilon,j})-\pi_>(\bar{u}_3^{\varepsilon,i_1}+\bar{u}_4^{\varepsilon,i_1},\bar{u}_1^{\varepsilon,j})$ for $i_1,j=1,2,3$:  Lemma 3.2 yields the following estimate
$$\aligned &\|\pi_>(u_3^{\varepsilon,i_1}+u_4^{\varepsilon,i_1},u_1^{\varepsilon,j})-\pi_>(\bar{u}_3^{\varepsilon,i_1}+\bar{u}_4^{\varepsilon,i_1},\bar{u}_1^{\varepsilon,j})
\|_{-2\delta-\kappa}\\\lesssim & (\|u_3^{\varepsilon,i_1}-\bar{u}_3^{\varepsilon,i_1}\|_{1/2-\delta}+\|u_4^{\varepsilon,i_1}-\bar{u}_4^{\varepsilon,i_1}\|_{1/2-\delta-\kappa})
\|u_1^{\varepsilon,j}\|_{-1/2-\delta/2}
+\|\bar{u}^\varepsilon_3+\bar{u}^\varepsilon_4\|_{1/2-\delta-\kappa}\|u_1^\varepsilon-\bar{u}_1^\varepsilon\|_{-1/2-\delta/2}\\\lesssim & \bigg(\|u_3^{\varepsilon,i_1}-\bar{u}_3^{\varepsilon,i_1}\|_{1/2-\delta}+\sum_{i_2,j_1=1}^3(\|u_3^{\varepsilon,i_2}+u_4^{\varepsilon,i_2}
-\bar{u}_3^{\varepsilon,i_2}-\bar{u}_4^{\varepsilon,i_2}\|_{1/2-\delta_0}\sup_{t\in[0,T]}\|u_1^{\varepsilon,j_1}\|_{-1/2-\delta/2}
\\&+\|\bar{u}_3^{\varepsilon,i_2}+\bar{u}_4^{\varepsilon,i_2}\|_{1/2-\delta_0}\sup_{t\in[0,T]}\|u_1^{\varepsilon,j_1}-\bar{u}_1^{\varepsilon,j_1}\|_{-1/2-\delta/2}
+\varepsilon^{\kappa/2}\|\bar{u}_3^{\varepsilon,i_2}+\bar{u}_4^{\varepsilon,i_2}\|_{1/2-\delta_0}\sup_{t\in[0,T]}\|u_1^{\varepsilon,j_1}\|_{-1/2-\delta/2})
\endaligned$$
$$\aligned&+\|u^{\varepsilon,\sharp,i_1}-\bar{u}^{\varepsilon,\sharp,i_1}\|_{1/2+\beta}\bigg)\|u_1^{\varepsilon,j}\|_{-1/2-\delta/2}
+\|\bar{u}^{\varepsilon,i_1}_3+\bar{u}^{\varepsilon,i_1}_4\|_{1/2-\delta-\kappa}\sup_{t\in[0,T]}\|u_1^{\varepsilon,j}-\bar{u}_1^{\varepsilon,j}\|_{-1/2-\delta/2}
\\\lesssim& [\delta C_W^{\varepsilon}+\varepsilon^{\kappa/2}(C_W^{\varepsilon}+\bar{C}^\varepsilon_W)](1+(C_W^{\varepsilon})^2+(\bar{C}^\varepsilon_W)^2)
(1+\|\bar{u}_4^{\varepsilon}\|_{\frac{1}{2}-\delta_0}+\|\bar{u}^{\varepsilon,\sharp}\|_{\frac{1}{2}+\beta}
+\|u_4^{\varepsilon}\|_{\frac{1}{2}-\delta_0})\\&+\|u_4^{\varepsilon}-\bar{u}_4^{\varepsilon}\|_{\frac{1}{2}-\delta_0}(C_W^{\varepsilon})^2
+\|u^{\varepsilon,\sharp}-\bar{u}^{\varepsilon,\sharp}\|_{\frac{1}{2}+\beta}C_W^{\varepsilon},\endaligned\eqno(3.34)$$
where in the last inequality we used (3.4) and (3.27).
Combining (3.33), (3.34), (3.5) and (3.23) we deduce that for $\beta_0>\frac{\kappa}{2}$
$$\aligned&\|\phi^{\varepsilon,\sharp,i}-\bar{\phi}^{\varepsilon,\sharp,i}\|_{-1-2\delta-2\kappa}\\\lesssim& \sum_{i_1,j_1,j=1}^3\bigg[\|u_2^{\varepsilon,i_1}\diamond u_2^{\varepsilon,j_1}-\bar{u}_2^{\varepsilon,i_1}\diamond \bar{u}_2^{\varepsilon,j_1}\|_{-\delta}+\|u_2^{\varepsilon,i_1}-\bar{u}_2^{\varepsilon,i_1}\|_{-\delta}\|u_3^{\varepsilon,j_1}+u_4^{\varepsilon,j_1}\|_{1/2-\delta_0}\\&+(\|\bar{u}_2^{\varepsilon,i_1}\|_{-\delta}+\|u_3^{\varepsilon,i_1}+u_4^{\varepsilon,i_1}\|_{\delta}+\|\bar{u}_3^{\varepsilon,i_1}+\bar{u}_4^{\varepsilon,i_1}\|_{\delta})\|u_3^{\varepsilon,j_1}+u_4^{\varepsilon,j_1}-(\bar{u}_3^{\varepsilon,j_1}+\bar{u}_4^{\varepsilon,j_1})\|_{\delta}
\\&+\|\pi_{0,\diamond}( u_3^{\varepsilon,i_1},u_1^{\varepsilon,j_1})-\pi_{0,\diamond}( \bar{u}_3^{\varepsilon,i_1},\bar{u}_1^{\varepsilon,j_1})\|_{-\delta}\bigg]
\\&+\varepsilon^{\kappa/2}\sum_{i_1,j_1=1}^3\bigg[\|\bar{u}_2^{\varepsilon,i_1}\diamond \bar{u}_2^{\varepsilon,j_1}\|_{-\delta}+\|\pi_{0,\diamond}( \bar{u}_3^{\varepsilon,i_1},\bar{u}_1^{\varepsilon,j_1})\|_{-\delta}+\|\bar{u}_2^{\varepsilon,i_1}\|_{-\delta}\|\bar{u}_3^{\varepsilon,j_1}+\bar{u}_4^{\varepsilon,j_1}\|_{1/2-\delta_0}+\|\bar{u}_3^{\varepsilon}+\bar{u}_4^{\varepsilon}\|^2_{\delta}
\\&+\|\bar{u}_3^{\varepsilon,i_1}\|_{1/2-\delta}
\|\bar{u}_1^{\varepsilon,j_1}\|_{-1/2-\delta/2}+\sum_{j_2=1}^3\|\bar{u}_3^{\varepsilon,i_1}+\bar{u}_4^{\varepsilon,i_1}\|_{1/2-\delta_0}
\|\bar{u}_1^{\varepsilon,j_1}\|_{-1/2-\delta/2}\bar{C}_W^\varepsilon\\&+\sum_{j,i_2=1}^3\|\bar{u}_3^{\varepsilon,i_1}+\bar{u}_4^{\varepsilon,i_1}\|_{1/2-\delta_0}
\|\pi_{0,\diamond}(P^{i_2i_1}D_{j_1}\bar{K}^{\varepsilon,j_1},\bar{u}_1^{\varepsilon,j})\|_{-\delta}
\\&+\sum_{j,i_2=1}^3\|\bar{u}_3^{\varepsilon,j_1}+\bar{u}_4^{\varepsilon,j_1}\|_{1/2-\delta_0}\|\pi_{0,\diamond}(P^{i_2i_1}D_{j_1}\bar{K}^{\varepsilon,i_1},
\bar{u}_1^{\varepsilon,j})\|_{-\delta}+\|\bar{u}^{\varepsilon,\sharp,i_1}\|_{1/2+\beta}\|\bar{u}_1^{\varepsilon,j_1}\|_{-1/2-\delta/2}\bigg]\\&+[\delta C_W^{\varepsilon}+\varepsilon^{\kappa/2}(C_W^{\varepsilon}+\bar{C}^\varepsilon_W+1)](1+(C_W^{\varepsilon})^2+(\bar{C}^\varepsilon_W)^2)
(1+\|\bar{u}_4^{\varepsilon}\|_{\frac{1}{2}-\delta_0}+\|\bar{u}^{\varepsilon,\sharp}\|_{\frac{1}{2}+\beta}
+\|u_4^{\varepsilon}\|_{\frac{1}{2}-\delta_0})\\&+\|u_4^{\varepsilon}-\bar{u}_4^{\varepsilon}\|_{1/2-\delta_0}[(C_W^{\varepsilon})^2+\bar{C}^\varepsilon_W]
+\|u^{\varepsilon,\sharp}-\bar{u}^{\varepsilon,\sharp}\|_{\frac{1}{2}+\beta}C_W^{\varepsilon}
\\\lesssim &(\delta C_W^{\varepsilon}+\varepsilon^{\kappa/2}(C_W^{\varepsilon}+\bar{C}^\varepsilon_W+1))(1+(C_W^{\varepsilon})^2+(\bar{C}^\varepsilon_W)^2)(1+\|\bar{u}^{\varepsilon,\sharp}\|_{1/2+\beta}+\|\bar{u}_4^{\varepsilon}\|_{1/2-\delta_0}+\|u_4^{\varepsilon}\|_{1/2-\delta_0})
\\&+\|u_4^\varepsilon-\bar{u}_4^\varepsilon\|_{1/2-\delta_0}(1+\bar{C}^\varepsilon_W
+(C_W^{\varepsilon})^2)+\|u_4^\varepsilon-\bar{u}_4^\varepsilon\|_{\delta}(\|u_4^{\varepsilon}\|_{\delta}+\|\bar{u}_4^{\varepsilon}\|_{\delta})
+\varepsilon^{\kappa/2}\|\bar{u}_4^\varepsilon\|_\delta^2\\&+C_W^{\varepsilon}\|u^{\varepsilon,\sharp}-\bar{u}^{\varepsilon,\sharp}\|_{1/2+\beta},\endaligned$$
where we used $\frac{1}{2}-\delta-\delta_0\geq0$.$\hfill\Box$
\vskip.10in
\no\textbf{ Estimate of $F^\varepsilon-\bar{F}^\varepsilon$}

Now we consider $F^\varepsilon-\bar{F}^\varepsilon$:
 $$\aligned&\|F^\varepsilon-\bar{F}^\varepsilon\|_{1/2+\beta}\\\lesssim& \sum_{i_1,j_1=1}^3\|\int_0^tP_{t-s}^\varepsilon P^{ii_1}D^\varepsilon_{j_1}\pi_<(u_3^{\varepsilon,i_1}(s)+u_4^{\varepsilon,i_1}(s)-(u_3^{\varepsilon,i_1}(t)+u_4^{\varepsilon,i_1}(t)),u_1^{\varepsilon,j_1}(s))ds
 \\&-\int_0^tP_{t-s}P^{ii_1}D_{j_1}\pi_<(\bar{u}_3^{\varepsilon,i_1}(s)+\bar{u}_4^{\varepsilon,i_1}(s)-(\bar{u}_3^{\varepsilon,i_1}(t)+\bar{u}_4^{\varepsilon,i_1}(t)),\bar{u}_1^{\varepsilon,j_1}(s))ds\|_{1/2+\beta}
\endaligned$$
$$\aligned&+\sum_{i_1,j_1=1}^3\|\int_0^tP_{t-s}^\varepsilon P^{ii_1}D^\varepsilon_{j_1}\pi_<(u_3^{\varepsilon,i_1}(t)+u_4^{\varepsilon,i_1}(t),u_1^{\varepsilon,j_1}(s))ds
\\&-\int_0^tP_{t-s}P^{ii_1}D_{j_1}\pi_<(\bar{u}_3^{\varepsilon,i_1}(t)+\bar{u}_4^{\varepsilon,i_1}(t),\bar{u}_1^{\varepsilon,j_1}(s))ds\\&-P^{ii_1}D^\varepsilon_{j_1}\pi_<(u_3^{\varepsilon,i_1}+u_4^{\varepsilon,i_1},\int_0^tP_{t-s}^\varepsilon u_1^{\varepsilon,j_1}ds)+P^{ii_1}D_{j_1}\pi_<(\bar{u}_3^{\varepsilon,i_1}+\bar{u}_4^{\varepsilon,i_1},\int_0^tP_{t-s}\bar{u}_1^{\varepsilon,j_1}ds)\|_{1/2+\beta}\\&+\sum_{i_1,j_1=1}^3\|\int_0^tP_{t-s}^\varepsilon P^{ii_1}D^\varepsilon_{j_1}\pi_<(u_3^{\varepsilon,j_1}(s)+u_4^{\varepsilon,j_1}(s)-(u_3^{\varepsilon,j_1}(t)+u_4^{\varepsilon,j_1}(t)),u_1^{\varepsilon,i_1}(s))ds
\\&-\int_0^tP_{t-s}P^{ii_1}D_{j_1}\pi_<(\bar{u}_3^{\varepsilon,j_1}(s)+\bar{u}_4^{\varepsilon,j_1}(s)-(\bar{u}_3^{\varepsilon,j_1}(t)+\bar{u}_4^{\varepsilon,j_1}
(t)),\bar{u}_1^{\varepsilon,i_1}(s))ds\|_{1/2+\beta}
\\&+\sum_{i_1,j_1=1}^3\|\int_0^tP_{t-s}^\varepsilon P^{ii_1}D^\varepsilon_{j_1}\pi_<(u_3^{\varepsilon,j_1}(t)+u_4^{\varepsilon,j_1}(t),u_1^{\varepsilon,i_1}(s))ds\\&-\int_0^tP_{t-s}P^{ii_1}D_{j_1}
\pi_<(\bar{u}_3^{\varepsilon,j_1}(t)+\bar{u}_4^{\varepsilon,j_1}(t),\bar{u}_1^{\varepsilon,i_1}(s))ds\\&-P^{ii_1}D^\varepsilon_{j_1}\pi_<(u_3^{\varepsilon,j_1}
+u_4^{\varepsilon,j_1},\int_0^tP_{t-s}^\varepsilon u_1^{\varepsilon,i_1}ds)+P^{ii_1}D_{j_1}\pi_<(\bar{u}_3^{\varepsilon,j_1}+\bar{u}_4^{\varepsilon,j_1},\int_0^tP_{t-s}\bar{u}_1^{\varepsilon,j_1}ds)\|_{1/2+\beta}
\\=&I_1+I_2+I_3+I_4.\endaligned$$
It is sufficient to estimate $I_1, I_2$ and the desired estimate for $I_3$ and $I_4$ can be obtained similarly.
Since by the assumption  on $h$ we know that $\textrm{supp}\mathcal{F}u_1^\varepsilon\subset \{\xi:|\xi|\leq L_0/2\varepsilon\}$, which by Lemmas 3.5-3.8 deduce that for $\frac{1}{2}>\beta+\frac{\delta}{2}+\delta_0+\kappa$ $$I_2\lesssim t^{\frac{1}{4}-\frac{\delta_0+\beta+\kappa+\frac{\delta}{2}}{2}}(\|u_3^\varepsilon+u_4^\varepsilon-\bar{u}_3^\varepsilon-\bar{u}_4^\varepsilon\|_{1/2-\delta_0}C_W^{\varepsilon}+\varepsilon^{\kappa/2}\|\bar{u}_3^\varepsilon+\bar{u}_4^\varepsilon\|_{1/2-\delta_0}C_W^{\varepsilon}+\|\bar{u}_3^\varepsilon+\bar{u}_4^\varepsilon\|_{1/2-\delta_0}\delta C_W^{\varepsilon}).\eqno(3.35)$$
For $I_1$ one has the following estimate
$$\aligned I_1\lesssim&\int_0^t\varepsilon^{\kappa/2}(t-s)^{-1-\frac{\delta/2+\beta+\kappa}{2}}\|u_1^\varepsilon(s)\|_{-1/2-\delta/2}\|u_3^\varepsilon(t)+u^\varepsilon_4(t)-u^\varepsilon_3(s)-u^\varepsilon_4(s)\|_{\kappa/2}ds
+\int_0^t(t-s)^{-1-\frac{\delta/2+\beta+\kappa}{2}}\\&\|u_1^\varepsilon(s)\|_{-1/2-\delta/2}\|u_3^\varepsilon(t)+u^\varepsilon_4(t)-u^\varepsilon_3(s)-u^\varepsilon_4(s)-(\bar{u}_3^\varepsilon(t)+\bar{u}^\varepsilon_4(t)-\bar{u}^\varepsilon_3(s)-\bar{u}^\varepsilon_4(s))\|_{\kappa/2}ds
\\&+\int_0^t(t-s)^{-1-\frac{\delta/2+\beta+\kappa}{2}}\|u_1^\varepsilon(s)-\bar{u}_1^\varepsilon(s)\|_{-1/2-\delta/2}\|\bar{u}_3^\varepsilon(t)
+\bar{u}^\varepsilon_4(t)-\bar{u}^\varepsilon_3(s)-\bar{u}^\varepsilon_4(s)\|_{\kappa/2}ds.\endaligned$$
Moreover Lemmas 3.2, 3.3 and 3.4 imply that for $t>s>0$
$$\aligned &\|u_3^\varepsilon(t)+u^\varepsilon_4(t)-u^\varepsilon_3(s)-u^\varepsilon_4(s)\|_{\kappa/2}\\\lesssim
&\|(P_t^\varepsilon-P_s^\varepsilon)(Pu_0^\varepsilon-u_1^\varepsilon(0))\|_{\kappa/2}+\|\int_0^s(P_{t-r}^\varepsilon-P_{s-r}^\varepsilon)
G^\varepsilon(r)dr\|_{\kappa/2}+\|\int_s^tP_{t-r}^\varepsilon G^\varepsilon(r)dr\|_{\kappa/2}
\\\lesssim &(t-s)^{b_0}s^{-(z+2\kappa+2b_0)/2}\|u_0^\varepsilon-u_1^\varepsilon(0)\|_{-z}+(t-s)^{b/2}\int_0^s(s-r)^{-(3/2+\delta/2+2\kappa+b)/2}\|G^\varepsilon(r)\|_{-3/2-\delta/2-\kappa/2}dr
\\&+(t-s)^{b_1}(\int_s^t\|P_{t-r}^\varepsilon G^\varepsilon(r)\|_{\kappa/2}^{\frac{1}{1-b_1}}dr)^{1-b_1}
\endaligned$$
$$\aligned\lesssim &(t-s)^{b_0}s^{-(z+2\kappa+2b_0)/2}\|u_0^\varepsilon-u_1^\varepsilon(0)\|_{-z}+(t-s)^{b/2}\int_0^s(s-r)^{-(3/2+\delta/2+2\kappa+b)/2}\|G^\varepsilon(r)\|_{-3/2-\delta/2-\kappa/2}dr
\\&+(t-s)^{b_1}(\int_s^t(t-r)^{-\frac{3/2+\delta/2+2\kappa}{2(1-b_1)}}\| G^\varepsilon(r)\|_{-3/2-\delta/2-\kappa/2}^{\frac{1}{1-b_1}}dr)^{1-b_1},\endaligned$$
and $$\aligned &\|u_3^\varepsilon(t)+u^\varepsilon_4(t)-u^\varepsilon_3(s)-u^\varepsilon_4(s)-(\bar{u}_3^\varepsilon(t)+\bar{u}^\varepsilon_4(t)-\bar{u}^\varepsilon_3(s)-\bar{u}^\varepsilon_4(s))\|_{\kappa/2}\\\lesssim
&\|P_t^\varepsilon-P_s^\varepsilon-(P_t-P_s)(Pu_0^\varepsilon-u_1^\varepsilon(0))\|_{\kappa/2}+\|(P_t-P_s)(Pu_0^\varepsilon-u_1^\varepsilon(0)-Pu_0+\bar{u}_1^\varepsilon(0))\|_{\kappa/2}\\&+\|\int_0^s[(P_{t-r}^\varepsilon-P_{s-r}^\varepsilon)G^\varepsilon(r)-(P_{t-r}-P_{s-r})\bar{G}^\varepsilon(r)]dr\|_{\kappa/2}+\|\int_s^t(P_{t-r}^\varepsilon G^\varepsilon(r)-P_{t-r}\bar{G}^\varepsilon(r))dr\|_{\kappa/2}
\\\lesssim &[((t-s)^{2b_0}s^{-(z+2\kappa+4b_0)/2})\wedge(\varepsilon^{\kappa}s^{-(z+2\kappa)/2})]\|u_0^\varepsilon-u_1^\varepsilon(0)\|_{-z}+\|(P_t-P_s)(Pu_0^\varepsilon-u_1^\varepsilon(0)-Pu_0+\bar{u}_1^\varepsilon(0))\|_{\kappa/2}\\&+\|\int_0^s[(P_{t-r}^\varepsilon-P_{s-r}^\varepsilon)G^\varepsilon(r)-(P_{t-r}-P_{s-r})\bar{G}^\varepsilon(r)]dr\|_{\kappa/2}+\|\int_s^t(P_{t-r}^\varepsilon G^\varepsilon(r)-P_{t-r}\bar{G}^\varepsilon(r))dr\|_{\kappa/2}
\\\lesssim &(t-s)^{b_0}s^{-(z+2\kappa+2b_0)/2}(\varepsilon^{\kappa/2}\|u_0^\varepsilon-u_1^\varepsilon(0)\|_{-z}+\|u_0^\varepsilon-u_0-u_1^\varepsilon(0)+\bar{u}_1^\varepsilon(0)\|_{-z})\\&+(t-s)^{b/2}\int_0^s(s-r)^{-(3/2+\delta/2+2\kappa+b)/2}(\varepsilon^{\kappa/2}\|G^\varepsilon(r)\|_{-3/2-\delta/2-\kappa/2}+\|G^\varepsilon(r)-\bar{G}^\varepsilon(r)\|_{-3/2-\delta/2-\kappa})dr
\\&+(t-s)^{b_1}(\int_s^t\|P_{t-r}^\varepsilon G^\varepsilon(r)-P_{t-r}\bar{G}^\varepsilon(r)\|_{\kappa/2}^{\frac{1}{1-b_1}}dr)^{1-b_1}
\\\lesssim &(t-s)^{b_0}s^{-(z+2\kappa+2b_0)/2}(\varepsilon^{\kappa/2}\|u_0^\varepsilon-u_1^\varepsilon(0)\|_{-z}+\|u_0^\varepsilon-u_0-u_1^\varepsilon(0)+\bar{u}_1^\varepsilon(0)\|_{-z})\\&+(t-s)^{b/2}\int_0^s(s-r)^{-(3/2+\delta/2+2\kappa+b)/2}(\varepsilon^{\kappa/2}\|G^\varepsilon(r)\|_{-3/2-\delta/2-\kappa/2}+\|G^\varepsilon(r)-\bar{G}^\varepsilon(r)\|_{-3/2-\delta/2-\kappa})dr
\\&+(t-s)^{b_1}(\int_s^t(t-r)^{-\frac{3/2+\delta/2+2\kappa}{2(1-b_1)}}(\varepsilon^{\kappa/2}\|G^\varepsilon(r)\|_{-3/2-\delta/2-\kappa/2}+\|G^\varepsilon(r)-\bar{G}^\varepsilon(r)\|_{-3/2-\delta/2-\kappa})^{\frac{1}{1-b_1}}dr)^{1-b_1},\endaligned$$
where $\delta/2+\beta+2\kappa<2b_0<2-z-2\kappa$, $\delta/2+\beta+2\kappa<b<1/2-2\kappa-\delta/2$, $\frac{1}{2}(\delta/2+\beta+2\kappa)<b_1<[1-(\delta+z+2\kappa)]\wedge \frac{1}{2}(1/2-\delta/2-2\kappa)$ and
$$\aligned G^\varepsilon=&\sum_{i_1,j=1}^3P^{ii_1}D_j^\varepsilon[u_1^{\varepsilon,i_1}\diamond u_2^{\varepsilon,j}+u_2^{\varepsilon,i_1}\diamond u_1^{\varepsilon,j}+u_1^{\varepsilon,i_1}\diamond (u_3^{\varepsilon,j}+u_4^{\varepsilon,j})+(u_3^{\varepsilon,i_1} +u_4^{\varepsilon,i_1})\diamond u_1^{\varepsilon,j}\\&+u_2^{\varepsilon,i_1}\diamond u_2^{\varepsilon,j}+u_2^{\varepsilon,i_1}(u_3^{\varepsilon,j}+u_4^{\varepsilon,j})+u_2^{\varepsilon,j}(u_3^{\varepsilon,i_1}+u_4^{\varepsilon,i_1})
+(u_3^{\varepsilon,i_1}+u_4^{\varepsilon,i_1})(u_3^{\varepsilon,j}+u_4^{\varepsilon,j})].\endaligned$$ Thus we obtain that
$$\aligned I_1\lesssim&(\varepsilon^{\kappa/2}C_W^{\varepsilon}+\delta C_W^{\varepsilon})\bigg{[}\int_0^t(t-s)^{-1-\frac{\delta/2+\beta+2\kappa}{2}+b_0}s^{-(z+2\kappa+2b_0)/2}ds(\|u_0^\varepsilon-u_1^\varepsilon(0)\|_{-z}+\|u_0-\bar{u}_1^\varepsilon(0)\|_{-z})
\\&+\int_0^t(t-s)^{-1-\frac{\delta/2+\beta+2\kappa}{2}+b/2}\int_0^s(s-r)^{-(3/2+\delta/2+2\kappa+b)/2}
(\|G^\varepsilon(r)\|_{-3/2-\delta/2-\kappa/2}\\&+\|\bar{G}^\varepsilon(r)\|_{-3/2-\delta/2-\kappa/2})drds
+\int_0^t(t-s)^{-1-\frac{\delta/2+\beta+2\kappa}{2}+b_1}(\int_s^t(t-r)^{-\frac{3/2+\delta/2+2\kappa}{2(1-b_1)}}(\| G^\varepsilon(r)\|_{-3/2+\delta/2+\kappa/2}^{\frac{1}{1-b_1}}\\&+\| \bar{G}^\varepsilon(r)\|_{-3/2-\delta/2-\kappa/2}^{\frac{1}{1-b_1}})dr)^{1-b_1}ds\bigg{]}
\endaligned$$
$$\aligned&+C_W^{\varepsilon}(\int_0^t(t-s)^{-1-\frac{\delta/2+\beta+2\kappa}{2}+b_0}s^{-(z+2\kappa+2b_0)/2}ds\|u_0^\varepsilon-u_1^\varepsilon(0)-(u_0-\bar{u}_1^\varepsilon(0))\|_{-z}\\
&+\int_0^t(t-s)^{-1-\frac{\delta/2+\beta+2\kappa}{2}+b/2}\int_0^s(s-r)^{-(3/2+\delta/2+2\kappa+b)/2}
\|G^\varepsilon(r)-\bar{G}^\varepsilon(r)\|_{-3/2-\delta/2-\kappa}dsdr
\\&+\int_0^t(t-s)^{-1-\frac{\delta/2+\beta+2\kappa}{2}+b_1}(\int_s^t(t-r)^{-\frac{3/2+\delta/2+2\kappa}{2(1-b_1)}}\| G^\varepsilon(r)- \bar{G}^\varepsilon(r)\|_{-3/2-\delta/2-\kappa}^{\frac{1}{1-b_1}}dr)^{1-b_1}ds)
\\\lesssim &(\varepsilon^{\kappa/2}C_W^{\varepsilon}+\delta C_W^{\varepsilon})\bigg[\int_0^t(t-s)^{-1-\frac{\delta/2+\beta+2\kappa}{2}+b_0}s^{-(z+2\kappa+2b_0)/2}ds(\|u_0^\varepsilon-u_1^\varepsilon(0)\|_{-z}
+\|u_0-\bar{u}_1^\varepsilon(0)\|_{-z})
\\&+\int_0^t(t-r)^{-\frac{3}{4}-\frac{\delta+\beta+4\kappa}{2}}\int_0^1(1-l)^{-1-\frac{\delta/2+\beta+2\kappa}{2}+b/2}l^{-(3/2+\delta/2+2\kappa+b)/2}dl
\\&(\|G^\varepsilon(r)\|_{-3/2-\delta/2-\kappa/2}+\|\bar{G}^\varepsilon(r)\|_{-3/2-\delta/2-\kappa/2})dr
\\&+(\int_0^t(t-s)^{-1-\frac{\delta/2+\beta+2\kappa}{2}+b_1}ds)^{b_1}(\int_0^t\int_0^r(t-s)^{-1-\frac{\delta/2+\beta+2\kappa}{2}+b_1}
(t-r)^{-\frac{3/2+\delta/2+2\kappa}{2(1-b_1)}}(\| G^\varepsilon(r)\|_{-3/2-\delta/2-\kappa/2}^{\frac{1}{1-b_1}}\\&+\| \bar{G}^\varepsilon(r)\|_{-3/2-\delta/2-\kappa/2}^{\frac{1}{1-b_1}})dsdr)^{1-b_1}\bigg]
\\&+C_W^{\varepsilon}(\int_0^t(t-s)^{-1-\frac{\delta/2+\beta+2\kappa}{2}+b_0}s^{-(z+2\kappa+2b_0)/2}ds\|u_0^\varepsilon-u_1^\varepsilon(0)-(u_0-\bar{u}_1^\varepsilon(0))\|_{-z}
\\
&+\int_0^t(t-r)^{-\frac{3}{4}-\frac{\delta+\beta+4\kappa}{2}}\int_0^1(1-l)^{-1-\frac{\delta/2+\beta+2\kappa}{2}+b/2}l^{-(3/2+\delta/2+2\kappa+b)/2}dl\|G^\varepsilon(r)-\bar{G}^\varepsilon(r)\|_{-3/2-\delta/2-\kappa}dr
\\&+(\int_0^t(t-s)^{-1-\frac{\delta/2+\beta+2\kappa}{2}+b_1}ds)^{b_1}(\int_0^t\int_0^r(t-s)^{-1-\frac{\delta/2+\beta+2\kappa}{2}+b_1}(t-r)^{-\frac{3/2+\delta/2+2\kappa}{2(1-b_1)}}\\&\| G^\varepsilon(r)- \bar{G}^\varepsilon(r)\|_{-3/2-\delta/2-\kappa}^{\frac{1}{1-b_1}}ds dr)^{1-b_1}).
\endaligned\eqno(3.36)$$
Hence (3.35) and (3.36) yield that
$$\aligned&\|F^\varepsilon(t)-\bar{F}^\varepsilon(t)\|_{1/2+\beta}\\\lesssim& t^{\frac{1}{4}-\frac{\delta_0+\beta+\kappa+\frac{\delta}{2}}{2}}[\|u_3^\varepsilon+u_4^\varepsilon-\bar{u}_3^\varepsilon-\bar{u}_4^\varepsilon\|_{1/2-\delta_0}C_W^{\varepsilon}+(\varepsilon^{\kappa/2}C_W^{\varepsilon}+\delta C_W^{\varepsilon})\|\bar{u}_3^\varepsilon+\bar{u}_4^\varepsilon\|_{1/2-\delta_0}]\\&+(\varepsilon^{\kappa/2}C_W^{\varepsilon}+\delta C_W^{\varepsilon})\bigg{[}t^{-\frac{\delta/2+\beta+z}{2}-2\kappa}(\|u_0^\varepsilon-u_1^\varepsilon(0)\|_{-z}+\|u_0-\bar{u}_1^\varepsilon(0)\|_{-z})
\\&+\int_0^t(t-r)^{-\frac{3}{4}-\frac{\delta+\beta+4\kappa}{2}}(\|G^\varepsilon(r)\|_{-3/2-\delta/2-\kappa/2}+\|\bar{G}^\varepsilon(r)\|_{-3/2-\delta/2-\kappa/2})dr
\\&+(\int_0^t(t-r)^{-\frac{3/2+\delta/2+2\kappa}{2(1-b_1)}}(\| G^\varepsilon(r)\|_{-3/2-\delta/2-\kappa/2}^{\frac{1}{1-b_1}}+\| \bar{G}^\varepsilon(r)\|_{-3/2-\delta/2-\kappa/2}^{\frac{1}{1-b_1}})dr)^{1-b_1}\bigg{]}
\\&+t^{-\frac{\delta/2+\beta+z}{2}-2\kappa}\|u_0^\varepsilon-u_1^\varepsilon(0)-(u_0-\bar{u}_1^\varepsilon(0))\|_{-z}C_W^{\varepsilon}+C_W^{\varepsilon}\int_0^t(t-r)^{-\frac{3}{4}-\frac{\delta+\beta+4\kappa}{2}}\|G^\varepsilon(r)-\bar{G}^\varepsilon(r)\|_{-3/2-\delta/2-\kappa}dr
\\&+(\int_0^t(t-r)^{-\frac{3/2+\delta/2+2\kappa}{2(1-b_1)}}\| G^\varepsilon(r)- \bar{G}^\varepsilon(r)\|_{-3/2-\delta/2-\kappa}^{\frac{1}{1-b_1}} dr)^{1-b_1}C_W^{\varepsilon},\endaligned\eqno(3.37)$$
and by a similar argument as (3.6), (3.28) we deduce that
$$\|G^\varepsilon\|_{-3/2-\delta/2-\kappa/2}\lesssim(1+({C}^\varepsilon_W)^3)(1+\|{u}^{\varepsilon,\sharp}\|_{1/2+\beta}+\|{u}_4^{\varepsilon}\|_{1/2-\delta_0})
+\|{u}_4^\varepsilon\|_\delta^2,$$
and
$$\aligned &\| G^\varepsilon(r)- \bar{G}^\varepsilon(r)\|_{-3/2-\delta/2-\kappa}\\\lesssim&(\delta C_W^{\varepsilon}+\varepsilon^{\kappa/2}(C_W^{\varepsilon}+\bar{C}^\varepsilon_W+1))(1+(C_W^{\varepsilon})^2+(\bar{C}^\varepsilon_W)^2)(1+\|\bar{u}^{\varepsilon,\sharp}\|_{1/2+\beta}+\|\bar{u}_4^{\varepsilon}\|_{1/2-\delta_0}+\|u_4^{\varepsilon}\|_{1/2-\delta_0})
\\&+\|u_4^\varepsilon-\bar{u}_4^\varepsilon\|_{1/2-\delta_0}(1+\bar{C}^\varepsilon_W+(C_W^{\varepsilon})^2)+\|u_4^\varepsilon-\bar{u}_4^\varepsilon\|_{\delta}(\|u_4^{\varepsilon}\|_{\delta}+\|\bar{u}_4^{\varepsilon}\|_{\delta})+\varepsilon^{\kappa/2}\|\bar{u}_4^\varepsilon\|_\delta^2\\&+C_W^{\varepsilon}\|u^{\varepsilon,\sharp}-\bar{u}^{\varepsilon,\sharp}\|_{1/2+\beta}.\endaligned$$
A similar argument as (3.37) deduce that
$$\aligned&\|F^\varepsilon(t)-\bar{F}^\varepsilon(t)\|_{\delta}\\\lesssim& t^{\frac{1}{4}-\frac{\delta}{4}-\frac{\kappa}{2}}(\|u_3^\varepsilon(t)+u_4^\varepsilon(t)-\bar{u}_3^\varepsilon(t)-\bar{u}_4^\varepsilon(t)\|_{\delta}
C_W^{\varepsilon}+(\varepsilon^{\kappa/2}C_W^{\varepsilon}+\delta C_W^{\varepsilon})\|\bar{u}_3^\varepsilon(t)+\bar{u}_4^\varepsilon(t)\|_{\delta})\\&+(\varepsilon^{\kappa/2}C_W^{\varepsilon}+\delta C_W^{\varepsilon})\bigg[t^{\frac{1}{4}-\frac{3\delta}{4}-\frac{z}{2}-\frac{\kappa}{2}}(\|u_0^\varepsilon-u_1^\varepsilon(0)\|_{-z}+\|u_0-\bar{u}_1^\varepsilon(0)\|_{-z})
\\&+\int_0^t(t-r)^{-\frac{1+2\delta+\kappa}{2}}(\|G^\varepsilon(r)\|_{-3/2-\delta/2-\kappa/2}+\|\bar{G}^\varepsilon(r)\|_{-3/2-\delta/2-\kappa/2})dr
\\&+(\int_0^t(t-r)^{-\frac{1+2\delta+2\kappa}{2}}(\| G^\varepsilon(r)\|_{-3/2-\delta/2-\kappa/2}+\| \bar{G}^\varepsilon(r)\|_{-3/2-\delta/2-\kappa/2})dr)\bigg]
\\&+t^{\frac{1}{4}-\frac{3\delta}{4}-\frac{z}{2}-\frac{\kappa}{2}}\|u_0^\varepsilon-u_1^\varepsilon(0)-(u_0-\bar{u}_1^\varepsilon(0))\|_{-z}
+\int_0^t(t-r)^{-\frac{1+2\delta+\kappa}{2}}\|G^\varepsilon(r)-\bar{G}^\varepsilon(r)\|_{-3/2-\delta/2-\kappa}dr
\\&+\int_0^t(t-r)^{-\frac{1+2\delta+2\kappa}{2}}\| G^\varepsilon(r)- \bar{G}^\varepsilon(r)\|_{-3/2-\delta/2-\kappa} dr,\endaligned\eqno(3.38)$$
where the only difference is that we use $\|u_3(t)+u_4(t)-u_3(s)-u_4(s)\|_{-1/2+3\delta/2+\kappa}$ to control $\|F^\varepsilon-\bar{F}^\varepsilon\|_\delta$.

\vskip.10in
\no\textbf{Estimate of $u_4^{\varepsilon,i}-\bar{u}_4^{ \varepsilon,i}$}
\vskip.10in

By paracontrolled ansatz and Lemma 3.2 we get that
$$\aligned\|u_4^{\varepsilon,i}-\bar{u}_4^{ \varepsilon,i}\|_{1/2-\delta_0}\lesssim& \sum_{i_1,j=1}^3(t^{\delta/4}\|u_3^{ \varepsilon,i_1}+u_4^{ \varepsilon,i_1}-\bar{u}_3^{ \varepsilon,i_1}-\bar{u}_4^{ \varepsilon,i_1}\|_{1/2-\delta_0}C_W^{\varepsilon}\\&+\|\bar{u}_3^{ \varepsilon,i_1}+\bar{u}_4^{ \varepsilon,i_1}\|_{1/2-\delta_0}\|K^{ \varepsilon,j}-\bar{K}^{ \varepsilon,j}\|_{3/2-\delta}\\&+\varepsilon^{\kappa/2}\|\bar{u}_3^{ \varepsilon,i_1}+\bar{u}_4^{ \varepsilon,i_1}\|_{1/2-\delta_0}\|\bar{K}^{ \varepsilon,j}\|_{3/2-\delta})+\|u^{ \varepsilon,\sharp,i}-\bar{u}^{ \varepsilon,\sharp,i}\|_{1/2-\delta_0},\endaligned$$
which by (3.26) shows that for $t$ small enough (only depending on $C_W^{\varepsilon}$)
$$\aligned\|u_4^{ \varepsilon,i}-\bar{u}_4^{ \varepsilon,i}\|_{1/2-\delta_0}\lesssim& \delta C_W^{\varepsilon}(C_W^{\varepsilon}+\bar{C}_W^{\varepsilon}+\|\bar{u}_4^{ \varepsilon}\|_{1/2-\delta_0})+\varepsilon^{\kappa/2}\bar{C}_W^{\varepsilon}(\bar{C}_W^{\varepsilon}+\|\bar{u}_4^{ \varepsilon}\|_{1/2-\delta_0})\\&+\|u^{ \varepsilon,\sharp,i}-\bar{u}^{ \varepsilon,\sharp,i}\|_{1/2-\delta_0}.\endaligned\eqno(3.39)$$
Similarly, one has for $t$ small enough (only depending on $C_W^{\varepsilon}$)
$$\aligned\|u_4^{ \varepsilon,i}-\bar{u}_4^{ \varepsilon,i}\|_{\delta}\lesssim& \delta C_W^{\varepsilon}(C_W^{\varepsilon}+\bar{C}_W^{\varepsilon}+\|\bar{u}_4^{ \varepsilon}\|_{\delta})+\varepsilon^{\kappa/2}\bar{C}_W^{\varepsilon}(\bar{C}_W^{\varepsilon}+\|\bar{u}_4^{ \varepsilon}\|_{\delta})+\|u^{ \varepsilon,\sharp,i}-\bar{u}^{ \varepsilon,\sharp,i}\|_{\delta}.\endaligned\eqno(3.40)$$
Then by (3.20), (3.21), (3.24) (3.25), (3.28) and (3.37) we get that for $\delta+z+\kappa<1$, $t$ small enough and $t\leq \tau_L\wedge \tau_{L_1}^\varepsilon\wedge \rho_{L_2}^\varepsilon\wedge\bar{\rho}_{L_3}^\varepsilon$ with $L, L_1, L_2, L_3\geq0$
$$\aligned &t^{\delta+z+\kappa}\| u^{ \varepsilon,\sharp}(t)-\bar{u}^{ \varepsilon,\sharp}(t)\|_{1/2+\beta}\\\lesssim & (\varepsilon^{\kappa/2}(C_W^{\varepsilon}+1)+\delta C_W^{\varepsilon})(\|Pu_0^\varepsilon-u_1^{ \varepsilon}(0)\|_{-z}+\|Pu_0-\bar{u}_1^{ \varepsilon}(0)\|_{-z})\\&+(C_W^{\varepsilon}+1)\|Pu_0^\varepsilon-u_1^{ \varepsilon}(0)-(Pu_0-\bar{u}_1^{ \varepsilon}(0))\|_{-z}\\&+t^{\delta+z+\kappa}\int_0^t(t-s)^{-3/4-\delta-\beta/2-3\kappa/2}s^{-(\delta+z+\kappa)}s^{\delta+z+\kappa}(\varepsilon^{\kappa/2}\|\bar{\phi}^{ \varepsilon,\sharp}\|_{-1-2\delta-2\kappa}+\|\phi^{ \varepsilon,\sharp}-\bar{\phi}^{ \varepsilon,\sharp}\|_{-1-2\delta-2\kappa})ds
\\&+C_W^{ \varepsilon} t^{\frac{1}{4}-\frac{\delta_0+\beta+\kappa+\frac{\delta}{2}}{2}}t^{\delta+z+\kappa}\|u_4^{ \varepsilon}-\bar{u}_4^{ \varepsilon}\|_{1/2-\delta_0}
+(\varepsilon^{\kappa/2}(C_W^{ \varepsilon}+\bar{C}_W^{ \varepsilon}+1)+\delta C_W^{ \varepsilon})C(C_W^{ \varepsilon},\bar{C}_W^{ \varepsilon},\|u_0^{ \varepsilon}\|_{-z},\|u_0\|_{-z})
\\&+C(C_W^{ \varepsilon},\bar{C}_W^{ \varepsilon})\int_0^t(t-r)^{-\frac{3}{4}-\frac{\delta+\beta+4\kappa}{2}}(\|u_4^{ \varepsilon}-\bar{u}_4^{ \varepsilon}\|_{1/2-\delta_0}+\|u_4^{ \varepsilon}-\bar{u}_4^{ \varepsilon}\|_{\delta}(\|u_4^{ \varepsilon}\|_{\delta}+\|\bar{u}_4^{ \varepsilon}\|_{\delta})\\&+\|u^{ \varepsilon,\sharp}-\bar{u}^{ \varepsilon,\sharp}\|_{1/2+\beta})dr
+C(C_W^{ \varepsilon},\bar{C}_W^{ \varepsilon})(\int_0^t(t-r)^{-\frac{3/2+\delta/2+2\kappa}{2(1-b_1)}}\| G^{ \varepsilon}(r)- \bar{G}^{ \varepsilon}(r)\|_{-3/2-\delta/2-\kappa}^{\frac{1}{1-b_1}}dr )^{1-b_1}
\\\lesssim & (\varepsilon^{\kappa/2}(L_2+1)+\delta C_W^{\varepsilon})(\|u_0^\varepsilon-u_1^{ \varepsilon}(0)\|_{-z}+\|u_0-\bar{u}_1^{ \varepsilon}(0)\|_{-z})+(L_2+1)\|u_0^\varepsilon-u_1^{ \varepsilon}(0)-(u_0-\bar{u}_1^{ \varepsilon}(0))\|_{-z}\\&+t^{\delta+z+\kappa}C(L_2,L_3)\int_0^t(t-s)^{-3/4-\delta-\beta/2-3\kappa/2}s^{-(\delta+z+\kappa)}s^{\delta+z+\kappa}(\|u^{ \varepsilon,\sharp}-\bar{u}^{ \varepsilon,\sharp}\|_{1/2+\beta}\\&+\|u^{ \varepsilon,\sharp}-\bar{u}^{ \varepsilon,\sharp}\|_{\delta}(\|u^{ \varepsilon,\sharp}\|_{\delta}+\|\bar{u}^{ \varepsilon,\sharp}\|_{\delta}))ds
\\&+L_2t^{\frac{1}{4}-\frac{\delta_0+\beta+\kappa+\frac{\delta}{2}}{2}}t^{\delta+z+\kappa} \|u^{ \varepsilon,\sharp}-\bar{u}^{ \varepsilon,\sharp}\|_{1/2+\beta}
+(\varepsilon^{\kappa/2}(L_2+L_3+1)+\delta C_W^{ \varepsilon})C(L,L_1,L_2,L_3)
\\&+t^{\delta+z+\kappa}C(L_2, L_3)\int_0^t(t-r)^{-\frac{3}{4}-\frac{\delta+\beta+4\kappa}{2}}(\|u^{ \varepsilon,\sharp}-\bar{u}^{ \varepsilon,\sharp}\|_{1/2+\beta}+\|u^{ \varepsilon,\sharp}-\bar{u}^{ \varepsilon,\sharp}\|_{\delta}(\|u^{ \varepsilon,\sharp}\|_{\delta}+\|\bar{u}^{ \varepsilon,\sharp}\|_{\delta}))dr
\\&+t^{\delta+z+\kappa}C(L_2,L_3)\big{(}\int_0^t(t-r)^{-\frac{3/2+\delta/2+2\kappa}{2(1-b_1)}}[\|u^{ \varepsilon,\sharp}-\bar{u}^{ \varepsilon,\sharp}\|_{1/2+\beta}\\&+\|u^{ \varepsilon,\sharp}-\bar{u}^{ \varepsilon,\sharp}\|_{\delta}(\|u^{ \varepsilon,\sharp}\|_{\delta}+\|\bar{u}^{ \varepsilon,\sharp}\|_{\delta})]^{\frac{1}{1-b_1}}dr \big{)}^{1-b_1},\endaligned\eqno(3.41)$$
where we used the condition on $\beta$ to deduce $-3/4-\delta-\beta/2-3\kappa/2>-1$, $\frac{1/2+\beta+z}{2}\leq\delta+z$ and $\delta+z+\kappa-\frac{\delta/2+\beta+z}{2}-2\kappa\geq0$.
Hence for $T_1$ small enough and $T_1\leq \tau_L\wedge \tau_{L_1}^\varepsilon\wedge \rho_{L_2}^\varepsilon\wedge\bar{\rho}_{L_3}^\varepsilon$ there exists some $\theta>0$ such that
$$\aligned &\sup_{t\in[0,T_1]}t^{\delta+z+\kappa}\| u^{ \varepsilon,\sharp}(t)-\bar{u}^{ \varepsilon,\sharp}(t)\|_{1/2+\beta}\\\lesssim & (\varepsilon^{\kappa/2}(L_2+1)+\delta C_W^{\varepsilon})(\|u_0^{ \varepsilon}-u_1^{ \varepsilon}(0)\|_{-z}+\|u_0^{ \varepsilon}-\bar{u}_1^{ \varepsilon}(0)\|_{-z})\\&+(L_2+1)\|u_0^\varepsilon-\bar{u}_1^\varepsilon(0)-(u_0-\bar{u}_1^\varepsilon(0))\|_{-z}
+(\varepsilon^{\kappa/2}(L_2+L_3+1)+\delta C_W^{ \varepsilon})C(L,L_1,L_2,L_3)
\\&+T_1^\theta C(L,L_1,L_2,L_3)(\sup_{t\in[0,T_1]}t^{\delta+z+\kappa}\|u^{ \varepsilon,\sharp}-\bar{u}^{ \varepsilon,\sharp}\|_{1/2+\beta}+\sup_{t\in[0,T_1]}t^{\frac{\delta+z+\kappa}{2}}\|u^{ \varepsilon,\sharp}-\bar{u}^{ \varepsilon,\sharp}\|_{\delta})
,\endaligned\eqno(3.42)$$
By a similar argument as (3.42) (3.16) and using  (3.28) (3.38)-(3.40) we also obtain that for $0<5\kappa<1-z-4\delta$
$$\aligned &\sup_{t\in[0,T]}t^{\frac{\delta+z+\kappa}{2}}\|u^{ \varepsilon,\sharp}(t)-\bar{u}^{ \varepsilon,\sharp}(t)\|_{\delta}\\\lesssim & (\varepsilon^{\kappa/2}(L_2+1)+\delta C_W^{\varepsilon})(\|u_0^{ \varepsilon}-u_1^{ \varepsilon}(0)\|_{-z}+\|u_0-\bar{u}_1^{ \varepsilon}(0)\|_{-z})\\&+(L_2+1)\|Pu_0^{ \varepsilon}-u_1^{ \varepsilon}(0)-(Pu_0-\bar{u}_1^{ \varepsilon}(0))\|_{-z}+
(\varepsilon^{\kappa/2}(L_2+L_3+1)+\delta C_W^{ \varepsilon})C(L,L_1,L_2,L_3)
\\&+T_1^\theta C(L,L_1,L_2,L_3)(\sup_{t\in[0,T_1]}t^{\delta+z+\kappa}\|u^{ \varepsilon,\sharp}-\bar{u}^{ \varepsilon,\sharp}\|_{1/2+\beta}+\sup_{t\in[0,T_1]}t^{\frac{\delta+z+\kappa}{2}}\|u^{ \varepsilon,\sharp}-\bar{u}^{ \varepsilon,\sharp}\|_{\delta})
.\endaligned\eqno(3.43)$$
Combining (3.42,3.43) we deduce that for $T_1$ small enough and $T_1\leq \tau_L\wedge \tau_{L_1}^\varepsilon\wedge \rho_{L_2}^\varepsilon\wedge\bar{\rho}_{L_3}^\varepsilon$
$$\aligned &\sup_{t\in[0,T_1]}(t^{\delta+z+\kappa}\|u^{ \varepsilon,\sharp}-\bar{u}^{ \varepsilon,\sharp}\|_{1/2+\beta}+t^{\frac{\delta+z+\kappa}{2}}\|u^{ \varepsilon,\sharp}(t)-\bar{u}^{ \varepsilon,\sharp}(t)\|_{\delta})\\\lesssim & (L_2+1)\|u_0^{ \varepsilon}-u_0\|_{-z}+
(\varepsilon^{\kappa/2}(L_2+L_3+1)+\delta C_W^{ \varepsilon})C(L,L_1,L_2,L_3)
,\endaligned\eqno(3.44)$$
which by (3.28) (3.39), (3.40) implies that
$$\aligned&\sup_{t\in[0,T_1]}t^{\delta+z+\kappa}\|\phi^{ \varepsilon,\sharp}-\bar{\phi}^{ \varepsilon,\sharp}\|_{-1-2\delta-2\kappa}\\\lesssim & C(L,L_1,L_2,L_3)\|u_0^{ \varepsilon}-u_0\|_{-z}+
(\varepsilon^{\kappa/2}(L_2+L_3+1)+\delta C_W^{ \varepsilon})C(L,L_1,L_2,L_3).\endaligned\eqno(3.45)$$
Moreover by paracontrolled ansatz and  Lemmas 3.7, 3.8 one also has that
$$\aligned\|u_4^{ \varepsilon,i}(t)-\bar{u}_4^{ \varepsilon,i}(t)\|_{-z}\lesssim& \sum_{i_1,j=1}^3(t^{\delta/4}\|u_3^{ \varepsilon,i_1}+u_4^{ \varepsilon,i_1}-\bar{u}_3^{ \varepsilon,i_1}-\bar{u}_4^{ \varepsilon,i_1}\|_{-z}C_W^{\varepsilon}
\\&+\|\bar{u}_3^{ \varepsilon,i_1}+\bar{u}_4^{ \varepsilon,i_1}\|_{-z}\|K^{ \varepsilon,j}-\bar{K}^{ \varepsilon,j}\|_{3/2-\delta}
\\&+\varepsilon^{\kappa/2}\|\bar{u}_3^{ \varepsilon,i_1}+\bar{u}_4^{ \varepsilon,i_1}\|_{-z}\|\bar{K}^{ \varepsilon,j}\|_{3/2-\delta})+\|u^{ \varepsilon,\sharp,i}-\bar{u}^{ \varepsilon,\sharp,i}\|_{-z},\endaligned$$
which combining with Lemmas 2.6, 3.2, 3.3, (3.17), (3.25), (3.44), (3.45)  implies that for $T_1$ small enough, $T_1\leq \tau_L\wedge \tau_{L_1}^\varepsilon\wedge \rho_{L_2}^\varepsilon\wedge\bar{\rho}_{L_3}^\varepsilon$ and $t\in [0,T_1]$
$$\aligned\|u_4^{ \varepsilon,i}(t)-\bar{u}_4^{ \varepsilon,i}(t)\|_{-z}\lesssim& (\delta C_W^{ \varepsilon}+\varepsilon^{\kappa/2}) C(L,L_1,L_2,L_3)+\|u^{ \varepsilon,\sharp,i}(t)-\bar{u}^{ \varepsilon,\sharp,i}(t)\|_{-z}
\\\lesssim &(\delta C_W^{ \varepsilon}+\varepsilon^{\kappa/2}) C(L,L_1,L_2,L_3)+\varepsilon^{\kappa/2}\|Pu_0^{ \varepsilon}-u_1^{ \varepsilon}(0)\|_{-z}\\&+\|Pu_0^{ \varepsilon}-u_1^{ \varepsilon}(0)-(Pu_0-\bar{u}_1^{ \varepsilon}(0))\|_{-z}\\&+\int_0^t(t-s)^{-\frac{1+2\delta+3\kappa-z}{2}}s^{-(\delta+z+\kappa)}s^{\delta+z+\kappa}(\varepsilon^{\kappa/2}\|\bar{\phi}^{ \varepsilon,\sharp}(s)\|_{-1-2\delta-2\kappa}\\&+\|\phi^{ \varepsilon,\sharp}-\bar{\phi}^{ \varepsilon,\sharp}\|_{-1-2\delta-2\kappa})ds+\|F^\varepsilon-\bar{F}^\varepsilon\|_{-z} \\\lesssim&(\delta C_W^{ \varepsilon}+\varepsilon^{\kappa/2}+\|u_0^\varepsilon-u_0\|_{-z})C(L,L_1,L_2,L_3)+L_2t^{\frac{1}{4}-\frac{\delta}{4}}  \|u_4^{ \varepsilon,i}(t)-\bar{u}_4^{ \varepsilon,i}(t)\|_{-z}.\endaligned$$
Here in the last inequality we used  $$\aligned&\|\bar{F}^{ \varepsilon}(t)-F^{ \varepsilon}(t)\|_{-z}\\\lesssim& (\varepsilon^{\kappa/2}+\delta C_W^\varepsilon)(C(C^{ \varepsilon}_W)\int_0^t(t-s)^{-\frac{3/2+\delta/2-z+\kappa}{2}}s^{-\frac{\delta+\kappa+z}{2}}ds\sup_{s\in[0,t]}s^{\frac{\delta+\kappa+z}{2}}\|u^{ \varepsilon}_3+u^{ \varepsilon}_4\|_{\delta}
+t^{\frac{1}{4}-\frac{\delta}{4}}\|u_3^{ \varepsilon}+u_4^{ \varepsilon}\|_{-z}C_W^{ \varepsilon})\\&+C(C^{ \varepsilon}_W,\bar{C}^{ \varepsilon}_W)\int_0^t(t-s)^{-\frac{3/2+\delta/2-z+\kappa}{2}}s^{-\frac{\delta+\kappa+z}{2}}ds\sup_{s\in[0,t]}s^{\frac{\delta+\kappa+z}{2}}\|u^{ \varepsilon}_3+u^{ \varepsilon}_4-(\bar{u}^{ \varepsilon}_3+\bar{u}^{ \varepsilon}_4)\|_{\delta}
\\&+t^{\frac{1}{4}-\frac{\delta}{4}}\|u_3^{ \varepsilon}+u_4^{ \varepsilon}-(\bar{u}_3^{ \varepsilon}+\bar{u}_4^{ \varepsilon})\|_{-z}C_W^{ \varepsilon}).\endaligned$$
Hence we obtain that for $T_1$ small enough and $T_1\leq \tau_L\wedge \tau_{L_1}^\varepsilon\wedge \rho_{L_2}^\varepsilon\wedge\bar{\rho}_{L_3}^\varepsilon$
$$\aligned\sup_{t\in[0, T_1]}\|u_4^{ \varepsilon,i}(t)-\bar{u}_4^{ \varepsilon,i}(t)\|_{-z}\lesssim&(\delta C_W^{ \varepsilon}+\varepsilon^{\kappa/2}+\|u_0^{ \varepsilon}-u_0\|_{-z})C(L,L_1,L_2,L_3).\endaligned$$

We can extend the time to $\tau_L\wedge \tau_{L_1}^\varepsilon\wedge \rho_{L_2}^\varepsilon\wedge\bar{\rho}_{L_3}^\varepsilon$ as we did in Subsection 3.2. By a similar argument as (3.19) and results in Section 4 we get that $\delta C_W^{\varepsilon}\rightarrow^P0$  as $\varepsilon\rightarrow0$. Since
 $$ \sup_{t\in[0,\tau_L\wedge\bar{\rho}_{L_3}^\varepsilon]}\|\bar{u}^{\varepsilon}(t)-u(t)\|_{-z}\rightarrow^P0,$$
 we obtain that if $\|u_0^\varepsilon-u_0\|_{-z}\rightarrow0$,
$$ \sup_{t\in[0,\tau_L\wedge \tau_{L_1}^\varepsilon\wedge \rho_{L_2}^\varepsilon\wedge\bar{\rho}_{L_3}^\varepsilon]}\|u^{\varepsilon}(t)-u(t)\|_{-z}\rightarrow^P0, \quad\varepsilon\rightarrow0$$

 \no\emph{Proof of Theorem 1.3} It is sufficient to show that for every $\epsilon>0$, $L>0$
$$\lim_{\varepsilon\rightarrow0}P(\sup_{t\in[0,\tau_L]}\|u^\varepsilon-u\|_{-z}>\epsilon)=0.$$
We have the following estimates:
$$\aligned &P(\sup_{t\in[0,\tau_L]}\|u^\varepsilon-u\|_{-z}>\epsilon)\\\leq& P(\sup_{t\in[0,\tau_L\wedge \tau_{L_1}^\varepsilon\wedge \rho_{L_2}^\varepsilon\wedge\bar{\rho}_{L_3}^\varepsilon]}\|u^\varepsilon-u\|_{-z}>\epsilon)+P(\tau_L\wedge\rho_{L_2}^\varepsilon\wedge\bar{\rho}_{L_3}^\varepsilon>\tau_{L_1}^\varepsilon)+P(\tau_L>\rho_{L_2}^\varepsilon)+P(\tau_L>\bar{\rho}_{L_3}^\varepsilon).\endaligned$$
The first term goes to zero by above proof. Also  for $L_1>L+\epsilon$
$$P(\tau_L\wedge\rho_{L_2}^\varepsilon\wedge\bar{\rho}_{L_3}^\varepsilon>\tau_{L_1}^\varepsilon)\leq P(\sup_{t\in[0,\tau_L\wedge \tau_{L_1}^\varepsilon\wedge \rho_{L_2}^\varepsilon\wedge\bar{\rho}_{L_3}^\varepsilon]}\|u^\varepsilon-u\|_{-z}>\epsilon),$$
which goes to zero by above proof.
The last two terms go to zero uniformly over $\varepsilon\in(0,1)$ as $L_2,L_3$ go to $\infty.$
Thus the result follows. $\hfill\Box$

\section{Convergence of renormalisation terms}
In the following we use notation $X,\bar{X}$ to represent $u_1, \bar{u}_1$ in the calculation respectively and $\hat{f}(k)= (2\pi)^{-\frac{3}{2}}\int_{\mathbb{T}^3} f(x)e^{\imath x\cdot k}dx$ for $k\in\mathbb{Z}^3$. To simplify the arguments below, we assume that $\hat{W}(0)=0$ and restrict ourselves to the flow of $\int_{\mathbb{T}^3} u(x)dx=0$. Then $$X_t^{\varepsilon,i}=\int_{-\infty}^t\sum_{i_1=1}^3P^{ii_1}P^\varepsilon_{t-s}H_\varepsilon dW^{i_1}_s =\sum _{k\in\mathbb{Z}^3\backslash\{0\}}\hat{X}_t^{\varepsilon,i}(k)e_k$$ is a centered Gaussian process with covariance function given by
$$E[\hat{X}_t^{\varepsilon,i}(k)\hat{X}_s^{\varepsilon,j}(k')]=1_{k+k'=0}\sum_{i_1=1}^3\frac{e^{-|k|^2f(\varepsilon k)|t-s|}h(\varepsilon k)^2}{2|k|^2f(\varepsilon k)}\hat{P}^{ii_1}(k)\hat{P}^{ji_1}(k),$$
and $\hat{X}_t(0)=0$, where $e_k(x)=(2\pi)^{-3/2}e^{\imath x\cdot k},x\in\mathbb{T}^3$ and $\hat{P}^{ii_1}(k)=\delta_{ii_1}-\frac{k^ik^{i_1}}{|k|^2}$ for $k\in\mathbb{Z}^3\backslash\{0\}$. Moreover, $\bar{X}_t^{\varepsilon,i}=\int_{-\infty}^t\sum_{i_1=1}^3P^{ii_1}P_{t-s}H_\varepsilon dW^{i_1} =\sum _{k\in\mathbb{Z}^3\backslash\{0\}}\hat{\bar{X}}_t^{\varepsilon,i}(k)e_k$ is a centered Gaussian process with covariance function given by
$$E[\hat{\bar{X}}_t^{\varepsilon,i}(k)\hat{\bar{X}}_s^{\varepsilon,j}(k')]=1_{k+k'=0}\sum_{i_1=1}^3\frac{e^{-|k|^2|t-s|}h(\varepsilon k)^2}{2|k|^2}\hat{P}^{ii_1}(k)\hat{P}^{ji_1}(k),$$
and $\hat{\bar{X}}_t(0)=0$. We also have for $t\leq s$
$$E[\hat{X}_t^{\varepsilon,i}(k)\hat{\bar{X}}_s^{\varepsilon,j}(k')]=1_{k+k'=0}\sum_{i_1=1}^3\frac{e^{-|k|^2(s-t)}h(\varepsilon k)^2}{|k|^2(f(\varepsilon k)+1)}\hat{P}^{ii_1}(k)\hat{P}^{ji_1}(k),$$
and for $t>s$
$$E[\hat{X}_t^{\varepsilon,i}(k)\hat{\bar{X}}_s^{\varepsilon,j}(k')]=1_{k+k'=0}\sum_{i_1=1}^3\frac{e^{-|k|^2(t-s)f(\varepsilon k)}h(\varepsilon k)^2}{|k|^2(f(\varepsilon k)+1)}\hat{P}^{ii_1}(k)\hat{P}^{ji_1}(k).$$
In this section we will prove that
 for $i,i_1,j,j_2=1,2,3,$ $u_1^{\varepsilon,i}-\bar{u}_1^{\varepsilon,i}\rightarrow 0$ in $C([0,T];\mathcal{C}^{-1/2-\delta/2})$, $u_1^{\varepsilon,i}\diamond u_1^{\varepsilon,j}-\bar{u}_1^{\varepsilon,i}\diamond \bar{u}_1^{\varepsilon,j}\rightarrow 0$ in $C([0,T];\mathcal{C}^{-1-\delta})$, $u_1^{\varepsilon,i}\diamond u_2^{\varepsilon,j}-\bar{u}_1^{\varepsilon,i}\diamond \bar{u}_2^{\varepsilon,j}\rightarrow 0$ in $C([0,T];\mathcal{C}^{-1/2-\delta})$, $u_2^{\varepsilon,i}\diamond u_2^{\varepsilon,j}-\bar{u}_2^{\varepsilon,i}\diamond \bar{u}_2^{\varepsilon,j}\rightarrow 0$ in $C([0,T];\mathcal{C}^{-\delta})$, $\pi_{0,\diamond}(u_3^{\varepsilon,i}, u_1^{\varepsilon,j})-\pi_{0,\diamond}(\bar{u}_3^{\varepsilon,i}, \bar{u}_1^{\varepsilon,j})\rightarrow 0$ in $C([0,T];\mathcal{C}^{-\delta})$, $\pi_{0,\diamond}(P^{ii_1}D^\varepsilon_jK^{\varepsilon,j}, u_1^{\varepsilon,j_2})-\pi_{0,\diamond}(P^{ii_1}D_j\bar{K}^{\varepsilon,j}, \bar{u}_1^{\varepsilon,j_2})\rightarrow 0$ in $C([0,T];\mathcal{C}^{-\delta})$ and $\pi_{0,\diamond}(P^{ii_1}D^\varepsilon_jK^{\varepsilon,i_1}, u_1^{\varepsilon,j_2})-\pi_{0,\diamond}(P^{ii_1}D_j\bar{K}^{\varepsilon,i_1}, \bar{u}_1^{\varepsilon,j_2})\rightarrow0$ in $C([0,T];\mathcal{C}^{-\delta})$, as $\varepsilon\rightarrow0$. 

\vskip.10in

 Now we introduce the following notations: $k_{1...n}=\sum_{i=1}^nk_i$. To obtain the results we first recall the following two lemmas from [ZZ14] for our later use:

\vskip.10in
\th{Lemma 4.1} ([ZZ14, Lemma 3.10]) Let $0<l,m<d,l+m-d>0$. Then we have
$$\sum_{k_1,k_2\in \mathbb{Z}^d\backslash\{0\},k_1+k_2=k}\frac{1}{|k_1|^{l}|k_2|^{m}}\lesssim \frac{1}{|k|^{l+m-d}}.$$

\vskip.10in
  \th{Lemma 4.2}([ZZ14, Lemma 3.11]) For any $0<\eta<1$, $i,j,l=1,2,3$ we have
  $$|e^{-|k_{12}|^2(t-s)}k_{12}^i\hat{P}^{jl}(k_{12})-e^{-|k_{2}|^2(t-s)}k_{2}^i\hat{P}^{jl}(k_{2})|\lesssim |k_1|^\eta|t-s|^{-(1-\eta)/2}.$$
  Here $\hat{P}^{ij}(x)=\delta_{ij}-\frac{x^i\otimes x^j}{|x|^2}$.
\vskip.10in
By a similar argument as the proof of Lemma 4.2 we have the following result.
\vskip.10in
  \th{Lemma 4.3} For any $0<\eta<1$, $i,j,l=1,2,3$ we have for $|\varepsilon k_{12}^i|\leq 3L_0, |\varepsilon k_2^i|\leq 3L_0$
  $$|e^{-|k_{12}|^2(t-s)\tilde{f}(\varepsilon k_{12})}k_{12}^ig(\varepsilon k_{12}^i)\hat{P}^{jl}(k_{12})-e^{-|k_{2}|^2(t-s)\tilde{f}(\varepsilon k_2)}k_{2}^ig(\varepsilon k_{2}^i)\hat{P}^{jl}(k_{2})|\lesssim |k_1|^\eta|t-s|^{-(1-\eta)/2}.$$
  Here $\hat{P}^{ij}(x)=\delta_{ij}-\frac{x^i\otimes x^j}{|x|^2}$.

\subsection{Convergence for $u_1^\varepsilon-\bar{u}_1^\varepsilon$}
In this subsection we consider the convergence of $u_1^\varepsilon-\bar{u}_1^\varepsilon$.

For $t_1<t_2$ we have
$$\aligned &E|\Delta_q[(u_1^{\varepsilon,i}(t_2)-\bar{u}_1^{\varepsilon,i}(t_2))-(u_1^{\varepsilon,i}(t_1)-\bar{u}_1^{\varepsilon,i}(t_1))]|^2
\\\lesssim &\sum_k\theta(2^qk)^2h(\varepsilon k)^2\bigg|\frac{1}{|k|^2f(\varepsilon k)}-\frac{4}{|k|^2(f(\varepsilon k)+1)}-\frac{e^{-|k|^2f(\varepsilon k)|t_2-t_1|}}{|k|^2f(\varepsilon k)}+\frac{2e^{-|k|^2f(\varepsilon k)|t_2-t_1|}}{|k|^2(f(\varepsilon k)+1)}\\&+\frac{1}{|k|^2}+\frac{2e^{-|k|^2|t_2-t_1|}}{|k|^2(f(\varepsilon k)+1)}-\frac{e^{-|k|^2|t_2-t_1|}}{|k|^2}\bigg|
\\\lesssim &\sum_k\theta(2^qk)h(\varepsilon k)^2(\frac{(\varepsilon|k|)^{2\eta}}{|k|^2}\wedge \frac{|k|^{2\eta}|t_2-t_1|^\eta}{|k|^2})\lesssim\varepsilon^\eta2^{q(2\eta+1)}|t_2-t_1|^{\eta/2}.\endaligned\eqno(4.1)$$
Here $\eta>0$ is small enough and in the  second inequality we used $|f(\varepsilon k)-1|\lesssim \varepsilon|k|,$ for $|\varepsilon k^i|\leq L_0$, $i=1,2,3$,  and $|f(\xi)|\geq c_f>0$.
Then by Gaussian hypercontractivity we have for $p>1$ that
$$\aligned &E\|\Delta_q[(u_1^{\varepsilon,i}(t_2)-\bar{u}_1^{\varepsilon,i}(t_2))-(u_1^{\varepsilon,i}(t_1)-\bar{u}_1^{\varepsilon,i}(t_1))]\|_{L^p}^p
\\\lesssim&\int_{\mathbb{T}^3} E(|\Delta_q[(u_1^{\varepsilon,i}(t_2)-\bar{u}_1^{\varepsilon,i}(t_2))-(u_1^{\varepsilon,i}(t_1)-\bar{u}_1^{\varepsilon,i}(t_1))](x)|^2)^{p/2}dx
\\\lesssim  &\varepsilon^{p\eta/2}2^{q(2\eta+1)p/2}|t_2-t_1|^{\eta p/4},\endaligned$$
which implies that for $\epsilon$ small enough
$$\aligned &E[\|(u_1^{\varepsilon,i}(t_2)-\bar{u}_1^{\varepsilon,i}(t_2))-(u_1^{\varepsilon,i}(t_1)-\bar{u}_1^{\varepsilon,i}(t_1))\|_{B^{-1/2-\eta-\epsilon}_{p,p}}^p]
\\\lesssim&\varepsilon^{p\eta/2}|t_2-t_1|^{\eta p/4},\endaligned$$
Then by Lemma 2.1 we obtain that for every $\delta>0, p>1$, $u_1^{\varepsilon,i}-\bar{u}_1^{\varepsilon,i}\rightarrow 0$ in $L^p(\Omega;C([0,T];\mathcal{C}^{-1/2-\delta/2}))$ as $\varepsilon\rightarrow0$.

 \subsection{Convergence for $u_1^{\varepsilon,i}\diamond u_1^{\varepsilon,j}-\bar{u}_1^{\varepsilon,i}\diamond \bar{u}_1^{\varepsilon,j}$}
 In this subsection we consider the convergence of $u_1^{\varepsilon,i}\diamond u_1^{\varepsilon,j}$. Recall that $u_1^{\varepsilon,i}\diamond u_1^{\varepsilon,j}=u_1^{\varepsilon,i} u_1^{\varepsilon,j}-C^{\varepsilon,ij}_0$ and $\bar{u}_1^{\varepsilon,i}\diamond \bar{u}_1^{\varepsilon,j}=\bar{u}_1^{\varepsilon,i} \bar{u}_1^{\varepsilon,j}-\bar{C}^{\varepsilon,ij}_0$.

 Take $$C^{\varepsilon,ij}_0=(2\pi)^{-3}\sum_k\frac{h(\varepsilon k)^2}{2|k|^2f(\varepsilon k)}\sum_{i_1=1}^3\hat{P}^{ii_1}(k)\hat{P}^{ji_1}(k)$$ and $$\bar{C}^{\varepsilon,ij}_0=(2\pi)^{-3}\sum_k\frac{h(\varepsilon k)^2}{2|k|^2}\sum_{i_1=1}^3\hat{P}^{ii_1}(k)\hat{P}^{ji_1}(k).$$
  For $t_1<t_2$
$$\aligned &E|\Delta_q[(u_1^{\varepsilon,i}\diamond u_1^{\varepsilon,j}(t_2)-\bar{u}_1^{\varepsilon,i}\diamond \bar{u}_1^{\varepsilon,j}(t_2))-(u_1^{\varepsilon,i}\diamond u_1^{\varepsilon,j}(t_1)-\bar{u}_1^{\varepsilon,i}\diamond \bar{u}_1^{\varepsilon,j}(t_1))]|^2
\\\lesssim &\sum_k\sum_{k_{12}=k}\theta(2^{-q}k)^2h(\varepsilon k_1)^2h(\varepsilon k_2)^2\bigg|\frac{1}{2|k_1|^2|k_2|^2f(\varepsilon k_1)f(\varepsilon k_2)}-\frac{4}{|k_1|^2|k_2|^2(f(\varepsilon k_1)+1)(f(\varepsilon k_2)+1)}\\&-\frac{e^{-(|k_1|^2f(\varepsilon k_1)+|k_2|^2f(\varepsilon k_2))(t_2-t_1)}}{2|k_1|^2|k_2|^2f(\varepsilon k_1)f(\varepsilon k_2)}+\frac{2e^{-(|k_1|^2+|k_2|^2)(t_2-t_1)}}{|k_1|^2|k_2|^2(f(\varepsilon k_1)+1)(f(\varepsilon k_2)+1)}+\frac{1}{2|k_1|^2|k_2|^2}\\&+\frac{2e^{-(|k_1|^2f(\varepsilon k_1)+|k_2|^2f(\varepsilon k_2))(t_2-t_1)}}{|k_1|^2|k_2|^2(1+f(\varepsilon k_1))(1+f(\varepsilon k_2))}-\frac{e^{-(|k_1|^2+|k_2|^2)(t_2-t_1)}}{2|k_1|^2|k_2|^2}\bigg|\\\lesssim&\sum_k\sum_{k_{12}=k}\theta(2^{-q}k)^2\bigg[\frac{(\varepsilon(|k_1|+|k_2|))^{2\eta}}{|k_1|^2|k_2|^2}\wedge \frac{(t_2-t_1)^\eta(|k_1|^{2\eta}+|k_2|^{2\eta})}{|k_1|^2|k_2|^2}\bigg]
\\\lesssim & \varepsilon^\eta(t_2-t_1)^{\eta/2}\sum_k\sum_{k_{12}=k}\theta(2^{-q}k)^2\frac{|k_1|^{2\eta}+|k_2|^{2\eta}}{|k_1|^2|k_2|^2}
\\\lesssim&\varepsilon^\eta(t_2-t_1)^{\eta/2}2^{(2\eta+2)q}.\endaligned$$
Here $\eta>0$ is small enough and in the last inequality we used Lemma 4.1. Then by Gaussian hypercontractivity and  Lemma 2.1 we obtain that for every $\delta>0, p>1,$ $u_1^{\varepsilon,i}\diamond u_1^{\varepsilon,j}-\bar{u}_1^{\varepsilon,i}\diamond \bar{u}_1^{\varepsilon,j}\rightarrow 0$ in $L^p(\Omega;C([0,T];\mathcal{C}^{-1-\delta}))$.

 \subsection{Convergence for $u_1^\varepsilon\diamond  u_2^\varepsilon-\bar{u}_1^\varepsilon\diamond  \bar{u}_2^\varepsilon$}
In this subsection we focus on $u_1^\varepsilon u_2^\varepsilon$ and prove that $u_1^{\varepsilon,i}\diamond u_2^{\varepsilon,j}-\bar{u}_1^{\varepsilon,i}\diamond  \bar{u}_2^{\varepsilon,j}\rightarrow 0$ in $C([0,T];\mathcal{C}^{-1/2-\delta})$ for $i,j=1,2,3$.
Recall that for $i,j=1,2,3$,
$$u_1^{\varepsilon,j}\diamond  u_2^{\varepsilon,i}= u_1^{\varepsilon,j} u_2^{\varepsilon,i}+\sum_{i_1=1}^3u_1^{\varepsilon,i_1} (C^{\varepsilon,i,i_1,j}(t)+\tilde{C}^{\varepsilon,i,i_1,j}(t)),$$
and
$$\bar{u}_1^{\varepsilon,j}\diamond  {\bar{u}}_2^{\varepsilon,i}= {\bar{u}}_1^{\varepsilon,j} \bar{u}_2^{\varepsilon,i}.$$
Now we have the following identity: for $t\in[0,T]$, $i,j=1,2,3$
 $$\aligned &\bar{u}_1^{\varepsilon,j} \bar{u}_2^{\varepsilon,i}(t)-u_1^{\varepsilon,j} u_2^{\varepsilon,i}(t)\\=&\frac{(2\pi)^{-3}}{2}\sum_{i_1,i_2=1}^3\sum_{k\in \mathbb{Z}^3\backslash\{0\}}\sum_{k_{123}=k}\bigg[\int_0^te^{-|k_{12}|^2(t-s)f(\varepsilon k_{12})}  k_{12}^{i_2}g(\varepsilon k_{12}^{i_2}):\hat{X}^{\varepsilon,i_1}_s(k_1)\hat{X}^{\varepsilon,i_2}_s(k_2)\hat{X}^{\varepsilon,j}_t(k_3):ds \hat{ P}^{ii_1}(k_{12})\\&-\int_0^te^{-|k_{12}|^2(t-s)} \imath k_{12}^{i_2}:\hat{\bar{X}}^{\varepsilon,i_1}_s(k_1)\hat{\bar{X}}^{\varepsilon,i_2}_s(k_2)\hat{\bar{X}}^{\varepsilon,j}_t(k_3):ds \hat{ P}^{ii_1}(k_{12})\bigg]e_k\\&+\frac{(2\pi)^{-3}}{2}\sum_{i_1,i_2,i_3=1}^3\sum_{k_1,k_2\in\mathbb{Z}^3\backslash\{0\}}\bigg{[}\int_0^te^{-f(\varepsilon k_{12})|k_{12}|^2(t-s)}g(\varepsilon k_{12}^{i_2}) k_{12}^{i_2}\hat{X}^{\varepsilon,i_1}_s(k_1)\frac{e^{-f(\varepsilon k_2)|k_2|^2(t-s)}h(\varepsilon k_2)^2}{2|k_2|^2f(\varepsilon k_2)}ds\endaligned$$
$$\aligned&-\int_0^te^{-|k_{12}|^2(t-s)}\imath k_{12}^{i_2}\hat{\bar{X}}^{\varepsilon,i_1}_s(k_1)\frac{e^{-|k_2|^2(t-s)}h(\varepsilon k_2)^2}{2|k_2|^2}ds\bigg{]}\hat{P}^{ii_1}(k_{12})\hat{P}^{i_2i_3}(k_{2})\hat{P}^{ji_3}(k_{2})e_{k_1}\\&+\frac{(2\pi)^{-3}}{2}\sum_{i_1,i_2,i_3=1}^3
 \sum_{k_1,k_2\in\mathbb{Z}^3\backslash\{0\}}\bigg{[}\int_0^te^{-f(\varepsilon k_{12})|k_{12}|^2(t-s)}g(\varepsilon k_{12}^{i_2}) k_{12}^{i_2}\hat{X}^{\varepsilon,i_2}_s(k_2)\frac{e^{-f(\varepsilon k_1)|k_1|^2(t-s)}h(\varepsilon k_1)^2}{2|k_1|^2f(\varepsilon k_1)}ds\\&-\int_0^te^{-|k_{12}|^2(t-s)}\imath k_{12}^{i_2}\hat{\bar{X}}^{\varepsilon,i_2}_s(k_2)\frac{e^{-|k_1|^2(t-s)}h(\varepsilon k_1)^2}{2|k_1|^2}ds\bigg{]}\hat{P}^{ii_1}(k_{12})\hat{P}^{i_1i_3}(k_{1})\hat{P}^{ji_3}(k_{1})e_{k_2}\\=&I_t^1+I_t^2+I_t^3.\endaligned$$

\no\textbf{Term in the first chaos:}
First we consider $I_t^2-\sum_{i_1=1}^3X_t^{\varepsilon,i_1} C^{\varepsilon,i,i_1,j}(t)$ and we have $$I_t^2-\sum_{i_1=1}^3X_t^{\varepsilon,i_1} C^{\varepsilon,i,i_1,j}(t)=I_t^2-\tilde{I}_t^2+\tilde{I}_t^2-\sum_{i_1=1}^3X_t^{\varepsilon,i_1} C^{\varepsilon,i,i_1,j}(t)+\sum_{i_1=1}^3\bar{X}^{\varepsilon,i_1}_t \bar{C}^{\varepsilon,i,i_1,j}(t),$$
 where $$\aligned\tilde{I}_t^2=&\frac{(2\pi)^{-3}}{2}\sum_{i_1,i_2,i_3=1}^3\sum_{k_1,k_2\in\mathbb{Z}^3\backslash\{0\}}\bigg{[}\hat{X}^{\varepsilon,i_1}_t(k_1)e_{k_1}
 \int_0^te^{-|k_{12}|^2(t-s)f(\varepsilon k_{12})} k_{12}^{i_2}g(\varepsilon k_{12}^{i_2})\frac{e^{-|k_2|^2(t-s)f(\varepsilon k_2)}}{2|k_2|^2f(\varepsilon k_{2})}h(\varepsilon k_2)^2ds\\&-\hat{\bar{X}}^{\varepsilon,i_1}_t(k_1)e_{k_1}
 \int_0^te^{-|k_{12}|^2(t-s)} \imath k_{12}^{i_2}\frac{e^{-|k_2|^2(t-s)}}{2|k_2|^2}h(\varepsilon k_2)^2ds\bigg{]}\hat{P}^{ii_1}(k_{12})\hat{P}^{i_2i_3}(k_{2})\hat{P}^{ji_3}(k_{2}),\endaligned$$
 and as introduction
  $$\aligned C^{\varepsilon,i,i_1,j}(t)=&\frac{(2\pi)^{-3}}{2}\sum_{i_2,i_3=1}^3\sum_{k_2\in\mathbb{Z}^3\backslash\{0\}}\int_0^te^{-2|k_2|^2(t-s)f(\varepsilon k_2)} k_2^{i_2}g(\varepsilon k_2^{i_2})\frac{h(\varepsilon k_2)^2}{2|k_2|^2f(\varepsilon k_2)}\hat{P}^{ii_1}(k_2)\hat{P}^{i_2i_3}(k_{2})\hat{P}^{ji_3}(k_{2})ds\\=&\frac{(2\pi)^{-3}}{4(a+b)\varepsilon}\sum_{i_2,i_3=1}^3\sum_{k_2\in\mathbb{Z}^3\backslash\{0\}}\int_0^te^{-2|k_2|^2(t-s)f(\varepsilon k_2)} (\cos(a\varepsilon k_2^{i_2})-\cos(b\varepsilon k_2^{i_2}))\frac{h(\varepsilon k_2)^2}{|k_2|^2f(\varepsilon k_2)}\\&\hat{P}^{ii_1}(k_2)\hat{P}^{i_2i_3}(k_{2})\hat{P}^{ji_3}(k_{2})ds\\=&\frac{(2\pi)^{-3}}{4(a+b)\varepsilon}\sum_{i_2,i_3=1}^3\sum_{k_2\in\mathbb{Z}^3\backslash\{0\}}(1-e^{-2|k_2|^2tf(\varepsilon k_2)}) (\cos(\varepsilon ak_2^{i_2})-\cos(\varepsilon bk_2^{i_2}))\frac{h(\varepsilon k_2)^2}{2|k_2|^4f(\varepsilon k_2)^2}\\&\hat{P}^{ii_1}(k_2)\hat{P}^{i_2i_3}(k_{2})\hat{P}^{ji_3}(k_{2})\\\rightarrow &\frac{(2\pi)^{-3}}{8(a+b)}\sum_{i_2,i_3=1}^3\int_{\mathbb{R}^3} (\cos (ax^{i_2})-\cos(b x^{i_2}))\frac{h(x)^2}{|x|^4f(x)^2}\hat{P}^{ii_1}(x)\hat{P}^{i_2i_3}(x)\hat{P}^{ji_3}(x)dx,\endaligned$$
  as $\varepsilon\rightarrow0$,
 and $$\bar{C}^{\varepsilon,i,i_1,j}(t)=(2\pi)^{-3}\sum_{i_2,i_3=1}^3\sum_{k_2\in\mathbb{Z}^3\backslash\{0\}}\int_0^te^{-2|k_2|^2(t-s)}\imath k_2^{i_2}\frac{h(\varepsilon k_2)^2}{4|k_2|^2}\hat{P}^{ii_1}(k_2)\hat{P}^{i_2i_3}(k_{2})\hat{P}^{ji_3}(k_{2})ds=0.$$
By a straightforward calculation we obtain that for $\eta>\epsilon>0$ small enough
$$\aligned& E[|\Delta_q(I_t^2-\tilde{I}_t^2)|^2]\\\lesssim& E\bigg[\bigg|\sum_{i_1,i_2,i_3=1}^3\int_0^t\sum_{k_1}\theta(2^{-q}k_1)e_{k_1}(a_{k_1}^{\varepsilon,i_1i_2i_3}
-\bar{a}_{k_1}^{\varepsilon,i_1i_2i_3})(t-s)(\hat{X}_s^{\varepsilon,i_1}(k_1)
-\hat{X}_t^{\varepsilon,i_1}(k_1))ds\bigg|^2\bigg]\\&+E\bigg[\bigg|\sum_{i_1,i_2,i_3=1}^3\int_0^t\sum_{k_1}\theta(2^{-q}k_1)e_{k_1}\bar{a}_{k_1}^{\varepsilon,i_1i_2i_3}(t-s)(\hat{X}_s^{\varepsilon,i_1}(k_1)
-\hat{X}_t^{\varepsilon,i_1}(k_1)-\hat{\bar{X}}_s^{\varepsilon,i_1}(k_1)
+\hat{\bar{X}}_t^{\varepsilon,i_1}(k_1))ds\bigg|^2\bigg]
\\\lesssim &\sum_{k_1}\theta(2^{-q}k_1)^2\bigg{[}\frac{1}{|k_1|^{2-2\eta}}\bigg(\int_0^t\sum_{k_2}\varepsilon^{\eta/2}(|k_{12}|^{\eta/2}+|k_{2}|^{\eta/2})|k_{12}|
\frac{e^{-\bar{c}_f|k_{12}|^2(t-s)-\bar{c}_f|k_2|^2(t-s)}}{|k_2|^2}(t-s)^{\eta/2}ds\bigg)^2
\\&+\frac{\varepsilon^{\eta}}{|k_1|^{2-2\eta}}\bigg(\int_0^t\sum_{k_2}|k_{12}|\frac{e^{-|k_{12}|^2(t-s)-|k_2|^2(t-s)}}{|k_2|^2}(t-s)^{\eta/4}ds\bigg)^2\bigg{]}
\\\lesssim &\varepsilon^\eta t^{(\eta-\epsilon)/2}2^{q(1+2\eta)}.\endaligned$$
Here $$a_{k_1}^{\varepsilon,i_1i_2i_3}(t-s)=\sum_{k_2}e^{-|k_{12}|^2(t-s)f(\varepsilon k_{12})} k_{12}^{i_2}g(\varepsilon k_{12}^{i_2})\frac{e^{-|k_2|^2(t-s)f(\varepsilon k_2)}h(\varepsilon k_2)^2}{|k_2|^2f(\varepsilon k_{2})}\hat{P}^{ii_1}(k_{12})\hat{P}^{i_2i_3}(k_{2})\hat{P}^{ji_3}(k_{2}),$$
$$\bar{a}_{k_1}^{\varepsilon,i_1i_2i_3}(t-s)=\sum_{k_2}e^{-|k_{12}|^2(t-s)} \imath k_{12}^{i_2}\frac{e^{-|k_2|^2(t-s)}h(\varepsilon k_2)^2}{|k_2|^2}\hat{P}^{ii_1}(k_{12})\hat{P}^{i_2i_3}(k_{2})\hat{P}^{ji_3}(k_{2}),$$
and in the  second inequality we used that for $\eta>0$ small enough
$$|e^{-|k_{12}|^2(t-s)f(\varepsilon k_{12})}-e^{-|k_{12}|^2(t-s)}|\lesssim e^{-|k_{12}|^2\bar{c}_f(t-s)}(1\wedge(t-s)^\eta|f(\varepsilon k_{12})-1|^\eta|k_{12}|^{2\eta})\lesssim  e^{-|k_{12}|^2\bar{c}_f(t-s)} |\varepsilon k_{12}|^{\eta/2},\eqno(4.2)$$
and
$$|g(\varepsilon k_{12}^{i_2})-\imath|\lesssim |\varepsilon k_{12}|^{\eta/2},\eqno(4.3)$$
which imply that
$$|(a_{k_1}^{\varepsilon,i_1i_2i_3}-\bar{a}_{k_1}^{\varepsilon,i_1i_2i_3})(t-s)|\lesssim \sum_{k_2}\varepsilon^{\eta/2}(|k_{12}|^{\eta/2}+|k_{2}|^{\eta/2})|k_{12}|\frac{e^{-|k_{12}|^2\bar{c}_f(t-s)-|k_2|^2\bar{c}_f(t-s)}}{|k_2|^2},$$
where $\bar{c}_f=c_f\wedge 1$.
 In the  second inequality we also used that for $\eta>0$ small enough
$$\aligned &E|(\hat{X}_s^{\varepsilon,i_1}
(k_1)-\hat{X}_t^{\varepsilon,i_1}(k_1))(\overline{\hat{X}_{\bar{s}}^{\varepsilon,i_1'}(k'_1)-\hat{X}_t^{\varepsilon,i_1'}(k'_1)})|
\\\leq& 1_{k_1=k_1'}(E|(\hat{X}_s^{\varepsilon,i_1}
(k_1)-\hat{X}_t^{\varepsilon,i_1}(k_1))|^2)^{1/2}(E|(\overline{\hat{X}_{\bar{s}}^{\varepsilon,i_1'}(k'_1)-\hat{X}_t^{\varepsilon,i_1'}(k'_1)})|^2)^{1/2}
\\\lesssim& 1_{k_1=k_1'}(\frac{h(\varepsilon k_1)^2}{|k_1|^2f(\varepsilon k_1)}(1-e^{-|k_1|^2f(\varepsilon k_1)(t-s)}))^{1/2}(\frac{h(\varepsilon k'_1)^2}{|k'_1|^2f(\varepsilon k'_1)}(1-e^{-|k'_1|^2f(\varepsilon k_1')(t-\bar{s})}))^{1/2}\\\lesssim&
\frac{h(\varepsilon k_1)^2}{|k_1|^2}|k_1|^{2\eta}|t-s|^{\eta/2}|t-\bar{s}|^{\eta/2},\endaligned$$
and by the estimates in Section 4.1
$$\aligned &E((\hat{X}_s^{\varepsilon,i_1}
(k_1)-\hat{X}_t^{\varepsilon,i_1}(k_1))-(\hat{\bar{X}}_{s}^{\varepsilon,i_1}(k_1)-\hat{\bar{X}}_t^{\varepsilon,i_1}(k_1)))(\overline{(\hat{X}_{\bar{s}}^{\varepsilon,i'_1}
(k'_1)-\hat{X}_t^{\varepsilon,i'_1}(k'_1))-(\hat{\bar{X}}_{\bar{s}}^{\varepsilon,i'_1}(k'_1)-\hat{\bar{X}}_t^{\varepsilon,i'_1}(k'_1))})
\\\lesssim&1_{k_1=k_1'} (E|(\hat{X}_s^{\varepsilon,i_1}
(k_1)-\hat{X}_t^{\varepsilon,i_1}(k_1))-(\hat{\bar{X}}_{s}^{\varepsilon,i_1}(k_1)-\hat{\bar{X}}_t^{\varepsilon,i_1}(k_1))|^2)^{1/2}\\&(E|(\hat{X}_{\bar{s}}^{\varepsilon,i'_1}
(k'_1)-\hat{X}_t^{\varepsilon,i'_1}(k'_1))-(\hat{\bar{X}}_{\bar{s}}^{\varepsilon,i'_1}(k'_1)-\hat{\bar{X}}_t^{\varepsilon,i'_1}(k'_1))|^2)^{1/2}
\\\lesssim&
\frac{h(\varepsilon k_1)^2}{|k_1|^2}\varepsilon^\eta|k_1|^{2\eta}|t-s|^{\eta/4}|t-\bar{s}|^{\eta/4},\endaligned$$
where in the last inequality we used (4.1).
Moreover, we obtain that
$$\aligned &E[|\Delta_q(\tilde{I}_t^2-\sum_{i_1=1}^3X_t^{\varepsilon,i_1} C^{\varepsilon,i,i_1,j}(t)+\sum_{i_1=1}^3\bar{X}_t^{\varepsilon,i_1} \bar{C}^{\varepsilon,i,i_1,j}(t))|^2]\\\lesssim &\sum_{k_1}\frac{h(\varepsilon k_1)^2h(\varepsilon k_2)^2}{2|k_1|^2f(\varepsilon k_1)}\theta(2^{-q}k_1)^2\bigg{[}\sum_{i_1,i_2,i_3=1}^3\sum_{k_2}\int_0^t\frac{e^{-|k_2|^2(t-s)f(\varepsilon k_2)}}{|k_2|^2f(\varepsilon k_2)}\\&\bigg(e^{-|k_{12}|^2(t-s)f(\varepsilon k_{12})}k_{12}^{i_2}g(\varepsilon k_{12}^{i_2})\hat{P}^{ii_1}(k_{12})-e^{-|k_{2}|^2f(\varepsilon k_2)(t-s)}k_{2}^{i_2}g(\varepsilon k_{2}^{i_2})
\hat{P}^{ii_1}(k_{2})\bigg)\hat{P}^{i_2i_3}(k_{2})\hat{P}^{ji_3}(k_{2})\\&-\frac{e^{-|k_2|^2(t-s)}}{|k_2|^2}\bigg(e^{-|k_{12}|^2(t-s)}\imath k_{12}^{i_2}\hat{P}^{ii_1}(k_{12})-e^{-|k_{2}|^2(t-s)}\imath k_{2}^{i_2}
\hat{P}^{ii_1}(k_{2})\bigg)\hat{P}^{i_2i_3}(k_{2})\hat{P}^{ji_3}(k_{2})ds\bigg{]}^2
\\&+\sum_{k_1}\theta(2^{-q}k_1)^2\frac{\varepsilon^\eta|k_1|^\eta}{|k_1|^2}\bigg{[}
\sum_{i_1,i_2,i_3=1}^3\sum_{k_2}\int_0^t\frac{e^{-|k_2|^2(t-s)}}{|k_2|^2}\bigg(e^{-|k_{12}|^2(t-s)}\imath k_{12}^{i_2}\hat{P}^{ii_1}(k_{12})\\&-e^{-|k_{2}|^2(t-s)}\imath k_{2}^{i_2}
\hat{P}^{ii_1}(k_{2})\bigg)\hat{P}^{i_2i_3}(k_{2})\hat{P}^{ji_3}(k_{2})ds\bigg{]}^2
=L_t^1+L_t^2,\endaligned\eqno(4.4)$$
where we used $$E|\hat{X}_t^{\varepsilon,i_1}(k_1)-\hat{\bar{X}}_t^{\varepsilon,i_1}(k_1)|^2\lesssim\frac{\varepsilon^\eta|k_1|^\eta}{|k_1|^2}.$$
Now we consider $$E_t=e^{-|k_{12}|^2(t-s)f(\varepsilon k_{12})}k_{12}^{i_2}g(\varepsilon k_{12}^{i_2})\hat{P}^{ii_1}(k_{12})-e^{-|k_{2}|^2f(\varepsilon k_2)(t-s)}k_{2}^{i_2}g(\varepsilon k_{2}^{i_2})
\hat{P}^{ii_1}(k_{2}):$$
By Lemma 4.3 we have that
if $|\varepsilon k_{12}^i|\leq L_0$ and $|\varepsilon k_{2}^i|\leq L_0$ for $i=1,2,3$, then for $0<\eta<1$
$$|E_t|\lesssim |k_1|^\eta(t-s)^{-(1-\eta)/2},$$
which combining with the condition on $h$, Lemma 4.2, (4.2) and (4.3) we obtain that for $\eta-\epsilon>0$ small enough
  $$\aligned|L_t^1|\lesssim& \sum_{k_1}\frac{1}{|k_1|^2}\theta(2^{-q}k_1)\bigg{(}\sum_{k_2}\int_0^t[(t-s)^{-(1-\eta)/2}|k_1|^\eta\wedge (t-s)^{-1/2}(\varepsilon^{\eta/2}|k_2|^{\eta/2}+\varepsilon^{\eta/2}|k_{12}|^{\eta/2} e^{-|k_{12}|^2(t-s)\bar{c}_f})]
\\&\frac{e^{-|k_2|^2(t-s)\bar{c}_f}}{|k_2|^2}ds\bigg{)}^2
\\\lesssim&\varepsilon^{\eta/2}t^{(\eta-\epsilon)/4}\sum_{k_1}\theta(2^{-q}k_1)\frac{1}{|k_1|^{2-\eta}}
\lesssim\varepsilon^{\eta/2}t^{(\eta-\epsilon)/4}2^{q(1+\eta)}.\endaligned\eqno(4.5)$$
 By Lemma 4.2 we have for $\eta-\epsilon>0$
$$|L_t^2|\lesssim \sum_{k_1}\frac{\varepsilon^\eta|k_1|^\eta}{|k_1|^2}\theta(2^{-q}k_1)^2(\sum_{k_2}\int_0^t\frac{e^{-|k_2|^2(t-s)}}{|k_2|^2}|k_1|^{\eta}(t-s)^{-(1-\eta)/2}ds)^2
\lesssim \varepsilon^\eta t^{\eta-\epsilon}2^{q(1+3\eta)},$$
which is the desired estimate.

Similarly, we obtain that
$$\aligned &E[|\Delta_q(I_t^3-\sum_{i_2=1}^3X_t^{\varepsilon,i_2} \tilde{C}^{\varepsilon,i,i_2,j}(t))|^2]\lesssim \varepsilon^\eta t^{\eta-\epsilon}2^{q(1+3\eta)},\endaligned$$
where   $$\aligned \tilde{C}^{\varepsilon,i,i_2,j}(t)=&\frac{(2\pi)^{-3}}{4}\sum_{i_1,i_3=1}^3\sum_{k_2\in\mathbb{Z}^3\backslash\{0\}}\int_0^te^{-2|k_2|^2(t-s)f(\varepsilon k_2)} k_2^{i_2}g(\varepsilon k_2^{i_2})\frac{h(\varepsilon k_2)^2}{|k_2|^2f(\varepsilon k_2)}\hat{P}^{ii_1}(k_2)\hat{P}^{i_1i_3}(k_{2})\hat{P}^{ji_3}(k_{2})ds\\
\rightarrow &(2\pi)^{-3}\sum_{i_1,i_3=1}^3\int_{\mathbb{R}^3} (\cos (ax^{i_2})-\cos(bx^{i_2}))\frac{h(x)^2}{8(a+b)|x|^4f(x)^2}\hat{P}^{ii_1}(x)\hat{P}^{i_1i_3}(x)\hat{P}^{ji_3}(x)dx,\endaligned$$

\textbf{Term in the third chaos:}
Now we focus on the bounds for $I_t^1$. We obtain the following inequalities:
$$\aligned &E|\Delta_q I_t^1|^2\\\lesssim&\sum_{i_1,i_2=1}^3\sum_{i_1',i_2'=1}^3\sum_k\theta(2^{-q}k)^2\sum_{k_{123}=k}\Pi_{i=1}^3\frac{h(\varepsilon k_i)^2}{|k_i|^2f(\varepsilon k_i)}\int_0^t\int_0^t\bigg|(e^{-|k_{12}|^2(t-s)f(\varepsilon k_{12})} k_{12}^{i_2}g(\varepsilon k_{12}^{i_2})-e^{-|k_{12}|^2(t-s)} \imath k_{12}^{i_2})\\&(e^{-|k_{12}|^2(t-\bar{s})f(\varepsilon k_{12})} k_{12}^{i'_2}g(\varepsilon k_{12}^{i'_2})-e^{-|k_{12}|^2(t-\bar{s})} \imath k_{12}^{i'_2})\hat{P}^{ii_1}(k_{12})\hat{P}^{ii'_1}(k_{12})\bigg|dsd\bar{s}
\\&+\sum_{i_1,i_2=1}^3\sum_{i_1',i_2'=1}^3\sum_k\theta(2^{-q}k)^2\sum_{k_{123}=k}\Pi_{i=1}^3\frac{h(\varepsilon k_i)^2}{|k_i|^2f(\varepsilon k_i)}\int_0^t\int_0^t\bigg|(e^{-|k_{12}|^2(t-s)f(\varepsilon k_{12})} k_{12}^{i_2}g(\varepsilon k_{12}^{i_2})-e^{-|k_{12}|^2(t-s)} \imath k_{12}^{i_2})\\&(e^{-|k_{23}|^2(t-\bar{s})f(\varepsilon k_{23})} k_{23}^{i'_2}g(\varepsilon k_{23}^{i'_2})-e^{-|k_{23}|^2(t-\bar{s})} \imath k_{23}^{i'_2})\hat{P}^{ii_1}(k_{12})\hat{P}^{ii'_1}(k_{23})\bigg|dsd\bar{s}\\&+\sum_{i_1,i_2=1}^3\sum_{i_1',i_2'=1}^3\sum_k\theta(2^{-q}k)^2
\sum_{k_{123}=k}\int_0^t\int_0^t\bigg|e^{-|k_{12}|^2(t-s) } k_{12}^{i_2}e^{-|k_{12}|^2(t-\bar{s})}  k_{12}^{i'_2}\hat{P}^{ii_1}(k_{12})\hat{P}^{ii'_1}(k_{12})\bigg|\\&E\bigg|[(:\hat{X}^{\varepsilon}_s(k_1)\hat{X}^{\varepsilon}_s(k_2)
\hat{X}^{\varepsilon}_t(k_3):-:\hat{\bar{X}}^{\varepsilon}_s(k_1)\hat{\bar{X}}^{\varepsilon}_s(k_2)\hat{\bar{X}}^{\varepsilon}_t(k_3):)
\\&(\overline{:\hat{X}^{\varepsilon}_{\bar{s}}(k_1)\hat{X}^{\varepsilon}_{\bar{s}}(k_2)\hat{X}^{\varepsilon}_t(k_3):
-:\hat{\bar{X}}^{\varepsilon}_{\bar{s}}(k_1)\hat{\bar{X}}^{\varepsilon}_{\bar{s}}(k_2)\hat{\bar{X}}^{\varepsilon}_t(k_3):})]\bigg|dsd\bar{s}
\endaligned$$
$$\aligned&+\sum_{i_1,i_2=1}^3\sum_{i_1',i_2'=1}^3\sum_k\theta(2^{-q}k)^2\sum_{k_{123}=k}\int_0^t\int_0^t\bigg|e^{-|k_{12}|^2(t-s) } k_{12}^{i_2}e^{-|k_{23}|^2(t-\bar{s})}  k_{23}^{i'_2}\hat{P}^{ii_1}(k_{12})\hat{P}^{ii'_1}(k_{12})\bigg|\\&E\bigg|(:\hat{X}^{\varepsilon}_s(k_1)\hat{X}^{\varepsilon}_s(k_2)\hat{X}^{\varepsilon}_t(k_3):-:\hat{\bar{X}}^{\varepsilon}_s(k_1)\hat{\bar{X}}^{\varepsilon}_s(k_2)\hat{\bar{X}}^{\varepsilon}_t(k_3):)
\\&(\overline{:\hat{X}^{\varepsilon}_{\bar{s}}(k_3)\hat{X}^{\varepsilon}_{\bar{s}}(k_2)\hat{X}^{\varepsilon}_t(k_1):-:\hat{\bar{X}}^{\varepsilon}_{\bar{s}}
(k_3)\hat{\bar{X}}^{\varepsilon}_{\bar{s}}(k_2)\hat{\bar{X}}^{\varepsilon}_t(k_1):})\bigg|dsd\bar{s}\\=&J_t^1+J_t^2+J_t^3+J_t^4,\endaligned$$
where we omit superscript of $\hat{X}^\varepsilon$ in $J^3, J^4$ for simplicity.
By (4.2) (4.3) we have for $\eta>0$ small enough
$$\aligned J_t^1\lesssim &\sum_k\theta(2^{-q}k)\sum_{k_{123}=k}\Pi_{i=1}^3\frac{1}{|k_i|^2}\frac{t^\eta|k_{12}|^\eta\varepsilon^\eta}{|k_{12}|^{2-2\eta}}\\\lesssim &\sum_k\theta(2^{-q}k)\sum_{k_{123}=k}\frac{t^\eta\varepsilon^\eta}{|k_3|^{2}|k_{12}|^{3-3\eta}}
\\\lesssim &t^\eta\varepsilon^\eta2^{q(1+3\eta)},\endaligned$$
and
$$\aligned J_t^2\lesssim &\sum_k\theta(2^{-q}k)\sum_{k_{123}=k}\frac{t^\eta\varepsilon^\eta}{|k_1|^{2}|k_{2}|^{2}|k_3|^{2}|k_{12}|^{1-3\eta/2}|k_{32}|^{1-3\eta/2}}\\\lesssim & \varepsilon^\eta\sum_k\theta(2^{-q}k)(\sum_{k_{123}=k}\frac{t^\eta}{|k_1|^{2}|k_2|^{2}|k_3|^{2}|k_{12}|^{2-3\eta}})^{1/2}(\sum_{k_{123}=k}\frac{t^\eta}{|k_1|^{2}|k_2|^{2}|k_3|^{2}|k_{32}|^{2-3\eta}})^{1/2}
\\\lesssim& \varepsilon^\eta t^\eta2^{q(1+3\eta)},\endaligned$$
where we used Lemma 4.1 in the last inequality.
Since
$$\aligned &E\big|(:\hat{X}^{\varepsilon}_s(k_1)\hat{X}^{\varepsilon}_s(k_2)\hat{X}^{\varepsilon}_t(k_3):-:\hat{\bar{X}}^{\varepsilon}_s(k_1)\hat{\bar{X}}^{\varepsilon}_s(k_2)\hat{\bar{X}}^{\varepsilon}_t(k_3):)
\\&(\overline{:\hat{X}^{\varepsilon}_{\bar{s}}(k_1)\hat{X}^{\varepsilon}_{\bar{s}}(k_2)\hat{X}^{\varepsilon}_t(k_3):-:\hat{\bar{X}}^{\varepsilon}_{\bar{s}}(k_1)
\hat{\bar{X}}^{\varepsilon}_{\bar{s}}(k_2)\hat{\bar{X}}^{\varepsilon}_t(k_3):})\big|
\\\lesssim&(E|:\hat{X}^{\varepsilon}_s(k_1)\hat{X}^{\varepsilon}_s(k_2)
\hat{X}^{\varepsilon}_t(k_3):-:\hat{\bar{X}}^{\varepsilon}_s(k_1)\hat{\bar{X}}^{\varepsilon}_s(k_2)\hat{\bar{X}}^{\varepsilon}_t(k_3):|^2)^{1/2}
\\&(E|:\hat{X}^{\varepsilon}_{\bar{s}}(k_1)\hat{X}^{\varepsilon}_{\bar{s}}(k_2)\hat{X}^{\varepsilon}_t(k_3):-:\hat{\bar{X}}^{\varepsilon}_{\bar{s}}(k_1)
\hat{\bar{X}}^{\varepsilon}_{\bar{s}}(k_2)\hat{\bar{X}}^{\varepsilon}_t(k_3):|^2)^{1/2}\\\lesssim&\bigg|\frac{1 }{8|k_1|^2|k_2|^2|k_3|^2f(\varepsilon k_1)f(\varepsilon k_2)f(\varepsilon k_3)}\\&-\frac{2 }{|k_1|^2|k_2|^2|k_3|^2(1+f(\varepsilon k_1))(1+f(\varepsilon k_2))(1+f(\varepsilon k_3))}+\frac{1 }{8|k_1|^2|k_2|^2|k_3|^2}\bigg|
\\\lesssim&\frac{(\varepsilon|k_1|)^\eta+(\varepsilon|k_2|)^\eta+(\varepsilon|k_3|)^\eta}{|k_1|^2|k_2|^2|k_3|^2},\endaligned\eqno(4.6)$$
we have $$\aligned J_t^3\lesssim &\sum_k\theta(2^{-q}k)\sum_{k_{123}=k}\Pi_{i=1}^3\frac{t^\eta }{|k_i|^2}\frac{(\varepsilon|k_1|)^\eta+(\varepsilon|k_2|)^\eta+(\varepsilon|k_3|)^\eta}{|k_{12}|^{2-2\eta}}\\\lesssim &\sum_k\theta(2^{-q}k)\sum_{k_{123}=k}\bigg(\frac{t^\eta\varepsilon^\eta}{|k_3|^{2}|k_{12}|^{3-3\eta}}
+\frac{t^\eta\varepsilon^\eta}{|k_3|^{2-\eta}|k_{12}|^{3-2\eta}}\bigg)
\lesssim t^\eta\varepsilon^\eta2^{q(1+3\eta)},\endaligned$$
where we used Lemma 4.1 in the second inequality.
By a similar argument as (4.6) and H\"{o}lder's inequality we have
$$\aligned J_t^4\lesssim &\sum_k\theta(2^{-q}k)\sum_{k_{123}=k}\frac{t^\eta\varepsilon^\eta(|k_1|^\eta+|k_2|^\eta+|k_3|^\eta)}{|k_1|^{2}|k_{2}|^{2}|k_3|^{2}|k_{12}|^{1-\eta}|k_{32}|^{1-\eta}}\\\lesssim & \varepsilon^\eta\sum_k\theta(2^{-q}k)(\sum_{k_{123}=k}\frac{t^\eta(|k_1|^\eta+|k_2|^\eta+|k_3|^\eta)}{|k_1|^{2}|k_2|^{2}|k_3|^{2}|k_{12}|^{2-2\eta}})^{1/2}(\sum_{k_{123}=k}\frac{t^\eta(|k_1|^\eta+|k_2|^\eta+|k_3|^\eta)}{|k_1|^{2}|k_2|^{2}|k_3|^{2}|k_{32}|^{2-2\eta}})^{1/2}
\\\lesssim& \varepsilon^\eta t^\eta2^{q(1+3\eta)},\endaligned$$
where we used Lemma 4.1 in the  last inequality. By a similar calculation as above  we also obtain that there exist $\eta,\epsilon,\gamma>0$ small enough such that
$$E[|\Delta_q(u_2^{\varepsilon,i}\diamond u_1^{\varepsilon,j}(t_1)-u_2^{\varepsilon,i}\diamond u_1^{\varepsilon,j}(t_2)-\bar{u}_2^{\varepsilon,i}\diamond \bar{u}_1^{\varepsilon,j}
(t_1)+\bar{u}_2^{\varepsilon,i}\diamond \bar{u}_1^{\varepsilon,j}(t_2))|^2]\lesssim \varepsilon^\gamma|t_1-t_2|^{\eta}2^{q(1+\epsilon)},$$
which by Gaussian hypercontractivity and Lemma 2.1 implies that for $\epsilon_0>0$ small enough and $p>1$
 $$\aligned &E[\|(u_2^{\varepsilon,i}\diamond u_1^{\varepsilon,j}(t_1)-u_2^{\varepsilon,i}\diamond u_1^{\varepsilon,j}(t_2)-\bar{u}_2^{\varepsilon,i}\diamond \bar{u}_1^{\varepsilon,j}
(t_1)+\bar{u}_2^{\varepsilon,i}\diamond \bar{u}_1^{\varepsilon,j}(t_2))\|^p_{\mathcal{C}^{-1/2-\epsilon/2-\epsilon_0-3/p}}]\\\lesssim &E[\|(u_2^{\varepsilon,i}\diamond u_1^{\varepsilon,j}(t_1)-u_2^{\varepsilon,i}\diamond u_1^{\varepsilon,j}(t_2)-\bar{u}_2^{\varepsilon,i}\diamond \bar{u}_1^{\varepsilon,j}
(t_1)+\bar{u}_2^{\varepsilon,i}\diamond \bar{u}_1^{\varepsilon,j}(t_2))\|^p_{B^{-1/2-\epsilon/2-\epsilon_0}_{p,p}}]\\\lesssim& \varepsilon^{p\gamma/2}|t_1-t_2|^{p\eta/2}.\endaligned$$
  Thus for every $i,j=1,2,3$, we choose $p$ large enough and deduce that for every $\delta>0, p>1$, $u_2^{\varepsilon,i} \diamond u_1^{\varepsilon,j}-\bar{u}_2^{\varepsilon,i} \diamond \bar{u}_1^{\varepsilon,j}\rightarrow0$ in $L^p(\Omega;C([0,T],\mathcal{C}^{-1/2-\delta/2}))$ .

\subsection{Convergence for $\pi_0(u_3^{\varepsilon,i_0},u_1^{\varepsilon,j_0})-\pi_0(\bar{u}_3^{\varepsilon,i_0},\bar{u}_1^{\varepsilon,j_0})$}
Now we treat $\pi_0(u_{31}^{\varepsilon, i_0},u_1^{\varepsilon, j_0})$ and the estimates for $\pi_0(u_3^{\varepsilon, i_0}-u_{31}^{\varepsilon, i_0},u_1^{\varepsilon, j_0})$ can be obtained similarly, where $(\partial_t-\Delta_\varepsilon)u_{31}^{\varepsilon,i}=-\frac{1}{2}\sum_{i_1=1}^3P^{ii_1}\sum_{j_1=1}^3D_{j_1}^\varepsilon(u_2^{\varepsilon,i_1}\diamond u_1^{\varepsilon, j_1}).$ We have the following identity:
$$\pi_0(u_{31}^{\varepsilon, i_0},u_1^{\varepsilon, j_0})-\pi_0(\bar{u}_{31}^{\varepsilon, i_0},\bar{u}_1^{\varepsilon, j_0})=\frac{1}{4}\bigg(\sum_{l=1}^7I_t^l-\bar{I}_t^5-\bar{I}_t^6\bigg),$$
where
$$\aligned I_t^1=&(2\pi)^{-9/2}\sum_{k\in \mathbb{Z}^3\backslash\{0\}}\sum_{|i-j|\leq1}\sum_{k_{1234}=k}\sum_{i_1,i_2,i_3,j_1=1}^3\theta(2^{-i} k_{123})\theta(2^{-j}k_4)\int_0^tds\bigg{[}e^{-|k_{123}|^2(t-s)f(\varepsilon k_{123})}\int_0^s:\hat{X}_\sigma^{\varepsilon ,i_2}(k_1)\\&\hat{X}^{\varepsilon, i_3}_\sigma(k_2)\hat{X}_s^{\varepsilon,j_1}(k_3)\hat{X}^{\varepsilon, j_0}_t(k_4):
e^{-|k_{12}|^2(s-\sigma)f(\varepsilon k_{12})}d\sigma  k_{12}^{i_3}g(\varepsilon k_{12}^{i_3}) k_{123}^{j_1}g(\varepsilon k_{123}^{j_1})-e^{-|k_{123}|^2(t-s)}\int_0^s:\hat{\bar{X}}_\sigma^{\varepsilon ,i_2}(k_1)\\&\hat{\bar{X}}^{\varepsilon, i_3}_\sigma(k_2)\hat{\bar{X}}_s^{\varepsilon,j_1}(k_3)\hat{\bar{X}}^{\varepsilon, j_0}_t(k_4):
e^{-|k_{12}|^2(s-\sigma)}d\sigma \imath k_{12}^{i_3}\imath k_{123}^{j_1}\bigg{]}\hat{P}^{i_1i_2}(k_{12})\hat{P}^{i_0i_1}(k_{123})e_k,\endaligned$$
 $$\aligned I_t^2=&(2\pi)^{-9/2}\sum_{k\in \mathbb{Z}^3\backslash\{0\}}\sum_{|i-j|\leq1}\sum_{k_{23}=k,k_1}\sum_{i_1,i_2,i_3,j_1=1}^3\theta(2^{-i} k_{123})\theta(2^{-j}k_1)\bigg{[}\int_0^tdse^{-|k_{123}|^2(t-s)f(\varepsilon k_{123})}\int_0^s:\hat{X}^{\varepsilon,i_3}_\sigma(k_2)\\&\hat{X}^{\varepsilon,j_1}_s(k_3)
:\frac{e^{-|k_1|^2(t-\sigma)f(\varepsilon k_1)}h(\varepsilon k_1)^2}{2|k_1|^2f(\varepsilon k_1)} k_{12}^{i_3}g(\varepsilon k_{12}^{i_3}) k_{123}^{j_1}g(\varepsilon k_{123}^{j_1})e^{-|k_{12}|^2(s-\sigma)f(\varepsilon k_{12})}d\sigma-\int_0^tdse^{-|k_{123}|^2(t-s)}\\&\int_0^s:\hat{\bar{X}}^{\varepsilon,i_3}_\sigma(k_2)\hat{\bar{X}}^{\varepsilon,j_1}_s(k_3)
:\frac{e^{-|k_1|^2(t-\sigma)}h(\varepsilon k_1)^2}{2|k_1|^2}\imath k_{12}^{i_3}\imath k_{123}^{j_1}e^{-|k_{12}|^2(s-\sigma)}d\sigma\bigg{]}\\&\sum_{i_4=1}^3\hat{P}^{i_2i_4}(k_1)\hat{P}^{j_0i_4}(k_1)\hat{P}^{i_1i_2}(k_{12})\hat{P}^{i_0i_1}(k_{123})e_k,\endaligned$$
$$\aligned I_t^3=&(2\pi)^{-9/2}\sum_{k\in \mathbb{Z}^3\backslash\{0\}}\sum_{|i-j|\leq1}\sum_{k_{23}=k,k_1}\sum_{i_1,i_2,i_3,j_1=1}^3\theta(2^{-i} k_{123})\theta(2^{-j}k_1)\bigg{[}\int_0^tdse^{-|k_{123}|^2(t-s)f(\varepsilon k_{123})}\int_0^s:\hat{X}^{\varepsilon,i_2}_\sigma(k_2)\\&\hat{X}^{\varepsilon,j_1}_s(k_3)
:\frac{e^{-|k_1|^2(t-\sigma)f(\varepsilon k_1)}h(\varepsilon k_1)^2}{2|k_1|^2f(\varepsilon k_1)} k_{12}^{i_3}g(\varepsilon k_{12}^{i_3}) k_{123}^{j_1}g(\varepsilon k_{123}^{j_1})e^{-|k_{12}|^2(s-\sigma)f(\varepsilon k_{12})}d\sigma-\int_0^tdse^{-|k_{123}|^2(t-s)}\\&\int_0^s:\hat{\bar{X}}^{\varepsilon,i_2}_\sigma(k_2)\hat{\bar{X}}^{\varepsilon,j_1}_s(k_3)
:\frac{e^{-|k_1|^2(t-\sigma)}h(\varepsilon k_1)^2}{2|k_1|^2}\imath k_{12}^{i_3}\imath k_{123}^{j_1}e^{-|k_{12}|^2(s-\sigma)}d\sigma\bigg{]}\\&\sum_{i_4=1}^3\hat{P}^{i_3i_4}(k_1)
\hat{P}^{j_0i_4}(k_1)\hat{P}^{i_1i_2}(k_{12})\hat{P}^{i_0i_1}(k_{123})e_k,\endaligned$$

$$\aligned I_t^4=&(2\pi)^{-9/2}\sum_{k\in \mathbb{Z}^3\backslash\{0\}}\sum_{|i-j|\leq1}\sum_{k_{12}=k,k_3}\sum_{i_1,i_2,i_3,j_1=1}^3\theta(2^{-i} k_{123})\theta(2^{-j}k_3)\bigg{[}\int_0^tdse^{-|k_{123}|^2(t-s)f(\varepsilon k_{123})}\int_0^s:\hat{X}_\sigma^{\varepsilon,i_2}(k_1)\\&\hat{X}_\sigma^{\varepsilon,i_3}(k_2):\frac{e^{-|k_3|^2(t-s)f(\varepsilon k_3)}h(\varepsilon k_3)^2}{2|k_3|^2f(\varepsilon k_3)}e^{-|k_{12}|^2(s-\sigma)f(\varepsilon k_{12})}d\sigma  k_{12}^{i_3}g(\varepsilon k_{12}^{i_3}) k_{123}^{j_1}g(\varepsilon k_{123}^{j_1}) -\int_0^tdse^{-|k_{123}|^2(t-s)}\\&\int_0^s:\hat{\bar{X}}_\sigma^{\varepsilon,i_2}(k_1)\hat{\bar{X}}_\sigma^{\varepsilon,i_3}(k_2):\frac{e^{-|k_3|^2(t-s)}h(\varepsilon k_3)^2}{2|k_3|^2}e^{-|k_{12}|^2(s-\sigma)}d\sigma \imath k_{12}^{i_3}\imath k_{123}^{j_1}\bigg{]}\\&\sum_{i_4=1}^3\hat{P}^{j_1i_4}(k_3)\hat{P}^{j_0i_4}(k_3)\hat{P}^{i_1i_2}(k_{12})\hat{P}^{i_0i_1}(k_{123})e_k,\endaligned$$
$$\aligned I_t^5=&(2\pi)^{-9/2}\sum_{k\in \mathbb{Z}^3\backslash\{0\}}\sum_{|i-j|\leq1}\sum_{k_{14}=k,k_2}\sum_{i_1,i_2,i_3,j_1=1}^3\theta(2^{-i} k_{1})\theta(2^{-j}k_4)\bigg{[}\int_0^tdse^{-|k_{1}|^2(t-s)f(\varepsilon k_1)}\int_0^s:\hat{X}^{\varepsilon,i_2}_\sigma(k_1)\\&\hat{X}^{\varepsilon,j_0}_t(k_4):\frac{e^{-|k_2|^2(s-\sigma)f(\varepsilon k_2)}h(\varepsilon k_2)^2}{2|k_2|^2f(\varepsilon k_2)}e^{-|k_{12}|^2(s-\sigma)f(\varepsilon k_{12})}d\sigma k_{12}^{i_3}g(\varepsilon k_{12}^{i_3}) k_{1}^{j_1}g(\varepsilon k_{1}^{j_1})\\&-\int_0^tdse^{-|k_{1}|^2(t-s)}\int_0^s:\hat{\bar{X}}^{\varepsilon,i_2}_\sigma(k_1)\hat{\bar{X}}^{\varepsilon,j_0}_t(k_4)
:\frac{e^{-|k_2|^2(s-\sigma)}h(\varepsilon k_2)^2}{2|k_2|^2}e^{-|k_{12}|^2(s-\sigma)}d\sigma \imath k_{12}^{i_3}\imath k_{1}^{j_1}\bigg{]}\\&\sum_{i_4=1}^3\hat{P}^{i_3i_4}(k_2)\hat{P}^{j_1i_4}(k_2)\hat{P}^{i_1i_2}(k_{12})\hat{P}^{i_0i_1}(k_{1})e_k,\endaligned$$
$$ \bar{I}_t^5=2\pi_0(\sum_{i_1,i_2,j_1=1}^3\int_0^tP_{t-s}^\varepsilon P^{i_0i_1} D_{j_1}^\varepsilon X_s^{\varepsilon,i_2} C^{\varepsilon,i_1,i_2,j_1}(s)ds, X^{\varepsilon,j_0}_t),$$
$$\aligned I_t^6=&(2\pi)^{-9/2}\sum_{k\in \mathbb{Z}^3\backslash\{0\}}\sum_{|i-j|\leq1}\sum_{k_{14}=k,k_2}\sum_{i_1,i_2,i_3,j_1=1}^3\theta(2^{-i} k_{1})\theta(2^{-j}k_4)\bigg{[}\int_0^tdse^{-|k_{1}|^2(t-s)f(\varepsilon k_1)}\int_0^s:\hat{X}^{\varepsilon,i_3}_\sigma(k_1)\\&\hat{X}^{\varepsilon,j_0}_t(k_4):\frac{e^{-|k_2|^2(s-\sigma)f(\varepsilon k_2)}h(\varepsilon k_2)^2}{2|k_2|^2f(\varepsilon k_2)}e^{-|k_{12}|^2(s-\sigma)f(\varepsilon k_{12})}d\sigma k_{12}^{i_3}g(\varepsilon k_{12}^{i_3}) k_{1}^{j_1}g(\varepsilon k_{1}^{j_1})\\&-\int_0^tdse^{-|k_{1}|^2(t-s)}\int_0^s:\hat{\bar{X}}^{\varepsilon,i_3}_\sigma(k_1)\hat{\bar{X}}^{\varepsilon,j_0}_t(k_4)
:\frac{e^{-|k_2|^2(s-\sigma)}h(\varepsilon k_2)^2}{2|k_2|^2}e^{-|k_{12}|^2(s-\sigma)}d\sigma \imath k_{12}^{i_3}\imath k_{1}^{j_1}\bigg{]}\\&\sum_{i_4=1}^3\hat{P}^{i_2i_4}(k_2)\hat{P}^{j_1i_4}(k_2)\hat{P}^{i_1i_2}(k_{12})\hat{P}^{i_0i_1}(k_{1})e_k,\endaligned$$
$$ \bar{I}_t^6=2\pi_0(\sum_{i_1,i_2,j_1=1}^3\int_0^tP_{t-s}^\varepsilon P^{i_0i_1} D_{j_1}^\varepsilon X_s^{\varepsilon,i_2} \tilde{C}^{\varepsilon,i_1,i_2,j_1}(s)ds, X^{\varepsilon,j_0}_t),$$
$$\aligned
I_t^7=&(2\pi)^{-6}\sum_{|i-j|\leq1}\sum_{k_1,k_2}\sum_{i_1,i_2,i_3,j_1=1}^3\theta(2^{-i} k_2)\theta(2^{-j}k_2)\bigg{[}\int_0^tdse^{-|k_2|^2(t-s)f(\varepsilon k_2)}\int_0^s\frac{h(\varepsilon k_1)^2h(\varepsilon k_2)^2}{4|k_1|^2|k_2|^2f(\varepsilon k_1)f(\varepsilon k_2)}\\&e^{-|k_{12}|^2(s-\sigma)f(\varepsilon k_{12}) -|k_1|^2(s-\sigma)f(\varepsilon k_1)-|k_2|^2(t-\sigma)f(\varepsilon k_2)}d\sigma g(\varepsilon k_{12}^{i_3}) k_{12}^{i_3} k_{2}^{j_1}g(\varepsilon k_{2}^{j_1})\\&-\int_0^tdse^{-|k_2|^2(t-s)}\int_0^s\frac{h(\varepsilon k_1)^2h(\varepsilon k_2)^2}{4|k_1|^2|k_2|^2}e^{-|k_{12}|^2(s-\sigma)-|k_1|^2(s-\sigma)-|k_2|^2(t-\sigma)}d\sigma \imath k_{12}^{i_3}\imath k_{2}^{j_1}\bigg{]}\\&\sum_{i_4,i_5=1}^3\big(\hat{P}^{i_3i_4}(k_1)\hat{P}^{j_1i_4}(k_1)\hat{P}^{i_2i_5}(k_2)\hat{P}^{j_0i_5}(k_2)
+\hat{P}^{i_2i_4}(k_1)\hat{P}^{j_1i_4}(k_1)\hat{P}^{i_3i_5}(k_2)\hat{P}^{j_0i_5}(k_2)\big)\hat{P}^{i_1i_2}
(k_{12})\hat{P}^{i_0i_1}(k_{2}).\endaligned$$

First we consider $I_t^7$: by simple calculations we have
$$\aligned I_t^7=&(2\pi)^{-6}\sum_{|i-j|\leq1}\sum_{k_1,k_2}\sum_{i_1,i_2,i_3,j_1=1}^3\theta(2^{-i} k_2)\theta(2^{-j}k_2)\hat{P}^{i_1i_2}(k_{12})\hat{P}^{i_0i_1}(k_{2})\\&\sum_{i_4,i_5=1}^3\big(\hat{P}^{i_3i_4}(k_1)\hat{P}^{j_1i_4}(k_1)\hat{P}^{i_2i_5}(k_2)\hat{P}^{j_0i_5}(k_2)
+\hat{P}^{i_2i_4}(k_1)\hat{P}^{j_1i_4}(k_1)\hat{P}^{i_3i_5}(k_2)\hat{P}^{j_0i_5}(k_2)\big)\\&\bigg{[}\frac{h(\varepsilon k_1)^2h(\varepsilon k_2)^2}{4|k_1|^2f(\varepsilon k_1)|k_2|^2f(\varepsilon k_2)(|k_1|^2f(\varepsilon k_1)+|k_2|^2f(\varepsilon k_2)+|k_{12}|^2f(\varepsilon k_{12}))} k_{12}^{i_3}g(\varepsilon k_{12}^{i_3}) k_{2}^{j_1}g(\varepsilon k_{2}^{j_1})\\&(\frac{1-e^{-2|k_2|^2f(\varepsilon k_2)t}}{2|k_2|^2f(\varepsilon k_2)}-\int_0^tdse^{-2|k_2|^2(t-s)f(\varepsilon k_2)}e^{-(|k_{12}|^2f(\varepsilon k_{12})+|k_1|^2f(\varepsilon k_1)+|k_2|^2f(\varepsilon k_2))s})\\&-\frac{h(\varepsilon k_1)^2h(\varepsilon k_2)^2}{4|k_1|^2|k_2|^2(|k_1|^2+|k_2|^2+|k_{12}|^2)}\imath k_{12}^{i_3}\imath k_{2}^{j_1}(\frac{1-e^{-2|k_2|^2t}}{2|k_2|^2}-\int_0^tdse^{-2|k_2|^2(t-s)}e^{-(|k_{12}|^2+|k_1|^2+|k_2|^2)s})\bigg{]}. \endaligned$$
Let $$\aligned &C_{11}^{\varepsilon,i_0,j_0}(t)\\=&(2\pi)^{-6}\sum_{|i-j|\leq1}\sum_{k_1,k_2}\sum_{i_1,i_2,i_3,j_1=1}^3\theta(2^{-i} k_2)\theta(2^{-j}k_2)k_{12}^{i_3}g(\varepsilon k_{12}^{i_3})  k_{2}^{j_1}g(\varepsilon  k_{2}^{j_1})\hat{P}^{i_1i_2}(k_{12})\hat{P}^{i_0i_1}(k_{2})\\&\frac{h(\varepsilon k_1)^2h(\varepsilon k_2)^2}{4|k_1|^2f(\varepsilon k_2)|k_2|^2f(\varepsilon k_2)(|k_1|^2f(\varepsilon k_1)+|k_2|^2f(\varepsilon k_2)+|k_{12}|^2f(\varepsilon k_{12}))}\frac{1-e^{-2|k_2|^2tf(\varepsilon k_{2})}}{2|k_2|^2f(\varepsilon k_{2})}\\&
\sum_{i_4,i_5=1}^3\big(\hat{P}^{i_3i_4}(k_1)\hat{P}^{j_1i_4}(k_1)\hat{P}^{i_2i_5}(k_2)\hat{P}^{j_0i_5}(k_2)
+\hat{P}^{i_2i_4}(k_1)\hat{P}^{j_1i_4}(k_1)\hat{P}^{i_3i_5}(k_2)\hat{P}^{j_0i_5}(k_2)\big),\endaligned$$
and $$\aligned &\bar{C}_{11}^{\varepsilon,i_0,j_0}(t)\\=&(2\pi)^{-6}\sum_{|i-j|\leq1}\sum_{k_1,k_2}\sum_{i_1,i_2,i_3,j_1=1}^3\theta(2^{-i} k_2)\theta(2^{-j}k_2)\imath k_{12}^{i_3}\imath k_{2}^{j_1}\hat{P}^{i_1i_2}(k_{12})\hat{P}^{i_0i_1}(k_{2})\frac{h(\varepsilon k_1)^2h(\varepsilon k_2)^2}{4|k_1|^2|k_2|^2(|k_1|^2+|k_2|^2+|k_{12}|^2)}\\&\frac{1-e^{-2|k_2|^2t}}{2|k_2|^2}
\sum_{i_4,i_5=1}^3\big(\hat{P}^{i_3i_4}(k_1)\hat{P}^{j_1i_4}(k_1)\hat{P}^{i_2i_5}(k_2)\hat{P}^{j_0i_5}(k_2)
+\hat{P}^{i_2i_4}(k_1)\hat{P}^{j_1i_4}(k_1)\hat{P}^{i_3i_5}(k_2)\hat{P}^{j_0i_5}(k_2)\big)\rightarrow\infty,\endaligned$$
as $\varepsilon\rightarrow0$.
Define $$\varphi^{\varepsilon,i_0j_0}_{11}-\bar{\varphi}^{\varepsilon,i_0j_0}_{11}:=I_t^7-C_{11}^{\varepsilon,i_0j_0}-\bar{C}_{11}^{\varepsilon,i_0j_0},$$
where $\varphi^\varepsilon, \bar{\varphi}^\varepsilon$ correspond to $u^\varepsilon, \bar{u}^\varepsilon$ respectively.
Then for every $\rho>0, \eta>0$ we deduce that  $$\aligned |\varphi^{\varepsilon,i_0,j_0}_{11}-\bar{\varphi}^{\varepsilon,i_0j_0}_{11}|\lesssim&\bigg|\sum_{k_1,k_2}(\frac{k_{12}^{i_3}k_2^{j_1}g(\varepsilon k_{12}^{i_3})g(\varepsilon k_2^{j_1}) }{|k_1|^2f(\varepsilon k_1)|k_2|^{2}f(\varepsilon k_2)(|k_1|^2f(\varepsilon k_1)+|k_2|^2f(\varepsilon k_2)+|k_{12}|^2f(\varepsilon k_{12}))}\\&\int_0^te^{-(|k_2|^2(2t-s)f(\varepsilon k_2)+|k_1|^2sf(\varepsilon k_1)+|k_{12}|^2sf(\varepsilon k_{12}))}ds-\frac{\imath k_{12}^{i_3}\imath k_2^{j_1}}{|k_1|^2|k_2|^{2}(|k_1|^2+|k_2|^2+|k_{12}|^2)}\\&\int_0^te^{-(|k_2|^2(2t-s)+|k_1|^2s+|k_{12}|^2s)}ds)\bigg|
\\\lesssim &t^{-\rho-\eta/2}\sum_{k_1,k_2}\frac{\varepsilon^\eta}{|k_1|^{3+r}|k_2|^{3+2\rho-r}}\lesssim \varepsilon^\eta t^{-\rho-\eta/2}.\endaligned$$ Here $r>0$ is small enough such that $2\rho>r>0$.
Similarly, we can also find similar $C_{12}^\varepsilon, \varphi_{12}^\varepsilon,\bar{\varphi}_{12}^\varepsilon, \bar{C}_{12}^\varepsilon$ for $u_3-u_{31}$  and satisfy similar estimates as $\varphi_{11}^\varepsilon$ and $\bar{\varphi}_{11}^\varepsilon$. Now define $C_1^\varepsilon=C_{11}^\varepsilon+C_{12}^\varepsilon$, $\varphi_1^\varepsilon= \varphi_{11}^\varepsilon+\varphi_{12}^\varepsilon$, $\bar{C}_1^\varepsilon=\bar{C}_{11}^\varepsilon+\bar{C}_{12}^\varepsilon$, $\bar{\varphi}_1^\varepsilon= \bar{\varphi}_{11}^\varepsilon+\bar{\varphi}_{12}^\varepsilon$ where we omit superscript for simplicity.

\textbf{Terms in the second chaos}: We come to $I_t^2$ and have the following calculations:
$$\aligned &E|\Delta_qI_t^2|^2\\\lesssim&\sum_{k\in \mathbb{Z}^3\backslash\{0\}}\sum_{|i-j|\leq1,|i'-j'|\leq1}\sum_{k_{23}=k,k_1,k_4}
\sum_{i_1,i_2,i_3,j_1=1}^3\sum_{i_1',i_2',i_3',j_1'=1}^3\theta(2^{-i} k_{123})\theta(2^{-j}k_1)
\theta(2^{-i'} k_{234})\theta(2^{-j'}k_4)\theta(2^{-q}k)^2\\&\bigg{[}\Pi_{i=1}^4\frac{h(\varepsilon k_i)^2}{|k_i|^2f(\varepsilon k_i)}\int_0^t\int_0^tdsd\bar{s}\big|(e^{-|k_{123}|^2(t-s)f(\varepsilon k_{123})}k_{123}^{j_1}g(\varepsilon k_{123}^{j_1})-e^{-|k_{123}|^2(t-s)}\imath k_{123}^{j_1})\\&(e^{-|k_{234}|^2f(\varepsilon k_{234})(t-\bar{s})}k_{234}^{j_1'}g(\varepsilon k_{234}^{j_1'})-e^{-|k_{234}|^2(t-\bar{s})}\imath k_{234}^{j_1'})\big|\int_0^s\int_0^{\bar{s}}d\sigma d\bar{\sigma}e^{-|k_1|^2f(\varepsilon k_1)(t-\sigma)-|k_4|^2f(\varepsilon k_4)(t-\bar{\sigma})}\\&(e^{-(|k_{12}|^2(s-\sigma)f(\varepsilon k_{12})+|k_{24}|^2(\bar{s}-\bar{\sigma}))f(\varepsilon k_{24}))} |k_{24}k_{12}|+e^{-(|k_{12}|^2(s-\sigma)f(\varepsilon k_{12})+|k_{34}|^2(\bar{s}-\bar{\sigma})f(\varepsilon k_{34}))} |k_{34}k_{12}|)
\\&+\frac{1}{|k_2|^2f(\varepsilon k_2)|k_3|^2f(\varepsilon k_3)}\int_0^t\int_0^tdsd\bar{s}e^{-|k_{123}|^2(t-s)}|k_{123}|e^{-|k_{234}|^2(t-\bar{s})} |k_{234}|\int_0^s\int_0^{\bar{s}}d\sigma d\bar{\sigma}\\&\bigg|(\frac{e^{-|k_1|^2f(\varepsilon k_1)(t-\sigma)}}{|k_1|^2f(\varepsilon k_1)}-\frac{e^{-|k_1|^2(t-\sigma)}}{|k_1|^2})(\frac{e^{-|k_4|^2f(\varepsilon k_4)(t-\bar{\sigma})}}{|k_4|^2f(\varepsilon k_4)}-\frac{e^{-|k_4|^2(t-\bar{\sigma})}}{|k_4|^2})\bigg|\\&(e^{-(|k_{12}|^2(s-\sigma)f(\varepsilon k_{12})+|k_{24}|^2(\bar{s}-\bar{\sigma}))f(\varepsilon k_{24}))} |k_{24}k_{12}|+e^{-(|k_{12}|^2(s-\sigma)f(\varepsilon k_{12})+|k_{34}|^2(\bar{s}-\bar{\sigma}))f(\varepsilon k_{34}))} |k_{34}k_{12}|)\\&+\Pi_{i=1}^4\frac{1}{|k_i|^2}\int_0^t\int_0^tdsd\bar{s}e^{-|k_{123}|^2(t-s)} |k_{123}|e^{-|k_{234}|^2(t-\bar{s})} |k_{234}|\int_0^s\int_0^{\bar{s}}d\sigma d\bar{\sigma}e^{-|k_1|^2(t-\sigma)-|k_4|^2(t-\bar{\sigma})}\\&|(e^{-|k_{12}|^2(s-\sigma)f(\varepsilon k_{12})}k_{12}^{i_3}g(\varepsilon k_{12}^{i_3})-\imath e^{-|k_{12}|^2(s-\sigma)}k_{12}^{i_3})\\&[(e^{-|k_{24}|^2(\bar{s}-\bar{\sigma})f(\varepsilon k_{24})}k_{24}^{i_3'}g(\varepsilon k_{24}^{i_3'})-\imath e^{-|k_{24}|^2(\bar{s}-\bar{\sigma})}k_{24}^{i_3'})+(e^{-|k_{34}|^2(\bar{s}-\bar{\sigma})f(\varepsilon k_{34})}k_{34}^{i_3'}g(\varepsilon k_{34}^{i_3'})-\imath e^{-|k_{34}|^2(\bar{s}-\bar{\sigma})}k_{34}^{i_3'})]|
\\&+\frac{1}{|k_1|^2|k_4|^2}\int_0^t\int_0^tdsd\bar{s}e^{-|k_{123}|^2(t-s)}|k_{123}|e^{-|k_{234}|^2(t-\bar{s})}|k_{234}|\int_0^s\int_0^{\bar{s}}d\sigma d\bar{\sigma}e^{-|k_1|^2(t-\sigma)-|k_4|^2(t-\bar{\sigma})}\\&\bigg(e^{-(|k_{12}|^2(s-\sigma)+|k_{24}|^2(\bar{s}-\bar{\sigma}))} |k_{24}k_{12}|E\bigg|(:\hat{X}^\varepsilon_\sigma(k_2)\hat{X}^\varepsilon_s(k_3):-:\hat{\bar{X}}^\varepsilon_\sigma(k_2)\hat{\bar{X}}^\varepsilon_s(k_3):)
\\&(\overline{:\hat{X}^\varepsilon_{\bar{\sigma}}(k_2)\hat{X}^\varepsilon_{\bar{s}}(k_3):
-:\hat{\bar{X}}^\varepsilon_{\bar{\sigma}}(k_2)\hat{\bar{X}}^\varepsilon_{\bar{s}}(k_3):})\bigg|
+e^{-(|k_{12}|^2(s-\sigma)+|k_{34}|^2(\bar{s}-\bar{\sigma}))} |k_{34}k_{12}|\\&E\bigg|(:\hat{X}^\varepsilon_\sigma(k_2)\hat{X}^\varepsilon_s(k_3):-:\hat{\bar{X}}^\varepsilon_\sigma(k_2)\hat{\bar{X}}^\varepsilon_s(k_3):)
(\overline{:\hat{X}^\varepsilon_{\bar{\sigma}}(k_3)\hat{X}^\varepsilon_{\bar{s}}(k_2):
-:\hat{\bar{X}}^\varepsilon_{\bar{\sigma}}(k_3)\hat{\bar{X}}^\varepsilon_{\bar{s}}(k_2):})\bigg|\bigg)\bigg{]}
\\\lesssim& \sum_{k\in \mathbb{Z}^3\backslash\{0\}}\sum_{|i-j|\leq1,|i'-j'|\leq1}\sum_{k_{23}=k,k_1,k_4}\theta(2^{-i} k_{123})\theta(2^{-j}k_1)\theta(2^{-i'} k_{234})\theta(2^{-j'}k_4)\theta(2^{-q}k)^2\\&\frac{t^\eta(|\varepsilon k_{123}|^{\eta/2}|\varepsilon k_{234}|^{\eta/2}+|\varepsilon k_{12}|^{\eta/2}|\varepsilon k_{24}|^{\eta/2}+(|\varepsilon k_{2}|^{\eta/2}+|
\varepsilon k_{3}|^{\eta/2})^2+|\varepsilon k_{1}|^{\eta/2}|
\varepsilon k_{4}|^{\eta/2})}{|k_2|^2|k_3|^2|k_1|^2(|k_1|+|k_{123}|)^{1-\eta}|k_4|^2
(|k_4|+|k_{234}|)^{1-\eta}(|k_4|+|k_{24}|)(|k_1|+|k_{12}|)}\\&+\frac{t^\eta(|\varepsilon k_{123}|^{\eta/2}|\varepsilon k_{234}|^{\eta/2}+|\varepsilon k_{12}|^{\eta/2}|\varepsilon k_{34}|^{\eta/2}+(|\varepsilon k_{2}|^{\eta/2}+|
\varepsilon k_{3}|^{\eta/2})^2+|\varepsilon k_{1}|^{\eta/2}|
\varepsilon k_{4}|^{\eta/2})}{|k_2|^2|k_3|^2|k_1|^2(|k_1|+|k_{123}|)^{1-\eta}|k_4|^2
(|k_4|+|k_{234}|)^{1-\eta}(|k_4|+|k_{34}|)(|k_1|+|k_{12}|)}
\\\lesssim& \varepsilon^\eta\sum_{k\in \mathbb{Z}^3\backslash\{0\}}\bigg(\sum_{q\lesssim i}2^{-(1-\eta-\epsilon)i}\sum_{q\lesssim i'}2^{-(1-\eta-\epsilon)i'}\sum_{k_{23}=k}\theta(2^{-q}k)^2\bigg(\frac{t^\eta}{|k_2|^{2-\eta}|k_3|^{2}}+\frac{t^\eta}{|k_2|^{2}|k_3|^{2-\eta}}\bigg)\\&+\sum_{q\lesssim i}2^{-(1-3\eta/2-\epsilon)i}\sum_{q\lesssim i'}2^{-(1-3\eta/2-\epsilon)i'}\sum_{k_{23}=k}\theta(2^{-q}k)^2\frac{t^\eta}{|k_2|^2|k_3|^{2}}\bigg)\lesssim \varepsilon^\eta t^\eta2^{q(3\eta+2\epsilon)},\endaligned$$
where $\eta,\epsilon>0$ are small enough. Here we used (4.2) (4.3) and the following in the  second inequality
$$\aligned &E|:\hat{X}^\varepsilon_\sigma(k_2)\hat{X}^\varepsilon_s(k_3):-:\hat{\bar{X}}^\varepsilon_\sigma(k_2)\hat{\bar{X}}^\varepsilon_s(k_3):|^2
\\\lesssim &\left|\frac{1}{4|k_2|^2|k_3|^2f(\varepsilon k_2)f(\varepsilon k_3)}-\frac{2}{|k_2|^2|k_3|^2(f(\varepsilon k_2)+1)(f(\varepsilon k_3)+1)}+\frac{1}{4|k_2|^2|k_3|^2}\right|\\\lesssim&\frac{(|\varepsilon k_2|^{\eta/2}+|\varepsilon k_3|^{\eta/2})^2}{|k_2|^2|k_3|^2},\endaligned\eqno(4.7)$$
and we used Lemma 4.1 in the last inequality
 and $q\lesssim i$ follows from $|k|\leq |k_{123}|+|k_1|\lesssim 2^{i}$ and $q\lesssim i'$ is similar.
For $I_t^3$ we have a similar estimate.

Now we deal with $I_t^4-\tilde{I}_t^4+\tilde{I}_t^4+4\sum_{i_1=1}^3u^{\varepsilon,i_1}_2C^{\varepsilon,i_0,i_1,j_0}(t)
-4\sum_{i_1=1}^3\bar{u}^{\varepsilon,i_1}_2\bar{C}^{\varepsilon,i_0,i_1,j_0}(t)$ with
$$\aligned \tilde{I}_t^4=&(2\pi)^{-\frac{9}{2}}\sum_{k\in \mathbb{Z}^3\backslash\{0\}}\sum_{|i-j|\leq1}\sum_{k_{12}=k,k_3}\sum_{i_1,i_2,i_3,j_1=1}^3\theta(2^{-i} k_{123})\theta(2^{-j}k_3)\bigg{[}\int_0^t:\hat{X}^{\varepsilon,i_2}_\sigma(k_1)
\hat{X}^{\varepsilon,i_3}_\sigma(k_2):e^{-|k_{12}|^2f(\varepsilon k_{12})(t-\sigma)} \\&k_{12}^{i_3}g(\varepsilon k_{12}^{i_3}) d\sigma\int_0^tdse^{-|k_{123}|^2(t-s)f(\varepsilon k_{123})}\frac{e^{-|k_3|^2(t-s)f(\varepsilon k_3)}h(\varepsilon k_3)^2}{2|k_3|^2f(\varepsilon k_3)} k_{123}^{j_1}g(\varepsilon k_{123}^{j_1})\\&-\int_0^t:\hat{\bar{X}}^{\varepsilon,i_2}_\sigma(k_1)
\hat{\bar{X}}^{\varepsilon,i_3}_\sigma(k_2):e^{-|k_{12}|^2(t-\sigma)}\imath k_{12}^{i_3} d\sigma\int_0^tdse^{-|k_{123}|^2(t-s)}\frac{e^{-|k_3|^2(t-s)}h(\varepsilon k_3)^2}{2|k_3|^2}\imath k_{123}^{j_1}\bigg{]}\\&\sum_{i_4}\hat{P}^{j_1i_4}(k_3)\hat{P}^{j_0i_4}(k_3)\hat{P}^{i_0i_1}(k_{123})\hat{P}^{i_1i_2}(k_{12})e_k ,\endaligned$$
and $$\aligned C^{\varepsilon,i_0,i_1,j_0}(t)=&\frac{(2\pi)^{-3}}{4}\sum_{|i-j|\leq1}\sum_{k_3}\sum_{j_1=1}^3\theta(2^{-i} k_{3})\theta(2^{-j}k_3)\int_0^tds\frac{e^{-2|k_3|^2(t-s)f(\varepsilon k_3)}h(\varepsilon k_3)^2}{|k_3|^2f(\varepsilon k_3)}\\& \sum_{i_4=1}^3\hat{P}^{j_1i_4}(k_3)\hat{P}^{j_0i_4}(k_3)k_{3}^{j_1}g(\varepsilon k_{3}^{j_1})\hat{P}^{i_0i_1}(k_{3})\\=&\frac{(2\pi)^{-3}}{4}\sum_{k_3}\sum_{j_1=1}^3\int_0^tds\frac{e^{-2|k_3|^2(t-s)f(\varepsilon k_3)}h(\varepsilon k_3)^2}{|k_3|^2f(\varepsilon k_3)} \sum_{i_4=1}^3\hat{P}^{j_1i_4}(k_3)\hat{P}^{j_0i_4}(k_3)k_{3}^{j_1}g(\varepsilon k_{3}^{j_1})\hat{P}^{i_0i_1}(k_{3})\\\rightarrow&(2\pi)^{-3}\sum_{j_1=1}^3\sum_{i_4=1}^3\int_{\mathbb{R}^3}\frac{(\cos (ax^{j_1})-\cos (bx^{j_1}))h(x)^2}{8|x|^4f(x)^2(a+b)} \hat{P}^{j_1i_4}(x)\hat{P}^{j_0i_4}(x)\hat{P}^{i_0i_1}(x)dx,\quad\varepsilon\rightarrow0,\endaligned$$
where in the second equality we used that $\theta(2^{-i}\cdot)\theta(2^{-j}\cdot)=0$ if $|i-j|>1$,
and $$\aligned \bar{C}^{\varepsilon,i_0,i_1,j_0}(t)=&(2\pi)^{-3}\sum_{|i-j|\leq1}\sum_{k_3}\sum_{j_1=1}^3\theta(2^{-i} k_{3})\theta(2^{-j}k_3)\int_0^tds\frac{e^{-2|k_3|^2(t-s)}h(\varepsilon k_3)^2}{4|k_3|^2}\\& \sum_{i_4=1}^3\hat{P}^{j_1i_4}(k_3)\hat{P}^{j_0i_4}(k_3)\imath k_{3}^{j_1}\hat{P}^{i_0i_1}(k_{3})\\=&0.\endaligned$$
Let $$c_{k_{123},k_3}^{\varepsilon,j_1}(t-s)=\sum_{i_1=1}^3e^{-|k_{123}|^2(t-s)f(\varepsilon k_{123})}\frac{e^{-|k_3|^2(t-s)f(\varepsilon k_3)}h(\varepsilon k_3)^2}{|k_3|^2f(\varepsilon k_3)} k_{123}^{j_1}g(\varepsilon k_{123}^{j_1})\hat{P}^{i_0i_1}(k_{123})$$
and $$\bar{c}_{k_{123},k_3}^{\varepsilon,j_1}(t-s)=\sum_{i_1=1}^3e^{-|k_{123}|^2(t-s)}\frac{e^{-|k_3|^2(t-s)}h(\varepsilon k_3)^2}{|k_3|^2} \imath k_{123}^{j_1}\hat{P}^{i_0i_1}(k_{123}).$$ Since by (4.2) and (4.3) for $t>s>\sigma$
$$\aligned &|e^{-|k_{12}|^2(s-\sigma)f(\varepsilon k_{12})}g(\varepsilon k_{12}^{i_3})k_{12}^{i_3}
-e^{-|k_{12}|^2(s-\sigma)}\imath k_{12}^{i_3}-e^{-|k_{12}|^2f(\varepsilon k_{12})(t-\sigma)}g(\varepsilon k_{12}^{i_3})k_{12}^{i_3}
+e^{-|k_{12}|^2(t-\sigma)}k_{12}^{i_3}\imath|\\\lesssim& (|\varepsilon k_{12}|^{2\eta}\wedge (t-s)^{1/2}|k_{12}|)e^{-|k_{12}|^2(s-\sigma)\bar{c}_f}|k_{12}|\\\lesssim& |\varepsilon k_{12}|^{\eta} (t-s)^{1/4}|k_{12}|^{3/2}e^{-|k_{12}|^2(s-\sigma)\bar{c}_f},\endaligned$$
 we obtain that for $\epsilon>0, \eta>0$ small enough,
$$\aligned& E|\Delta_q(I_t^4-\tilde{I}_t^4)|^2\\\lesssim &\sum_{k\in \mathbb{Z}^3\backslash\{0\}}\sum_{|i-j|\leq1,|i'-j'|\leq1}\sum_{k_{12}=k,k_3,k_4}\theta(2^{-q}k)^2\theta(2^{-i} k_{123})\theta(2^{-j}k_3)\theta(2^{-i'} k_{124})\theta(2^{-j'}k_4)\\& \bigg{[}\bigg|\int_0^tds\int_0^td\bar{s}\sum_{j_1,j_1',i_3,i_3'=1}^3(c^{\varepsilon,j_1}_{k_{123},k_3}-\bar{c}^{\varepsilon,j_1}_{k_{123},k_3})
(t-s)(\bar{c}^{\varepsilon,j_1'}_{k_{124},k_4}-\bar{c}^{\varepsilon,j_1'}_{k_{124},k_4})(t-\bar{s})
\\&[\int_0^sd\sigma\int_0^{\bar{s}}d\bar{\sigma}|(e^{-|k_{12}|^2(s-\sigma)f(\varepsilon k_{12})}
-e^{-|k_{12}|^2f(\varepsilon k_{12})(t-\sigma)})(e^{-|k_{12}|^2f(\varepsilon k_{12})(\bar{s}-\bar{\sigma})}-e^{-|k_{12}|^2f(\varepsilon k_{12})(t-\bar{\sigma})})||k_{12}|^2\\&\frac{1}{|k_1|^2|k_2|^2}+\int_s^td\sigma\int_{\bar{s}}^td\bar{\sigma}e^{-|k_{12}|^2f(\varepsilon k_{12})(t-\sigma)}e^{-|k_{12}|^2f(\varepsilon k_{12})(t-\bar{\sigma})}|k_{12}|^2\frac{1}{|k_1|^2|k_2|^2}]\bigg|
\\&+\bigg|\int_0^tds\int_0^td\bar{s}\sum_{j_1,j_1'=1}^3\bar{c}^{\varepsilon,j_1}_{k_{123},k_3}(t-s)\bar{c}^{\varepsilon,j_1'}_{k_{124},k_4}(t-\bar{s})
\\&[\int_0^sd\sigma\int_0^{\bar{s}}d\bar{\sigma}|(e^{-|k_{12}|^2(s-\sigma)f(\varepsilon k_{12})}g(\varepsilon k_{12}^{i_3})k_{12}^{i_3}
-e^{-|k_{12}|^2(s-\sigma)}\imath k_{12}^{i_3}-e^{-|k_{12}|^2f(\varepsilon k_{12})(t-\sigma)}g(\varepsilon k_{12}^{i_3})k_{12}^{i_3}\\&
+e^{-|k_{12}|^2(t-\sigma)}k_{12}^{i_3}\imath)(e^{-|k_{12}|^2(\bar{s}-\bar{\sigma})f(\varepsilon k_{12})}g(\varepsilon k_{12}^{i_3'})k_{12}^{i'_3}
-e^{-|k_{12}|^2(\bar{s}-\bar{\sigma})}\imath k_{12}^{i_3'}-e^{-|k_{12}|^2f(\varepsilon k_{12})(t-\bar{\sigma})}g(\varepsilon k_{12}^{i_3'})k_{12}^{i_3'}\\&+e^{-|k_{12}|^2(t-\bar{\sigma})}k_{12}^{i_3'}\imath)|\frac{1}{|k_1|^2|k_2|^2}
+\int_s^td\sigma\int_{\bar{s}}^td\bar{\sigma}|(e^{-|k_{12}|^2f(\varepsilon k_{12})(t-\sigma)}g(\varepsilon k_{12}^{i_3})-e^{-|k_{12}|^2(t-\sigma)}\imath)\\&(e^{-|k_{12}|^2f(\varepsilon k_{12})(t-\bar{\sigma})}g(\varepsilon k_{12}^{i_3})-e^{-|k_{12}|^2(t-\bar{\sigma})}\imath)||k_{12}|^2\frac{1}{|k_1|^2|k_2|^2}
\\&+\int_s^td\sigma\int_{\bar{s}}^td\bar{\sigma}e^{-|k_{12}|^2(t-\sigma)}e^{-|k_{12}|^2(t-\bar{\sigma})}|k_{12}|^2\frac{|\varepsilon k_1|^{\eta}
+|\varepsilon k_2|^{\eta}}{|k_1|^2|k_2|^2}\\&+\int_0^sd\sigma\int^{\bar{s}}_0d\bar{\sigma}|(e^{-|k_{12}|^2(t-\sigma)}
-e^{-|k_{12}|^2(s-\sigma)})(e^{-|k_{12}|^2(t-\bar{\sigma})}-e^{-|k_{12}|^2(s-\bar{\sigma})})||k_{12}|^2\frac{|\varepsilon k_1|^{\eta}+|\varepsilon k_2|^{\eta}}{|k_1|^2|k_2|^2}]\bigg|\bigg{]}
\endaligned$$
$$\aligned\lesssim &\sum_{k\in \mathbb{Z}^3\backslash\{0\}}\sum_{|i-j|\leq1,|i'-j'|\leq1}\sum_{k_{12}=k,k_3,k_4}\theta(2^{-q}k)^2\theta(2^{-i} k_{123})\theta(2^{-j}k_3)\theta(2^{-i'} k_{124})\theta(2^{-j'}k_4)\\&\frac{t^{2\epsilon}((|\varepsilon k_{123}|^{\eta/2}+|\varepsilon k_3|^{\eta/2})(|\varepsilon k_{124}|^{\eta/2}+|\varepsilon k_4|^{\eta/2})+|\varepsilon k_{12}|^{\eta}+|\varepsilon k_1|^{\eta}+|\varepsilon k_2|^{\eta})}{|k_{12}||k_1|^2|k_2|^2|k_3|^2|k_4|^2(|k_{123}|^2+|k_3|^2)^{3/4-\epsilon}(|k_{124}|^2+|k_4|^2)^{3/4-\epsilon}}
\\\lesssim&\varepsilon^\eta t^{2\epsilon}\bigg(\sum_{q\lesssim i}\sum_{q\lesssim i'}2^{-(i+i')(1/2-3\epsilon-\eta/2)}\sum_k\sum_{k_{12}=k}\theta(2^{-q}k)\frac{1}{|k_{12}||k_1|^2|k_2|^2}\\&+\sum_{q\lesssim i}\sum_{q\lesssim i'}2^{-(i+i')(1/2-3\epsilon)}\sum_k\sum_{k_{12}=k}\theta(2^{-q}k)\frac{|k_1|^{\eta}+|k_2|^{\eta}+|k_{12}|^\eta}{|k_{12}||k_1|^2|k_2|^2}\bigg)
\\\lesssim& \varepsilon ^\eta t^{2\epsilon}2^{q(6\epsilon+\eta)}.\endaligned$$
Here we used (4.2) (4.3) and (4.7) in the second inequality and we used Lemma 4.1 in the last inequality.
Moreover, by similar argument as (4.4) and (4.5) we obtain for $\eta,\epsilon>0$ small enough
$$\aligned& E[|\Delta_q(\tilde{I}_t^4+4\sum_{i_1=1}^3u^{\varepsilon,i_1}_2(t) C^{\varepsilon,i_0,i_1,j_0}_t-4\sum_{i_1=1}^3\bar{u}^{\varepsilon,i_1}_2(t) \bar{C}^{\varepsilon,i_0,i_1,j_0}_t))|^2]\\\lesssim&\sum_{k}\sum_{k_{12}=k}\frac{h(\varepsilon k_1)^2h(\varepsilon k_2)^2}{|k_1|^2|k_2|^2|k_{12}|^2}\theta(2^{-q}k)^2\bigg{[}\sum_{i_1,j_1=1}^3\sum_{|i-j|\leq1}\sum_{k_3}\theta(2^{-j}k_3)\int_0^t\frac{e^{-|k_3|^2(t-s)f(\varepsilon k_3)}h(\varepsilon k_3)^2}{|k_3|^2f(\varepsilon k_3)}\\&(\theta(2^{-i}k_{123})e^{-|k_{123}|^2(t-s)f(\varepsilon k_{123})}k_{123}^{j_1}g(\varepsilon k_{123}^{j_1})\hat{P}^{i_0i_1}(k_{123})-\theta(2^{-i}k_3)e^{-|k_{3}|^2(t-s)f(\varepsilon k_3)}k^{j_1}_{3}g(\varepsilon k^{j_1}_{3})\hat{P}^{i_0i_1}(k_{3}))ds\\&-\sum_{|i-j|\leq1}\sum_{i_1,j_1=1}^3\sum_{k_3}\theta(2^{-j}k_3)\int_0^t\frac{e^{-|k_3|^2(t-s)}h(\varepsilon k_3)^2}{|k_3|^2}\\&(\theta(2^{-i}k_{123})e^{-|k_{123}|^2(t-s)}k_{123}^{j_1}\hat{P}^{i_0i_1}(k_{123})-\theta(2^{-i}k_3)e^{-|k_{3}|^2(t-s)}k^{j_1}_{3}
\hat{P}^{i_0i_1}(k_{3}))ds\bigg{]}^2
\\&+\sum_{k}\sum_{k_{12}=k}\frac{|\varepsilon k_{12}|^\eta+|\varepsilon k_1|^{\eta}+|\varepsilon k_2|^{\eta}}{|k_1|^2|k_2|^2|k_{12}|^2}\theta(2^{-q}k)^2\bigg{[}\sum_{|i-j|\leq1}\sum_{k_3}\theta(2^{-j}k_3)\int_0^t\frac{e^{-|k_3|^2(t-s)
f(\varepsilon k_3)}h(\varepsilon k_3)^2}{|k_3|^2f(\varepsilon k_3)}\\&(\theta(2^{-i}k_{123})e^{-|k_{123}|^2(t-s)}k_{123}^{j_1}\hat{P}^{i_0i_1}(k_{123})-\theta(2^{-i}k_3)e^{-|k_{3}|^2(t-s)}k^{j_1}_{3}
\hat{P}^{i_0i_1}(k_{3}))ds\bigg{]}^2\\\lesssim&\sum_{k}\sum_{k_{12}=k}\frac{h(\varepsilon k_1)^2h(\varepsilon k_2)^2}{|k_1|^2|k_2|^2|k_{12}|^2}\theta(2^{-q}k)^2\bigg{[}\sum_{|i-j|\leq1}\sum_{k_3}\theta(2^{-j}k_3)\int_0^t\frac{e^{-|k_3|^2(t-s)\bar{c}_f}h(\varepsilon k_3)^2}{|k_3|^2f(\varepsilon k_3)}(|k_{12}|^{2\eta} \\&(t-s)^{-\frac{1-2\eta}{2}})\wedge (t-s)^{-\frac{1}{2}}(\varepsilon^{\eta}|k_{123}|^{\eta}e^{-\bar{c}_f|k_{123}|^2(t-s)}+\varepsilon^{\eta}|k_{3}|^{\eta}) ds\bigg{]}^2
\\&+\sum_{k}\sum_{k_{12}=k}\frac{|\varepsilon k_{12}|^\eta+|\varepsilon k_1|^{\eta}+|\varepsilon k_2|^{\eta}}{|k_1|^2|k_2|^2|k_{12}|^2}\theta(2^{-q}k)^2\bigg{[}\sum_{|i-j|\leq1}\sum_{k_3}\theta(2^{-j}k_3)\int_0^t\frac{e^{-|k_3|^2(t-s)
f(\varepsilon k_3)}h(\varepsilon k_3)^2}{|k_3|^2f(\varepsilon k_3)}\\&|k_{12}|^\eta(t-s)^{-\frac{1-\eta}{2}} ds\bigg{]}^2 \\\lesssim& \varepsilon^{\eta}t^\epsilon2^{3q\eta}.\endaligned$$
Here in the first inequality we used (4.7) and in the second inequality we used (4.2) (4.3) and by a similar argument as Lemmas 4.2 and 4.3 we obtain similar estimate for the corresponding terms.

Similarly, we have for $I_t^4$ of $\pi_0(u^{\varepsilon,i_0}_3-u^{\varepsilon,i_0}_{31},u_1^{\varepsilon,j_0})-\pi_0(\bar{u}^{\varepsilon,i_0}_3-\bar{u}^{\varepsilon,i_0}_{31},\bar{u}_1^{\varepsilon,j_0})$
$$\aligned& E[|\Delta_q(I_t^4+4\sum_{j_1=1}^3u^{\varepsilon,j_1}_2(t) \tilde{C}^{\varepsilon,i_0,j_1,j_0}_t)|^2] \lesssim \varepsilon^{\eta}t^\epsilon2^{3q\eta},\endaligned$$
where $$\aligned \tilde{C}^{\varepsilon,i_0,j_1,j_0}_t=&\frac{(2\pi)^{-3}}{4}\sum_{|i-j|\leq1}\sum_{k_3}\sum_{i_1=1}^3\theta(2^{-i} k_{3})\theta(2^{-j}k_3)\int_0^tds\frac{e^{-2|k_3|^2(t-s)f(\varepsilon k_3)}h(\varepsilon k_3)^2}{|k_3|^2f(\varepsilon k_3)}\\& \sum_{i_4}\hat{P}^{i_1i_4}(k_3)\hat{P}^{j_0i_4}(k_3)k_{3}^{j_1}g(\varepsilon k_{3}^{j_1})\hat{P}^{i_0i_1}(k_{3})\\=&\frac{(2\pi)^{-3}}{4}\sum_{k_3}\sum_{i_1=1}^3\int_0^tds\frac{e^{-2|k_3|^2(t-s)f(\varepsilon k_3)}h(\varepsilon k_3)^2}{|k_3|^2f(\varepsilon k_3)} \sum_{i_4=1}^3\hat{P}^{i_1i_4}(k_3)\hat{P}^{j_0i_4}(k_3)k_{3}^{j_1}g(\varepsilon k_{3}^{j_1})\hat{P}^{i_0i_1}(k_{3})\\\rightarrow&(2\pi)^{-3}\sum_{i_1=1}^3\sum_{i_4=1}^3\int_{\mathbb{R}^3}\frac{(\cos (ax^{j_1})-\cos (bx^{j_1}))h(x)^2}{8|x|^4f(x)^2(a+b)} \hat{P}^{i_1i_4}(x)\hat{P}^{j_0i_4}(x)\hat{P}^{i_0i_1}(x)dx,\quad\varepsilon\rightarrow0,\endaligned$$.

Now we consider $I_t^5-\tilde{I}_t^5+\tilde{I}_t^5-\bar{I}_t^5$ with
$$\aligned \tilde{I}_t^5=&(2\pi)^{-9/2}\sum_{k\in \mathbb{Z}^3\backslash\{0\}}\sum_{|i-j|\leq1}\sum_{k_{14}=k,k_2}\sum_{i_1,i_2,i_3,j_1=1}^3\theta(2^{-i} k_{1})\theta(2^{-j}k_4)\bigg{[}\int_0^t:\hat{X}^{\varepsilon,i_2}_s(k_1)\hat{X}^{\varepsilon,j_0}_t(k_4):e^{-|k_{1}|^2(t-s)f(\varepsilon k_1)}\\& k_{1}^{j_1}g(\varepsilon k_{1}^{j_1}) ds\int_0^sd\sigma e^{-|k_{12}|^2(s-\sigma)f(\varepsilon k_{12})}\frac{e^{-|k_2|^2(s-\sigma)f(\varepsilon k_2)}h(\varepsilon k_2)^2}{2|k_2|^2f(\varepsilon k_2)} k_{12}^{i_3}g(\varepsilon k_{12}^{i_3})-\int_0^t:\hat{\bar{X}}^{\varepsilon,i_2}_s(k_1)\hat{\bar{X}}^{\varepsilon,j_0}_t(k_4):\\&e^{-|k_{1}|^2(t-s)}\imath k_{1}^{j_1} ds\int_0^sd\sigma e^{-|k_{12}|^2(s-\sigma)}\frac{e^{-|k_2|^2(s-\sigma)}h(\varepsilon k_2)^2}{2|k_2|^2} \imath k_{12}^{i_3}\bigg{]}\hat{P}^{i_1i_2}(k_{12})\sum_{i_4=1}^3\hat{P}^{i_3i_4}(k_2)\hat{P}^{j_1i_4}(k_2)\hat{P}^{i_0i_1}(k_{1})e_k,\endaligned$$
and $$\aligned \bar{I}_t^5=&(2\pi)^{-9/2}\sum_{k\in \mathbb{Z}^3\backslash\{0\}}\sum_{|i-j|\leq1}\sum_{k_{14}=k,k_2}\sum_{i_1,i_2,i_3,j_1=1}^3\theta(2^{-i} k_{1})\theta(2^{-j}k_4)\bigg{[}\int_0^t:\hat{X}^{\varepsilon,i_2}_s(k_1)\hat{X}^{\varepsilon,j_0}_t(k_4):e^{-|k_{1}|^2(t-s)f(\varepsilon k_1)}\\& k_{1}^{j_1}g(\varepsilon k_{1}^{j_1}) ds\int_0^sd\sigma e^{-|k_{2}|^2(s-\sigma)f(\varepsilon k_2)}\frac{e^{-|k_2|^2(s-\sigma)f(\varepsilon k_2)}h(\varepsilon k_2)^2}{2|k_2|^2f(\varepsilon k_2)} k_{2}^{i_3}g(\varepsilon k_{2}^{i_3})-\int_0^t:\hat{\bar{X}}^{\varepsilon,i_2}_s(k_1)\hat{\bar{X}}^{\varepsilon,j_0}_t(k_4):\\&e^{-|k_{1}|^2(t-s)}\imath k_{1}^{j_1} ds\int_0^sd\sigma e^{-|k_{2}|^2(s-\sigma)}\frac{e^{-|k_2|^2(s-\sigma)}h(\varepsilon k_2)^2}{2|k_2|^2} \imath k_{2}^{i_3}\bigg{]}\hat{P}^{i_1i_2}(k_{2})\sum_{i_4=1}^3\hat{P}^{i_3i_4}(k_2)\hat{P}^{j_1i_4}(k_2)\hat{P}^{i_0i_1}(k_{1})e_k
\endaligned$$
$$\aligned=&(2\pi)^{-9/2}\sum_{k\in \mathbb{Z}^3\backslash\{0\}}\sum_{|i-j|\leq1}\sum_{k_{14}=k,k_2}\sum_{i_1,i_2,i_3,j_1=1}^3\theta(2^{-i} k_{1})\theta(2^{-j}k_4)\int_0^t:\hat{X}^{\varepsilon,i_2}_s(k_1)\hat{X}^{\varepsilon,j_0}_t(k_4):e^{-|k_{1}|^2(t-s)f(\varepsilon k_1)}\\& k_{1}^{j_1}g(\varepsilon k_{1}^{j_1}) ds\int_0^sd\sigma e^{-|k_{2}|^2(s-\sigma)f(\varepsilon k_2)}\frac{e^{-|k_2|^2(s-\sigma)f(\varepsilon k_2)}h(\varepsilon k_2)^2}{2|k_2|^2f(\varepsilon k_2)} k_{2}^{i_3}g(\varepsilon k_{2}^{i_3})\\&\hat{P}^{i_1i_2}(k_{2})\sum_{i_4=1}^3\hat{P}^{i_3i_4}(k_2)\hat{P}^{j_1i_4}(k_2)\hat{P}^{i_0i_1}(k_{1})e_k\\=&2(2\pi)^{-3/2}\sum_{k\in \mathbb{Z}^3\backslash\{0\}}\sum_{|i-j|\leq1}\sum_{k_{14}=k}\sum_{i_1,i_2,j_1=1}^3\theta(2^{-i} k_{1})\theta(2^{-j}k_4)\int_0^t:\hat{X}^{\varepsilon,i_2}_s(k_1)\hat{X}^{\varepsilon,j_0}_t(k_4):e^{-|k_{1}|^2(t-s)f(\varepsilon k_1)}\\& k_{1}^{j_1}g(\varepsilon k_{1}^{j_1})C_s^{\varepsilon,i_1,i_2,j_1}\hat{P}^{i_0i_1}(k_{1})e_kds.\endaligned$$
Let $$d_{k_{12},k_2}^\varepsilon(s-\sigma)=\sum_{i_2,i_3=1}^3e^{-|k_{12}|^2(s-\sigma)f(\varepsilon k_{12})}\frac{e^{-|k_2|^2(s-\sigma)f(\varepsilon k_2)}h(\varepsilon k_2)^2}{|k_2|^2f(\varepsilon k_2)} k_{12}^{i_3}g(\varepsilon k_{12}^{i_3})\hat{P}^{i_1i_2}(k_{12})$$
and $$\bar{d}_{k_{12},k_2}^\varepsilon(s-\sigma)=\sum_{i_2,i_3=1}^3e^{-|k_{12}|^2(s-\sigma)}\frac{e^{-|k_2|^2(s-\sigma)}h(\varepsilon k_2)^2}{|k_2|^2} k_{12}^{i_3}\imath\hat{P}^{i_1i_2}(k_{12}).$$
Since by H\"{o}lder's inequality we obtain that for $\eta>0$ small enough
$$\aligned &E(:\hat{X}^{\varepsilon,i_2}_s(k_1)\hat{X}^{\varepsilon,j_0}_t(k_4):-:\hat{X}^{\varepsilon,i_2}_\sigma(k_1)\hat{X}^{\varepsilon,j_0}_t(k_4):)
(\overline{:\hat{X}^{\varepsilon,i_2}_{\bar{s}}(k'_1)\hat{X}^{\varepsilon,j_0}_t(k'_4):
-:\hat{X}^{\varepsilon,i_2}_{\bar{\sigma}}(k'_1)\hat{X}^{\varepsilon,j_0}_t(k'_4):})
\\\lesssim&(1_{k_1=k'_1}1_{k_4=k'_4}+1_{k_1=k'_4}1_{k_4=k'_1})
\bigg(\frac{1-e^{-|k_1|^2|s-\sigma|f(\varepsilon k_1)}}{|k_1|^2|k_4|^2f(\varepsilon k_1)}\bigg)^{1/2}\bigg(\frac{1-e^{-|k'_1|^2|\bar{s}-\bar{\sigma}|f(\varepsilon k_1')}}{|k'_1|^2|k'_4|^2f(\varepsilon k_1')}\bigg)^{1/2}
\\\lesssim&(1_{k_1=k'_1}1_{k_4=k'_4}+1_{k_1=k'_4}1_{k_4=k'_1})\frac{|k_1|^\eta|k'_1|^\eta}{|k_1||k'_1||k_4||k'_4|}|s-\sigma|^{\eta/2}|\bar{s}-\bar{\sigma}|^{\eta/2},\endaligned$$
and by (4.7) $$\aligned &E(:\hat{X}^{\varepsilon,i_2}_s(k_1)\hat{X}^{\varepsilon,j_0}_t(k_4):-:\hat{X}^{\varepsilon,i_2}_\sigma(k_1)\hat{X}^{\varepsilon,j_0}_t(k_4):-:\hat{\bar{X}}^{\varepsilon,i_2}_s(k_1)\hat{\bar{X}}^{\varepsilon,j_0}_t(k_4):+:\hat{\bar{X}}^{\varepsilon,i_2}_\sigma(k_1)\hat{\bar{X}}^{\varepsilon,j_0}_t(k_4):)\\&
(\overline{:\hat{X}^{\varepsilon,i_2}_s(k'_1)\hat{X}^{\varepsilon,j_0}_t(k'_4):-:\hat{X}^{\varepsilon,i_2}_\sigma(k'_1)\hat{X}^{\varepsilon,j_0}_t(k'_4):-:\hat{\bar{X}}^{\varepsilon,i_2}_s(k'_1)\hat{\bar{X}}^{\varepsilon,j_0}_t(k'_4):+:\hat{\bar{X}}^{\varepsilon,i_2}_\sigma(k'_1)\hat{\bar{X}}^{\varepsilon,j_0}_t(k'_4):)}
\\\lesssim&(1_{k_1=k'_1}1_{k_4=k'_4}+1_{k_1=k'_4}1_{k_4=k'_1})
\bigg(\frac{(\varepsilon^{2\eta}|k_1|^{2\eta}+\varepsilon^{2\eta}|k_4|^{2\eta})\wedge(|s-\sigma|^{\eta}|k_1|^{2\eta})}{|k_1|^2|k_4|^2}\bigg)^{1/2}\\&\bigg(
\frac{(\varepsilon^{2\eta}|k'_1|^{2\eta}+\varepsilon^{2\eta}|k'_4|^{2\eta})\wedge(|s-\sigma|^{\eta}|k'_1|^{2\eta})}{|k'_1|^2|k'_4|^2}\bigg)^{1/2}
\\\lesssim&(1_{k_1=k'_1}1_{k_4=k'_4}+1_{k_1=k'_4}1_{k_4=k'_1})
\frac{\varepsilon^{\eta}|k_1|^{\frac{\eta}{2}}(|k_1|^{\frac{\eta}{2}}+|k_4|^{\frac{\eta}{2}})
|k_1'|^{\frac{\eta}{2}}(|k_1'|^{\frac{\eta}{2}}+|k_4'|^{\frac{\eta}{2}})}{|k_1||k'_1||k_4||k'_4|}|s-\sigma|^{\eta/4}|\bar{s}-\bar{\sigma}|^{\eta/4},\endaligned$$
it follows that for $\eta,\epsilon>0$ small enough
$$\aligned &E|\Delta_q(I_t^5-\tilde{I}_t^5)|^2\\\lesssim &\sum_{k\in \mathbb{Z}^3\backslash\{0\}}\sum_{|i-j|\leq1,|i'-j'|\leq1}\sum_{k_{14}=k,k_{14}'=k,k_2,k'_2}\sum_{i_1,i_1',j_1,j_1'=1}^3\theta(2^{-q}k)^2\theta(2^{-i} k_{1})\theta(2^{-j}k_4)\theta(2^{-i'} k'_{1})\theta(2^{-j'}k'_4)\\& \int_0^tds\int_0^td\bar{s}\int_0^sd\sigma\int_0^{\bar{s}}d\bar{\sigma}\bigg[|
(e^{-|k_{1}|^2(t-s)f(\varepsilon k_1)}k_1^{j_1}g(\varepsilon k_1^{j_1})\hat{P}^{i_0i_1}(k_1)-e^{-|k_{1}|^2(t-s)}\imath k_1^{j_1}\hat{P}^{i_0i_1}(k_1))|\\&|e^{-|k'_{1}|^2(t-s)f(\varepsilon k'_1)}k_1^{'j_1'}g(\varepsilon k_1^{'j_1'})\hat{P}^{i_0'i_1'}(k'_1)-e^{-|k'_{1}|^2(t-s)}\imath k_1^{'j_1'}\hat{P}^{i_0'i_1'}(k'_1)|\\&(1_{k_1=k'_1}1_{k_4=k'_4}+1_{k_1=k'_4}1_{k_4=k'_1})
\frac{|k_1|^\eta|k'_1|^\eta}{|k_1||k'_1||k_4||k'_4|}|s-\sigma|^{\eta/2}|\bar{s}-\bar{\sigma}|^{\eta/2}|d^\varepsilon_{k_{12},k_2}(s-\sigma)
d^\varepsilon_{k'_{12},k'_2}(\bar{s}-\bar{\sigma})|
\\&+e^{-|k_{1}|^2(t-s)}e^{-|k'_{1}|^2(t-\bar{s})}|k_{1}||k'_1|(1_{k_1=k'_1}1_{k_4=k'_4}
+1_{k_1=k'_4}1_{k_4=k'_1})\frac{|k_1|^\eta|k'_1|^\eta}{|k_1||k'_1||k_4||k'_4|}
\\&|s-\sigma|^{\eta/2}|\bar{s}-\bar{\sigma}|^{\eta/2}|(d^\varepsilon_{k_{12},k_2}(s-\sigma)-\bar{d}^\varepsilon_{k_{12},k_2}(s-\sigma))
(d^\varepsilon_{k'_{12},k'_2}(\bar{s}-\bar{\sigma})-\bar{d}^\varepsilon_{k'_{12},k'_2}(\bar{s}-\bar{\sigma}))|
\\&+e^{-|k_{1}|^2(t-s)}e^{-|k'_{1}|^2(t-\bar{s})}|k_{1}||k'_1||\bar{d}^\varepsilon_{k_{12},k_2}(s-\sigma)\bar{d}^\varepsilon_{k_{12}',k_2'}
(\bar{s}-\bar{\sigma})|(1_{k_1=k'_1}1_{k_4=k'_4}+1_{k_1=k'_4}1_{k_4=k'_1})
\\&\frac{\varepsilon^{\eta}|k_1|^{\frac{\eta}{2}}|k'_1|^{\frac{\eta}{2}}(|k_1|^{\eta}+|k_4|^{\eta})}{|k_1||k'_1||k_4||k'_4|}|s-\sigma|^{\eta/4}|\bar{s}-\bar{\sigma}|^{\eta/4}\bigg]
\\\lesssim &\sum_{k\in \mathbb{Z}^3\backslash\{0\}}\sum_{|i-j|\leq1,|i'-j'|\leq1}\sum_{k_{14}=k}\theta(2^{-q}k)^2\theta(2^{-i} k_{1})\theta(2^{-j}k_4)\theta(2^{-i'} k_{1})\theta(2^{-j'}k_4)\\&\bigg(\frac{t^{\epsilon}\varepsilon ^\eta(|k_1|^{\eta}+|k_4|^{\eta})}{|k_1|^{4-2\eta-2\epsilon}|k_4|^2}+\frac{t^{\epsilon}\varepsilon ^\eta|k_1|^{\eta/2}|k_4|^{\eta/2}}{|k_1|^{3-\eta-\epsilon}|k_4|^{3-\eta-\epsilon}}\bigg)
\\\lesssim&t^{\epsilon}\sum_k\sum_{k_{14}=k}\theta(2^{-q}k)\sum_{q\lesssim i}2^{-i}\frac{\varepsilon ^\eta(|k_1|^{\eta}+|k_4|^\eta)}{|k_1|^{3-2\eta-2\epsilon}|k_4|^2}+t^{\epsilon}\sum_k\sum_{k_{14}=k}\theta(2^{-q}k)\sum_{q\lesssim j}2^{-j\epsilon}\frac{\varepsilon ^\eta|k_1|^{\eta/2}|k_4|^{\eta/2}}{|k_1|^{3-\eta-2\epsilon}|k_4|^{3-\eta-\epsilon}}
\\\lesssim& \varepsilon ^\eta t^{\epsilon}2^{q(2\epsilon+3\eta)},\endaligned$$where in the second inequality we used (4.2) and (4.3) and in the last inequality we used Lemma 4.1 and  $q\lesssim i$ follows from $|k|\leq |k_1|+|k_4|\lesssim 2^i$.

Moreover, it follows  that for $\eta,\epsilon>0$ small enough
$$\aligned &E[|\Delta_q(\tilde{I}_t^5-\bar{I}_t^5)|^2]\\\lesssim&\sum_{k\in \mathbb{Z}^3\backslash\{0\}}\sum_{|i-j|\leq1,|i'-j'|\leq1}\sum_{k_{14}=k,k_{14}'=k, k_2,k'_2}\sum_{i_1,i_2,i_3,i_1',i_2',i_3'=1}^3\theta(2^{-q}k)^2\theta(2^{-i} k_{1})\theta(2^{-j}k_4)\theta(2^{-i'} k'_{1})\theta(2^{-j'}k'_4)\\&(1_{k_1=k'_1}1_{k_4=k'_4}+1_{k_1=k'_4}1_{k_4=k'_1})
\bigg(\int_0^t\int_0^t|k_1||k'_1|e^{-|k_1|^2f(\varepsilon k_1)(t-s)-|k'_1|^2(t-\bar{s})f(\varepsilon k_1')}\frac{h(\varepsilon k_1)h(\varepsilon k_1')h(\varepsilon k_4)h(\varepsilon k_4')}{|k_1||k'_1||k_4||k'_4|}\\&\int_0^s\int_0^{\bar{s}}\bigg|\bigg[\frac{e^{-|k_2|^2(s-\sigma)
f(\varepsilon k_2)}h(\varepsilon k_2)^2}{|k_2|^2f(\varepsilon k_2)}(e^{-|k_{12}|^2f(\varepsilon k_{12})(s-\sigma)}k_{12}^{i_3}g(\varepsilon k_{12}^{i_3})\hat{P}^{i_1i_2}(k_{12})-e^{-|k_2|^2f(\varepsilon k_2)(s-\sigma)}k_2^{i_3}g(\varepsilon k_2^{i_3})\\&\hat{P}^{i_1i_2}(k_2))-\frac{e^{-|k_2|^2(s-\sigma)}h(\varepsilon k_2)^2}{|k_2|^2}(e^{-|k_{12}|^2(s-\sigma)}\imath k_{12}^{i_3}\hat{P}^{i_1i_2}(k_{12})-e^{-|k_2|^2(s-\sigma)}\imath k_2^{i_3}\hat{P}^{i_1i_2}(k_2))\bigg]\bigg|\endaligned$$
$$\aligned&\bigg|\bigg[\frac{e^{-|k'_2|^2(\bar{s}-\bar{\sigma})f(\varepsilon k'_2)}h(\varepsilon k_2')^2}{|k'_2|^2f(\varepsilon k'_2)}(e^{-|k'_{12}|^2f(\varepsilon k'_{12})(\bar{s}-\bar{\sigma})}k^{'i_3'}_{12}g(\varepsilon k^{'i_3'}_{12})\hat{P}^{i_1'i_2'}(k^{'i_3'}_{12})-e^{-|k'_2|^2f(\varepsilon k'_2)(\bar{s}-\bar{\sigma})}k^{'i_3'}_2g(\varepsilon k^{'i_3'}_2)\hat{P}^{i_1'i_2'}(k'_2))\\&-\frac{e^{-|k'_2|^2(\bar{s}-\bar{\sigma})}h(\varepsilon k_2')^2}{|k'_2|^2}(e^{-|k'_{12}|^2(\bar{s}-\bar{\sigma})}\imath k^{'i_3'}_{12}\hat{P}^{i_1'i_2'}(k'_{12})-e^{-|k'_2|^2(\bar{s}-\bar{\sigma})}\imath k^{'i_3'}_2\hat{P}^{i_1'i_2'}(k'_2))\bigg]\bigg|dsd\bar{s}d\sigma d\bar{\sigma}\\+&\int_0^t\int_0^t|k_1||k'_1|e^{-\bar{c}_f|k_1|^2(t-s)-\bar{c}_f|k'_1|^2(t-\bar{s})}\frac{|\varepsilon k_1|^{\eta}+|\varepsilon k_4|^{\eta}+|\varepsilon k_1|^{\eta/2}|\varepsilon k'_1|^{\eta/2}}{|k_1||k'_1||k_4||k'_4|}\int_0^s\int_0^{\bar{s}}\\&\frac{e^{-|k_2|^2(s-\sigma)}}{|k_2|^2}|e^{-|k_{12}|^2(s-\sigma)}\imath k_{12}^{i_3}\hat{P}^{i_1i_2}(k_{12})-e^{-|k_2|^2(s-\sigma)}\imath k_2^{i_3}\hat{P}^{i_1i_2}(k_2)|\\&\frac{e^{-|k'_2|^2(\bar{s}-\bar{\sigma})}}{|k'_2|^2}|(e^{-|k'_{12}|^2(\bar{s}-\bar{\sigma})}\imath k_{12}^{'i_3'}\hat{P}^{i_1'i_2'}(k_{12}')-e^{-|k'_2|^2(\bar{s}-\bar{\sigma})}\imath k^{'i_3'}_2\hat{P}^{i_1'i_2'}(k'_2))|dsd\bar{s}d\sigma d\bar{\sigma}\bigg)\\\lesssim&\sum_{k\in \mathbb{Z}^3\backslash\{0\}}\sum_{|i-j|\leq1,|i'-j'|\leq1}\sum_{k_{14}=k,k_{14}'=k, k_2,k'_2}\sum_{i_1,i_2,i_3,i_1',i_2',i_3'=1}^3\theta(2^{-q}k)^2\theta(2^{-i} k_{1})\theta(2^{-j}k_4)\theta(2^{-i'} k'_{1})\theta(2^{-j'}k'_4)\\&(1_{k_1=k'_1}1_{k_4=k'_4}+1_{k_1=k'_4}1_{k_4=k'_1})
\bigg(\int_0^t\int_0^t|k_1||k'_1|e^{-|k_1|^2f(\varepsilon k_1)(t-s)-|k'_1|^2(t-\bar{s})f(\varepsilon k_1')}\frac{h(\varepsilon k_1)h(\varepsilon k_1')h(\varepsilon k_4)h(\varepsilon k_4')}{|k_1||k'_1||k_4||k'_4|}\\&\int_0^s\int_0^{\bar{s}}\frac{e^{-|k_2|^2\bar{c}_f(s-\sigma)}}{|k_2|^2}\big[((s-\sigma)^{-\frac{(1-2\eta)}{2}}
|k_1|^{2\eta})\wedge ((s-\sigma)^{-\frac{1}{2}}(\varepsilon^\eta|k_2|^\eta+\varepsilon^\eta|k_{12}|^\eta e^{-|k_{12}|^2(s-\sigma)\bar{c}_f}))\big]\\&\frac{e^{-|k_2'|^2\bar{c}_f(s-\sigma)}}{|k_2'|^2}\big[((s-\sigma)^{-\frac{(1-2\eta)}{2}}
|k_1'|^{2\eta})\wedge ((s-\sigma)^{-\frac{1}{2}}(\varepsilon^\eta|k_2'|^\eta+\varepsilon^\eta|k_{12}'|^\eta e^{-|k_{12}'|^2\bar{c}_f(s-\sigma)}))\big]dsd\bar{s}d\sigma d\bar{\sigma}\\+&\int_0^t\int_0^t|k_1||k'_1|e^{-\bar{c}_f|k_1|^2(t-s)-\bar{c}_f|k'_1|^2(t-\bar{s})}\frac{|\varepsilon k_1|^{\eta}+|\varepsilon k_4|^{\eta}+|\varepsilon k_1|^{\eta/2}|\varepsilon k'_1|^{\eta/2}}{|k_1||k'_1||k_4||k'_4|}\int_0^s\int_0^{\bar{s}}\\&\frac{e^{-|k_2|^2(s-\sigma)}}{|k_2|^2}|k_{1}|^\eta|s-\sigma|^{-\frac{1-\eta}{2}}|
\frac{e^{-|k'_2|^2(\bar{s}-\bar{\sigma})}}{|k'_2|^2}|k_{1}'|^\eta|\bar{s}-\bar{\sigma}|^{-\frac{1-\eta}{2}}|dsd\bar{s}d\sigma d\bar{\sigma}\bigg)\\\lesssim&\varepsilon^\eta t^{\epsilon} \bigg(\sum_k\sum_{k_{14}=k}\theta(2^{-q}k)\sum_{q\lesssim i}2^{-i}(\frac{1}{|k_1|^{3-2\eta-2\epsilon}|k_4|^{2-\eta}}+\frac{1}{|k_1|^{3-3\eta-2\epsilon}|k_4|^2})
\\&+\sum_k\sum_{k_{14}=k}\theta(2^{-q}k)\sum_{q\lesssim i}2^{-i\epsilon}\frac{1}{|k_1|^{3-2\eta-2\epsilon}|k_4|^{3-\eta-\epsilon}}+\frac{1}{|k_1|^{3-\eta-2\epsilon}|k_4|^{3-2\eta-\epsilon}}\bigg)
\\\lesssim& \varepsilon^{\eta}t^{\epsilon}2^{q(3\epsilon+3\eta)},\endaligned$$
where we used (4.7) in the first inequality and a similar argument as (4.4) (4.5) in the second inequality  and in the last inequality we used Lemma 4.1 and $q\lesssim i$ follows from $|k|\leq |k_1|+|k_4|\lesssim 2^i$.

\textbf{Terms in the fouth chaos}: Now for $I_t^1$ we have the following calculations:
$$\aligned &E[|\Delta_qI_t^1|^2]\\\lesssim&\sum_{k\in \mathbb{Z}^3\backslash\{0\}}\sum_{|i-j|\leq1,|i'-j'|\leq1}\sum_{k_{1234}=k,k'_{1234}=k}\sum_{i_1,i_2,i_3,j_1,i_1',i_2',i_3',j_1'=1}^3\theta(2^{-q}k)^2\theta(2^{-i} k_{123})\theta(2^{-j}k_4)\theta(2^{-i'} k'_{123})\theta(2^{-j'}k'_4)\\&(1_{k_1=k'_1,k_2=k'_2,k_3=k'_3,k_4=k'_4}+1_{k_1=k'_4,k_2=k'_2,k_3=k'_3,k_4=k'_1}
+1_{k_1=k'_1,k_2=k'_2,k_3=k'_4,k_4=k'_3}+1_{k_1=k'_3,k_2=k'_4,k_3=k'_1,k_4=k'_2}\\&+1_{k_1=k'_1,k_2=k'_3,k_3=k'_2,k_4=k'_4}+1_{k_1=k'_3,k_2=k'_2,k_3=k'_4,k_4=k'_1}
+1_{k_1=k'_4,k_2=k'_2,k_3=k'_1,k_4=k'_3})\endaligned$$
$$\aligned&\bigg[\int_0^tds\int_0^td\bar{s}|e^{-|k_{123}|^2(t-s)f(\varepsilon k_{123})}k_{123}^{j_1}g(\varepsilon k_{123}^{j_1})-e^{-|k_{123}|^2(t-s)}\imath k_{123}^{j_1}||e^{-|k'_{123}|^2(t-s)f(\varepsilon k'_{123})}k^{'j_1'}_{123}g(\varepsilon k^{'j_1'}_{123})\\&-e^{-|k'_{123}|^2(t-s)}\imath k'_{123}|\int_0^s\int_0^{\bar{s}}\frac{1}{|k_1|^2f(\varepsilon k_1)|k_2|^2f(\varepsilon k_2)|k_3|^2f(\varepsilon k_3)|k_4|^2f(\varepsilon k_4)}\\&e^{-(|k_{12}|^2(s-\sigma)f(\varepsilon k_{12})+|k'_{12}|^2f(\varepsilon k'_{12})(\bar{s}-\bar{\sigma}))}d\sigma d\bar{\sigma} |k_{12}k_{12}'|\\&+\int_0^tds\int_0^td\bar{s}e^{-|k_{123}|^2(t-s)} |k_{123}^{j_1}e^{-|k'_{123}|^2(t-s)} k_{123}^{'j_1'}|\int_0^s\int_0^{\bar{s}}\frac{1}{|k_1|^2f(\varepsilon k_1)|k_2|^2f(\varepsilon k_2)|k_3|^2f(\varepsilon k_3)|k_4|^2f(\varepsilon k_4)}\\&|e^{-|k_{12}|^2(s-\sigma)f(\varepsilon k_{12})}k_{12}^{i_4}g(\varepsilon k_{12}^{i_4})-e^{-|k_{12}|^2(s-\sigma)}k_{12}^{i_4}\imath|\cdot|e^{-|k'_{12}|^2(\bar{s}-\bar{\sigma})f(\varepsilon k'_{12})}k_{12}^{'i_4'}g(\varepsilon k_{12}^{'i_4'})-e^{-|k'_{12}|^2(\bar{s}-\bar{\sigma})}k_{12}^{'i_4'}\imath|\\&(1_{\{i_4=i_3,i_4'=i_3'\}}+1_{\{i_4=i_2,i_4'=i_2'\}})d\sigma d\bar{\sigma} \\&+\int_0^tds\int_0^td\bar{s}e^{-|k_{123}|^2(t-s)} |k_{123}^{j_1}|e^{-|k'_{123}|^2(t-s)} |k^{'j_1'}_{123}|\int_0^s\int_0^{\bar{s}}e^{-|k_{12}|^2(s-\sigma)}|k_{12}| e^{-|k'_{12}|^2(\bar{s}-\bar{\sigma})}|k_{12}'|\\&(E|:\hat{X}_\sigma^{\varepsilon,i_2}(k_1)\hat{X}^{\varepsilon,i_3}_\sigma(k_2)
\hat{X}_s^{\varepsilon,j_1}(k_3)\hat{X}_t^{\varepsilon,j_0}(k_4):-:\hat{\bar{X}}_\sigma^{\varepsilon,i_2}(k_1)
\hat{\bar{X}}_\sigma^{\varepsilon,i_3}(k_2)\hat{\bar{X}}_s^{\varepsilon,j_1}(k_3)\hat{\bar{X}}_t^{\varepsilon,j_0}(k_4):|^2)^{1/2}
\\&(E|:\hat{X}_{\bar{\sigma}}^{\varepsilon,i_2'}(k'_1)\hat{X}_{\bar{\sigma}}^{\varepsilon,i_2'}(k'_2)\hat{X}_{\bar{s}}^{\varepsilon,j_1'}(k'_3)\hat{X}_{\bar{t}}
^{\varepsilon,j_0'}(k'_4):-:\hat{\bar{X}}_\sigma^{\varepsilon,i_2'}(k'_1)\hat{\bar{X}}_\sigma^{\varepsilon,i_3'}(k'_2)\hat{\bar{X}}_s^{\varepsilon,j_1'}(k'_3)
\hat{\bar{X}}_t^{\varepsilon,j_0'}(k'_4):|^2)^{1/2}d\sigma d\bar{\sigma}\bigg]\\=&E_t^1+E_t^2+E_t^3+E_t^4+E_t^5+E_t^6+E_t^7. \endaligned.$$
Here $E_t^i$ means the term corresponding to each characteristic function. Since
$$\aligned &E|:\hat{X}_\sigma^{\varepsilon,i_2}(k_1)\hat{X}^{\varepsilon,i_3}_\sigma(k_2)
\hat{X}_s^{\varepsilon,j_1}(k_3)\hat{X}_t^{\varepsilon,j_0}(k_4):-:\hat{\bar{X}}_\sigma^{\varepsilon,i_2}(k_1)
\hat{\bar{X}}_\sigma^{\varepsilon,i_3}(k_2)\hat{\bar{X}}_s^{\varepsilon,j_1}(k_3)\hat{\bar{X}}_t^{\varepsilon,j_0}(k_4):|^2
\\\lesssim&|\frac{1}{16|k_1|^2|k_2|^2|k_3|^2|k_4|^2f(\varepsilon k_1)f(\varepsilon k_2)f(\varepsilon k_3)f(\varepsilon k_4)}+\frac{1}{16|k_1|^2|k_2|^2|k_3|^2|k_4|^2}\\&-\frac{2}{|k_1|^2|k_2|^2|k_3|^2|k_4|^2(f(\varepsilon k_1)+1)(f(\varepsilon k_2)+1)(f(\varepsilon k_3)+1)(f(\varepsilon k_4)+1)}|\\\lesssim&\frac{\sum_{i=1}^4|\varepsilon k_i|^\eta}{|k_1|^2|k_2|^2|k_3|^2|k_4|^2},\endaligned\eqno(4.8)$$
for $\epsilon,\eta>0$ small enough  we have
$$\aligned E_t^1\lesssim&\varepsilon^\eta\sum_{k\in \mathbb{Z}^3\backslash\{0\}}\sum_{|i-j|\leq1,|i'-j'|\leq1}\sum_{k_{1234}=k}\frac{\theta(2^{-q}k)^2\theta(2^{-i} k_{123})\theta(2^{-j}k_4)\theta(2^{-i'} k_{123})\theta(2^{-j'}k_4)t^\eta}{|k_1|^2|k_2|^2|k_3|^2|k_4|^2|k_{12}|^2|k_{123}|^{2-2\eta}}\\&(|k_{123}|^\eta+|k_{12}|^\eta+\sum_{i=1}^4|k_i|^{\eta})
\\\lesssim& \sum_{k\in \mathbb{Z}^3\backslash\{0\}}\sum_{|i-j|\leq1,|i'-j'|\leq1}\sum_{k_{1234}=k}\theta(2^{-q}k)^2\theta(2^{-i} k_{123})\theta(2^{-j}k_4)\theta(2^{-i'} k_{123})\theta(2^{-j'}k_4)\bigg(\frac{t^\eta \varepsilon^\eta}{|k_4|^2|k_{123}|^{4-3\eta-\epsilon}}\\&+\frac{t^\eta \varepsilon^\eta}{|k_4|^{2-\eta}|k_{123}|^{4-2\eta-\epsilon}}\bigg)
\\\lesssim& \sum_{k\in \mathbb{Z}^3\backslash\{0\}}\sum_{q\lesssim i}2^{-(2-3\eta-\epsilon)i}\theta(2^{-q}k)^2\frac{t^\eta \varepsilon^\eta}{|k|}
\lesssim\varepsilon^\eta2^{q(3\eta+\epsilon)}t^\eta,\endaligned$$
where we used (4.2), (4.3) in the first inequality, Lemma 4.1 in the second inequality and $|k|\leq |k_{123}|+|k_4|\lesssim 2^i$ in the third inequality. Moreover, we have for $\epsilon,\eta>0$ small enough $$\aligned E_t^2\lesssim&\varepsilon^\eta\sum_{k\in \mathbb{Z}^3\backslash\{0\}}\sum_{|i-j|\leq1,|i'-j'|\leq1}\sum_{k_{1234}=k}\frac{\theta(2^{-q}k)^2\theta(2^{-i} k_{123})\theta(2^{-j}k_4)\theta(2^{-i'} k_{234})\theta(2^{-j'}k_1)t^\eta}{|k_1|^2|k_2|^2|k_3|^2|k_4|^2|k_{12}||k_{24}||k_{123}|^{1-\eta}|k_{234}|^{1-\eta}}
\\&((|k_{123}||k_{234}|)^{\eta/2}+(|k_{12}||k_{24}|)^{\eta/2}+\sum_{i=1}^4|k_i|^{\eta})
\\\lesssim& \sum_{k\in \mathbb{Z}^3\backslash\{0\}}\sum_{k_{1234}=k}\frac{\theta(2^{-q}k)^2t^\eta2^{-q(2-4\eta)}
((|k_{123}||k_{234}|)^{\eta/2}+(|k_{12}||k_{24}|)^{\eta/2}+\sum_{i=1}^4|k_i|^{\eta})}{|k_1|^{1+2\eta}|k_2|^2|k_3|^2|k_4|^{1+2\eta}|k_{12}||k_{24}||k_{123}|^{1-\eta}|k_{234}|^{1-\eta}}
\\\lesssim& \varepsilon^\eta\sum_{k\in \mathbb{Z}^3\backslash\{0\}}\bigg[(\sum_{k_{1234}=k}\frac{\theta(2^{-q}k)^2t^\eta2^{-q(2-4\eta)}}{|k_1|^{1+2\eta}|k_2|^2|k_3|^2|k_4|^{1+2\eta}
|k_{12}|^2|k_{123}|^{2-3\eta}})^{1/2}\\&(\sum_{k_{1234}=k}\frac{\theta(2^{-q}k)^2t^\eta2^{-q(2-4\eta)}}{|k_1|^{1+2\eta}|k_2|^2|k_3|^2|k_4|^{1+2\eta}
|k_{24}|^2|k_{234}|^{2-3\eta}})^{1/2}
\\&+(\sum_{k_{1234}=k}\frac{\theta(2^{-q}k)^2t^\eta2^{-q(2-4\eta)}}{|k_1|^{1+2\eta}|k_2|^2|k_3|^2|k_4|^{1+2\eta}|k_{12}|^{2-\eta}|k_{123}|^{2-2\eta}})^{1/2}
\\&(\sum_{k_{1234}=k}\frac{\theta(2^{-q}k)^2t^\eta2^{-q(2-4\eta)}}{|k_1|^{1+2\eta}|k_2|^2|k_3|^2|k_4|^{1+2\eta}|k_{24}|^{2-\eta}|k_{234}|^{2-2\eta}})^{1/2}
\\&+(\sum_{k_{1234}=k}\frac{\theta(2^{-q}k)^2t^\eta2^{-q(2-4\eta)}\sum_{i=1}^4|k_i|^{\eta}}{|k_1|^{1+2\eta}|k_2|^2|k_3|^2|k_4|^{1+2\eta}|k_{12}|^2|k_{123}
|^{2-2\eta}})^{1/2}\\&(\sum_{k_{1234}=k}\frac{\theta(2^{-q}k)^2t^\eta2^{-q(2-4\eta)}\sum_{i=1}^4|k_i|^{\eta}}{|k_1|^{1+2\eta}|k_2|^2|k_3|^2|k_4|^{1+2\eta}
|k_{24}|^2|k_{234}|^{2-2\eta}})^{1/2}\bigg]
\\\lesssim& \sum_{k\in \mathbb{Z}^3\backslash\{0\}}2^{-(2-4\eta)q}\frac{t^\eta\varepsilon^\eta}{|k|^{1+\eta-\epsilon}}
\lesssim2^{q(3\eta+\epsilon)}t^\eta\varepsilon^\eta,\endaligned$$
where we used Lemma 4.1 in the fourth inequality.
By a similar argument we can also obtain the same bounds for $E_t^3, E_t^4$, $E_t^5$ $E_t^6, E_t^7$,  which implies that for $\eta,\epsilon>0$ small enough
$$E[|\Delta_qI_t^1|^2]\lesssim2^{q(3\eta+\epsilon)}t^\eta\varepsilon^\eta.$$
By a similar calculation as above we also get that for $\eta,\epsilon, \gamma>0$ small enough
$$\aligned &\sum_{i_0,j_0=1}^3E[|\Delta_q(\pi_{0,\diamond}(u_3^{\varepsilon,i_0},u_1^{\varepsilon,j_0})(t_1)-\pi_{0,\diamond}(u_3^{\varepsilon,i_0},u_1^{\varepsilon,j_0})(t_2)
-\pi_{0,\diamond}(\bar{u}_3^{\varepsilon,i_0},\bar{u}_1^{\varepsilon,j_0})(t_1)\\&+\pi_{0,\diamond}(\bar{u}_3^{\varepsilon,i_0},\bar{u}_1^{\varepsilon,j_0})(t_2))|^2]\\\lesssim & \varepsilon^\gamma|t_1-t_2|^\eta2^{q\epsilon},\endaligned$$
 which by Gaussian hypercontractivity and Lemma 2.1  implies that for $i_0,j_0=1,2,3,$ and $\delta>0$ small enough, $p>1$
$$\pi_{0,\diamond}(u_3^{\varepsilon,i_0},u_1^{\varepsilon,j_0})-\pi_{0,\diamond}(\bar{u}_3^{\varepsilon,i_0},\bar{u}_1^{\varepsilon,j_0})\rightarrow 0\textrm{ in } L^p(\Omega;C([0,T],\mathcal{C}^{-\delta})).$$

\subsection{Convergence of $\pi_0(P^{i_1i_2}D_{j_0}^\varepsilon K^{\varepsilon,j_0},u_1^{\varepsilon,j_1})-\pi_0(P^{i_1i_2}D_{j_0}\bar{K}^{\varepsilon,j_0},\bar{u}_1^{\varepsilon,j_1})$ and
$\pi_0(P^{i_1i_2}$ $D_{j_0}K^{\varepsilon,i_2},u_1^{\varepsilon,j_1})-\pi_0(P^{i_1i_2}D_{j_0}K^{\varepsilon,i_2},u_1^{\varepsilon,j_1})$}
In this subsection for $i_1,i_2,j_0,j_1=1,2,3$,  we consider $\pi_0(P^{i_1i_2}D_{j_0}K^{\varepsilon,j_0},u_1^{\varepsilon,j_1})-\pi_0(P^{i_1i_2}D_{j_0}\bar{K}^{\varepsilon,j_0},\bar{u}_1^{\varepsilon,j_1})$. Similar results for $\pi_0(P^{i_1i_2}D_{j_0}K^{\varepsilon,i_2},u_1^{\varepsilon,j_1})-\pi_0(P^{i_1i_2}D_{j_0}K^{\varepsilon,i_2},u_1^{\varepsilon,j_1})$ can be deduced . We have the following identity:
$$\aligned &\pi_0(P^{i_1i_2}D^\varepsilon_{j_0}K^{\varepsilon,j_0},u_1^{\varepsilon,j_1})-\pi_0(P^{i_1i_2}D_{j_0}\bar{K}^{\varepsilon,j_0},\bar{u}_1^{\varepsilon,j_1}):=I^1+I^2
,\endaligned$$
where
$$\aligned I^1=&(2\pi)^{-\frac{3}{2}}\sum_{k\in\mathbb{Z}^3\backslash\{0\}}\sum_{|i-j|\leq1}
\sum_{k_{12}=k}\theta(2^{-i}k_1)\theta(2^{-j}k_2)\bigg[\int_0^te^{-(t-s)f(\varepsilon k_1)|k_1|^2}
g(\varepsilon k_1^{j_0}) k_1^{j_0}:\hat{X}^{\varepsilon,j_0}_s(k_1)\hat{X}^{\varepsilon,j_1}_t(k_2):ds
\\&-\int_0^te^{-(t-s)|k_1|^2}
\imath k_1^{j_0}:\hat{\bar{X}}^{\varepsilon,j_0}_s(k_1)\hat{\bar{X}}^{\varepsilon,j_1}_t(k_2):ds\bigg]e_k
\hat{P}^{i_1i_2}(k_1),\endaligned$$
and
$$\aligned I^2=&(2\pi)^{-3}\sum_{|i-j|\leq1}\sum_{k_{1}}\theta(2^{-i}k_1)\theta(2^{-j}k_1)\bigg[\int_0^te^{-2(t-s)f(\varepsilon k_1)|k_1|^2}g(\varepsilon k_1^{j_0}) k_1^{j_0}\frac{h(\varepsilon k_1)^2}{2|k_1|^2f(\varepsilon k_1)}ds\\&-\int_0^te^{-2(t-s)|k_1|^2}\imath k_1^{j_0}\frac{h(\varepsilon k_1)^2}{2|k_1|^2}ds\bigg]\hat{P}^{i_1i_2}(k_1)\sum_{j_2=1}^3\hat{P}^{j_0j_2}(k_1)\hat{P}^{j_1j_2}(k_1).\endaligned$$
It is easy to get that the second term in the right hand side of the above equality equals zero and
$$\aligned I_2=&C_3^{\varepsilon,i_1,i_2,j_1,j_0}=(2\pi)^{-3}\sum_{|i-j|\leq1}\sum_{k_{1}}\theta(2^{-i}k_1)\theta(2^{-j}k_1)\int_0^te^{-2(t-s)f(\varepsilon k_1)|k_1|^2}g(\varepsilon k_1^{j_0}) k_1^{j_0}\frac{h(\varepsilon k_1)^2}{2|k_1|^2f(\varepsilon k_1)}ds\\&\hat{P}^{i_1i_2}(k_1)\sum_{j_2=1}^3\hat{P}^{j_0j_2}(k_1)\hat{P}^{j_1j_2}(k_1).\endaligned$$
Here we have
$$\sum_{j_0=1}^3C_3^{\varepsilon,i_1,i_2,j_1,j_0}=2C^{\varepsilon,i_1,i_2,j_1}.$$
Moreover, we have
$$\aligned &E|\Delta_qI^1|^2\\\lesssim&\sum_{k\in\mathbb{Z}^3\backslash\{0\}}\sum_{|i-j|\leq1,|i'-j'|\leq1}\sum_{k_{12}=k,k_{12}'=k}\theta(2^{-q}k)^2\theta(2^{-i}k_1)
\theta(2^{-j}k_2)\theta(2^{-i'}k_1')\theta(2^{-j'}k_2')\\&
\int_0^t\int_0^t|e^{-(t-s)f(\varepsilon k_1)|k_1|^2}k_1^{j_0}g(\varepsilon k_1^{j_0})-e^{-(t-s)|k_1|^2}\imath k_1^{j_0}|E|:\hat{X}_s^{\varepsilon,j_0}(k_1)\hat{X}_t^{\varepsilon,j_1}(k_2)::\hat{X}_{\bar{s}}^{\varepsilon,j_0}(k_1')\hat{X}_t^{\varepsilon,j_1}(k_2'):|
\\&|e^{-(t-\bar{s})f(\varepsilon k_1')|k_1'|^2}k_1^{'j_0}g(\varepsilon k_1^{'j_0})-e^{-(t-\bar{s})|k_1'|^2}k_1^{'j_0}|
dsd\bar{s}+\int_0^t\int_0^te^{-(t-s)|k_1|^2}|k_1|e^{-(t-\bar{s})|k_1'|^2}|k_1'|
\\&E|(:\hat{X}_s^{\varepsilon,j_0}(k_1)\hat{X}_t^{\varepsilon,j_1}(k_2):-:\hat{\bar{X}}_s^{\varepsilon,j_0}(k_1)\hat{\bar{X}}_t^{\varepsilon,j_1}(k_2):)
(:\hat{X}_{\bar{s}}^{\varepsilon,j_0}(k_1')\hat{X}_t^{\varepsilon,j_1}(k_2'):-:\hat{\bar{X}}_{\bar{s}}^{\varepsilon,j_0}(k_1')
\hat{\bar{X}}_t^{\varepsilon,j_1}(k_2'):)|dsd\bar{s}\endaligned$$
$$\aligned\lesssim &t^{\epsilon}\sum_k\sum_{q\lesssim i}\sum_{k_{12}=k}\theta(2^{-q}k)\theta(2^{-i}k_1)\big(\frac{1}{|k_1|^{4-2\epsilon}|k_2|^2}+\frac{1}{|k_1|^{3-\epsilon}|k_2|^{3-\epsilon}}\big)(|\varepsilon k_1|^\eta+|k_2|^{\eta}\varepsilon^\eta)
\\\lesssim& \varepsilon^{\eta}t^\epsilon2^{q(2\epsilon+\eta)},\endaligned$$
where we used (4.2), (4.3) and (4.7) in the second inequality and in the last inequality we used Lemma 4.1.
By a similar calculation we also get that for $\epsilon,\eta>0$ small enough
$$\aligned &E[|\Delta_q(\pi_{0,\diamond}(P^{i_1i_2}D_{j_0}^\varepsilon K^{\varepsilon,j_2},u_1^{\varepsilon,j_1})(t_1)
-\pi_{0,\diamond}(P^{i_1i_2}D_{j_0}^\varepsilon K^{\varepsilon,j_2},u_1^{\varepsilon,j_1})(t_2)\\&-\pi_{0,\diamond}(P^{i_1i_2}D_{j_0}\bar{K}^{\varepsilon,j_2},\bar{u}_1^{\varepsilon,j_1}
)(t_1)+\pi_{0,\diamond}(P^{i_1i_2}D_{j_0}\bar{K}^{\varepsilon,j_2},\bar{u}_1^{\varepsilon,j_1})(t_2))|^2]\\\lesssim & \varepsilon^\eta|t_1-t_2|^\eta2^{q(\epsilon+3\eta)},\endaligned$$
which by Gaussian hypercontractivity and Lemma 2.1  implies that for $i_1,i_2,j_0,j_1=1,2,3$ such that
$$\pi_{0,\diamond}(P^{i_1i_2}D_{j_0}^\varepsilon K^{\varepsilon,j_2},u_1^{\varepsilon,j_1})-\pi_{0,\diamond}(P^{i_1i_2}D_{j_0}\bar{K}^{\varepsilon,j_2},\bar{u}_1^{\varepsilon,j_1})\rightarrow 0\textrm{ in } C([0,T],\mathcal{C}^{-\delta}).$$
By a similar argument we also obtain that for $i_1,i_2,j_0,j_1=1,2,3$ such that
$$\pi_{0}(P^{i_1i_2}D_{j_0}^\varepsilon K^{\varepsilon,i_2},u_1^{\varepsilon,j_1})-\tilde{C}^{\varepsilon,i_1,j_0,j_1,i_2}_3-\pi_{0,\diamond}(P^{i_1i_2}D_{j_0} \bar{K}^{\varepsilon,i_2},\bar{u}_1^{\varepsilon,j_1})\rightarrow 0\textrm{ in } L^p(\Omega;C([0,T],\mathcal{C}^{-\delta})),$$
where
$$\aligned &\tilde{C}^{\varepsilon,i_1,j_0,j_1,i_2}_3\\=&(2\pi)^{-3}\sum_{i_3=1}^3\sum_{k_2\in\mathbb{Z}^3\backslash\{0\}}\int_0^te^{-2|k_2|^2(t-s)f(\varepsilon k_2)} k_2^{j_0}g(\varepsilon k_2^{j_0})\frac{h(\varepsilon k_2)^2}{2|k_2|^2f(\varepsilon k_2)}\hat{P}^{i_1i_2}(k_2)\hat{P}^{i_2i_3}(k_{2})\hat{P}^{j_1i_3}(k_{2})ds,\endaligned$$
$$\sum_{i_2=1}^3\tilde{C}^{\varepsilon,i_1,j_0,j_1,i_2}_3=2\tilde{C}^{\varepsilon,i_1,j_0,j_1}.$$

\subsection{Convergence of  $u_2^{\varepsilon,i}u_2^{\varepsilon,j}-\bar{u}_2^{\varepsilon,i}\bar{u}_2^{\varepsilon,j}$}
In this subsection for $i,j=1,2,3$, we deal with $u_2^{\varepsilon,i}u_2^{\varepsilon,j}$ and prove that $u_2^{\varepsilon,i}\diamond u_2^{\varepsilon,j}-\bar{u}_2^{\varepsilon,i}\diamond\bar{u}_2^{\varepsilon,j}\rightarrow 0$ in $L^p(\Omega;C([0,T];\mathcal{C}^{-\delta}))$.
Recall that for $i,j=1,2,3$
$${u}_2^{\varepsilon,i}\diamond {u}_2^{\varepsilon,j}:={u}_2^{\varepsilon,i} {u}_2^{\varepsilon,j}-{\varphi}_2^{\varepsilon,ij}(t)-{C}_2^{\varepsilon,ij},$$
$$\bar{u}_2^{\varepsilon,i}\diamond \bar{u}_2^{\varepsilon,j}:=\bar{u}_2^{\varepsilon,i} \bar{u}_2^{\varepsilon,j}-\bar{\varphi}_2^{\varepsilon,ij}(t)-\bar{C}_2^{\varepsilon,ij}.$$
We have the following identities:
$$\aligned &u_2^{\varepsilon,i}u_2^{\varepsilon,j}-\bar{u}_2^{\varepsilon,i}\bar{u}_2^{\varepsilon,j}:=L^1+L^2+L^3,\endaligned$$
where
$$\aligned L_t^1=&(2\pi)^{-9/2}\bigg[\sum_{i_1,i_2,j_1,j_2=1}^3
\sum_{k_{1234}=k}\bigg(\int_0^t\int_0^te^{-|k_{12}|^2f(\varepsilon k_{12})(t-s)-|k_{34}|^2f(\varepsilon k_{34})(t-\bar{s})}k_{12}^{i_2}g(\varepsilon k_{12}^{i_2})k_{34}^{j_2}g(\varepsilon k_{34}^{j_2})
\\&:\hat{X}_s^{\varepsilon,i_1}(k_1)\hat{X}^{\varepsilon,i_2}_s(k_2)\hat{X}^{\varepsilon,j_1}_{\bar{s}}(k_3)\hat{X}^{\varepsilon,j_2}_{\bar{s}}(k_4):
dsd\bar{s}e_k\\&-\int_0^t\int_0^te^{-|k_{12}|^2(t-s)-|k_{34}|^2(t-\bar{s})}
:\hat{\bar{X}}_s^{\varepsilon,i_1}(k_1)\hat{\bar{X}}^{\varepsilon,i_2}_s(k_2)\hat{\bar{X}}^{\varepsilon,j_1}_{\bar{s}}(k_3)\hat{\bar{X}}^{\varepsilon,j_2}_{\bar{s}}(k_4):dsd\bar{s}e_k
\\&\imath k^{i_2}_{12}
\imath k^{j_2}_{34}\bigg)\hat{P}^{ii_1}(k_{12})\hat{P}^{jj_1}(k_{34})\bigg]
\endaligned$$
$$L^2=\sum_{i=1}^4{I^i},$$
with
$$\aligned I^1_t=&(2\pi)^{-9/2}\sum_{i_1,i_2,j_1,j_2=1}^3\sum_{k_{24}=k,k_1}\bigg(\int_0^t\int_0^te^{-|k_{12}|^2f(\varepsilon k_{12})(t-s)-|k_{4}-k_1|^2f(\varepsilon (k_4-k_1))(t-\bar{s})}\frac{h(\varepsilon k_1)^2e^{-|k_1|^2f(\varepsilon k_{1})|s-\bar{s}|}}{2|k_1|^2f(\varepsilon k_{1})}\\&:\hat{X}^{\varepsilon,i_2}_s(k_2)\hat{X}^{\varepsilon,j_2}_{\bar{s}}(k_4):dsd\bar{s}e_k\hat{P}^{ii_1}
(k_{12})g(\varepsilon(k_{12}^{i_2})) k^{i_2}_{12}\hat{P}^{jj_1}(k_4-k_1)g(\varepsilon(k_4^{j_2}-k_1^{j_2}))(k^{j_2}
_{4}-k_1^{j_2})\\&\sum_{j_3=1}^3\hat{P}^{i_1j_3}(k_1)\hat{P}^{j_1j_3}(k_1)-\int_0^t\int_0^te^{-|k_{12}|^2(t-s)-|k_{4}-k_1|^2(t-\bar{s})}\frac{h(\varepsilon k_1)^2e^{-|k_1|^2|s-\bar{s}|}}{2|k_1|^2}:\hat{\bar{X}}^{\varepsilon,i_2}_s(k_2)\hat{\bar{X}}^{\varepsilon,j_2}_{\bar{s}}(k_4):dsd\bar{s}e_k\\&\hat{P}^{ii_1}
(k_{12})\imath k^{i_2}_{12}\imath\hat{P}^{jj_1}(k_{4}-k_1)(k^{j_2}
_{4}-k_1^{j_2})\sum_{j_3=1}^3\hat{P}^{i_1j_3}(k_1)\hat{P}^{j_1j_3}(k_1)\bigg),\endaligned$$
$$\aligned I^2_t=&(2\pi)^{-9/2}\sum_{i_1,i_2,j_1,j_2=1}^3\sum_{k_{24}=k,k_1}\bigg(\int_0^t\int_0^te^{-|k_{12}|^2f(\varepsilon k_{12})(t-s)-|k_{4}-k_1|^2f(\varepsilon (k_4-k_1))(t-\bar{s})}\frac{h(\varepsilon k_1)^2e^{-|k_1|^2f(\varepsilon k_{1})|s-\bar{s}|}}{2|k_1|^2f(\varepsilon k_{1})}\\&:\hat{X}^{\varepsilon,i_1}_s(k_2)\hat{X}^{\varepsilon,j_2}_{\bar{s}}(k_4):dsd\bar{s}e_k\hat{P}^{ii_1}
(k_{12})g(\varepsilon(k_{12}^{i_2})) k^{i_2}_{12}\hat{P}^{jj_1}(k_4-k_1)g(\varepsilon(k_4^{j_2}-k_1^{j_2}))(k^{j_2}
_{4}-k_1^{j_2})\\&\sum_{j_3=1}^3\hat{P}^{i_2j_3}(k_1)\hat{P}^{j_1j_3}(k_1)-\int_0^t\int_0^te^{-|k_{12}|^2(t-s)-|k_{4}-k_1|^2(t-\bar{s})}\frac{h(\varepsilon k_1)^2e^{-|k_1|^2|s-\bar{s}|}}{2|k_1|^2}:\hat{\bar{X}}^{\varepsilon,i_1}_s(k_2)\hat{\bar{X}}^{\varepsilon,j_2}_{\bar{s}}(k_4):dsd\bar{s}e_k\\&\hat{P}^{ii_1}
(k_{12})\imath k^{i_2}_{12}\imath\hat{P}^{jj_1}(k_{4}-k_1)(k^{j_2}
_{4}-k_1^{j_2})\sum_{j_3=1}^3\hat{P}^{i_2j_3}(k_1)\hat{P}^{j_1j_3}(k_1)\bigg),\endaligned$$
$$\aligned I^3_t=&(2\pi)^{-9/2}\sum_{i_1,i_2,j_1,j_2=1}^3\sum_{k_{24}=k,k_1}\bigg(\int_0^t\int_0^te^{-|k_{12}|^2f(\varepsilon k_{12})(t-s)-|k_{4}-k_1|^2f(\varepsilon (k_4-k_1))(t-\bar{s})}\frac{h(\varepsilon k_1)^2e^{-|k_1|^2f(\varepsilon k_{1})|s-\bar{s}|}}{2|k_1|^2f(\varepsilon k_{1})}\\&:\hat{X}^{\varepsilon,i_2}_s(k_2)\hat{X}^{\varepsilon,j_1}_{\bar{s}}(k_4):dsd\bar{s}e_k\hat{P}^{ii_1}
(k_{12})g(\varepsilon(k_{12}^{i_2})) k^{i_2}_{12}\hat{P}^{jj_1}(k_4-k_1)g(\varepsilon(k_4^{j_2}-k_1^{j_2}))(k^{j_2}
_{4}-k_1^{j_2})\\&\sum_{j_3=1}^3\hat{P}^{i_1j_3}(k_1)\hat{P}^{j_2j_3}(k_1)-\int_0^t\int_0^te^{-|k_{12}|^2(t-s)-|k_{4}-k_1|^2(t-\bar{s})}\frac{h(\varepsilon k_1)^2e^{-|k_1|^2|s-\bar{s}|}}{2|k_1|^2}:\hat{\bar{X}}^{\varepsilon,i_2}_s(k_2)\hat{\bar{X}}^{\varepsilon,j_1}_{\bar{s}}(k_4):dsd\bar{s}e_k\\&\hat{P}^{ii_1}
(k_{12})\imath k^{i_2}_{12}\imath\hat{P}^{jj_1}(k_{4}-k_1)(k^{j_2}
_{4}-k_1^{j_2})\sum_{j_3=1}^3\hat{P}^{i_1j_3}(k_1)\hat{P}^{j_2j_3}(k_1)\bigg),\endaligned$$
$$\aligned I^4_t=&(2\pi)^{-9/2}\sum_{i_1,i_2,j_1,j_2=1}^3\sum_{k_{24}=k,k_1}\bigg(\int_0^t\int_0^te^{-|k_{12}|^2f(\varepsilon k_{12})(t-s)-|k_{4}-k_1|^2f(\varepsilon (k_4-k_1))(t-\bar{s})}\frac{h(\varepsilon k_1)^2e^{-|k_1|^2f(\varepsilon k_{1})|s-\bar{s}|}}{2|k_1|^2f(\varepsilon k_{1})}\\&:\hat{X}^{\varepsilon,i_1}_s(k_2)\hat{X}^{\varepsilon,j_1}_{\bar{s}}(k_4):dsd\bar{s}e_k\hat{P}^{ii_1}
(k_{12})g(\varepsilon(k_{12}^{i_2})) k^{i_2}_{12}\hat{P}^{jj_1}(k_4-k_1)g(\varepsilon(k_4^{j_2}-k_1^{j_2}))(k^{j_2}
_{4}-k_1^{j_2})\\&\sum_{j_3=1}^3\hat{P}^{i_2j_3}(k_1)\hat{P}^{j_2j_3}(k_1)-\int_0^t\int_0^te^{-|k_{12}|^2(t-s)-|k_{4}-k_1|^2(t-\bar{s})}\frac{h(\varepsilon k_1)^2e^{-|k_1|^2|s-\bar{s}|}}{2|k_1|^2}:\hat{\bar{X}}^{\varepsilon,i_1}_s(k_2)\hat{\bar{X}}^{\varepsilon,j_1}_{\bar{s}}(k_4):dsd\bar{s}e_k\\&\hat{P}^{ii_1}
(k_{12})\imath k^{i_2}_{12}\imath\hat{P}^{jj_1}(k_{4}-k_1)(k^{j_2}
_{4}-k_1^{j_2})\sum_{j_3=1}^3\hat{P}^{i_2j_3}(k_1)\hat{P}^{j_2j_3}(k_1)\bigg),\endaligned$$
and
$$\aligned L^3_t=&(2\pi)^{-6}\sum_{i_1,i_2,j_1,j_2=1}^3\sum_{k_1,k_2}\bigg[\int_0^t\int_0^t
e^{-|k_{12}|^2f(\varepsilon k_{12})(t-s+t-\bar{s})}\frac{h(\varepsilon k_1)^2h(\varepsilon k_2)^2e^{-(|k_1|^2f(\varepsilon k_{1})+|k_2|^2f(\varepsilon k_{2}))|s-\bar{s}|}}{4|k_1|^2|k_2|^2f(\varepsilon k_{1})f(\varepsilon k_{2})}dsd\bar{s}\\&\hat{P}^{ii_1}(k_{12})\hat{P}^{jj_1}(k_{12})g(\varepsilon k_{12}^{i_2}) k_{12}^{i_2}(-g(-\varepsilon k_{12}^{j_2}) k_{12}^{j_2})\\&-\int_0^t\int_0^t
e^{-|k_{12}|^2(t-s+t-\bar{s})}\frac{h(\varepsilon k_1)^2h(\varepsilon k_2)^2e^{-(|k_1|^2+|k_2|^2)|s-\bar{s}|}}{4|k_1|^2|k_2|^2}dsd\bar{s}\hat{P}^{ii_1}(k_{12})\hat{P}^{jj_1}(k_{12})\\&\imath k_{12}^{i_2}(-\imath k_{12}^{j_2})\bigg]\sum_{j_3,j_4=1}^3
\big(\hat{P}^{i_1j_3}(k_1)\hat{P}^{j_1j_3}(k_1)\hat{P}^{i_2j_4}(k_2)\hat{P}^{j_2j_4}(k_2)
+\hat{P}^{i_1j_3}(k_1)\hat{P}^{j_2j_3}(k_1)\hat{P}^{i_2j_4}(k_2)\hat{P}^{j_1j_4}(k_2)\big).\endaligned$$
By a easy computation we obtain that
$$\aligned L_t^3=&(2\pi)^{-6}\sum_{i_1,i_2,j_1,j_2=1}^3\sum_{k_1,k_2}h(\varepsilon k_1)^2h(\varepsilon k_2)^2\hat{P}^{ii_1}(k_{12})\hat{P}^{jj_1}(k_{12})\sum_{j_3,j_4=1}^3
\big(\hat{P}^{i_1j_3}(k_1)\hat{P}^{j_2j_3}(k_1)\hat{P}^{i_2j_4}(k_2)\\&\hat{P}^{j_1j_4}(k_2)+\hat{P}^{i_1j_3}(k_1)\hat{P}^{j_1j_3}(k_1)\hat{P}^{i_2j_4}(k_2)
\hat{P}^{j_2j_4}(k_2)\big)\\&
\bigg[-k_{12}^{i_2}g(\varepsilon k_{12}^{i_2})k_{12}^{j_2}g(-\varepsilon k_{12}^{j_2})\frac{1}{2|k_1|^2f(\varepsilon k_1)|k_2|^2f(\varepsilon k_2)(|k_1|^2f(\varepsilon k_1)+|k_2|^2f(\varepsilon k_2)+|k_{12}|^2f(\varepsilon k_{12}))}\\&(\frac{1-e^{-2|k_{12}|^2tf(\varepsilon k_{12})}}{2|k_{12}|^2f(\varepsilon k_{12})}-
\int_0^te^{-2|k_{12}|^2f(\varepsilon k_{12})(t-s)-(|k_1|^2f(\varepsilon k_{1})+|k_2|^2f(\varepsilon k_{2})+|k_{12}|^2f(\varepsilon k_{12}))s}ds)\\&-k_{12}^{i_2}k_{12}^{j_2}\frac{1}{2|k_1|^2|k_2|^2(|k_1|^2+|k_2|^2+|k_{12}|^2)}(\frac{1-e^{-2|k_{12}|^2t}}{2|k_{12}|^2}-
\int_0^te^{-2|k_{12}|^2(t-s)-(|k_1|^2+|k_2|^2+|k_{12}|^2)s}ds)\bigg].\endaligned$$
Let
$$\aligned C_2^{\varepsilon,ij}=&(2\pi)^{-6}\sum_{i_1,i_2,j_1,j_2=1}^3\sum_{k_1,k_2}h(\varepsilon k_1)^2h(\varepsilon k_2)^2\hat{P}^{ii_1}(k_{12})\hat{P}^{jj_1}(k_{12})k_{12}^{i_2}g(\varepsilon k_{12}^{i_2})(-k_{12}^{j_2})g(-\varepsilon k_{12}^{j_2})\sum_{j_3,j_4=1}^3\bigg(
\hat{P}^{i_1j_3}(k_1)\\&\hat{P}^{j_1j_3}(k_1)\hat{P}^{i_2j_4}(k_2)\hat{P}^{j_2j_4}(k_2)+\hat{P}^{i_1j_3}(k_1)\hat{P}^{j_2j_3}(k_1)\hat{P}^{i_2j_4}(k_2)
\hat{P}^{j_1j_4}(k_2)\bigg)
\\&\frac{1}{4|k_1|^2f(\varepsilon k_1)|k_2|^2f(\varepsilon k_2)(|k_1|^2f(\varepsilon k_1)+|k_2|^2f(\varepsilon k_2)+|k_{12}|^2f(\varepsilon k_{12}))}\frac{1}{|k_{12}|^2f(\varepsilon k_{12})}.\endaligned$$

$$\aligned \bar{C}_2^{\varepsilon,ij}=&(2\pi)^{-6}\sum_{i_1,i_2,j_1,j_2=1}^3\sum_{k_1,k_2}h(\varepsilon k_1)^2h(\varepsilon k_2)^2\hat{P}^{ii_1}(k_{12})\hat{P}^{jj_1}(k_{12})k_{12}^{i_2}k_{12}^{j_2}\sum_{j_3,j_4=1}^3
\bigg(\hat{P}^{i_1j_3}(k_1)\hat{P}^{j_1j_3}(k_1)\\&\hat{P}^{i_2j_4}(k_2)\hat{P}^{j_2j_4}(k_2)+\hat{P}^{i_1j_3}(k_1)\hat{P}^{j_2j_3}(k_1)
\hat{P}^{i_2j_4}(k_2)\hat{P}^{j_1j_4}(k_2)\bigg)\\&
\frac{1}{4|k_1|^2|k_2|^2(|k_1|^2+|k_2|^2+|k_{12}|^2)}\frac{1}{|k_{12}|^2}.\endaligned$$
Define $$\varphi_2^{\varepsilon,ij}-\bar{\varphi}_2^{\varepsilon,ij}=L_t^3-C_2^{\varepsilon,ij}+\bar{C}_2^{\varepsilon,ij},$$
where $\varphi_2^{\varepsilon}, \bar{\varphi}_2^{\varepsilon}$ corresponds to $u^\varepsilon, \bar{u}^\varepsilon$ respectively.
Then for $\rho>0$ we have
$$\aligned |\varphi_2^\varepsilon-\bar{\varphi}_2^\varepsilon|\lesssim& \bigg|\sum_{k_1,k_2}\frac{-k_{12}^{i_2}k_{12}^{j_2}g(\varepsilon k_{12}^{i_2})g(-\varepsilon k_{12}^{j_2})}{|k_1|^2f(\varepsilon k_1)|k_2|^2f(\varepsilon k_2)(|k_1|^2f(\varepsilon k_1)+|k_2|^2f(\varepsilon k_2)+|k_{12}|^2f(\varepsilon k_{12}))}(\frac{e^{-2|k_{12}|^2f(\varepsilon k_{12})t}}{2|k_{12}|^2f(\varepsilon k_{12})}\\&+\int_0^te^{-2|k_{12}|^2f(\varepsilon k_{12})(t-s)-(|k_1|^2f(\varepsilon k_1)+|k_2|^2f(\varepsilon k_2)+|k_{12}|^2f(\varepsilon k_{12}))s}ds)\\&+ \sum_{k_1,k_2}k_{12}^{i_2}k_{12}^{j_2}\frac{1}{|k_1|^2|k_2|^2(|k_1|^2+|k_2|^2+|k_{12}|^2)}(\frac{e^{-2|k_{12}|^2t}}{2|k_{12}|^2}
+\int_0^te^{-2|k_{12}|^2(t-s)-(|k_1|^2+|k_2|^2+|k_{12}|^2)s}ds)\bigg|
\\\lesssim &t^{-\rho}\sum_{k_1,k_2}\frac{\varepsilon^\eta(|k_1|^\eta+|k_2|^\eta+|k_{12}|^\eta)}{|k_1|^2|k_2|^2|k_{12}|^{2+2\rho}}\lesssim t^{-\rho}\varepsilon^\eta,\endaligned$$
where $2\rho>\eta>0$.
Then $\varphi^\varepsilon_{2}-\bar{\varphi}^\varepsilon_{2}$ converges to  $0$ with respect to $\|\varphi\|=\sup_{t\in[0,T]}t^\rho|\varphi(t)|$ for any $\rho>0$.

\textbf{Terms in the second chaos}:
Now we come to $L_t^2$: it is sufficient to consider $I_t^1$ and the desired estimates for other terms can be obtained similarly.  For $\epsilon>0$ small enough we have the following inequalities
$$\aligned &E|\Delta_qI_t^i|^2\\\lesssim &\sum_k\sum_{k_{24}=k,k_{24}'=k,k_1,k_1'}\sum_{i_2,j_2,i_2',j_2'=1}^3\theta(2^{-q}k)^2\bigg[\int_0^t\int_0^t\int_0^t\int_0^t\bigg|\bigg(e^{-|k_{12}|^2f(\varepsilon k_{12})(t-s)-|k_{4}-k_1|^2f(\varepsilon (k_{4}-k_1))(t-\bar{s})}k_{12}^{i_2}g(\varepsilon k_{12}^{i_2})\\&(k_4^{j_2}-k_1^{j_2})g(\varepsilon (k_4^{j_2}-k_1^{j_2}))\frac{e^{-|k_1|^2f(\varepsilon k_1)|s-\bar{s}|}}{|k_1|^2f(\varepsilon k_1)}-e^{-|k_{12}|^2(t-s)-|k_{4}-k_1|^2(t-\bar{s})}\imath \imath k_{12}^{i_2}(k_4^{j_2}-k_1^{j_2})\frac{e^{-|k_1|^2|s-\bar{s}|}}{|k_1|^2}\bigg)\\&\bigg(\frac{e^{-|k_1'|^2f(\varepsilon k_1')|s-\bar{s}|}}{|k_1'|^2f(\varepsilon k_1')}e^{-|k_{12}'|^2f(\varepsilon k_{12}')(t-\sigma)-|k_{4}'-k_1'|^2f(\varepsilon (k_{4}'-k_1'))(t-\bar{\sigma})}k_{12}^{'i_2'}g(\varepsilon k_{12}^{'i_2'})(k_4^{'j_2'}-k_1^{'j_2'})g(\varepsilon (k_4^{j_2'}-k_1^{j_2'}))\\&-e^{-|k_{12}'|^2(t-\sigma)-|k_{4}'-k_1'|^2(t-\bar{\sigma})}\frac{e^{-|k_1'|^2|s-\bar{s}|}}{|k_1'|^2}\imath \imath k_{12}^{'i_2'}(k_4^{'j_2'}-k_1^{'j_2'})\bigg)\frac{1_{\{k_2=k_2',k_4=k_4'\}}+1_{\{k_2=k_4',k_4=k_2'\}}}{|k_2|^2|k_4|^2}\bigg|dsd\bar{s}d\sigma d\bar{\sigma}\\+&\int_0^t\int_0^t\int_0^t\int_0^t|e^{-|k_{12}|^2(t-\sigma)-|k_{4}-k_1|^2(t-\bar{\sigma})}\imath \imath k_{12}^{i_2}(k_4^{j_2}-k_1^{j_2})e^{-|k_{12}'|^2(t-s)-|k_{4}'-k_1'|^2(t-\bar{s})}\imath \imath k_{12}^{'i_2'}(k_4^{'j_2'}-k_1^{'j_2'})\\&\frac{|\varepsilon k_2|^{\eta}+|\varepsilon k_4|^{\eta}}{|k_1|^2|k_1'|^2|k_2|^2|k_4|^2}1_{\{k_2=k_2',k_4=k_4'\}}+1_{\{k_2=k_4',k_4=k_2'\}}\bigg|dsd\bar{s}d\sigma d\bar{\sigma}\bigg]
\\\lesssim & t^\epsilon\sum_k \sum_{k_{24}=k,k_1,k_3}\bigg{(}\frac{\theta(2^{-q}k)^2\bigg[(|\varepsilon k_{12}|^{\frac{\eta}{2}}+|\varepsilon(k_4-k_1)|^{\frac{\eta}{2}}+|\varepsilon k_1|^{\frac{\eta}{2}})(|\varepsilon k_{23}|^{\frac{\eta}{2}}+|\varepsilon(k_4-k_3)|^{\frac{\eta}{2}}+|\varepsilon k_3|^{\frac{\eta}{2}})}{|k_1|^2|k_2|^2|k_3|^2|k_4|^2|k_1-k_4|^{1-\epsilon}|k_4-k_3||k_{12}|^{1-\epsilon}}\\&\frac{+|\varepsilon k_2|^{\eta}+|\varepsilon k_4|^{\eta}\bigg]}{|k_{23}|}+\frac{\theta(2^{-q}k)^2\bigg[(|\varepsilon k_{12}|^{\frac{\eta}{2}}+|\varepsilon(k_4-k_1)|^{\frac{\eta}{2}}+|\varepsilon k_1|^{\frac{\eta}{2}})(|\varepsilon k_{34}|^{\frac{\eta}{2}}+|\varepsilon(k_2-k_3)|^{\frac{\eta}{2}}+|\varepsilon k_3|^{\frac{\eta}{2}})}{|k_1|^2|k_2|^2|k_3|^2|k_4|^2|k_1-k_4|^{1-\epsilon}|k_2-k_3||k_{12}|^{1-\epsilon}}\\&\frac{+|\varepsilon k_2|^{\eta}+|\varepsilon k_4|^{\eta}\bigg]}{|k_{34}|}\bigg{)},\endaligned$$
where in the first inequality we used (4.7).
Now in the following we only estimate the first term on the right hand side of the inequality and the second term can be estimated similarly:
$$\aligned E|\Delta_qI_t^i|^2
\lesssim & t^\epsilon\sum_k\sum_{k_{24}=k}\frac{\theta(2^{-q}k)^2}{|k_2|^2|k_4|^2}\sum_{k_1}\frac{1}{|k_1-k_4|^{1-\epsilon}|k_1|^2|k_{12}|^{1-\epsilon}}
\sum_{k_3}\frac{1}{|k_3-k_4||k_3|^2|k_{23}|}\bigg[(|\varepsilon k_{12}|^{\frac{\eta}{2}}\\&+|\varepsilon(k_4-k_1)|^{\frac{\eta}{2}}+|\varepsilon k_1|^{\frac{\eta}{2}})(|\varepsilon k_{23}|^{\frac{\eta}{2}}+|\varepsilon(k_4-k_3)|^{\frac{\eta}{2}}+|\varepsilon k_3|^{\frac{\eta}{2}})+|\varepsilon k_2|^{\eta}+|\varepsilon k_4|^{\eta}\bigg]\\\lesssim & t^\epsilon\sum_k\sum_{k_{24}=k}\frac{\theta(2^{-q}k)^2}{|k_2|^2|k_4|^2}(\sum_{k_1}\frac{1}{|k_1-k_4|^{2-2\epsilon}|k_1|^2})^{1/2}
(\sum_{k_1}\frac{1}{|k_{12}|^{2-2\epsilon}|k_1|^2})^{1/2}(\sum_{k_3}\frac{1}{|k_3-k_4|^2|k_3|^2})^{1/2}\\&(\sum_{k_3}\frac{1}{|k_{23}|^2|k_3|^2})^{1/2}
\bigg[(|\varepsilon k_{12}|^{\frac{\eta}{2}}+|\varepsilon(k_4-k_1)|^{\frac{\eta}{2}}+|\varepsilon k_1|^{\frac{\eta}{2}})(|\varepsilon k_{23}|^{\frac{\eta}{2}}+|\varepsilon(k_4-k_3)|^{\frac{\eta}{2}}+|\varepsilon k_3|^{\frac{\eta}{2}})\\&+|\varepsilon k_2|^{\eta}+|\varepsilon k_4|^{\eta}\bigg]\\\lesssim &t^\epsilon\varepsilon^\eta\sum_k\sum_{k_{24}=k}\bigg(\frac{\theta(2^{-q}k)^2}{|k_2|^{3-\epsilon-\eta}|k_4|^{3-\epsilon}}+\frac{\theta(2^{-q}k)^2}{|k_2|^{3-\epsilon}
|k_4|^{3-\epsilon-\eta}}\bigg)\lesssim \varepsilon^\eta2^{2\epsilon+\eta}t^\epsilon,\endaligned$$
where  in the last two inequalities we used Lemma 4.1.

\textbf{Terms in the fourth chaos}:

Now we consider $L_t^1$: For $\epsilon,\eta>0$ small enough we have the following calculations:
$$\aligned &E|\Delta_qL_t^1|^2\\\lesssim&\sum_k\sum_{k_{1234}=k}\sum_{i_1,i_2,j_1,j_2=1}^3\sum_{i_1',i_2',j_1',j_2'=1}^3\theta(2^{-q}k)^2\bigg(\int_0^t\int_0^t
\int_0^t\int_0^t\bigg|(e^{-f(\varepsilon k_{12})|k_{12}|^2(t-s)
-|k_{34}|^2f(\varepsilon k_{34})(t-\bar{s})}\\&k_{12}^{i_2}g(\varepsilon k_{12}^{i_2})k_{34}^{j_2}g(\varepsilon k_{34}^{j_2})-e^{-|k_{12}|^2(t-s)
-|k_{34}|^2(t-\bar{s})}k_{12}^{i_2}\imath\imath k_{34}^{i_2})\\&(e^{-f(\varepsilon k_{12}')|k_{12}'|^2(t-\sigma)
-|k_{34}'|^2f(\varepsilon k_{34}')(t-\bar{\sigma})}k_{12}^{'i_2'}g(\varepsilon k_{12}^{'i_2'})k_{34}^{'j_2'}
g(\varepsilon k_{34}^{'j_2'})-e^{-|k_{12}'|^2(t-s)
-|k_{34}'|^2(t-\bar{s})}k_{12}^{'i_2'}\imath\imath k_{34}^{'i_2'})\bigg|\\&(1_{\{k_1=k_1',k_2=k_2',k_3=k_3',k_4=k_4'\}}
+1_{\{k_1=k_3',k_2=k_2',k_3=k_1',k_4=k_4'\}})\frac{1}{|k_1|^2|k_2|^2|k_3|^2|k_4|^2}dsd\bar{s}d\sigma d\bar{\sigma}\\&+\int_0^t\int_0^t\int_0^t\int_0^t\bigg|e^{-|k_{12}|^2(t-s)-|k_{34}|^2(t-\bar{s})}k_{12}^{i_2}k_{34}^{j_2}\imath\imath
e^{-|k_{12}'|^2(t-s)-|k_{34}'|^2(t-\bar{s})}k_{12}^{'i_2'}k_{34}^{'j_2'}\imath\imath\bigg|\\&
\frac{\sum_{i=1}^4|\varepsilon k_i|^\eta}{|k_1|^2|k_2|^2|k_3|^2|k_4|^2}(1_{\{k_1=k_1',k_2=k_2',k_3=k_3',k_4=k_4'\}}
+1_{\{k_1=k_3',k_2=k_2',k_3=k_1',k_4=k_4'\}})dsd\bar{s}d\sigma d\bar{\sigma}\bigg)\\\lesssim& t^\epsilon\sum_k\sum_{k_{1234}=k}\theta(2^{-q}k)^2\bigg(\frac{|\varepsilon k_{12}|^\eta+|\varepsilon k_{34}|^\eta+\sum_{i=1}^4|\varepsilon k_i|^\eta}{|k_1|^2|k_2|^2|k_3|^2|k_4|^2|k_{12}|^{2-\epsilon}|k_{34}|^{2-\epsilon}}\\&+\frac{(|\varepsilon k_{12}|^{\frac{\eta}{2}}+|\varepsilon k_{34}|^{\frac{\eta}{2}})(|\varepsilon k_{14}|^{\frac{\eta}{2}}+|\varepsilon k_{23}|^{\frac{\eta}{2}})+\sum_{i=1}^4|\varepsilon k_i|^\eta}{|k_1|^2|k_2|^2|k_3|^2|k_4|^2|k_{12}|^{1-\epsilon/2}|k_{34}|^{1-\epsilon/2}|k_{14}|^{1-\epsilon/2}|k_{23}|^{1-\epsilon/2}}\bigg)
\endaligned$$
$$\aligned\lesssim&\varepsilon^\eta t^\epsilon\bigg(\sum_k\bigg[(\sum_{k_{1234}=k}\frac{\theta(2^{-q}k)^2(| k_{12}|^{\eta}+|k_{34}|^{\eta}) }{|k_1|^2|k_2|^2|k_3|^2|k_4|^2|k_{12}|^{2-\epsilon}|k_{34}|^{2-\epsilon}})^{1/2}(\sum_{k_{1234}=k}\frac{\theta(2^{-q}k)^2 (| k_{14}|^{\eta}+|k_{23}|^{\eta}) }{|k_1|^2|k_2|^2|k_3|^2|k_4|^2|k_{14}|^{2-\epsilon}|k_{23}|^{2-\epsilon}})^{1/2}\\&+(\sum_{k_{1234}=k}\frac{\theta(2^{-q}k)^2\sum_{i=1}^4|k_i|^\eta }{|k_1|^2|k_2|^2|k_3|^2|k_4|^2|k_{12}|^{2-\epsilon}|k_{34}|^{2-\epsilon}})^{1/2}(\sum_{k_{1234}=k}\frac{\theta(2^{-q}k)^2 \sum_{i=1}^4| k_i|^\eta }{|k_1|^2|k_2|^2|k_3|^2|k_4|^2|k_{14}|^{2-\epsilon}|k_{23}|^{2-\epsilon}})^{1/2}\bigg]\\&+2^{q(2\epsilon+\eta)}\bigg)\\\lesssim&\varepsilon^\eta2^{q(2\epsilon+\eta)}t^\epsilon,\endaligned$$
where we used (4.8) in the first inequality and  Lemma 4.1 in the last inequality.
By a similar calculation we also get that for $\epsilon,\eta>0$ small enough
$$\aligned &E[|\Delta_q(u_2^{\varepsilon,i}\diamond u_2^{\varepsilon,j}(t_1)-u_2^{\varepsilon,i}\diamond u_2^{\varepsilon,j}(t_2)-\bar{u}_2^{\varepsilon,i}\diamond \bar{u}_2^{\varepsilon_2,j}(t_1)+\bar{u}_2^{\varepsilon_2,i}\diamond \bar{u}_2^{\varepsilon_2,j})(t_2))|^2]\\\lesssim & \varepsilon^\eta|t_1-t_2|^\eta2^{q(\epsilon+3\eta)},\endaligned$$
 which by Gaussian hypercontractivity and Lemma 2.1  implies that for $i,j=1,2,3$ and every $\delta>0$
$$u_2^{\varepsilon,i}\diamond u_2^{\varepsilon,j}-\bar{u}_2^{\varepsilon,i}\diamond \bar{u}_2^{\varepsilon,j}\rightarrow 0\textrm{ in } L^p(\Omega; C([0,T],\mathcal{C}^{-\delta})).$$


\begin{thebibliography}{99}
\bibitem[BCD11]{} H. Bahouri, J.-Y. Chemin, R. Danchin,  Fourier analysis and nonlinear
partial differential equations, vol. 343 of Grundlehren der Mathematischen
Wissenschaften [Fundamental Principles of Mathematical Sciences]. Springer, Heidelberg,
2011.
\bibitem[BG97]{} L. Bertini, G. Giacomin,  Stochastic Burgers and KPZ equations from particle
systems. Comm. Math. Phys. 183, no. 3, (1997), 571-607.
\bibitem[Bon81]{} J.-M. Bony,  Calcul symbolique et propagation des singularit\'{e}s pour les \'{e}quations
aux d\'{e}riv\'{e}es partielles non lin\'{e}aires. Ann. Sci. \'{E}cole Norm. Sup. (4) 14, no. 2, (1981),
209-246.
\bibitem[CC13]{}R\'{e}mi Catellier, Khalil Chouk, Paracontrolled Distributions and the 3-dimensional Stochastic Quantization Equation, arXiv:1310.6869
  \bibitem[De13]{}A. Debussche, Ergodicity Results for the Stochastic Navier-Stokes Equations: An Introduction, Topics in Mathematical Fluid Mechanics
Lecture Notes in Mathematics 2013, 23-108
\bibitem[Cor11]{} I. Corwin. The Kardar-Parisi-Zhang equation and universality class. ArXiv 1106.1596
(2011). arXiv:1106.1596.
 \bibitem[DD02]{}G. Da Prato, A. Debussche,  2D-Navier-Stokes equations driven by
a space-time white noise. J. Funct. Anal. 2002, 196 (1), 180-210.
 \bibitem[DD03]{}   G. Da Prato and A. Debussche, Ergodicity for the 3D stochastic Navier-Stokes
equations, \emph{J. Math. Pures Appl.} (9) 82 (2003), no. 8, 877-947.
\bibitem[DDT94]{} G. Da Prato, A. Debussche, R. Temam,  Stochastic Burgers equation. NoDEA Nonlinear Differential Equations Appl. 389-402
        (1994)
        \bibitem[DG01]{} A. M. Davie and J. G. Gaines. Convergence of numerical schemes for the solution
of parabolic stochastic partial differential equations. Math. Comp. 70, no. 233, (2001),
121-134.
       \bibitem[FG95]{} F. Flandoli, D. Gatarek, Martingale and stationary solutions for stochastic Navier-Stokes equations, \emph{Probability Theory and Related Fields} \textbf{102} (1995), 367-391
\bibitem[FR08] {} F. Flandoli, M. Romito, Markov selections for the 3D stochastic Navier-Stokes equations, \emph{Probab. Theory Relat. Fields} \textbf{140} (2008), 407-458
    \bibitem[Gub04]{} M. Gubinelli, Controlling rough paths. J. Funct. Anal. 216, no. 1, (2004), 86-140.
     \bibitem[GIP13]{} M. Gubinelli, P. Imkeller, N. Perkowski, Paracontrolled distributions and singular PDEs, arXiv:1210.2684
    \bibitem[Gy98]{} I. Gy\"{o}ngy. Lattice approximations for stochastic quasi-linear parabolic partial differential
equations driven by space-time white noise. I. Potential Anal. 9, no. 1, (1998),
1-25.
    \bibitem[Gy99]{} I. Gy\"{o}ngy. Lattice approximations for stochastic quasi-linear parabolic partial differential
equations driven by space-time white noise. II. Potential Anal. 11, no. 1, (1999),
1-37.
    \bibitem[Hai11]{} M. Hairer, Rough stochastic PDEs. Comm. Pure Appl. Math. 64, no. 11, (2011), 1547-1585.
doi:10.1002/cpa.20383.
\bibitem[Hai13]{} M. Hairer, Solving the KPZ equation. Ann. of Math. (2) 178, no. 2, (2013), 559-664.
\bibitem[Hai14]{} M. Hairer, A theory of regularity structures. Invent. Math. (2014).
\bibitem[HM06]{} M. Hairer, J. C. Mattingly, Ergodicity of the 2D Navier-Stokes equations with degenerate stochastic forcing  \emph{Annals of Math.}, \textbf{164} (2006), 993-1032
    \bibitem[HM12]{} M. Hairer and J. Maas, A spatial version of the Ito-Stratonovich correction. ˆ Ann.
Probab. 40, no. 4, (2012), 1675-1714
   \bibitem[HMW14]{} M. Hairer, J. Maas, H. Weber, Approximating Rough Stochastic PDEs, Comm. Pure Appl. Math, 67(5): 776-870, 2014
    \bibitem[HW13]{} M. Hairer,  H. Weber, Rough Burgers-like equations with multiplicative noise, Prob. Theory and Rel. Fields 155(1-2): 71-126, 2013
    \bibitem[KPZ86]{} M. Kardar, G. Parisi, Y.-C. Zhang, Dynamic scaling of growing interfaces.
Phys. Rev. Lett. 56, no. 9, (1986), 889-892.
\bibitem[KT01]{}H. Koch, D. Tataru,
Well posedness for the Navier–Stokes equations
Adv. Math., 157 (1) (2001),  22-35
   \bibitem[Lyo98]{} T. J. Lyons, Differential equations driven by rough signals. Rev. Mat. Iberoamericana
14, no. 2, (1998), 215-310.
 \bibitem[RZZ14]{}M. R\"{o}ckner, R.-C. Zhu, X.-C. Zhu, Local existence and non-explosion of solutions for  stochastic fractional partial differential equations driven by multiplicative noise, Stochastic Processes and their Applications 124 (2014) 1974-2002
     \bibitem[S85]{}W. Sickel, Periodic spaces and relations to strong summability of multiple Fourier series. Math. Nachr.
124, 15-44 (1985)

\bibitem[SW71]{}E. M. Stein, G. L. Weiss, Introduction to Fourier Analysis on Euclidean Spaces,
Princeton University Press, 1971
\bibitem[Te84]{} R. Temam, Navier-Stokes Equations, North-Holland, Amsterdam,(1984)
\bibitem[Tri78]{} H. Triebel,  Interpolation theory, function spaces, differential operators. North-Holland Mathematical
Library 18. North-Holland Publishing Co. Amsterdam-New York 1978.
\bibitem[ZZ14]{} Rongchan Zhu, Xiangchan Zhu, Three-dimensional Navier-Stokes equations driven by space-time white noise, arXiv:1406.0047

\end{thebibliography}
\end{document}